\documentclass[a4paper,reqno,10pt]{article}
 \usepackage[utf8]{inputenc}
\usepackage{pifont}
\usepackage[margin=1in]{geometry}
\usepackage[english]{babel}
\usepackage[T1]{fontenc}
\usepackage{times}
\usepackage{amsmath}
\usepackage{amsthm}
\usepackage{amsfonts,mathrsfs,relsize,bigints,stmaryrd,esint}
\usepackage{amssymb}
\usepackage{hyperref}
\usepackage{appendix}
\usepackage{mathtools}
\usepackage{mathabx}
\usepackage{graphicx}
\usepackage{xcolor}
\usepackage{enumitem}
\usepackage{ulem}
\usepackage{math}
\usepackage{dsfont} 
\usepackage{braket}
\mathtoolsset{showonlyrefs}
\makeatletter  \@addtoreset{equation}{section} \makeatother

\newtheorem{defi}{Definition}[section]
\newtheorem{lem}{Lemma}[section]
\newtheorem{theo}{Theorem}[section]
\newtheorem{cor}{Corollary}[section]
\newtheorem{pro}{Proposition}[section]

\newtheorem{rem}{Remark}[section]

\newcommand{\pa}{\partial}

\newcommand{\wrta}{\,\di {\rm A}}

\newcommand{\OpM}{{\rm O}{\rm P}{\rm M}}

\newcommand{\sL}{\mathscr{L}}
\newcommand{\sE}{\mathscr{E}}

\newcommand{\sD}{\mathscr{D}}
\newcommand{\ii }{{\rm i} }

\newcommand{\fv}{\mathfrak{v}}

\newcommand{\LL}{\mathbb{L}}

\newcommand{\bomega}{{\boldsymbol{\omega}}}

\allowdisplaybreaks 

\date{}
\title
{\bf Time-periodic vortices near translating symmetric dipole patches}
\author{Luca Franzoi\footnote{
	Dipartimento di Matematica, Universit\`a degli Studi di Roma Tor Vergata, Via della Ricerca Scientifica 1, 00133 Roma, Italy. \textit{Email:} \texttt{franzoi@mat.uniroma2.it};}, \quad 
    Taoufik Hmidi\footnote{New York University Abu Dhabi, Abu Dhabi, United Arab Emirates.
	 	\textit{Email:} \texttt{th2644@nyu.edu}; 
	}, \quad  Riccardo Montalto\footnote{Dipartimento di Matematica ``Federigo Enriques'', Universit\`a degli Studi di Milano, Via Cesare Saldini 50, 20133 Milano, Italy. \textit{Email:} \texttt{riccardo.montalto@unimi.it}}}

\begin{document}

\maketitle


\begin{abstract}
 We prove the existence of time-periodic solutions of the two-dimensional incompressible Euler equations bifurcating from a translating vortex pair. The reference configuration consists of two symmetric vortex patches of equal strength and opposite sign traveling at constant speed. In a regime of large separation between the vortices, the dynamics may be viewed as a small perturbation of an integrable system.
Working in a co-moving frame and using the contour dynamics formulation, we reduce the problem to a nonlinear transport equation for the vortex boundaries. The linearized operator exhibits degeneracies associated with symmetries and transport effects. By combining a Lyapunov–Schmidt reduction, Nash–Moser scheme, and spectral analysis with sharp asymptotic expansions of the eigenvalues in order to overcome the degeneracy, we construct families of non-rigid, time-periodic vortex patch solutions for a large Cantor set of parameters. The analysis reveals that translating dipoles possess a surprisingly rich nearby dynamics, far beyond the classical rigid paradigm usually associated with vortex patch motion. More generally, the approach developed in this work is flexible and robust, and is expected to extend to a broader class of nonlocal PDEs from Fluid Mechanics with degenerate behaviours.
\end{abstract}

 \medskip



\tableofcontents

\section{Introduction and main results}

\subsection{General overview and formal result}
Time-periodic structures    are among the most striking coherent phenomena arising in inviscid fluid dynamics. In the two-dimensional incompressible Euler equations, such structures often arise from a subtle balance between nonlinearity, symmetry, and long-range interactions, and they play a fundamental role in shaping  vortex patterns. Moreover, it is expected that in 2D inviscid fluids, time-periodic and quasi-periodic structures are the building blocks to determine the long time dynamics for ``generic'' initial data (see the review \cite{DE23}). While rigidly rotating or translating vortex configurations have been extensively studied, much less is known about genuinely time-periodic motions that occur near traveling equilibria, especially when the vortices possess internal degrees of freedom.

In this work, we investigate the existence of small-amplitude time-periodic solutions of the Euler equations bifurcating from a translating vortex pair. In  vorticity formulation, the two-dimensional Euler equations read as
\begin{align}\label{eq:Euler}
\begin{cases}
  \partial_t\omega+(u\cdot \nabla)\omega=0,\quad & (t,{x})\in (0,\infty)\times \R^2\\
  \Delta u=\nabla^{\perp}\omega,\quad & \nabla^\perp=(-\partial_2,\partial_1) \\
  \omega({x},0)=\omega_0({x}),\quad & {x}\in \R^2,
\end{cases}
\end{align}
where $\omega:\R\times \R^2 \to \R$ is the scalar vorticity.
Our starting point is a {\it symmetric translating dipole}, consisting of two simply connected vortex patches of equal strength and opposite sign. The patches translate rigidly with constant velocity, and their centers of mass move along parallel straight-line trajectories. Such translating dipoles are classical objects in vortex dynamics and are typically regarded as rigid coherent structures. Our goal is to describe the dynamics in a neighborhood of these configurations and to establish the existence of small-amplitude time-periodic solutions bifurcating from them. 
From a physical perspective, the problem captures the interplay between the self-induced motion of each patch and the weak interaction exerted by the distant companion vortex. 
This regime  suggests a perturbative approach, but the infinite-dimensional nature of the problem, combined with the presence of symmetries and neutral directions, leads to substantial analytical difficulties. 

From a mathematical point of view, the dynamics of the vortex boundaries is governed by a nonlinear, nonlocal evolution equation derived from contour dynamics, see for instance \cite{BertozziMajda,chemin}. 
After imposing symmetry conditions and working in a moving frame adapted to the translation of the dipole, the problem reduces to finding  solutions of a functional equation where the space domain is the periodic torus, coupling the shape deformation of the patch to the unknown traveling speed.
The resulting equation is Hamiltonian, quasilinear, reversible, and features a degenerate linearization.

In the main result of the present paper, we prove the existence of families of time-periodic, non-rigid vortex patches bifurcating from symmetric translating dipoles. 
 Our analysis highlights the subtle structure of vortex interactions in the Euler equations and provides a framework for studying time-periodic phenomena near traveling equilibria in other nonlocal equations arising in geophysical fluid dynamics.

Our approach is based on a Lyapunov–Schmidt reduction that separates the dynamics into a finite-dimensional bifurcation equation and an infinite-dimensional range equation. The bifurcation equation governs the tangential modes associated with phase and translation invariance, and it is solved using action–angle coordinates and a fixed point argument. The range equation captures the normal dynamics and requires a delicate Nash–Moser iteration to overcome the loss of derivatives caused by small divisors and transport effects.

A central difficulty lies in the analysis of the linearized operator around the traveling vortex pair. Although the operator exhibits a leading-order transport structure, its lower-order terms are smoothing and can be controlled through a careful normal form reduction. This allows us to construct an approximate inverse satisfying tame estimates, provided suitable non-resonance conditions are imposed on the temporal frequency. Measure estimates then show that these conditions are satisfied on a large Cantor-like subset of parameters. A central difficulty is the presence of degenerate modes that requires a very careful perturbative expansion with respect to the distance between the two interacting vortices. This is crucial in order to verify the usual non-resonance conditions which arise in the inversion of the linearized operator. 

\subsubsection{Translating dipoles}
Translating vortex dipoles are among the most classical coherent structures in the two-dimensional incompressible Euler equations. They play a central role in vortex dynamics and arise in a variety of applications, including the modeling of trailing vortices behind aircraft wings and vortex--wall interactions. These configurations consist of pairs of vortices with opposite circulations that propagate at constant speed while rigidly preserving their overall shape. As such, they provide fundamental building blocks for more complex phenomena in geophysical flows, plasma physics, and turbulence.

The simplest explicit example is given by a pair of point vortices with opposite circulation, separated by a fixed distance, which translate uniformly with constant speed, see for \mbox{instance \cite{Lamb}.} This configuration is too singular: it solves the point-vortex system, but not the Euler equations.
Beyond this singular setting, Chaplygin \& Lamb \cite{Chaplygin,Lamb} constructed a smooth vortex pair with compactly supported, non-uniform vorticity, expressed in terms of Bessel functions. These solutions are  commonly referred to as Lamb--Chaplygin dipoles.

In the context of vortex patches, where the vorticity is uniformly distributed over a compact set, a substantial body of numerical and theoretical work has been devoted to the study of translating dipoles, notably through contour dynamics methods. In particular, numerical experiments by Deem \& Zabusky \cite{DZ} and Pierrehumbert \cite{humbert} revealed continuous families of steadily translating symmetric vortex patches and provided evidence for the existence of a bifurcation curve that, in an appropriate asymptotic regime, connects to the point-vortex system.

On the analytical side, variational methods initiated by Turkington \cite{Tur} and further developed by Keady \cite{keady} establish the existence of such rigidly translating vortex pairs by maximizing the excess kinetic energy under suitable constraints. While this approach yields existence and asymptotic estimates for physically relevant quantities such as the propagation speed, it provides only limited insight into the qualitative and topological structure of the solutions. This gap was addressed by Hmidi \& Mateu \cite{Hmidi}, who developed an alternative approach based on contour dynamics and the implicit function theorem in an infinite-dimensional setting. Further developments in this direction, including numerical investigations, the study of limiting profiles, and the exploration of multi-vortex configurations, can be found in \cite{drit,GH23,Hassain-Whee,Hassain1,Saffm,Wu}

Beyond existence of vortex dipoles, substantial progress has been made in the analysis of the stability of translating dipoles. In particular, a number of works have established orbital stability results for these coherent structures by relying on variational techniques inspired by Arnold’s stability theory and rearrangement methods. The central idea in these approaches is to characterize steady dipoles as extrema of suitable energy functionals under appropriate constraints, and to exploit the associated variational structure and symmetry invariance to control the evolution of perturbations. We refer, for instance, to \cite{Abe1, Abe2, Burton, Burton2, wang} for rigorous results in this direction.

Despite these advances, the vast majority of available results concern rigid or steadily translating configurations. By contrast, much less is known about genuinely time-dependent motions occurring in the vicinity of a translating dipole, even in the simplest setting of vortex patches. In this framework, the existence of time-periodic solutions describing non-rigid oscillations of the vortex boundaries around a traveling dipole remains largely unexplored.
From a mathematical viewpoint, this problem is highly challenging due to the infinite-dimensional nature of the dynamics. The analysis is further complicated by strong degeneracies arising from symmetries, as well as the presence of small divisors and space-time resonances.
The present work addresses this gap by constructing families of time-periodic solutions bifurcating from a translating vortex pair. Further details will be provided in Section \ref{time-periodic section}.

\subsubsection{Kida vortices in linear shear and non-rigid periodic dynamics}

Elliptic vortex patches placed in a linear background shear give rise to a remarkable class of exact solutions of the two-dimensional incompressible Euler equations, known as Kida vortices \cite{Kida}. Indeed,  when the ambient velocity field is linear in space, an initially elliptic patch of uniform vorticity remains elliptic for all times. This invariance reduces the infinite-dimensional contour dynamics to a finite-dimensional system of ordinary differential equations governing the semi-axes and orientation of the ellipse. 
The dynamics of the Kida vortex exhibit a rich variety of regimes governed by the interplay between self-induced rotation and external strain. When rotation dominates, the elliptical vortex remains close to a rigidly rotating configuration, akin to Kirchhoff ellipses, undergoing only weak periodic deformation. However, when the strain is comparable to the self-induced rotation, the system enters an oscillatory regime: the ellipse periodically deforms, with its major and minor axes alternating, while its orientation evolves non-uniformly. In contrast, when the external strain dominates, the vortex undergoes strong deformation, with the aspect ratio increasing significantly and leading to highly elongated configurations.

A particularly interesting regime  is  when the initial ellipse is close to the disc. In a frame translating with its center of mass, the motion becomes periodic but not rigid: the boundary oscillates while remaining within the elliptic class. In the absence of shear, this regime reduces to the classical Kirchhoff ellipse. Moreover, it is known from \cite{Burbea,DZ,HMV} that these ellipses constitute the first nontrivial branch in a countable family of rigidly rotating $m$-fold vortex patches bifurcating from the disc.
A natural question is how such rigid periodic motions are affected by the presence of a small external shear. Does the rigidity persist, or does the shear induce genuine deformation? While the first branch, corresponding to the Kida ellipse, is well understood, the behavior of the higher $m$-fold branches remains largely unexplored.

In this work, we address this question in the more realistic setting of interacting vortex pairs. Indeed, the interaction between two well-separated vortices can, to leading order, be approximated by a single vortex subject to an external shear. We show that, under such perturbations, the rigidly rotating branches are transformed into non-rigid time-periodic motions. A detailed analysis of this phenomenon is provided in Section \ref{time-periodic section}.

\subsubsection{Small divisors problems in fluids}
In the last decade, there has been a  surge of works proving the existence of time quasi-periodic waves for PDEs arising in fluid mechanics.
Most of these results in literature are proved by means of hard normal form spectral techniques, microlocal analysis and Nash-Moser methods to deal with the presence of small divisors issues and consequent losses of regularity.
For the two dimensional water waves equations, we mention
Berti \& Montalto \cite{BM}, Baldi, Berti, Haus \& Montalto \cite{BBHM} for time quasi-periodic standing waves and Berti, Franzoi \& Maspero \cite{BFM24,BFM21}, Feola \& Giuliani \cite{FG} for time quasi-periodic traveling wave solutions. The latter works have recently been extended to the much more difficult 3D case by Feola, Montalto \& Terracina \cite{FMT25}.

By means of the contour dynamics, in the last years the existence of time quasi-periodic solutions was proved also for vortex patches in active scalar equations. 
We mention Berti, Hassainia \& Masmoudi \cite{BHM23} for vortex patches of
the Euler equations close to Kirchhoff ellipses, Hmidi \& Roulley \cite{HR21} for the 
quasi-geostrophic shallow water equations,  Hassainia, Hmidi \& Masmoudi \cite{HHM21} for generalized surface quasi-geostrophic equations, Roulley \cite{R22} for Euler-$\alpha$ flows, Hassainia \& Roulley \cite{HR21-1} for Euler equations in the unit disk close to Rankine vortices and Hassainia, Hmidi \& Roulley \cite{HHR23} for 2D Euler annular vortex patches. Very recently, patch-type solutions with a traveling quasi-periodic structure have been constructed for the one dimensional Vlasov-Poisson equation \cite{FR26}.

Another striking application of KAM theory and Nash–Moser methods concerns the classical leapfrogging phenomenon. Hassainia, Hmidi, and Masmoudi established the existence of time-periodic leapfrogging vortex patches for the two-dimensional Euler equations \cite{HHM24}. More recently, Garcia, Hassainia, and Hmidi extended these results to the axisymmetric three-dimensional Euler equations by constructing periodically leapfrogging vortex rings \cite{GHH26}.

Time quasi-periodic solutions were also constructed for the 3D Euler equations with time quasi-periodic external force \cite{BM20} and for the forced 2D Navier-Stokes
equations \cite{FrMo} approaching in the zero viscosity limit time quasi-periodic solutions of the 2D Euler equations for all times. Other recent extensions include space quasi-periodic steady solutions (which become quasi-periodic traveling waves in a moving frame) for the Euler equations in the bounded channel \cite{FMM24}, and large amplitude traveling waves for the $\beta$-plane approximation of forced rotating fluids \cite{BFMT25} with their nonlinear stability \cite{FFM26}. 
At last, we mention some other constructions of time quasi-periodic and almost-periodic solutions to Euler equations \cite{CF13,EPSL23,FM24} on the periodic torus in dimension $3$ and even dimensions:   we remark that these latter solutions are engineered so that there are no small divisors issues to deal with, with consequently much easier proofs and a drawback of not having information on the potential stability of the solutions.

\subsection{Our main contributions}\label{time-periodic section}

As emphasized above, the vortex dynamics  in the vicinity of translating dipoles remain poorly understood. 
In this section, we detail the main results of this paper partially bridge this gap. The first one is Theorem \ref{Theo-1}, which states the existence of steadily traveling dipoles consisting of well-separated patches of order-one size. The second result is Theorem \ref{Theo-2}, which establishes the existence of non-rigid time-periodic solutions in a neighborhood of the traveling pairs of Theorem \ref{Theo-1}.

To properly present the main results, we shall first describe the underlying rigid configuration from which the bifurcation is performed, followed by a detailed scheme of their proofs. 
Our construction differs slightly from the classical desingularization approach developed in \cite{Hmidi}, where the dipole is obtained by smoothing a pair of point vortices of opposite circulation into concentrated vortex patches. In that setting, each patch has diameter of \mbox{order $\varepsilon$,} whereas the distance between them remains of order one.
In contrast, we work in a complementary regime in which the vortex patches have diameter of order one, whereas their separation is taken to be large. In the formal limit where the separation tends to infinity, the two patches become dynamically independent, and the system reduces to that of  a single isolated patch, which is stationary in the circular case.

\subsubsection{Rigid Symmetric Dipoles}\label{Sec:RSD}

\begin{figure}[ht]
\label{figure 1}

\centering

\includegraphics[width=0.42\textwidth]{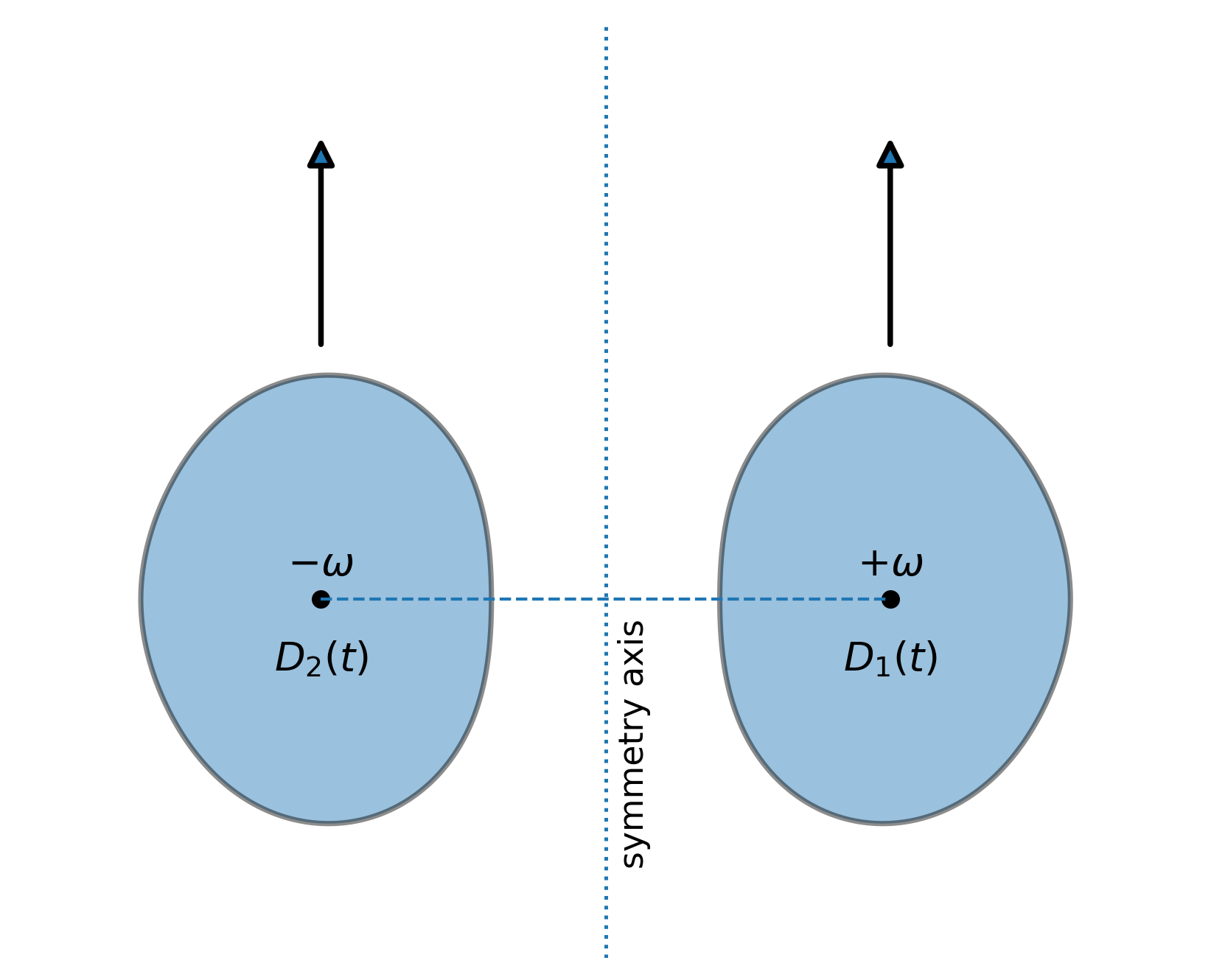}

\caption{Symmetric translating vortex patches with opposite vorticity. }
\end{figure}

We start with the bifurcation of the steadily translating dipole in Theorem \ref{Theo-1}. We consider a symmetric dipole configuration consisting of two bounded vortex patches of opposite sign as described in Figure \ref{figure 1}.
More precisely, let $D_1, D_2 \subset  \C \simeq \R^2 $ be two bounded simply connected domains such that 
\begin{equation}
    D_1 \subset \{ z \in \mathbb{C} \, : \,  {\rm Re} (z) > 0 \}\,,
\qquad
D_2 = \{ -\overline{z}  \,: \, z \in D_1 \}\,,
\end{equation}
so that $D_2$ is obtained from $D_1$ by reflection across the vertical axis. 
We prescribe the initial vorticity
\begin{equation}
   \omega_0(x) = \mathds{1}_{D_1}(x) - \mathds{1}_{D_2}(x) \in L^1(\R^2) \cap L^\infty(\R^2) \,,
\end{equation}
corresponding to two symmetric patches with opposite circulation. 
Following the classical theory of Yudovi\v{c} \cite{Yudo}, the evolution of $\omega_{0}(x)$ remains in the class $L^1(\R^2) \cap L^\infty(\R^2)$, transported through the flow of the Euler equations, and, thanks to the symmetry of the initial data, the configuration is preserved for all times, retaining the form
\begin{equation}\label{struct-vort}
\omega(t,x) = \mathds{1}_{D_1(t)}(x) - \mathds{1}_{D_2(t)}(x) \,,
\qquad
D_2(t) = \{ -\overline{z} : z \in D_1(t) \} \,.
\end{equation}
We first seek special solutions in which each patch evolves by rigid translation, without deformation of its internal shape. 
That is, we consider configurations of the form
\begin{equation}\label{stru-vort1}
D_1(t) = z_1(t) + O\,, \quad D_2(t) = z_2(t) -\overline{O}\,,
\end{equation}
where $O$ is a time-independent  bounded domain (a small perturbation of the unit disc, centered at the origin) 
and $z_k(t) \in \mathbb{C}$ denotes the center of mass of the $k$-th patch, $k=1,2$. The boundary of the steady domain $O$ takes the form
\begin{equation}\label{param.gamma01}
	\gamma(\theta) := \sqrt{1 + 2r(\theta)} e^{\im \theta} \,,  \quad \theta\in\T\,,
\end{equation}
where $r:\T\to\R$ being a smooth  periodic function. In the regime of large separation, the interaction between the two patches is weak and, at leading order, the induced velocity field is accurately described by a modified point-vortex system. 
Accordingly, we assume that the centers satisfy
\begin{align}\label{point.dynX}
\begin{cases}
  -\im \,\overline{\dot{z}_1} (t) & = \frac{2}{z_1(t) - z_2(t)} + \tU  \,, \\
 		-\im\,   \overline{\dot{z}_2} (t) & =   \frac{2}{z_1(t) - z_2(t)} + \tU \,,
\end{cases}
\end{align}
where the function $\tU$ will play the role of a Lagrange multiplier, and will be adjusted with respect to the patch shape. 
From \eqref{point.dynX}, it follows that the centers move along vertical lines while maintaining a constant separation:
\begin{equation}\label{def.alpha.intro}
   |z_1(t) - z_2(t)| = |z_1(0) - z_2(0)|:=\tfrac{1}{2\alpha}  \,.
\end{equation}
Our first main result establishes the existence of rigidly translating vortex patches with analytic boundary.
These solutions can be viewed as perturbations of two independent stationary vortices located at infinite separation.
They provide the fundamental building block for the subsequent construction of non-rigid time-periodic configurations. 
We introduce the analytic function space that captures the regularity of the patch boundary: for $\rho>0$ and $\delta>0$, we define
\begin{equation}
		  B_\delta^{\rho}:=\bigg\{ f=\sum_{j\in\mathbb{Z}} f_j e^{\ii \theta}\, :\,  \| f\|_{X^{\rho}}  := \sum_{j\in \Z} e^{\rho|j|}  |f_{j}| \leqslant \delta  \bigg\}\,.   
\end{equation}
The exponential weight ensures analyticity in a strip of width $\rho$.
Our first main result reads as follows.

\begin{theo}\label{Theo-1}
Let $k\in\mathbb{N}$. There exist $\alpha_0 \in (0,1)$ small enough and $\delta>0$ such that the following holds.
There exist two $C^k$ curves
\[
\alpha\in[-\alpha_0,\alpha_0]\longmapsto r_\alpha\in B_\delta^{\rho}\,,
\qquad
\alpha\in[-\alpha_0,\alpha_0]\longmapsto \tU_\alpha\in\mathbb{R}\,,
\]
for which equation \eqref{eq:Euler} admits a solution of the form \eqref{struct-vort}, \eqref{stru-vort1}, and \eqref{param.gamma01}.
Moreover, as $\alpha\to 0$, the following asymptotic expansions hold
\begin{align}
		r_\alpha(\theta) = - \alpha^2 \cos(2\theta)  +\tfrac12 \alpha^3 \cos(3\theta)+\tfrac23 \alpha^4 \cos(4\theta)  + O(\alpha^5)\quad\hbox{and}\quad \tU_\alpha=O(\alpha^4)\,.  \label{r.alpha.th1}
   \end{align}
\end{theo}

\begin{rem} Some remarks on Theorem \ref{Theo-1} are in order:
\begin{enumerate}[label=(\roman*)]
    \item The parameter $\alpha$ encodes the distance between the two vortices and is therefore naturally positive. However, as will become apparent from the formulation, the equation governing the boundary admits a natural extension to small negative values of $\alpha$ as well.
    \item The boundary of the traveling vortex pairs is analytic. It is worth noting that in \cite{Hmidi}, the corresponding solutions constructed near point vortices were shown to possess only H\"older regularity. This regularity was subsequently improved in \cite{GH23}, where analyticity of the boundary was established by reformulating the problem as an elliptic free boundary problem and invoking the classical regularity theory of Kinderlehrer, Nirenberg \& Spruck \cite{Kin1}. In contrast, our approach provides a more direct and self-contained route to analyticity, relying on the implementation of the implicit function theorem within an analytic functional framework.
\end{enumerate}
    
\end{rem}
We quickly describe the strategy of the proof of Theorem \ref{Theo-1}. It relies on a perturbative construction around the trivial configuration corresponding to two independent circular vortices at infinite separation. In this limiting regime, the interaction between the patches vanishes and the problem reduces to a stationary solution of the Euler equations. The main idea is to treat the finite but large separation case as a small nonlinear perturbation of this integrable configuration.

We parametrize the boundary of each patch as a normal graph over the unit circle through a function $r(\theta)$ belonging to the analytic space $B_\delta^{\rho}$. Within the contour dynamics formulation of the vortex motion, this leads to a nonlinear functional equation coupling the boundary deformation and the translation speed \( \tU \).
The construction proceeds by expanding the nonlinear functional with respect to the small parameter $\alpha$, which measures the inverse separation distance. At leading order, the linearized operator around the circular configuration is diagonal in Fourier variables. Its kernel corresponds to the translational invariance and is handled by introducing the speed $\tU$ as an additional unknown. After fixing this degeneracy, the linearized operator becomes invertible on the orthogonal complement.

The core of the proof is an application of the implicit function theorem in the analytic Banach space $B_\delta^{\rho}$. A careful analysis of the nonlinear terms shows that the functional is smooth in both $r$ and $\alpha$, and that the linearized operator is an isomorphism for sufficiently small $\alpha$. This yields the existence of $C^k$ curves $\alpha \mapsto r_\alpha$ and $\alpha \mapsto \tU_\alpha$ solving the nonlinear equation.

Finally, the asymptotic expansions are obtained by computing successive orders in $\alpha$ in the Fourier representation of the equation. The quadratic interaction generates the $\cos(2\theta)$ mode at order $\alpha^2$, while higher-order nonlinear interactions produce the $\cos(3\theta)$ and $\cos(4\theta)$ contributions. The translation speed appears at order $\alpha^4$ due to symmetry cancellations at lower orders.

\subsubsection{Non rigid time-periodic solutions near the traveling  dipoles}
We now turn to the construction of non-rigid time-periodic vortex patches in a neighborhood of the rigidly translating symmetric dipoles obtained in Theorem \ref{Theo-1}. While the dipoles constructed above move by pure translation with fixed shape, our goal here is to exhibit nearby solutions whose boundaries undergo genuine time-dependent deformations, yet remain time-periodic in a suitable moving frame.
To capture these non-rigid dynamics, we introduce a slight modification of the ansatz discussed in Section~\ref{Sec:RSD}. We continue to look for solutions of the form \eqref{struct-vort}, but we now allow the reference domain $O$ in \eqref{stru-vort1} to deform  over time. More precisely, we parametrize its  boundary as
\begin{align}\label{param.gamma02}
\gamma(t,\theta)
=
\sqrt{1 + 2 r(t,\theta)}\, e^{\im \theta}\,,
\qquad
t \in \mathbb{R}\,, \quad \theta \in \T\,,
\end{align}
where the deformation function $r(t,\theta)$ is now time-dependent.
Regarding the motion of the centers, we retain the reduced system \eqref{point.dynX}, except that the translation speed is no longer constant: the parameter $\tU$ is now replaced by a time-dependent function $\tU=\tU(t)$. 
We use this speed as a modulation parameter in order to remove the partial degeneracy associated with the first sine mode, which stems from the translation invariance along the vertical axis. The function $\tU(t)$ will be determined uniquely by imposing an appropriate orthogonality condition, and therefore depends in a nonlinear way on the evolving shape of the patch.
\\
These time-periodic solutions are constructed in Sobolev spaces, defined as follows
\begin{equation}\label{sob.space.intro}
	H^s:= H^s(\T^2) := \bigg\{\, f:\T^2 \to \R   \, : \, \| f \|_{H^s(\T^2)}^2 :=  \| f \|_{s}^2  := \sum_{(\ell,j)\in \Z^2} \braket{\ell,j}^{2s} |f_{\ell,j}|^2 < \infty  \bigg\}\,,  
\end{equation}
where $\braket{\ell,j}:= \max\{ 1, |\ell|, |j| \}$. We also need to define the following subspaces
\begin{equation}
     H_\circ^{s}:=H_{\circ}^{s}(\mathbb{T}^{2}):=\Big\{ h\in H^{s}(\mathbb{T}^{2})  \, : \, \int_{\T} h(\vf,\theta) \wrt\theta = \int_{\T}h(\varphi,\theta)\sin(\theta)\wrt\theta=0,\, \ \forall\,\varphi\in\T \Big\}\,. \label{Hcirc.def.intro}
\end{equation}
The space $H_\circ^{s}$ plays a crucial role in our construction, as we are going to explain in Remark \ref{remark.theo2} below.

We are now in a position to state our second main result. 

\begin{theo}\label{Theo-2}
There exist $\alpha_0>0$ small enoughe, $s \in \N$ large enough, and a family $\alpha\in[0,\alpha_0]\mapsto \cM_{\alpha}$ of  smooth, invertible, reversibility preserving linear operators $\cM_{\alpha}:H_{\circ}^{s}(\T^2)\to H_{\circ}^{s}(\T^2)$, with $H_{\circ}^{s}(\T^2)$ as in \eqref{Hcirc.def.intro},  such that,  for any $0<\alpha_1<\alpha_2 < \alpha_{0}$ and for any integer $\tJ\geq 2$, the following holds.   There exists $0<\varepsilon_0 \ll \alpha_1$  such that, for any $\varepsilon \in (0,\varepsilon_0)$, there exists a Cantor-like set $\mathtt{C}_\varepsilon\subset [\alpha_{1},\alpha_{2}] $ of asymptotically full measure, namely
\begin{equation}\label{cantor.limit.theo2}
    \lim_{\varepsilon\to 0}|\mathtt{C}_\varepsilon|= \alpha_{2}-\alpha_{1}\,,
\end{equation}
such that, for each $\alpha \in\mathtt{C}_\varepsilon$, there exists a $\frac{2\pi}{\omega_{\mathtt{J,\varepsilon}}}$-periodic in time solution to the  equation \eqref{eq:Euler} supplemented with  \eqref{struct-vort}, \eqref{stru-vort1},  \eqref{param.gamma02}, of the form
\begin{equation}\label{sol.theo2}
     r(t,\theta)=r_\alpha(\theta)+\varepsilon\cM_{\alpha}[\cos (\mathtt{J}\cdot - \omega_{\mathtt{J,\varepsilon}}t)](\theta)+\varepsilon^2P(\omega_{\mathtt{J,\varepsilon}} t,\theta)\,,\quad 
\tU(t)= \tU_\alpha+\varepsilon\mathtt{W}(\omega_{\mathtt{J},\varepsilon} t)\,,
\end{equation}
with $(r_\alpha(\theta),\tU_{\alpha})$ as in Theorem \ref{Theo-1}, for some functions $P\in H_\circ^s(\T^2)$ even in $(\vf,\theta)$,  and $\mathtt{W}\in H^s(\T)$, oscillating in time with frequency $\omega_{\mathtt{J},\varepsilon}>0$ satisfying
\begin{equation}
    \omega_{\mathtt{J},\varepsilon}
		=\omega_{\tJ,0}+O(\varepsilon)\,,
		\qquad
		\omega_{\tJ,0}=
 \big(\tfrac12 -  \alpha^4 +O(\alpha^5) \big)\tJ - \tfrac12 + O(\alpha^5)\,. 
\end{equation}
\end{theo}

\begin{rem}\label{remark.theo2}
Some remarks on Theorem \ref{Theo-2} are in order:
\begin{enumerate}[label=(\roman*)]
    \item We construct time-periodic solutions bifurcating from a fixed  steady solution $(r_\alpha,\tU_{\alpha})$ of Theorem \ref{Theo-1} that is even in $(\vf,\theta)$ and belongs to the space $H_\circ^s(\T^2)$ defined in \eqref{Hcirc.def.intro}. The requirement to have zero average in space follows from the conservation of the mass of the system, which can be easily seen in \eqref{Func.G-ini} later. On the other hand, the requirement of not having the mode $\sin(\theta)$ in the Fourier expansion is needed to deal with the important degeneracy of the mode 1 that we have already come across in the construction of the steady solutions of Theorem \ref{Theo-1}. The unknown $\tU(t)$ will play indeed the role of Lagrange multiplier to ultimately guarantee that the  time-periodic solution  in \eqref{sol.theo2} lives in the space $H_\circ^s$, overcoming the issue generated by the degeneracy;
\item At the first leading order in $\varepsilon$, the time-periodic perturbation behaves as a solution of the linearized equation at the steady solution localized on the Fourier mode $\tJ\geq 2$ in space. In particular, the mode $\tJ=2$ is included: this is different compared to the (quasi-)periodic solutions constructed in \cite{BHM23} close to Kirchhoff ellipses. Indeed, in this latter work, the steady ellipses have a 2-fold symmetry, essentially leading to a degeneracy of the mode 2 which has to be excluded by the choice of the tangential modes. In our result, the family of rigidly traveling equilibria has only a 1-fold symmetry, which can be read from the asymptotic expansion \eqref{r.alpha.th1} in Theorem \ref{Theo-1}. \\
The case $\tJ=2$ is also intimately related to the dynamics of Kida vortices \cite{Kida} as described before: in the regime of a weak external linear shear, time-periodic deformations of the elliptic patches do not break the 2-fold symmetry, with small oscillations of the aspect ratio or fluctuations of the angle of the mutating axis. This cannot happen to our solutions even if the leading order perturbation is 2-fold, since the long-distance interaction between the two patches brings contributions to all the modes. As a pure mathematical model, if we truncate the nonlinear interaction between the two patches at its quadratic leading order, one should expect to recover the 2-fold symmetry and might be able to construct, still in the long-distance regime, a pair of translating ellipses that oscillate in time as Kida vortices. \\
On the other hand, the mode $\tJ=1$ is excluded as it is degenerate. Nevertheless, the existence of solutions that are mode-1 time-periodic perturbations of steady solutions might be investigated without preserving the translation invariance of the dipole and by performing a more refined analysis close to the degeneracy;
\item The map $\cM_{\alpha}$ is a time independent linear operator that we will construct with the aim of diagonalizing the linearized equation at the nontrivial equilibrium $r_\alpha(\theta)$. This is the key step in order to extract the spectral information around the equilibrium that we need when searching for time-periodic solutions. Furthermore, we will transform the full nonlinear system by means of this linear map to ensure that, when linearizing at any approximate solution of the Nash-Moser iteration, the leading order operator is already diagonal; 
\item The present work opens several directions for future investigation. 
A natural problem is the construction of quasi-periodic vortex patches and higher-dimensional invariant tori near translating dipoles. 
Another challenging question concerns the nonlinear stability and long-time dynamics of the periodic solutions constructed here. 
More broadly, the methods developed in this paper suggest that KAM-type phenomena may arise in a much wider class of nonlocal active scalar equations.
\end{enumerate}
\end{rem}

\subsection{Scheme of the proof of Theorem \ref{Theo-2}}

We now describe in more technical detail the main mechanisms underlying the proof 
of Theorem \ref{Theo-2}. The construction combines perturbative analysis around the rigid 
dipole of Theorem \ref{Theo-1}, a Nash-Moser scheme and spectral analysis of the linearized operators adapted to a 
quasilinear, nonlocal Hamiltonian PDEs.
\paragraph{Reformulation around the translating dipole.}
We begin by working in a frame translating with the symmetric dipole constructed in 
Theorem \ref{Theo-1}. In this frame, the vortex boundaries are parametrized as small time-dependent 
deformations of the rigid profile. The evolution equation for the boundary deformation 
is written as a nonlinear, nonlocal equation on the bi-dimensional torus $(\varphi,\theta)\in\mathbb{T}^2$, 
where $\varphi=\omega t$ is the time angle. The translation speed is promoted to a 
time-dependent modulation parameter in order to remove the degeneracy associated with 
vertical translation invariance.
As shown in Section \ref{model.derivation} and Section \ref{Ref-VP}, after passing to a co-moving frame and introducing the time angle 
$\varphi=\omega t \in \T$, the boundary deformation radius $r(\varphi,\theta)$ satisfies a 
quasilinear equation of the form
\begin{align}\label{Func.G-ini}
		\cG(\tU,r;\alpha):= \omega \pa_{\vf} r + \pa_{\theta} \Big( \tU(\vf) \sqrt{1+2 r} \cos(\theta) + F[r;\alpha] \Big)=0 \,, 
         \end{align}
       where $\alpha$ is as in \eqref{def.alpha.intro}, $F[r;\alpha]$ is a nonlinear nonlocal operator induced by the 
Biot–Savart law,
         \begin{align}
         \label{stream.F25.intro}
		\nonumber F[r;\alpha](\vf,\theta) 
		& =   \frac{1}{2\pi} \int_{O(\vf)} \log \big|  \sqrt{1+2r(\vf,\theta)}e^{\im\theta} -\xi \big| \wrta(\xi) 
		\\
		& \quad - \frac{1}{2\pi} \int_{O(\vf)} K_\alpha\big(\sqrt{1+2 r(\vf,\theta)}e^{\im\theta},\xi\big) \wrta(\xi) \,,
	  \end{align}
      with the kernel $K_{\alpha}$ explicitly defined in \eqref{Kalpha.kernel} and 
     the domain $O(\varphi)$  defined, for each $\vf\in\T$, by
      \begin{equation}\label{O.intro}
 	O(\vf):= \Big\{ \xi = \rho e^{\im\eta} \,:\, \eta\in[0,2\pi], \ \  0\leq \rho  \leq \sqrt{1+2r(\vf,\eta)} \Big\}\,.
 \end{equation}
By \eqref{Func.G-ini}, \eqref{stream.F25.intro}, \eqref{O.intro}, we note that
 \begin{equation}
     \cG : H_{\rm even}^s(\T) \times H_{\rm even}^s(\T^2) \times [-\alpha_0,\alpha_0] \to H_{\rm odd}^{s-1}(\T^2) \,,
 \end{equation}
 with $ H_{\rm even}^s(\T^2)$, $ H_{\rm odd}^s(\T^2)$ as in \eqref{FS-even}, \eqref{FS-odd}.
   In Theorem \ref{Theo-1}, we  construct smooth curves of functions $\alpha\in[-\alpha_0,\alpha_0]\longmapsto \big(r_\alpha(\theta), \tU_\alpha\big)$ such that
   $$
   \forall \, \alpha\in[-\alpha_0,\alpha_0]\,,\quad\cG(\tU_\alpha,r_\alpha;\alpha)=0\,.
   $$
   Our goal is to construct {reversible} time-periodic solutions of the form
\begin{equation}\label{perturbed.solution0}
	\tU(\vf) = \tU_{\alpha} + \varepsilon \tW(\vf) \,, \quad r(\vf,\theta) = r_\alpha(\theta) + \varepsilon \rho(\vf,\theta)\,,
\end{equation}
where $\rho\in H^{s}_{\circ,\textnormal{even}}$ (with this space  defined in \eqref{Hs-circ}) admits the decomposition 
\begin{equation}
	\rho(\vf,\theta) =\rho_{1}(\vf) \cos(\theta) + \rho_{\geqslant 2}(\vf,\theta) \,, \quad \rho_{\geqslant 2}= \Pi_{\geqslant 2} \rho_{\geqslant 2}\,,
\end{equation}
with  $\Pi_{\geqslant 2}$ being the projector on the Fourier modes $|j|\geq 2$, see \eqref{proj.geq2}.
Then, solving equation \eqref{Func.G-ini} in a neighborhood of the traveling pairs is equivalent to solving the following problem in terms of $(\tW, \rho)$,
\begin{equation}\label{full.nonlin.eq.pert0}
	\wt\cG (\tW,\rho;\alpha,\varepsilon) := \varepsilon^{-1} \cG(\tU_{\alpha}+ \varepsilon \tW, r_\alpha+ \varepsilon \rho;\alpha)=0\,.
\end{equation}

\paragraph{Trivial degeneracy and speed modulation.}
In the limit system at $\alpha=0$, the Rankine disk vortex, namely $r\equiv 0$, is an exact trivial equilibrium of the system. In the proof of Theorem \ref{Theo-1}, to bifurcate the rigid dipole equilibrium $r_\alpha$ for $|\alpha|\ll 1$ small by an implicit function argument, we are led to invert the linear operator $\di_{r} \pa_{\theta} F[0;0]$, which acts on the exponential basis as (see Proposition \ref{lemma.constrain1}-$(iv)$)
\begin{equation}
   \forall\, j\in \Z\setminus\{0\}\,, \quad \di_{r} \pa_{\theta} F[0;0] \big[ e^{\im j \theta} \big] =  \tfrac{\im \,j}{2|j|} (|j|-1) e^{\im j \theta} \,.
\end{equation}
It is clear that such operator cannot be inverted when $|j|=1$. 
This degeneracy is the source of the main difficulties of the problem. 
We solve this issue by promoting the translation speed $\tU(\varphi)$ to a modulation parameter. Morally speaking, we use it as a ``Lagrange multiplier'' to ensure that the map $\cG$ in \eqref{Func.G-ini} acts only on Sobolev spaces $H_\circ^s(\T^2)$ as in \eqref{Hcirc.def.intro}: we determine $\tU \equiv \tU[r]$ by imposing the orthogonality condition
\begin{equation}\label{ortho.circ.cond}
    \int_\T \cG(\tU,r;\alpha)(\vf,\eta) \sin(\eta) \wrt\eta = 0 \,.
\end{equation}
When $\varepsilon=0$, the speed $\tU_{\alpha}\equiv \tU_{\alpha}[r_\alpha]$ in Theorem \ref{Theo-1} has been determined by the same condition. This implies that, recalling \eqref{perturbed.solution0}, the condition  \eqref{ortho.circ.cond} actually fixes $\tW\equiv \tW[\rho]$:
it is not an independent unknown but is uniquely determined by the shape 
deformation $\rho$. In other words, the translation speed is fully fixed 
through the geometry of the patch.
After fixing $\tW$ in this way in Proposition \ref{prop.U1}, the original problem \eqref{full.nonlin.eq.pert0} reduces to a single equation 
for the deformation $\rho$, that is
\begin{align}\label{equation.rho.full0intr}
\tG(\rho;\alpha,\varepsilon)
:=
\wt\cG(\tW[\rho],\rho;\alpha,\varepsilon)
=
0 \,,
\end{align}
and, thanks to the orthogonality condition \eqref{ortho.circ.cond}, we ensure that
\begin{equation}
    \tG(\cdot;\alpha,\varepsilon):
H^{s}_{\circ,\mathrm{even}}
\longrightarrow
H^{s-1}_{\circ,\mathrm{odd}} \,,
\end{equation}
where the spaces $H^s_{\circ}$, as defined through \eqref{Hcirc.def.intro}, incorporate the reversibility symmetry and 
the restriction excluding the first sine mode. We point out that the translating dipole is the equilibrium corresponding to $\rho\equiv 0$, that is,
\begin{equation}
    \forall\, \varepsilon\,, \  \forall \,|\alpha|\leqslant\alpha_0\,, \quad \tG(0;\alpha,\varepsilon)=0 \,.
\end{equation}
This reformulation in \eqref{equation.rho.full0intr} transforms the problem into an infinite-dimensional equation free of neutral directions. It is this reduced equation that will be analyzed in the following sections.


\paragraph{Linearization at the translating dipoles and reducibility.}
A key step in the proof is the spectral study of the linearized operator at the rigidly
translating symmetric dipole constructed in \mbox{Theorem \ref{Theo-1}.} 
This analysis, carried out in Section \ref{sect.linear.eq}, reveals a fundamental structural feature:
at the linear level, the dynamics at  the traveling dipole is completely reducible 
to a Fourier multiplier by means of {time-independent} transformations. 
In particular, no Cantor-type restriction or Melnikov condition is required at this stage.
\\
More precisely, the linearized operator admits the asymptotic structure
\begin{equation}\label{equi.linear.intro}
    \mathcal{L}_{0,\alpha}
:= \partial_\rho \tG(0;\alpha,\varepsilon)
=
\omega \partial_\varphi
+
\partial_\theta \circ \tV_{0,\alpha}
-\tfrac{1}{2}\mathcal{H}
+
\mathcal{R}_\alpha\,,
\end{equation}
where $\tV_{0,\alpha}\equiv \tV_{0,\alpha}(\theta)$ is a function,  $\mathcal{H}$ denotes the toroidal Hilbert transform, 
and $\mathcal{R}_\alpha$ is an operator of order $-\infty$ with coefficients depending only on $\theta$. We already know that the unperturbed translation speed satisfies $\tU_{\alpha}=O(\alpha^4)$.
A detailed asymptotic expansion of $\mathcal{R}_\alpha$ is given in Lemma \ref{lemma.R.expand}. 
Furthermore, Lemma~\ref{speed.V.expan} shows that the transport coefficient 
$\tV_{0,\alpha}$ admits the expansion
$$
 \tV_{0,\alpha}(\theta)= \tfrac12 + \tU_{\alpha}\cos(\theta) +\alpha^2\cos(2\theta) -\tfrac34 \alpha^3 \cos(3\theta)  + \tfrac{2}{3}\alpha^4\cos(4\theta)  + O(\alpha^5) \,.
$$
The diagonalization of the operator \eqref{equi.linear.intro} is carried out in two steps:
\begin{itemize}
\item[1)] The first step consists in reducing the transport operator, which is developed in Section \ref{Redu-transpo}. 
Using a time-independent, space-periodic change of variables 
$\theta \mapsto \theta + \beta(\theta)$, implemented by a bounded symplectic 
isomorphism, we conjugate the transport part to a constant-coefficient 
transport operator. 
We remark once more that this reduction does not involve any 
small divisor argument. More precisely, in Proposition~\ref{diffeo.conj.no.proj}
we construct an invertible operator
\[
\cS_{\rm ph} : H_\circ^s(\T^2) \longrightarrow H_\circ^s(\T^2)
\]
such that the conjugation of the  operator $\cL_{0,\alpha}$ takes the form
        $$
	\cL_{0,\alpha}^{(1)}  := 	\cS_{\rm ph}^{-1} \cL_{0,\alpha} \cS_{\rm ph}=\omega\,\pa_{\vf} + \cD_{0,\alpha}+ \cE_{0,\alpha},
	$$
where the operator $\cD_{0,\alpha}$ is diagonal in the Fourier basis of $H_\circ^s(\T)$ and is given by 
	\begin{equation*}
		\begin{aligned}
\cD_{0,\alpha} :=  {\rm diag}\big\{ \lambda_{j}(\alpha) \,: \,\ j \in \Z_{\rm ph} \big\}\,, \quad \lambda_{j}(\alpha)  &:= \begin{cases}
            0 & j=1\,, \\
			    \im \big[j\big(\tfrac12 - \alpha^4 + O(\alpha^5)\big)
- \tfrac12 {\rm sgn}(j) \big]\,, &  |j|\geq 2 \,,
			\end{cases} 
		\end{aligned}
	\end{equation*}
   Here, with a slight abuse of notation, we set
    \begin{equation}
        \Z_{\rm ph} := \Z \setminus \{0,-1\} \,,
    \end{equation}
    where the index \(j=1\) corresponds to the Fourier mode
\(\cos\theta\).
    The remainder $\cE_{0,\alpha} \in \OpM^{-\infty}$ is a smoothing operator
of order $-\infty$, whose dependence on $\alpha$ admits an explicit expansion
up to order $O(\alpha^5)$;
\item[2)]  The next step, carried out in Section \ref{section-Full-redu} (see also Proposition \ref{prop.diag.red.equi}), is to complete the reduction of the linearized operator by fully diagonalizing $\cL_{0,\alpha}^{(1)}$. Remarkably, this diagonalization is also  achieved without imposing any small divisor condition, by solving in a single step a nonlinear homological equation with formal power series expansions. We construct  a reversibility-preserving transformation  $\mathcal{U} = \mathrm{Id} + \mathcal{Y}$, with $\mathcal{Y}$ being a smoothing operator with small size in $\alpha$, such that
\begin{align}
		&  \cU^{-1} \cL_{0,\alpha}^{(1)} \cU = \omega\,\pa_{\vf} + \cD_{0,\alpha}^{(\infty)} :H_{\circ,{\rm even}}^s(\T^2) \to H_{\circ,{\rm odd}}^{s-1}(\T^2) \,,
	\end{align}
	where $\cD_{0,\alpha}^{(\infty)} $ is the diagonal operator
	\begin{align*}
		\cD_{0,\alpha}^{(\infty)} & := \cD_{0,\alpha} + \cZ_{0,\alpha} =  {\rm diag}\big\{ \lambda_{j}^{(\infty)}(\alpha) \, : \, j\in\Z_{\rm ph} \big\}
        \end{align*}
        and the eigenvalues take the form
        \begin{align}\label{frequency.intro}
		\lambda_{j}^{(\infty)}(\alpha) &:= \begin{cases}
			0\,, & j= 1 \\
			\lambda_{j}(\alpha) + \tz_{j}(\alpha) \,, & |j|\geq 2,
		\end{cases} 
	\end{align}
	with 
    \begin{equation}\label{mode.2.asympt.intro}
        \tz_{\pm 2}(\alpha) = \pm \tfrac32 \alpha^4 + O(\alpha^5) \,, \quad \tz_{j}(\alpha) = O(\alpha^5) \quad \forall \, |j|\geq 3\,,
    \end{equation}
    and 
		\begin{align}
			&	\lambda_{-j}^{(\infty)}(\alpha)=	-\,\overline{\lambda_{j}^{(\infty)}}(\alpha) \,, 
			\quad |j|\geq 2\,.
		\end{align} 
\end{itemize}
    Hence, the linearized dynamics at the translating dipole is spectrally integrable. 
     Combining back to the previous two reduction steps (see Proposition \ref{prop.Malpha}), and, defining
\begin{equation*}
	\cM_{\alpha} := \cS_{\rm ph} \cU \,,
\end{equation*}
we obtain the complete diagonalization
\begin{equation*}
	\cM_{\alpha}^{-1} \cL_{0,\alpha} \cM_{\alpha}
	= \omega\,\pa_{\vf} + \cD_{0,\alpha}^{(\infty)}
	: H_{\circ,{\rm even}}^{s}(\T^2)
	\to H_{\circ,{\rm odd}}^{s-1} (\T^2)\,.
\end{equation*}
Having completely reduced the linearized dynamics around the translating dipole, we now conjugate the nonlinear equation \eqref{equation.rho.full0intr} through the transformation $\cM_{\alpha}$, namely,
\begin{equation}\label{Eq-New1}
	\begin{aligned}
		\bG(\fu;\omega,\alpha,\varepsilon)
		:= \cM_{\alpha}^{-1}\big[ \tG(\cM_{\alpha}[\fu];\alpha,\varepsilon)\big]=0 \,.
	\end{aligned}
\end{equation}
This reformulation transfers the problem into coordinates adapted to the spectral structure of the equilibrium, where the principal linear dynamics is completely diagonalized. In particular, the trivial configuration remains an equilibrium state, namely
\begin{equation}\label{bG.equizero}
	\forall \,  \omega\in\R\,, \  \ \varepsilon\in\R\,, \  \ |\alpha|\leqslant\alpha_0\,,
	\qquad
	\bG(0;\omega,\alpha,\varepsilon)=0 \,.
\end{equation}
Moreover, the linearization of $\bG$ at the origin is now given by the diagonal operator
\begin{equation}\label{Normal-lin-equili}
\partial_{\fu} \bG(0;\omega,\alpha,\varepsilon)
=
\omega\,\pa_{\vf} + \cD_{0,\alpha}^{(\infty)}\,,
\end{equation}
which provides the fundamental starting point for the nonlinear bifurcation analysis carried out in the subsequent sections.

\paragraph{Lyapunov--Schmidt reduction and solution to the bifurcation equation.}
This step will be carried out in Section \ref{Luap-Schm}. The main idea is to perform a Lyapunov–Schmidt reduction to the equation \eqref{Eq-New1} by decomposing the dynamics into two components: a finite-dimensional part generated by the distinguished low Fourier modes, and an infinite-dimensional part corresponding to the normal directions. This decomposition isolates the potentially resonant modes responsible for the main bifurcation mechanism, while the remaining components are treated through the analysis of the associated range equation.
More precisely, we look for solutions of \eqref{Eq-New1} of the form
\begin{equation*}
	\fu(\vf,\vartheta)
	=
	\fu_{\intercal}(\vf,\vartheta)
	+
	\fu_{\perp}(\vf,\vartheta)\,,
\end{equation*}
where the tangential component contains the degenerate first Fourier mode together with the selected resonant mode indexed by $\tJ\geqslant2$. More explicitly,
\begin{equation*}
	\fu_{\intercal}(\vf,\vartheta)
	:=\Pi_{\intercal}\fu(\vf,\vartheta)=
	\fc_{1}(\vf)\cos(\vartheta)
	+
	\fa_{\tJ}(\vf)\cos(\tJ\vartheta)
	+
	\fb_{\tJ}(\vf)\sin(\tJ\vartheta)\,,
\end{equation*}
where the functions 
\[
\fc_{1},\fa_{\tJ},\fb_{\tJ}:\T\to\R
\]
are smooth and satisfy the following  symmetry conditions  reflecting the reversibility structure of the
equation:
\begin{equation*}
\fc_{1}(\vf), \, \fa_{\tJ}(\vf)\quad \text{to be even in } \vf\,,
\qquad
\fb_{\tJ}(\vf)\quad \text{to be odd in } \vf\,.
\end{equation*}
On the other hand, the normal component $\fu_{\perp}$  belongs to an infinite-dimensional  subspace and is defined through the projection onto the complementary Fourier modes
\begin{equation*}
	\fu_{\perp}(\vf,\vartheta)
	=
	\Pi_{\perp}\fu(\vf,\vartheta)
	:=
	\sum_{j\in\Z_{\perp}}
	\wh\fu_{\perp}(\vf,j)\te_j(\vartheta),\quad \Z_{\perp}:=\big\{j\in\Z \,:\,  |j|\geqslant 2, \ j\neq\pm \tJ  \big\}\,,
\end{equation*}
and is assumed to be even in $(\vf,\vartheta).$ 
 As a consequence, solving the
nonlinear equation \eqref{Eq-New1}  becomes equivalent to solving the coupled system
\begin{equation}
\begin{cases}\label{LyaSchint}
	\bG_{\intercal}(\fu_{\intercal},\fu_{\perp};\omega,\alpha,\varepsilon)
	:=
	\Pi_{\intercal}\bG(\fu_{\intercal}+\fu_{\perp};\omega,\alpha,\varepsilon)
	=0\,, & \textnormal{[\textit{bifurcation} equation]}  \\
	\bG_{\perp}(\fu_{\intercal},\fu_{\perp};\omega,\alpha,\varepsilon)
	:=
	\Pi_{\perp}\bG(\fu_{\intercal}+\fu_{\perp};\omega,\alpha,\varepsilon)
	=0\,. & \textnormal{[\textit{range} equation]} 
\end{cases}
\end{equation}
Compared to the standard approach, we choose to solve the system \eqref{LyaSchint} in the following order:
\begin{itemize}
	\item[1)] For a fixed normal component $\fu_{\perp}$, we solve the finite dimensional
	{\it{bifurcation equation}}
	with respect to $\fu_{\intercal}=\fu_{\intercal}(\vf,\vartheta;\fu_{\perp})$;
	\item[2)] Next, we insert $\fu_{\intercal}= \fu_{\intercal}(\fu_\perp)$ into the infinite dimensional {\it{range equation}}
	and we solve it with respect to $\fu_{\perp}=\fu_{\perp}(\vf,\theta)$;
	\item[3)]
    Finally, we come back to the bifurcation equation, both for determining $\fu_{\intercal}= \fu_{\intercal}(\fu_\perp)$ and the final frequency of oscillation
	\begin{equation}
		\omega_{\varepsilon}
		=\omega_{\varepsilon}(\fu_{\perp})
		=\omega_{0}+O(\varepsilon) \,, 
		\qquad
		\omega_{0}=-\im\,\lambda_{\tJ}^{(\infty)}(\alpha) \,,
	\end{equation}
    with $\lambda_{\tJ}^{(\infty)}(\alpha)$ the unperturbed frequency in \eqref{frequency.intro} relative to the tangential mode $\tJ\geq 2$.
\end{itemize}
According to this plan,
for a fixed normal component \(\fu_{\perp}\), we first solve the finite-dimensional
bifurcation equation
\[
\bG_{\intercal}(\fu_{\intercal},\fu_{\perp};\omega,\alpha,\varepsilon)=0
\]
with respect to the tangential unknown \(\fu_{\intercal}\). This step relies on
a fixed-point argument combined with suitable action--angle coordinates allowing
the inversion of the restricted linearized operator near equilibrium. Particular
attention is devoted to tracking the dependence of the tangential solution on
the normal component, which acts at this stage as a parameter. This analysis is
developed in Section~\ref{section.bifurcation}. At this stage, it is important to emphasize that we do not solve the bifurcation equation directly. Instead, we introduce and analyze the modified system \eqref{mod.aa.sys}, which involves an additional parameter $\mathtt{w}$. The role of this parameter is to enforce the zero-average condition, a constraint that is crucial in the periodic setting under consideration. This reformulation allows us to work within a functional framework compatible with the periodic structure of the problem. The parameter $\mathtt{w}$ is then uniquely determined by the solvability conditions
\begin{equation}
    \tw = \omega  + \im \lambda_{\tJ}^{(\infty)}(\alpha) \,,
\end{equation}
and its value ultimately provides, in the last step of the plan above, the oscillation frequency $\omega\in\R$ of the resulting solutions. 
The introduction of this extra parameter, sometimes referred as the \textit{hyphotetical conjugation argument}, has the advantage of momentarily turning the unknown oscillation frequency $\omega$ into a parameter, which is useful when imposing non-resonance conditions in the following step.

\paragraph{Solution to the range equation.}
Once the tangential component $\fu_\intercal = \fu_{\intercal}(\fu_\perp)$ has been determined, we substitute it into the
normal equation and obtain the reduced range equation
\begin{equation}\label{bifurcation.eqint}
	\wt\bG_{\perp}(\fu_{\perp};\omega,\alpha,\varepsilon)
	:=
	\bG_{\perp}(\fu_{\intercal}(\fu_{\perp}),\fu_{\perp};\omega,\alpha,\varepsilon)
	=0 \,.
\end{equation}
In Section \ref{section.normal}, we shall deal with   the resolution of the infinite-dimensional range equation \eqref{bifurcation.eqint}. 
The normal component is determined via a nonlinear Nash--Moser iteration. This requires constructing a right inverse of the linearized operator at each approximate solution. The main difficulty comes from the quasilinear nature of the problem and the small divisor phenomena generated by the interaction between the temporal frequency and the spatial spectrum of the linearized dynamics.  
The resolution of the range equation \eqref{bifurcation.eqint} consists of three main points:
\begin{itemize}
    \item [(i)] {\it Linearized operator at an approximate solution and reduction of the transport.}
    In Section \ref{section.linear.eps}, we start to study the  linearization of the equation around an approximate solution obtained along the iteration scheme, which is explicitly given by
    \begin{equation}\label{linear.pert.intro}
         \sL_{\perp} := \Pi_{\perp} \big( \omega\pa_{\vf} + \pa_{\vartheta}\circ\big(\tfrac12 - \alpha^4 +O(\alpha^4) + \fv_{\varepsilon}(\vf,\vartheta )  \big) + \pa_{\vartheta} \cR(0;0) + \cZ_{0,\alpha} \big)  \Pi_{\perp} +\sE_{0}  \,.
    \end{equation}
    The resulting operator is viewed as a quasi-linear perturbation of the diagonal operator obtained in the reducibility analysis around the traveling dipole, described through \eqref{Normal-lin-equili}. Its leading part is governed by a transport operator with time-dependent variable coefficients, while the remainder term $\sE_{0}$ is smoothing in space, see Lemma \ref{expand.sL.pert}. In Section \ref{sect.redu.eps} we will reduce the transport part through reversible and symplectic changes of variables that flatten the highest-order coefficient. This conjugation transforms the full operator into 
    \begin{equation}\label{linear.pertredu.intro}
        \sL_{\perp,{\rm red}} = \omega\pa_{\vf} + \sD_{\perp,{\rm red}} + \sE_{\perp,{\rm red}}\,,
    \end{equation}
     where $\sD_{\perp,{\rm red}}$  is  diagonal Fourier multiplier and $\sE_{\perp,{\rm red}}$ is a $\varepsilon$-small, smoothing in space perturbation. In this procedure, which has been implemented in several other papers (for instance, see \cite{FMT25,HR21}), small divisors appear and we are required to impose non-resonance condition of the form
    \begin{equation}
        |\omega\ell + \tm_{1,\varepsilon} j| \geq 4\upsilon |j|^{-\tau} \,, \quad \forall \, (\ell,j) \in \Z^2\setminus\{0\} \,,
    \end{equation}
    for some $\upsilon\in (0,1)$, $\tau\gg1$ large enough and where $\tm_{1,\varepsilon} = \tfrac12 -\alpha^4 + O(\alpha^5) + O(\varepsilon)$ is the reduced transport constant. We remark that, since we are searching for time-periodic solutions, we can impose a lower bound on these small divisors that loses derivative only in space.
    For more details, we refer to Lemma~\ref{almost.straight.lemma} and Lemma~\ref{diffeo.conj.pert}.
    \item [(ii)] {\it Right inverse of the linearized operator.} In Section \ref{sect.almost.inv}, we then construct a right inverse of the operator $\sL_{\perp}$ in \eqref{linear.pert.intro}, satisfying tame estimates,  by means of a perturbative argument. 
    As a preliminary step, we first invert the reduced operator $\sL_{\perp,{\rm red}}$ in \eqref{linear.pertredu.intro}. The diagonal operator $\sD_{\perp,{\rm red}}$ is fully invertible,
    provided that its eigenvalues are controlled by first Melnikov non-resonance conditions
    \begin{equation}
        |\omega\ell - \im \mu_{j,\varepsilon}(\omega,\alpha)| \geq 2\upsilon |j|^{-\tau} \,, \quad \forall\, \ell \in \Z\,, \ \ j \in \Z_\perp \,,
    \end{equation}
    for some $\upsilon\in (0,1)$, $\tau\gg1$ large enough and where each $\mu_{j,\varepsilon}(\omega,\alpha)\in \im \R$ is an eigenvalue of $\sD_{\perp,{\rm red}}$. At this level, too, the small divisor loses derivatives only in space: combined with the smoothing of the remainder $\sE_{\perp,{\rm red}}$, it allows us to find a right inverse of the full operator $\sL_{\perp,{\rm red}}$ by a Neumann series argument. Finally, being the map that conjugates $\sL_{\perp,{\rm red}}$ back to $\sL_{\perp}$ $\varepsilon$-close to the identity and satisfying tame estimates, we are able to find a right inverse of the operator $\sL_{\perp}$ as well.
    For more details, we refer to Proposition \ref{prop-perp.full} and Proposition \ref{prop.full.inv}.
    \item [(iii)] {\it Nash--Moser iteration.}
    This discussion will be developed in detail in Sections~\ref{section.NASH} and \ref{Concusion}. It concerns the implementation of the Nash--Moser iteration scheme and the completion of the proof of the range equation. The purpose of this iterative procedure is to solve the nonlinear equation despite the loss of derivatives produced by the quasilinear structure of the problem and the presence of small divisors in the inversion of the linearized operators. The scheme is implemented in a rather classical way, following the strategy developed in several previous works on quasilinear Hamiltonian PDEs and KAM theory, for instance see \cite{BBHM,HHM24}.  

Starting from a sufficiently accurate approximate solution, we construct a sequence of corrections by solving linearized equations at each step of the iteration. The key ingredient is the right inverse obtained in the previous sections through the reducibility and normal form analysis. Thanks to the tame estimates satisfied by this inverse, the corrections can be controlled in high Sobolev norms while preserving the rapid convergence at lower regularity levels.  
The reversibility and symmetry properties of the problem are preserved throughout the iteration.  
Finally, combining the Nash--Moser convergence scheme with the non-resonance conditions yields the existence of solutions for parameters belonging to a Cantor-like set of asymptotically full measure. This completes the resolution of the range equation and therefore the construction of non-rigid time-periodic vortex patches near the translating dipoles.
\end{itemize}

\paragraph{Measure of the final Cantor set of parameters.}
The last point to check is that the measure of the final Cantor set of parameters arising from the Nash--Moser scheme, see Theorem \ref{meas.est.thm} is actually of full asymptotic measure. This will be done in
Section \ref{section-measures}. 

The construction proceeds by decomposing the resonant set into elementary bad regions associated with the small divisors appearing in the homological equations. Using the asymptotic expansions of the eigenvalues together with suitable transversality properties, we obtain quantitative estimates on the size of each resonant region, which in the end sum up to a convergent and small quantity.  The smoothing structure of the perturbative terms and the separation properties of the spectrum play a crucial role in this analysis.  We finally remark that we have to be careful when dealing with  the mode $|j|=2$, depending on having $\pm 2 \in \Z_{\perp}$ in the normal mode or being $\tJ=2$ the tangential mode, due to its distinguished asymptotics in \eqref{mode.2.asympt.intro}.

\section{Contour dynamics reformulation for  dipoles}\label{model.derivation}

In this section, we derive the contour dynamics equation governing the evolution of the vortex patch boundaries in the general case. Starting from the Euler equations written in vorticity form, we first describe the dynamics of interacting vortex patches and their associated point-vortex approximation: we present this description in both cases of translating and co-rotating dipoles, without assuming symmetries on the two patches. We then introduce a moving frame adapted to the translation of the dipole and exploit the symmetry properties of the configuration in order to reduce the system to a single boundary equation.

The derivation proceeds through several steps. We first reformulate the boundary dynamics in terms of contour equations associated with the Biot--Savart law. We then introduce a suitable parametrization of the patch boundary and compute the induced evolution equation for the deformation variable. Particular attention is devoted to the role of the translation speed modulation and to the contribution of the interaction kernel generated by the distant companion vortex.

We consider the dynamics of two simply connected planar vortex patches 
$D_1(t), D_2(t) \subset \R^2$ with constant vorticity evolving in the fluid domain 
$\R^2 \simeq \C$. Denoting by $\omega(t,x)$ the scalar vorticity of the flow, we assume that
\begin{equation}
	\omega(t,x) = {\mathds{1}}_{D_1(t)}(x) + \epsilon\, {\mathds{1}}_{D_2(t)}(x) \,, \quad \epsilon= \pm 1\,.
\end{equation}
The evolution of the vortex patch configuration is globally well-posed within the Yudovich class for the two-dimensional incompressible Euler equations. The two patches are assumed to have the form
\begin{equation}\label{patches.form}
	D_{k}(t) = z_k(t) + O_{k}(t) \,, \quad k=1,2\,,
\end{equation}
where $O_k(t)$ denotes the shape of the $k$-th patch in a frame centered at its center of mass, and where the centers of mass $(z_1(t),z_2(t)) \in \C^2$ evolve according to the following {\it modified} point vortex dynamics
\begin{equation}\label{point.dyn.app}
	\begin{aligned}
		 \im \overline{\dot{z}_1} (t) & = \, \epsilon\,  \Big(\frac{1}{Z(t)} + \tU(t) +\epsilon \,\t\Omega'(t) \overline{z_1}(t) \Big) \,, \\
			\im   \overline{\dot{z}_2} (t) & =   - \Big(\frac{1}{Z(t)} + \tU(t) -  \t\Omega'(t) \overline{z_2}(t)\Big)\,,
	\end{aligned} \quad \quad  Z(t) :=\frac12 \big( z_{1}(t) - z_{2}(t) \big)\,.
\end{equation}
The functions $\tU:\R\to\C$ and $\t\Omega:\R\to\R$ are unknown and will serve as modulation parameters. The time evolution of the two boundaries $\partial D_k(t)$ is governed by the system of classical {\it contour dynamics} equations, as for instance in \cite{BertozziMajda,chemin},
\begin{equation}\label{contour.dyn.app}
	\pa_{t} w_{k}(t,\theta) \cdot \bn(\gamma_{k}(t,\theta))  + \pa_{\theta} \big( \psi(t,w_{k}(t,\theta)) \big) =0 \,, \quad t\in \R\,, \ \theta\in \T\,, \ k=1,2 \,,
\end{equation}
where $w_k(t,\,\cdot\,):\T\to \partial D_k(t)$ denotes a parametrization of the boundary of $D_k(t)$ of the form
\begin{equation}\label{contour.form}
	w_k(t,\theta) := z_{k}(t) + \gamma_{k}(t,\theta) \,,
\end{equation}
 $\gamma_{k}(t,\,\cdot\,):\T\to\pa O_k(t)$ is a parametrization of the boundary of $O_{k}(t)$, and $\bn(\gamma_{k})$ is the  outward normal vector to $\pa O_{k}(t)$ at the point $\gamma_{k}(t,\theta)$, which reads as
\begin{equation}\label{normal.gamma.app}
	\bn(\gamma_{k}(t,\theta)) = -\im\, \pa_{\theta}\gamma_{k}(t,\theta) \,.
\end{equation}
Here, the real scalar product between two complex numbers is given by
\begin{equation}\label{real.scalar.prod.app}
	z_a \cdot z_b := {\rm Re} (  \overline{z_a}\,z_b) = {\rm Im} (\im \,\overline{z_a}\, z_b) \,.
\end{equation}
 The stream function $\psi(t,\,\cdot\,)$ associated with the vorticity distribution is given by
 \begin{equation}\label{stream.gen.app}
 	\begin{aligned}
 		\psi(t,w) & := \frac{1}{2\pi} \int_{O_{1}(t)} \log |w - z_{1}(t) -\xi_{1} | \wrta(\xi_{1}) \\
 		& \   + \frac{\epsilon}{2\pi} \int_{O_{2}(t)} \log |w - z_{2}(t) - \xi_{2} | \wrta(\xi_{2}) \,,
 	\end{aligned}
 \end{equation}
 which nonlinearly couples the evolutions of the two boundaries.

 \subsection{Modulated point vortex system}
We now explore the modulated version of the classical point vortex dynamics in \eqref{point.dyn.app} associated with two interacting vortices. This finite-dimensional model serves as an effective approximation for the motion of the centers of the two vortex patches considered throughout the paper. The modulation is designed to capture the influence of the boundary deformations on the translational/rotational dynamics of the pair.

The classical dynamics of two point vortices with circulations $\Gamma_1$ and $\Gamma_2$, with centers $z_1(t), z_2(t)\in\C$, is described by the standard system
 \begin{equation}\label{hiroshi1}
 	 \overline{\dot{z}_1}(t) = \frac{\Gamma_{2}}{2\pi \im} \frac{1}{z_1(t)-z_2(t)} \,, \quad \overline{\dot{z}_2}(t) = \frac{\Gamma_{1}}{2\pi \im} \frac{1}{z_2(t)-z_1(t)} \,,\quad\hbox{where}\quad \dot{z}:=\frac{\di z}{\di t}.
 \end{equation}
We modify the dynamics of the two point vortices by means of smooth modulation functions in the following way
 \begin{equation}\label{hiroshi1.mod}
 \begin{aligned}
 		 \overline{\dot{z}_1}(t) & = \frac{\Gamma_{2}}{2\pi \im}\Big( \frac{1}{z_1(t)-z_2(t)}  +\frac{\tU(t)}{2}  + \frac{2\pi\,\dot{\t{\Omega}}(t)}{\Gamma_{2}}  \overline{z_1}(t) \Big)   \,, \\
 		 \overline{\dot{z}_2}(t) & = \frac{\Gamma_{1}}{2\pi \im} \Big(  \frac{1}{z_2(t)-z_1(t)} -\frac{\tU(t)}{2} + \frac{2\pi\,\dot{\t\Omega}(t)}{\Gamma_{1}}  \overline{z_2}(t)\Big)\,,
 \end{aligned}
 \end{equation}
 where $\tU:\R\to\C$ is a translational modulation and   $\t\Omega:\R\to\R$ is a rotational modulation.
 The  first parameter $\tU$ accounts for corrections to the propagation speed generated by the deformation of the patches, while the second one $\t\Omega$ compensates for the rotational drift induced by the moving frame. They are uniquely determined through suitable orthogonality conditions removing the neutral directions generated by the symmetries of the Euler equations during the interaction between the point vortex dynamics and the boundary deformation of the patches.
 In the case of the present paper, the modulation $\tU$ will be crucially used, for instance in Proposition \ref{lemma.constrain} and Proposition \ref{prop.U1}, to eliminate the first sine mode from the nonlinear equations.
 
   Throughout this paper, we use the normalization 
 \begin{equation}\label{circulations}
 	\Gamma_{1}:= 4\pi\,, \quad \Gamma_{2} := 4\pi \epsilon\,.
 \end{equation}
 We define the center of mass $C(t)$ and the displacement $Z(t)$ from the center of mass by
 \begin{equation}\label{hiroshi2}
 	C(t) := \tfrac12\big( z_1(t) + z_2(t)  \big) \,, \quad Z(t) := \tfrac12\big( z_1(t) - z_2(t)  \big)\,.
 \end{equation}
 Then, by \eqref{circulations} and \eqref{hiroshi2}, the equations \eqref{hiroshi1.mod} become
 \begin{equation}\label{hiroshi3}
 		\im \overline{\dot{z}_1} (t) =  \epsilon\Big(\frac{1}{Z(t)} + \tU(t) +\epsilon\, \dot{\t\Omega}(t) \overline{{z}_1}(t) \Big) \,, \quad  	\im   \overline{z_2} (t) =   -\Big(\frac{1}{Z(t)} + \tU(t) - \dot{\t\Omega}(t) \overline{z_2}(t) \Big)\,, 
 \end{equation}
 which coincide with \eqref{point.dyn.app}. Moreover, the evolution of $(C(t),Z(t))$ is given by
 \begin{equation}\label{hiroshi4}
 \begin{aligned}
 		 \overline{\dot{C}(t)} & = \frac{\epsilon-1}{2\,\im} \Big(\frac{1}{Z(t)} + \tU(t) \Big) + \frac{1}{\im} \dot{\t\Omega}(t) \overline{C(t)}  \,, \\
 		\overline{\dot{Z}(t)} & = \frac{\epsilon+1}{2\,\im} \Big(\frac{1}{Z(t)} + \tU(t) \Big) + \frac{1}{\im} \dot{\t\Omega}(t) \overline{Z(t)} \,.
 \end{aligned}
 \end{equation}
 We obtain that, when $\epsilon= +1$, the center of mass of the two points is constant up to a time-dependent phase shift $e^{\im\,\t\Omega(t)}$ , whereas, when $\epsilon = -1$, it is the displacement $Z(t)$ that becomes a constant of the motion, still up to a time-dependent phase shift $e^{\im\,\t\Omega(t)}$ . 
 \begin{rem}\label{remark.distance}
 	 If the displacement evolves as $Z(t)=e^{\im\,\t\Omega(t)} Z(0)$ (that is, when $\epsilon=-1$), then the distance $|z_1(t)- z_2(t)|$ is a constant of the motion. On the other hand, by \eqref{hiroshi4}, we have that
 	\begin{equation}
 		\frac14 \frac{\di}{\di t} |z_1(t) - z_2(t) |^2 = \dot{Z}(t) \, \overline{Z(t)} + Z(t)\,  \overline{\dot{Z}(t)} = (\epsilon+1) {\rm Im} (\tU(t) Z(t)) \,.
 	\end{equation}
 	We obtain that the distance is a constant of motion in the  case $\epsilon= 1$ only if $\tU\equiv 0$.
 \end{rem}
 \noindent Integrating the equations in \eqref{hiroshi4}, we are ready to describe the orbits of \eqref{hiroshi3}:
\\[1mm]
 \noindent $\bullet$ When $\epsilon =1$ and $\tU\equiv 0$, writing $Z(0)=e^{\im \vf_{0}}|Z(0)|$ for some $\vf_{0}\in\T,$   the trajectories are
 	\begin{equation}
 		\begin{aligned}
 				z_1(t) &= C(t) + Z(t) = e^{\im\,\t\Omega(t)}\Big(C(0) + e^{\im\big(\frac{t}{|Z(0)|^2} +\vf_{0}\big)} |Z(0)| \Big)\,, \\
 				 z_2(t) &= C(t) - Z(t) = e^{\im\,\t\Omega(t)}\Big(C(0) - e^{\im\big(\frac{t}{|Z(0)|^2} +\vf_{0}\big)} |Z(0)| \Big)\,. 
 		\end{aligned}
\end{equation}
The points $z_1(t)$ and $z_2(t)$ remain diametrically opposite with respect to the center of mass $e^{\im\,\t\Omega(t)}C(0)$ and rotate around it. In particular, when $C(0)=0$, one has
$z_2(t) = - z_1(t)$;
 \\[1mm]
 \noindent $\bullet$
 When $\epsilon = -1$ and $\t\Omega\equiv0,$  by writing $Z(0)=e^{\im\vf_{0}}|Z(0)|$ and choosing $\tU(t)=e^{-\im\vf_{0}}|\tU(t)|$ for some $\vf_{0}\in\T$, the trajectories are given by, for $k=1,2$,
 	\begin{equation}\label{point.traj}
\begin{aligned}
		& z_k(t) = z_k(0) -\im\Big( \frac{t}{\overline{Z(0)}} + \int_{0}^t \overline{\tU(\tau)} \wrt\tau  \Big)=  z_k(0)- \im e^{\im\vf_{0}} \Big( \frac{t}{|Z(0)|}  +\int_{0}^t |\tU(\tau)|\wrt \tau  \Big) \,, \\
		& C(t)=C(0)   -  \im e^{\im\vf_{0}}\Big( \frac{t}{|Z(0)|}  +\int_{0}^t |\tU(\tau)| \wrt \tau  \Big) \,,
\end{aligned}
\end{equation}
 	that is, $z_1(t)$ and $z_2(t)$ translate  together  in the direction $-\im e^{\im \vf_{0}}$. In particular, when $\vf_{0}=0 \mod \pi$,  the translation is parallel to the vertical axis, whereas, when $\vf_{0}= \frac{\pi}{2} \mod \pi$, the translation is parallel to the horizontal axis.

\subsection{Symmetry reduction in the translating case}
We now exploit the symmetry properties of the translating dipole configuration in order to reduce the contour dynamics system to a single scalar equation. More precisely, we require that the two vortex patches are symmetric with respect to the imaginary axis passing through the center of mass. This symmetry is preserved by the dynamics and allows us to express the evolution of one patch entirely in terms of the other.
From now on we choose to work only with   the configuration  of the two patches with opposite vorticity, that is, we fix
\begin{equation}
	\epsilon = -1\,.
\end{equation}
We further restrict ourselves to the case $\t\Omega(t)\equiv 0$, so that the corresponding point vortex system reduces to \eqref{point.traj}. In this setting, the two points $(z_1(t),z_2(t))$ travel parallel to a fixed axis with constant separation, moving at a speed inversely proportional to the distance $|Z(t)|=|Z(0)|$ from the center of mass.
 Without loss of generality, we also fix $C(0)=0$ and $\vf_{0}=\pi$, so that the motion is parallel to the vertical direction $\big(\begin{smallmatrix}
	0 \\ 1
\end{smallmatrix}\big)\simeq \im$ with the center of mass moving on the vertical axis. That is, we have
\begin{equation}\label{point.sym.app}
	z_2(t) = - \overline{z_1(t)} \,, \quad Z(t) = Z(0) = |Z(0)| \in \R\,, \quad \tU(t)=|\tU(t)|\in\R\,.
\end{equation}
We shall reflect this symmetry at the level of the vortex patches by requiring that the two patches are symmetric with respect to the imaginary axis passing through the center of mass.
Recalling \eqref{patches.form}, \eqref{contour.form},  this is equivalent to requiring that
  \begin{equation}\label{symmetry.cond.app}
  	D_2(t) = - \overline{D_1(t)} \,, \quad O_2(t) = -\overline{O_1(t)}\,, \quad \gamma_{2}(t,\theta) = - \overline{\gamma_{1}(t,\theta)} \,.
   \end{equation}
   We plan to further explore the contour dynamics equation \eqref{contour.dyn.app} under these symmetry conditions. First, by \eqref{contour.form}, \eqref{normal.gamma.app}, \eqref{real.scalar.prod.app}, \eqref{point.sym.app} and \eqref{symmetry.cond.app}, we have
   \begin{equation}\label{dot.sym}
   	\begin{aligned}
   		\pa_{t} w_2 \cdot \bn(\gamma_{2}) &= {\rm Im} \big( \im \pa_{t} \overline{w_2} \,(-\im \pa_{\theta} \gamma_{2}) \big) \\
   		&= {\rm Im} \big( \pa_{t} w_1 \,\pa_{\theta}\overline{\gamma_{1}} \big) = - {\rm Im} \big( \pa_{t} \overline{w_1} \, \pa_{\theta} \gamma_{1} \big) = -  \pa_{t} w_1 \cdot \bn(\gamma_{1}) \,.
    	\end{aligned}
   \end{equation}
   We now look at the stream function \eqref{stream.gen.app}. By \eqref{contour.form}, \eqref{point.sym.app} and \eqref{symmetry.cond.app}, we note that
      \begin{equation}\label{stream.sym}
   	\begin{aligned}
   		\psi(t,w_2) & = \frac{1}{2\pi} \int_{O_{1}(t)} \log | z_2 - z_{1} + \gamma_{2} -\xi_{1} | \wrta(\xi_{1})   - \frac{1}{2\pi} \int_{O_{2}(t)} \log | \gamma_{2}- \xi_{2} | \wrta(\xi_{2}) \\
   		& = \frac{1}{2\pi} \int_{O_{1}(t)} \log | -2Z(0)- \overline{\gamma_{1}} -\xi_{1} | \wrta(\xi_{1})  - \frac{1}{2\pi} \int_{O_{2}(t)} \log | -\overline{\gamma_{1}}- \xi_{2} | \wrta(\xi_{2}) \\
   		& = \frac{1}{2\pi} \int_{O_{2}(t)} \log | 2Z(0)+ \overline{\gamma_{1}} -\overline{\xi_{2}} | \wrta(\xi_{2})   - \frac{1}{2\pi} \int_{O_{1}(t)} \log | \overline{\gamma_{1}}- \overline{\xi_{1}} | \wrta(\xi_{1}) \\
   		& = - \psi(t,w_1)\,,
   	\end{aligned}
   \end{equation}
  and
   \begin{equation}\label{stream.only1}
   	\begin{aligned}
   			\psi(t,w_1) & = \frac{1}{2\pi} \int_{O_{1}(t)} \log |  \gamma_{1} -\xi_{1} | \wrta(\xi_{1})   - \frac{1}{2\pi} \int_{O_{2}(t)} \log | 2Z(0)+ \gamma_{1}- \xi_{2} | \wrta(\xi_{2}) \\
   			& = \frac{1}{2\pi} \int_{O_{1}(t)} \log |  \gamma_{1} -\xi_{1} | \wrta(\xi_{1})   - \frac{1}{2\pi} \int_{O_{1}(t)} \log | 2Z(0)+ \gamma_{1} +\overline{\xi_{1}} | \wrta(\xi_{1}) \,.
   	\end{aligned}
   \end{equation}
Combining \eqref{dot.sym} and \eqref{stream.sym}, we conclude that if $w_1(t,\theta)$ solves \eqref{contour.dyn.app}, with the constraint
$$
w_2(t,\theta)=-\overline{w_1(t,\theta)} \,,
$$
then, $w_2$ automatically solves \eqref{contour.dyn.app} as well. Therefore, the two contour dynamics equations reduce to a single scalar equation. Consequently, it is sufficient to study \eqref{contour.dyn.app} with \eqref{stream.only1} for $k=1$, where the whole dynamics is encoded through the evolution of the single patch
$$
D_1(t)=z_1(t)+O_1(t) \,.
$$

\subsection{Polar coordinates parametrization and  dynamics}
In the previous section, we performed the symmetry reduction and we are led to study the scalar equation
\begin{equation}\label{contour.dyn.only1}
	\pa_{t} w_{1}(t,\theta) \cdot \bn(\gamma_{1}(t,\theta))  + \pa_{\theta} \big( \psi(t,w_{1}(t,\theta)) \big) =0 \,, \quad t\in \R\,, \ \theta\in \T\,,
\end{equation}
where $w_{1}(t,\,\cdot\,):\T\to \pa D_1(t)$ is of the form \eqref{contour.form}.
We now introduce the parametrization of the vortex boundary that will be used throughout the paper. The boundary of the reference patch is described as a perturbation of the unit circle through a scalar deformation radius depending on time and the angular variable. Substituting this representation into the contour dynamics formulation, under the symmetry reduction,  allows us to derive the evolution equation satisfied by the deformation function. This reformulation transforms the geometric motion of the vortex boundary into a nonlinear nonlocal equation more suitable for the subsequent analytical study.
\\
We parametrize the boundary $\partial O_1(t)$ using polar coordinates as follows
\begin{equation}\label{param.gamma}
	\gamma_{1}(t,\theta) = R(t,\theta) e^{\im\theta} := \sqrt{1+ 2 r(t,\theta)}\, e^{\im\theta} \,, \qquad t\in \R\,, \quad \theta \in \T\,,
\end{equation}
and therefore
\begin{equation*}
	O_1(t)=\big\{ \rho e^{\im\eta} \,: \, \eta \in \T\,, \ 0\leq \rho \leq \sqrt{1+2r(t,\eta)} \big\}\,.
\end{equation*}
Here,  the time-dependent deformation radius $ r$, assumed to be small enough, is one of the unknowns of the problem, together with the modulation parameter $\tU(t)$. In this way, we seek vortex patches that can be viewed as deformations of the unit disk $\D$.
We note that, by \eqref{param.gamma}, 
\begin{equation}
	\pa_{t} \gamma_{1}(t,\theta) = \pa_{t}R(t,\theta) e^{\im\theta} \,, \quad \pa_{\theta}\gamma_{1}(t,\theta) = \big( \pa_{\theta}R(t,\theta) +\im\,R(t,\theta) \big)e^{\im\theta}\,,
\end{equation}
and using  \eqref{normal.gamma.app}, \eqref{real.scalar.prod.app} and \eqref{param.gamma}, we infer
	\begin{align}
		\pa_{t} w_1 \cdot \bn(\gamma_{1}) & = {\rm Im} \big( \overline{\dot{z}_1} \, \pa_{\theta} \gamma_{1} + \pa_{t}\overline{\gamma_{1}} \, \pa_{\theta}\gamma_{1} \big) \\
		& = {\rm Im}\big(  \overline{\dot{z}_1} \,\pa_{\theta}\gamma_{1} \big) + {\rm Im}\big( \pa_{t}R \,\pa_{\theta} R + \im R\, \pa_{t} R \big) \\
		& = {\rm Im}\big(  \overline{\dot{z}_1}\, \pa_{\theta}\gamma_{1} \big)  + R \,\pa_{t} R  = {\rm Im}\big(  \overline{\dot{z}_1}\, \pa_{\theta}\gamma_{1} \big)  + \pa_{t}r \,.\label{time.part}
	\end{align}
To compute the terms in \eqref{stream.only1}, we use the following expansion.
\begin{lem}\label{loga.expa}
	For any $|z|<1$, we have that
	\begin{equation}
		\log |1+z| = {\rm Re}(z) - \tfrac{1}{2} {\rm Re}(z^2) + \sum_{n\geq 3} \tfrac{(-1)^{n+1}}{n} {\rm Re}(z^n)\,.
	\end{equation}
\end{lem}
\begin{proof}
	Note that $\log|1+z| = \frac12\big( \log(1+z) + \log(1+\bar z) \big)$. Then, it suffices to apply the standard power series expansion for $\log(1+z)$, with $|z|<1$.
\end{proof}
We recall the parameter $\alpha$ introduced in \eqref{def.alpha.intro}. For our purpose, we shall work in the regime where $\alpha$ is sufficiently small, corresponding to configurations in which the two concentrated vortex patches are sufficiently far apart.
Then, by using  Lemma \ref{loga.expa}, we find that  for any $\xi_{1}\in O_1(t)$,
	\begin{align}
	\log |2Z(0) &+ \gamma_{1} +\overline{\xi_{1}} |  =	\log |\alpha^{-1} + \gamma_{1} +\overline{\xi_{1}} |\\
	 & = \log|\alpha^{-1}| + \log\big| 1  + \alpha(\gamma_{1}+\overline{\xi_1}) \big|  \\
	 & = \log|\alpha^{-1}| + \alpha {\rm Re}( \gamma_{1}+\overline{\xi_{1}} )  + \sum_{n\geq 2} (-1)^{n+1}\tfrac{\alpha^n}{n}{\rm Re}\big( (\gamma_{1}+\overline{\xi_{1}})^n \big) \,.
	\label{stream.part}
\end{align}
We now impose, following the mass conservation, the constraint
\begin{equation*}
	\int_{O_{1}(t)} \wrta(\xi_{1}) = 4\pi  \quad \forall\, t\in\R\,.
\end{equation*}
By inserting \eqref{time.part} and \eqref{stream.part} into \eqref{contour.dyn.app}, together with   \eqref{param.gamma} and \eqref{point.dyn.app}, 
we conclude that
\begin{equation}\label{quasi.final.eq}
	\pa_{t} r(t,\theta)  + {\rm Im} \Big(\big(\,  \overline{\dot{z}_1} -\tfrac{\im}{Z(0)} \big) \pa_{\theta}\gamma_{1}(t,\theta) \Big) + \pa_{\theta} \big( F_1[\gamma_{1};\alpha](t,\theta) \big) = 0 \,,
\end{equation}
where
\begin{equation}\label{F1}
	\begin{aligned}
		F_1[\gamma_{1};\alpha](t,\theta)  & :=  \frac{1}{2\pi} \int_{O_{1}(t)} \log |  \gamma_{1}(t,\theta) -\xi_{1} | \wrta(\xi_{1})  - \frac{1}{2\pi} \int_{O_{1}(t)} K_{\alpha}(\gamma_{1}(t,\theta),\xi_{1}) \wrta(\xi_{1})\,,
	\end{aligned}
\end{equation}
and
\begin{equation}\label{Kalpha.kernel}
\begin{aligned}
		K_\alpha(\gamma_{1},\xi_{1}) & := \log \big| 1+\alpha( \gamma_{1}+ \overline{\xi_{1}} )\big| - \alpha {\rm Re}(\gamma_{1}+\overline{\xi_{1}})\\
		& = \sum_{n\geq 2}(-1)^{n+1} \tfrac{\alpha^{n}}{n}{\rm Re}\big(( \gamma_{1}(t,\theta)+\overline{\xi_{1}} )^n\big) \,.
\end{aligned}
\end{equation}
In the present setting, where $\epsilon=-1$ and $\t\Omega\equiv0$, and under the symmetry condition \eqref{point.sym.app}, the system \eqref{point.dyn.app} takes the form
\begin{equation}\label{point.shift.dyn}
	 \overline{\dot{z}_1}(t) = \im \Big( \tfrac{1}{Z(0)} + \tU(t)  \Big) \,, \quad  \overline{\dot{z}_2}(t) = \im \Big( \tfrac{1}{Z(0)} + \tU(t)  \Big) \,.
\end{equation}
Then, by \eqref{point.shift.dyn} and \eqref{param.gamma}, we find that
\begin{equation}\label{ImReU}
	\begin{aligned}
		 {\rm Im} \Big( \big( \,\overline{\dot{z}_1} -\tfrac{\im}{Z(0)} \big) \pa_{\theta}\gamma_{1}(t,\theta) \Big) & =  {\rm Im} \big(\im\, \tU(t)\, \pa_{\theta}\gamma_{1}(t,\theta) \big) =  {\rm Re} \big( \tU(t)\, \pa_{\theta}\gamma_{1}(t,\theta) \big) \\
		 &  =\pa_{\theta}\Big( {\rm Re} \big( \tU(t)\, \gamma_{1}(t,\theta) \big)\Big) \\
		 &= \pa_{\theta}\big( \tU(t) \, \sqrt{1+2r(t,\theta)} \cos(\theta) \big)\,.
	\end{aligned}
\end{equation}
Plugging \eqref{ImReU} into \eqref{quasi.final.eq}, we finally conclude that the contour dynamics for the patch $D_1(t)$ is determined by the equation
\begin{equation}\label{the.equation}
	\begin{aligned}
 \pa_{t} r(t,\theta) +   \pa_{\theta}\Big(  \tU(t)  \sqrt{1+2 r(t,\theta)}\cos(\theta ) + F[ r ;\alpha ](t,\theta) \Big) = 0 \,,
	\end{aligned}
\end{equation}
where
\begin{equation}\label{F.app}
	\begin{aligned}
		F[r;\alpha](t,\theta) & := F_1[\gamma_{1};\alpha](t,\theta) \\
& =   \frac{1}{2\pi} \int_{O_{1}(t)} \log \big|  \sqrt{1+2r(t,\theta)}e^{\im\theta} -\xi \big| \wrta(\xi) 
\\
		& \quad - \frac{1}{2\pi} \int_{O_{1}(t)} K_\alpha\big(\sqrt{1+2 r(t,\theta)}e^{\im\theta},\xi\big) \wrta(\xi)\\       
        &:= F_{\textnormal{self}}[r;\alpha](t,\theta)-F_{\textnormal{int}}[r;\alpha](t,\theta)\,,
	\end{aligned}
\end{equation}
with  the kernel $K_\alpha(\gamma_{1},\xi_{1})$ as in \eqref{Kalpha.kernel}. 
The goal in the upcoming sections is to construct periodic solutions in time with an a priori unknown  frequency $\omega\in\R$. More precisely, we look for solutions of the form
$$ r(t,\theta)=\check{r}(\vf,\theta)|_{\vf=\omega t}\quad\hbox{and}\quad \tU(t)=\check{\tU}(\vf)|_{\vf=\omega t} \,,
$$
where  $\check{r}:(\vf,\theta)\in\T^2\to \check{r}(\vf,\theta)\in\R$ and $\tU: \vf\in\T\to \tU(\vf)\in\R$ solve the equation, see also \eqref{Func.G-ini},
\begin{equation}\label{equation.to.solve}
		\omega\,\pa_{\vf} \check{r}(\vf,\theta) +  \pa_{\theta}  \Big( \check{\tU}(\vf) \sqrt{1+2\check{r}(\vf,\theta)} \cos(\theta)  + F[\check{r};\alpha](\vf,\theta) \Big) =0\,.
\end{equation}
For the sake of notational simplicity, in the remainder of the paper we shall work exclusively with equation
 \eqref{equation.to.solve} relabeling $\check{r}$ and $\check{\tU}$ back to $r$ and $\tU$, but keeping them as functions of the angle $\vf\in\T$.

\smallskip
 
 To conclude this section, we derive some boundary integral representations for $\pa_{\theta}F[r;\alpha]$, with $F[r;\alpha]$ as  in \eqref{F.app}. For simplicity of notation, we suppress the time dependence throughout this discussion. 
 We compute separately the contributions in \eqref{F.app}.
 \\[1mm]
 \noindent $\bullet$ We start with $\pa_{\theta}F_{\rm self}[r;\alpha]$, with $F_{\rm self}[r,;\alpha]$ as in \eqref{F.app}.
 From straightforward computations, recalling the parametrization in \eqref{param.gamma}, we obtain
\begin{align}
		\partial_\theta F_{\textnormal{self}}[r;\alpha](\theta)&=
\frac{1}{2\pi} \partial_\theta\int_{O_{1}} \log |  \gamma_1(\theta) -\xi | \wrta(\xi)\\
&=-\frac{1}{\pi}{\rm Re} \int_{O_{1}} \partial_{\overline{\xi}}\log |  \gamma_1(\theta) -\xi |\overline{\partial_\theta\gamma(\theta)} \wrta(\xi)\,,
\end{align}
where we have used that
\begin{equation}\label{partiaz-log}
    \pa_{\bar z} \log |z-\xi| = - \pa_{\bar \xi} \log |z-\xi| \,.
\end{equation}
We now employ Stokes Theorem in its complex form, which can be stated  as
\begin{equation}\label{stokes.complex.app}
	2 \im \, \int_{D} \pa_{\bar \xi} f(\xi) \wrta(\xi) = \int_{\pa D} f(\xi) \wrt\xi\,,
\end{equation}
for any $\cC^1$ function $f:\overline D \to \C$  on a simply connected bounded domain $D$ with boundary $\pa D$,  oriented counterclockwise. It follows that
\begin{align}\label{Fself-1}
	\nonumber 	\partial_\theta F_{\textnormal{self}}[r;\alpha](\theta)
&=-{\rm Re} \bigg(\frac{1}{2\im \pi}\int_{\partial O_{1}} \log |  \gamma_1(\theta) -\xi |\overline{\partial_\theta\gamma_1(\theta)} \wrt \xi\bigg)\\
&=-{\rm Im} \bigg(\frac{1}{2\pi}\int_{\partial O_{1}} \log |  \gamma_1(\theta) -\xi |\overline{\partial_\theta\gamma_1(\theta)} \wrt \xi\bigg)\,.
\end{align}
Using polar coordinates $\xi=\gamma_{1}(\eta)$ we infer
\begin{align*}
{\rm Im}\big(\overline{\partial_\theta\gamma_1(\theta)} \wrt \xi\big)&=\partial_\theta\partial_\eta {\rm Im}\big(\overline{\gamma_1(\theta)}\,\gamma_1(\eta) \big) \wrt \eta\\
&=\partial_\theta\partial_\eta {\rm Im}\big(R(\theta) R(\eta) e^{\im (\eta-\theta)} \big) \wrt \eta\\
&=\partial_\theta\partial_\eta \big(R(\theta) R(\eta) \sin( \eta-\theta) \big) \wrt \eta\,.
\end{align*}
Plugging this identity into \eqref{Fself-1} yields
\begin{align}\label{Fself-01}
	 	\partial_\theta F_{\textnormal{self}}[r;\alpha](\theta)
&= \fint_{\mathbb{T}} \log |  \gamma_1(\theta) -\gamma_1(\eta) | \partial_\theta\partial_\eta \Big(R(\theta) R(\eta) {\sin( \theta-\eta)} \Big) \wrt \eta\,,
\end{align}
where we have used the notation \eqref{average-f}.
\\[1mm]
\noindent $\bullet$
Let us now turn to the second contribution, namely $\partial_\theta F_{\mathrm{int}}[r;\alpha]$, with $F_{\rm int}[r;\alpha]$ as in \eqref{F.app}.
Using the general  identity
\begin{align}\label{diff-gener}
 \partial_\theta\big[ H(\gamma_1(\theta),\xi)\big]&=(\nabla_z H)(\gamma_1(\theta),\xi)\cdot\partial_\theta\gamma_1 \\
&=2{\rm Re}\big( (\partial_z H)(\gamma_1(\theta),\xi)\,\partial_\theta\gamma_1\big) 
\end{align}
together with the specific one  
\begin{align}\label{specific-K}
\pa_{z} K_\alpha(z,\xi) = \pa_{\bar \xi} K_\alpha(z,\xi)\,,
\end{align}
we find
\begin{align*}
\partial_\theta F_{\textnormal{int}}[r;\alpha](\theta)&=\frac{1}{\pi} {\rm Re}\int_{O_{1}} \pa_{\bar \xi} K_\alpha(\gamma_1(\theta),\xi) \partial_\theta\gamma_1(\theta) \wrta(\xi) \,.
    \end{align*}
Applying once again the identity \eqref{stokes.complex.app} yields
\begin{align*}
\partial_\theta F_{\textnormal{int}}[r;\alpha](\theta)&= {\rm Re}\bigg(\frac{1}{2\im \pi} \int_{\partial O_{1}}  K_\alpha(\gamma_1(\theta),\xi) \partial_\theta\gamma_1(\theta) \wrt \xi\bigg)\\
&= {\rm Im}\bigg(\frac{1}{2\pi} \int_{\partial O_{1}}  K_\alpha(\gamma_1(\theta),\xi) \partial_\theta\gamma_1(\theta) \wrt \xi\bigg)\,.
    \end{align*}
    Direct computations imply
    \begin{align*}
{\rm Im}({\partial_\theta\gamma_1(\theta)} \wrt \xi)&=\partial_\theta\partial_\eta {\rm Im}({\gamma_1(\theta)}\,\gamma_1(\eta) ) \wrt \eta\\
&=\partial_\theta\partial_\eta {\rm Im}\big(R(\theta) R(\eta) e^{\im (\eta+\theta)} \big) \wrt \eta\\
&=\partial_\theta\partial_\eta \big(R(\theta) R(\eta) \sin( \eta+\theta) \big)\wrt \eta
\end{align*}
Therefore
\begin{align}
\partial_\theta F_{\textnormal{int}}[r;\alpha](\theta)&= \fint_{\mathbb{T}}  K_\alpha(\gamma_1(\theta),\gamma_1(\eta)) \partial_\theta\partial_\eta \Big(R(\theta) R(\eta) \sin( \eta+\theta) \Big)\wrt \eta\,. \label{F.int-01}
    \end{align}
    
Finally, we rewrite the equation \eqref{equation.to.solve} in the form
    \begin{align}\label{Func.G}
		\cG(\tU,r;\alpha):= \omega \pa_{\vf} r + \pa_{\theta} \Big( \tU(\vf) \sqrt{1+2 r} \cos(\theta) + F[r;\alpha] \Big)=0 \,, 
         \end{align}
         where, using \eqref{F.app}, \eqref{Fself-01} and \eqref{F.int-01}, we have the representation
         \begin{align}\label{F-compact}
\nonumber \pa_\theta F[r;\alpha] (\vf,\theta)
		& =    \fint_{\T} \log | R(\vf,\theta) e^{\im\theta} - R(\vf,\eta) e^{\im\eta} | \pa_{\theta,\eta}^2 \big[ R(\vf,\theta) R(\vf,\eta) \sin(\theta-\eta) \big] \wrt \eta 
		\\
		& \quad -  \fint_{\T} K_\alpha ( R(\vf,\theta) e^{\im\theta}, R(\vf,\eta) e^{\im\eta} )|\pa_{\theta,\eta}^2 \big[ R(\vf,\theta) R(\vf,\eta) \sin(\eta+\theta) \big] \wrt \eta\,.
\end{align}

\section{Functional setting}\label{Section-Func-sett}
This section  is devoted to the introduction of  the functional framework used throughout the paper. We define the standard Sobolev and analytic Sobolev spaces on the torus that are adapted to our problem, together with subspaces encoding the reversibility and symmetry properties of the problem. 
We also describe the class of symplectic transformations employed to simplify the transport part of the linearized operator. Finally, we introduce a matrix representation of linear operators and suitable decay norms, which will be crucial for the normal form analysis used in the reducibility of the linearized operator at the traveling waves.  
\medskip
 
\paragraph{\bf Notations.} Throughout this paper we shall use the following notations:
\\[1mm]
\noindent $\blacktriangleright$
Let $\N:=\{1,2,\ldots\}$ and $\N_{0}:=\{0\}\cup \N$;
\\[1mm]
\noindent $\blacktriangleright$ For a $2\pi$-periodic and integrable function $f:\T\to\R$, we define its average by
\begin{align}\label{average-f}
\langle f\rangle_\theta:=\fint_{\T} f(\theta)\,\wrt \theta := \frac{1}{2\pi}\int_{\T} f(\theta)\,\wrt \theta\,;
\end{align}
\noindent $\blacktriangleright$ We denote the H\"ormander derivative
\begin{equation}
    {\rm D}:= - \im \pa_{\theta} \quad \textnormal{ and } \quad |{\rm D}|e^{\im j \theta} := |j|e^{\im j \theta} \,;
\end{equation}
\noindent $\blacktriangleright$ We denote  $a \lesssim_{s} b$, to  mean that $0\leqslant a \le C(s)\, b$ for some constant $C(s)>0$;
\\[1mm]
\noindent $\blacktriangleright$
Let   $\gamma\mapsto f(\gamma),g(\gamma), h(\gamma)$ be functions taking values in a Banach space $X$, with $h(\gamma)\neq 0$. We write  $$f(\gamma)=g(\gamma)+O(h(\gamma)) \quad \textnormal{as}\quad  \gamma\to 0
$$
to  mean that the quotient 
$
\frac{\|f(\gamma)-g(\gamma)\|_{X}}{\|h(\gamma)\|_{X}}
$
 is bounded for $\gamma$ small enough.

\subsection{Sobolev-analytic functional spaces}

We shall work with functional spaces defined through Fourier expansions. More precisely, for any periodic function $f:\T^2 \to \mathbb{R}$, we write
\begin{equation}
f(\varphi,\theta) = \sum_{(\ell,j)\in \Z^2} f_{\ell,j} \, e^{\im(\ell \varphi + j \theta)} \,,  
\qquad 
f_{\ell,j}:= \frac{1}{(2\pi)^{2}}\int_{\T^2} f(\phi,\eta)\, e^{-\im(\ell\phi + j\eta)} \wrt \phi \wrt \eta \,.
\end{equation}
Depending on the analytical framework under consideration, we shall use the following classes of spaces:
\\[1mm]
\noindent $\blacktriangleright$ {\it Sobolev-analytic spaces}: for $\rho\geq 0$ and $s\in\mathbb{R}$, we define the anisotropic Banach spaces
\begin{equation}\label{X.space}
		 X^{\rho,s}:= X^{\rho,s}(\T^2):=\bigg\{ f:\T^2 \to \R \, :\,  \| f\|_{X^{\rho,s}} :=  \| f \|_{\rho,s}  := \sum_{(\ell,j)\in \Z^2} e^{\rho|j|} \braket{\ell,j}^s |f_{\ell,j}| < \infty  \bigg\}\,,   
\end{equation}
where $\braket{\ell,j}:={\rm max}\{1,|\ell|,|j|\}$. These spaces will be used in  Section \ref{Traveling pairs} to construct the stationary solutions, with respect to  a co-moving frame, stated in Theorem \ref{Theo-1};
\\[1mm]
\noindent $\blacktriangleright$ {\it Sobolev spaces}: for $s\in \R$, we define the isotropic Banach space (see also \eqref{sob.space.intro})
\begin{equation}\label{sob.space}
	H^s:= H^s(\T^2) := \bigg\{\, f:\T^2 \to \R   \, : \, \| f \|_{H^s(\T^2)}^2 :=  \| f \|_{s}^2  := \sum_{(\ell,j)\in \Z^2} \braket{\ell,j}^{2s} |f_{\ell,j}|^2 < \infty  \bigg\}\,.  
\end{equation}
These are the classical Sobolev spaces that will be used from Section \ref{sect.linear.eq} onward in the construction of the time-periodic solutions stated in Theorem \ref{Theo-2}, viewed as small-amplitude perturbations of the stationary solutions obtained in Theorem \ref{Theo-1}.
\begin{rem}
We have, following the standard inclusion $\ell^1(\Z^2) \subset \ell ^2(\Z^2)$, that $X^{\rho,s}\subset H^s$, with
\begin{equation}
    \| f \|_{s} \leq \| f \|_{0,s} \leq \| f \|_{\rho,s} \,.
\end{equation}
In order to construct time-periodic solutions as Sobolev-regular perturbations of Sobolev-analytic traveling waves, we implement a Nash--Moser scheme, which requires a sufficiently high level of regularity for the patch boundary. In fact, our construction yields a stronger conclusion, namely that the boundary of the traveling waves is analytic. Questions concerning the boundary regularity of translating and rotating vortex patches for active scalar equations have been extensively investigated in several previous works, see for instance \cite{Casrro,GH23,HMW20,HMV}.
\end{rem}
Another important class of subspaces, which encodes the reversibility structure of the system and plays a fundamental role in the construction of time-periodic solutions, is given by the spaces of even and odd functions defined as follows:
\begin{align}
	 H^s_{\textnormal{even}}:=  H^s_{\textnormal{even}}(\T^2)& =\Big\{h\in H^{s}(\mathbb{T}^{2}) \, : \,  h(-\varphi,-\theta)=h(\varphi,\theta),\;\, \forall \, (\varphi,\theta)\in\T^{2}\Big\} \,,  \label{FS-even}\\
	H^s_{\textnormal{odd}} :=  H^s_{\textnormal{odd}}(\T^2) & =\Big\{h\in H^{s}(\mathbb{T}^{2})  \, : \,  h(-\varphi,-\theta)=-h(\varphi,\theta), \;\, \forall \,(\varphi,\theta)\in\T^{2}\Big\} \,.\label{FS-odd}
\end{align}
The following subspaces are also fundamental in the construction of the phase spaces associated with the time-periodic solutions,
\begin{align}
    H_0^s & := H_0^s(\T^2) := \Big\{ h\in H^s(\T^2)  \, : \,\int_{\T}h(\vf,\theta) \wrt\theta = 0,\, \, \forall\,\vf\in\T \Big\}\,, \\
    H_\circ^{s}&:=H_{\circ}^{s}(\mathbb{T}^{2}):=\Big\{ h\in H^{s}_0(\mathbb{T}^{2})  \, : \, \int_{\T}h(\varphi,\theta)\sin(\theta)\wrt\theta=0,\, \ \forall\,\varphi\in\T \Big\}\,. \label{Hcirc.def}
\end{align}
The space $H_0^s$  encodes the conservation of mass in the Euler equations, which, in the framework of patch dynamics, translates into the zero-average condition for the boundary perturbation. 
The space $H_\circ^s$, instead, is related to the invariance of the Euler equations under vertical translations in the case of symmetric vortex pairs. This symmetry is reflected in the degeneracy of the first sine mode in the patch formulation.
Finally, we also define the subspaces
\begin{align}\label{Hs-circ}
H^{s}_{\circ,\textnormal{even}}=H^{s}_{\textnormal{even}}\cap H^s_{\circ}\quad\hbox{and}\quad H^{s}_{\circ,\textnormal{odd}}=H^{s}_{\textnormal{odd}}\cap H^s_{\circ} \,.
\end{align}
By setting
\begin{align}\label{basis-1}
\bc_{1}(\theta):= \sqrt{2}\cos(\theta)\,,\quad \bs_{1}(\theta):=\sqrt{2}\sin(\theta)\,,\quad \be_{j}(\theta) := e^{\im j\theta}\,, \quad |j|\geq 2\,,
\end{align}
any element $h\in H^{s}_{\circ}(\T^2)$ decomposes as
\begin{align}
h(\varphi,\theta)
&=h_1(\varphi) {\bf{c}}_1(\theta)+\sum_{|j|\geqslant 2}h_j(\varphi) \be_{j}(\theta)\,,
\end{align}
where
\begin{equation}\label{coeff.ph}
    h_{1}(\vf) := \frac{1}{2\pi} \int_{\T} h(\vf,\eta) {\bf{c}}_1(\eta) \wrt\eta \,, \quad h_{j}(\vf):= \frac{1}{2\pi} \int_{\T} h(\vf,\eta) \be_{j}(-\eta) \wrt\eta\,, \quad |j|\geq 2 \,.
\end{equation}
We also introduce the following projectors: for $h\in H_0^s(\T^2)$,
\begin{align}
 \Pi_{1,\bc}h(\vf,\theta) &:=\frac{1}{2\pi}\int_{\T}h(\varphi,\eta) {\bf{c}}_1(\eta)\wrt \eta\,\, {\bf{c}}_1(\theta) \,, \label{proj.1c} \\
 \Pi_{1,\bs}h(\vf,\theta) & :=\frac{1}{2\pi}\int_{\T}h(\varphi,\eta) {\bf{s}}_1(\eta)\wrt\eta\,\, {\bf{s}}_1(\theta) \,, \label{proj.1s} \\
 \Pi_{\geq 2} h(\vf,\theta) & := \sum_{|j|\geq 2} h_{j}(\vf) e^{\im j \theta} \,.\label{proj.geq2}
 \end{align}
 In particular, we have in $L_0^{2}(\T)$ the decomposition
 \begin{equation}
     {\rm Id} = \Pi_{1,\bs_{1}} + \Pi_{1,\bc_{1}} + \Pi_{\geq 2} \,,
 \end{equation}
 where the three projections are orthogonal with respect to the scalar product
 \begin{equation}\label{scalar.prod.L2}
     \braket{f,g}_{L^2(\T)}:= \fint_{\T} f(\eta) \overline{g(\eta)} \wrt\eta \,.
 \end{equation}
Furthermore, we define
\begin{equation}\label{proj.ph}
    \Pi_{\rm ph} := \Pi_{1,\bc} + \Pi_{\geq 2} \,.
\end{equation}
We want to control also the regularity with respect to parameters $\lambda =(\omega,\alpha)\in  \R \times  [\alpha_{1},\alpha_{2}]$. 
We fix once and for all
\begin{equation}\label{s0-embed}
    s_0 \geq  7 \,.
\end{equation}
We introduce the following norm.
\begin{defi}\label{W.space.def}
	Let $\gamma\in (0,1)$,  $\Lambda \subset \R\times [\alpha_{1},\alpha_{2}]$ a closed set and  $s\geq s_0$, with $s_0$ as in \eqref{s0-embed}. We denote by $W_\gamma^5(\Lambda; H^s(\T^2))$ the space of functions $f:\Lambda \to H^s(\T^2)$, $\lambda \mapsto f(\lambda)$ that are $5$-times differentiable with respect to $\lambda$ and such that their weighted norm, defined by
	\begin{equation}\label{norm.weighted.def}
        \| f\|_{s}^{5,\gamma,\Lambda} := \| f \|_{s}^{5,\gamma}  := \sum_{k=0}^{5} \gamma^{k} \sup_{\lambda\in\Lambda} \|  \pa_{\lambda}^k f(\lambda)\|_{H^{s-k}(\T^2)} \,,
	\end{equation}
	is finite. For simplicity, 
    if $f(\vf,\theta;\lambda) \equiv f(\lambda)\in \C$, we denote $|f|^{5,\gamma,\Lambda}:= \| f \|_{0}^{5,\gamma,\Lambda}$. Also, we shall omit the dependence on the parameter set $\Lambda$ unless it is needed.
\end{defi}

The following properties are standard.

\begin{lem}\label{lem-productlaw}
Let $\gamma \in (0,1)$. The following results hold:
\\[1mm]
\noindent $(i)$
    Let $\rho > 0$ and $s\geq 0$. Then the space  $ X^{\rho,s}$ is an algebra and in particular we have  the estimates
		\begin{align}
			&\| f g \|_{\rho,s}\lesssim_{s} \| f \|_{\rho,s}\| g \|_{\rho,0}+\| f \|_{\rho,0}\| g \|_{\rho,s} \,;
		\end{align}
\noindent $(ii)$
    Let $s\geq s_0$. Then the space  $ W_\gamma^5(\Lambda;H^s(\T^2))$ is an algebra with the tame estimates
		\begin{align}
			&\| f g \|_{s}^{5,\gamma} \lesssim_{s} \| f\|_{s}^{5,\gamma} \| g\|_{s_0}^{5,\gamma} +  \| f\|_{s_0}^{5,\gamma} \| g\|_{s}^{5,\gamma} \lesssim_{s} 2 \| f \|_{s}^{5,\gamma} \| g \|_{s}^{5,\gamma} \,;
		\end{align}
\noindent $(iii)$ For any $n,m\in \N_{0}$, $s\in\mathbb{R}$, we have
    $$
    \|\partial_\varphi^m\partial_\theta^n f\|_{s-m-n}^{5,\gamma}\lesssim \| f\|_{s}^{5,\gamma} \,, \quad \|\partial_\varphi^m |{\rm D}|^n f\|_{s-m-n}^{5,\gamma}\lesssim \| f\|_{s}^{5,\gamma} \,;
    $$
\noindent $(iv)$ If $F:\R\to\R$ is analytic around the origin, there exists $\delta_0>0$ such that, if  $\| f\|_{0}^{5,\gamma}\leq \delta_0 $, then, for any $s\geq s_0$,
    \begin{align}
    	 & \|F(f)-F(0)\|_{s}^{5,\gamma}\lesssim_{s,F} \| f\|_{s}^{5,\gamma} \,, \\
    	 &  \|F(f)-F(0) - F'(0)f\|_{s}^{5,\gamma}\lesssim_{s,F} \| f\|_{s_0}^{5,\gamma}\| f\|_{s}^{5,\gamma} \lesssim \big(\| f\|_{s}^{5,\gamma} \big)^2\,.
    \end{align}
\end{lem}

\subsection{Symplectic change of the space variable}

In Section \ref{sect.linear.eq} and Section \ref{section.normal}, we will consider bounded maps of the form
\begin{equation}\label{map.compo.sympl}
	\cS: H^{s}(\T^2) \mapsto H^{s}(\T^2) \,, \quad (\cS h)(\vf,\theta) := (1+ \pa_{\theta}\beta(\vf,\theta)) (\cB h)(\vf,\theta)\,,
\end{equation}
where $\cB$ is the map
\begin{equation}\label{map.compo}
	(\cB h)(\vf,\theta) := h(\vf,\theta+\beta(\vf,\theta))\,,
\end{equation}
for some function $\beta:\T^2 \to \R$. We have the following result.
\begin{lem}\label{lemma.diffeo.gen}
	Let $\| \beta\|_{s_0}^{5,\gamma}\leq \delta(s_0)$ sufficiently small. The map $\cB$ in \eqref{map.compo} and the  map $\cS$ in \eqref{map.compo.sympl} satisfy the following tame estimates: for any $s\geq s_0$,
	\begin{equation}
	\begin{aligned}
			\| \cB h \|_{s}^{5,\gamma}+	\| \cS h \|_{s}^{5,\gamma}  &  \lesssim_{s} \| h \|_{s}^{5,\gamma} + \| \beta \|_{s}^{5,\gamma} \| h \|_{s_0}^{5,\gamma} \,.
	\end{aligned}
	\end{equation}
	Moreover, the maps $\cB$ and $\cE$ are invertible, with inverses given by
	\begin{equation}
		(\cB h)^{-1}(\vf,\vartheta) = h(\vf,\theta+\breve{\beta}(\vf,\vartheta)) \,, \quad \cS^{-1} = \cB^{-1} \circ (1+\pa_{\theta} \beta(\vf,\theta))^{-1} = \cB^{*} \,,
 	\end{equation}
for some function $\breve{\beta}$ defined by
\[
\theta=\vartheta+\breve{\beta}(\vf,\vartheta) \Longleftrightarrow \vartheta=\theta+\beta(\vf,\theta)\,.
\]
Furthermore, for $s\geq s_0$, the following estimates hold:
 	\begin{equation}
 		\begin{aligned}
 		&	\| \breve{\beta} \|_{s}^{5,\gamma} \lesssim_{s} \| \beta \|_{s}^{5,\gamma} \,, \\
 			&	\| \cB^{-1} h \|_{s}^{5,\gamma} +\| \cS^{-1} h \|_{s}^{5,\gamma} \lesssim_{s} \| h \|_{s}^{5,\gamma} + \| \beta \|_{s}^{5,\gamma} \| h \|_{s_0}^{5,\gamma} \,, \\
 			& 	\| (\cS^{\pm 1} -{\rm Id} ) h \|_{s}^{5,\gamma} \lesssim_{s}  \| \beta\|_{s_0+1}^{5,\gamma} \| h \|_{s+1}^{5,\gamma} + \| \beta \|_{s+1}^{5,\gamma} \| h \|_{s_0+1}^{5,\gamma} \,.
 		\end{aligned}
 	\end{equation}
\end{lem}
\begin{proof}
	It follows as a direct application of Lemma 2.3-$(ii)$ in \cite{BM20}, Lemma 6.1-
\end{proof}
We point out that  $\cS$ and its inverse preserve the spatial average. More precisely, the following lemma holds.
\begin{lem}
	We have that $\braket{\cS h}_{\theta}=\braket{h}_{\theta}$ and $\braket{\cS^{-1}h}_{\vartheta} = \braket{h}_{\vartheta}$. As a consequence, we have that $\cS^{\pm 1}:H_0^{s}(\T^2)\to H_0^{s}(\T^2)$.
\end{lem}
\begin{proof}
	Using a change of variables together with $\wrt\vartheta = (1+\pa_{\theta}\beta(\vf,\theta)) \wrt\theta$,
	we get
	\begin{equation}\label{avera}
		\begin{aligned}
			\braket{\cS h}_{\theta} & =  \fint_{\T} (1+\pa_{\theta}\beta (\vf,\theta)) h(\theta+\beta(\vf,\theta)) \wrt\theta=  \fint_{\T} h(\vartheta)  \wrt\vartheta = \braket{h}_{\theta} \,, \\
			\braket{\cS^{-1} h}_{\theta} & =  \fint_{\T} \frac{h(\vartheta+\breve{\beta}(\vf,\vartheta))}{1+\pa_{\theta}\beta(\vf,\vartheta+\breve{\beta}(\vf,\vartheta))} \wrt \vartheta =  \fint_{\T} h(\theta) \wrt\theta = \braket{h}_{\theta} \,,
		\end{aligned}
	\end{equation}
	from which the claims follow.
\end{proof}

\subsection{Matrix representation of linear operators}
We introduce now the matrix representation of $\vf$-independent linear operators that will be used in Section \ref{sect.linear.eq} to study the linear dynamics around the equilibrium. Recalling \eqref{proj.ph}, let
\begin{align}
	L_{\rm ph}^2(\T) & := \Big\{ h \in L^2(\T) \, : \, h(\theta) = h_{1}\cos(\theta)+ \sum_{|j|\geq 2} h_{j} e^{\im j\theta} =: \Pi_{\rm ph}h(\theta) \Big\} \subset H_{\circ}^0(\T^2) \,.
\end{align}
By direct computations, each $h\in H_\circ^0(\T^2)$, recalling \eqref{Hcirc.def} and \eqref{coeff.ph}, 
 expands as
 \begin{equation}
     h(\vf,\theta) = h_{1}(\vf)\cos(\theta) + \sum_{|j|\geq 2} h_{j}(\vf) e^{\im j\theta}\,,
 \end{equation}
and satisfies the following characterizations:
\begin{align}
    h \quad \hbox{is real} \  & \Longleftrightarrow \quad h_{1}(\vf) = \overline{h_{1}(\vf)} \ \quad \hbox{and}\quad  \ h_{j}(\vf) = \overline{h_{-j}(\vf)} , \ |j|\geq 2 \,; \label{real.h}\\
    h\quad \hbox{is even} \  & \Longleftrightarrow \quad h_{1}(\vf) = h_{1}(-\vf) \ \quad \hbox{and}\quad \ h_{j}(\vf) = h_{-j}(-\vf), \ |j|\geq 2  \,; \label{even.h} \\
     h\quad \hbox{is odd} \ & \Longleftrightarrow \quad h_{1}(-\vf) = - h_{1}(\vf) \ \quad \hbox{and}\quad \ h_{-j}(-\vf) = -h_{j}(\vf), \ |j|\geq 2  \,. \label{odd.h}
\end{align}
Let $\cA: L_{\rm ph}^2(\T)\to L_{\rm ph}^2(\T)$ be a linear operator. It can be represented as
\begin{equation}\label{matrix.repr.ph}
	\cA [h(\theta)] = \Big( \cA_{1}^{1} h_{1}+ \sum_{|j'|\geq 2} \cA_{1}^{j'} h_{j'}  \Big)\cos(\theta) + \sum_{|j|\geq 2}\Big( \cA_{j}^{1} h_{1} +\sum_{|j'|\geq 2} \cA_{j}^{j'} h_{j'} \Big)e^{\im j\theta} \,,
\end{equation}
where the matrix elements $(\cA_{j}^{j'})$ are defined by, for $|j|, |j'|\geq 2$,
\begin{equation}
	\begin{aligned}
		\cA_{1}^{1}&:= \frac{1}{\pi} \int_{\T} \cA[\cos(\eta)] \cos(\eta) \wrt\eta \,, \quad	&	\cA_{1}^{j'}  := \frac{1}{\pi} \int_{\T}\cA[e^{\im\,j'\eta}] \cos(\eta) \wrt\eta \,, \\
		\cA_{j}^{1} & := \frac{1}{2\pi} \int_{\T} \cA[\cos(\eta)] e^{-\im\,j\eta} \wrt\eta \,, \quad&	\cA_{j}^{j'}  := \frac{1}{2\pi} \int_{\T} \cA[e^{\im\,j'\eta}] e^{-\im\,j\eta} \wrt\eta \,.
	\end{aligned}
\end{equation}
To unify the notation (and with a slight conceptual abuse), we define, recalling also \eqref{basis-1},
\begin{equation}\label{basis.ph}
	\Z_{\rm ph} := \Z\setminus\{0,-1\}\,, \quad \textnormal{and} \quad  \be_{j}(\theta):=\begin{cases}
		\bc_{1}(\theta) =   \sqrt{2}\cos(\theta)\,, & j= 1\,, \\
		\be_{j}(\theta) = e^{\im j\theta}\,, & |j|\geq 2\,,
	\end{cases} \quad j \in\Z_{\rm ph}\,.
\end{equation}
In this way, we represent a linear operator $\cA: L_{\rm ph}^2(\T)\to L_{\rm ph}^2(\T)$ as
\begin{equation}\label{matrix.rep.compact}
	h(\theta)= \sum_{j\in\Z_{\rm ph}} h_{j} \be_{j}(\theta) \quad \mapsto \quad \cA[h(\theta)]=\sum_{j,j'\in\Z_{\rm ph}}\cA_{j}^{j'} h_{j'} \be_{j}(\theta) \,,
\end{equation}
where the matrix elements $(\cA_{j}^{j'})_{j,j'\in\Z_{\rm ph}}$ are given by
\begin{equation}\label{matrix.elem.compact}
	\cA_{j}^{j'} := \frac{1}{2\pi} \int_{\T} \cA[\be_{j'}(\eta)] \overline{\be_{j}(\eta)} \wrt\eta \,, \quad j,j'\in\Z_{\rm ph}\,.
\end{equation}
For the implementation of the normal form analysis of the linearized operator developed in Section \ref{sect.linear.eq}, we introduce the following norms. They are designed to measure the off-diagonal decay of the matrix coefficients associated with linear operators acting on $L_{\rm ph}^2(\T)$.
\begin{defi}\label{block norm}
{\bf (Matrix decay norm and the space $\OpM_{s}^{m}$).} Let $m\in\R$, $s\geq s_0$ and $\cA$ be a linear operator in $L_{\rm ph}^2(\T)$ represented as in \eqref{matrix.rep.compact}, \eqref{matrix.elem.compact}. We say that $\cA$ belongs to the class $\OpM_{s}^m$ if we have that
\begin{equation}\label{def decay norm}
		|\cA|_{m,s} := \sup_{j'\in\Z_{\rm ph}} \braket{j'}^{-m}\sum_{j\in\Z_{\rm ph}} \braket{j-j'}^{s}|\cA_{j}^{j'}|  <\infty \,.
\end{equation}  
If a family of linear operators $\cA=(\cA(\lambda))_{\lambda\in\Lambda}$ if $5$-times differentiable with respect to the parameter $\lambda\in\Lambda \subset \R\times [\alpha_{1},\alpha_{2}]$, then, for $\gamma\in(0,1)$, we define
\begin{equation}
	|\cA|_{m,s}^{5,\gamma} := |\cA|_{m,s}^{5,\gamma,\Lambda} := \sum_{k=0}^{5} \gamma^{k} \sup_{\lambda\in\Lambda} |\pa_{\lambda}^k\cA(\lambda)|_{m,s-k}\,. 
\end{equation}
Finally, we define
\begin{equation}\label{OPM-infty}
    \OpM_{s}^{-\infty} := \bigcap_{m\in\N_{0}} \OpM_{s}^{-m} \,.
\end{equation}
\end{defi}
It follows directly from the definition that
\begin{align}
    m\leq m' \ & \Rightarrow \ \OpM_{s}^{m} \subseteq \OpM_{s}^{m'} \ \text{ and }\ |\,\cdot\,|_{m',s}^{5,\gamma} \leq |\,\cdot\,|_{m,s}^{5,\gamma}  \,, \label{monotone.OPM}\\
     s\leq s' \ & \Rightarrow \ \OpM_{s'}^{m} \subseteq \OpM_{s}^{m} \ \text{ and } \ |\,\cdot\,|_{m,s}^{5,\gamma} \leq |\,\cdot\,|_{m,s'}^{5,\gamma}  \,.
\end{align}
We state some standard properties of the decay norms, see for instance Lemma 2.6 in \cite{FrMo}. 
\begin{lem}\label{standard prop decay norm}
$(i)$ Let $s \geq s_0$, $m\geq 0$ and $\cR \in \OpM_{s+m}^{-m}$. 
If $\| u \|_s^{5,\gamma} < \infty$, then
\begin{equation}\label{action.matrix.decay.est}
    \| |{\rm D}|^{m}\cR u \|_{s}^{5,\gamma} \lesssim_{s} |\cR|_{-m, s+m}^{5,\gamma} \| u \|_s^{5,\gamma}\,; 
\end{equation}
\noindent $(ii)$ Let $s \geq s_0$, $m, m' \in \R$, and let $\cR \in\OpM_{s}^m$, $\cQ \in \OpM_{s + |m|}^{m'}$. 
Then $\cR \cQ \in \OpM_{s}^{m+m'}$ and there exist constants $C(s),C(s_{0})>0$ independent of $m,m'$ such that
\begin{equation}
   |\cR \cQ|_{m + m', s}^{5,\gamma} \leq C(s_{0}) |\cR|_{m, s}^{5,\gamma} |\cQ|_{m', s_{0} + |m|}^{5,\gamma} + C(s)|\cR|_{m, s_{0}}^{5,\gamma} |\cQ|_{m', s + |m|}^{5,\gamma} \,. 
\end{equation}
In particular, when $m=0$ and $\cR\in\OpM_{s}^{0}$, $\cQ \in \OpM_{s}^{m'}$. 
then $\cR \cQ \in \OpM_{s}^{m'}$ and there exist constants $C(s),C(s_{0})>0$ independent of $m'$ such that
\begin{equation}
    |\cR \cQ|_{m', s}^{5,\gamma} \leq C(s_{0})|\cR|_{0, s}^{5,\gamma} |\cQ|_{m',s_{0}}^{5,\gamma} + C(s) |\cR|_{0,s_{0}}^{5,\gamma} |\cQ|_{m', s}^{5,\gamma} \,;
\end{equation}
\noindent $(iii)$ Let $s \geq s_{0}$ and ${\mathcal R} \in \OpM_{s}^0$. 
Then, for any integer $n \geq 1$, ${\mathcal R}^n \in \OpM_{s}^0$ and there exist constants $C(s_0),C(s) > 0$, independent of $n$, such that 
\begin{equation*}
	\begin{aligned}
		& |{\mathcal R}^n|_{0, s_0}^{5,\gamma} \leq C(s_0)^{n - 1} \big(|{\mathcal R}|_{0, s_0}^{5,\gamma}\big)^{n} \,, \\
		& |{\mathcal R}^n|_{0, s}^{5,\gamma} \leq \,C(s)^{n - 1} \big(C(s_0)|{\mathcal R}|_{0, s_0}^{5,\gamma}\big)^{n - 1} |{\mathcal R}|_{0, s}^{5,\gamma}\,;
	\end{aligned}
\end{equation*}
\noindent $(iv)$ Let $s \geq s_{0}$, $m \geq 0$ and ${\mathcal R} \in \OpM_{s}^{-m}$. 
Then there exists $\delta(s) \in (0, 1)$ small enough such that, 
if $|{\mathcal R}|_{- m, s_0}^{5,\gamma} \leq \delta(s)$, 
then the map $\Phi = {\rm Id} + {\mathcal R}$ is invertible 
and the inverse satisfies the estimate 
\[
|\Phi^{- 1} - {\rm Id}|_{- m, s}^{5,\gamma} \lesssim_{s, m} |{\mathcal R}|_{- m, s}^{5,\gamma} \,;
\]
\noindent
$(v)$ Let $s \geq s_0$, $m \in \R$ and ${\mathcal R} \in \OpM_{s}^m$. 
Then 
\[
| \mathcal R_j^j |^{5,\gamma} \lesssim \langle j \rangle^m|{\mathcal R}|_{s_0}^{5,\gamma}\,.
\] 
\end{lem}

\begin{proof}
We first prove item $(i)$. We start by proving the claim when $\cR \in \OpM_{s}^{0}$ and we assume that $\cR$ does not depend on parameters. We have
\begin{align}
   \| \cR u \|_{s}^2 &   =\sum_{j\in \Z_{\rm ph}}  \Big| \sum_{j'\in \Z_{\rm ph}} \braket{j}^{s} \cR_{j}^{j'} u_{j'} \Big|^2  \leq ({\rm I}) + ({\rm II})  \,,
\end{align}
where
\begin{align}
    ({\rm I}) &  :=   \sum_{j\in \Z_{\rm ph}}  \bigg| \sum_{j'\in \Z_{\rm ph} \atop  \braket{j}\leq 2 \braket{j'} } \frac{\braket{j}^s}{\braket{j-j'}^{s_0} \braket{j'}^s} \braket{j-j'}^{s_0}\cR_{j}^{j'} \braket{j'}^{s} u_{j'} \bigg|^2  \,,  \\
    ({\rm II}) &  :=   \sum_{j\in \Z_{\rm ph}}  \bigg| \sum_{j'\in \Z_{\rm ph}\atop \braket{j} >  2 \braket{j'}} \frac{\braket{j}^s}{\braket{j-j'}^{s} \braket{j'}^{s_0}} \braket{j-j'}^{s}\cR_{j}^{j'} \braket{j'}^{s_0} u_{j'} \bigg|^2 \,. \label{bas}
\end{align}
We first estimate $({\rm I})$.  By Cauchy-Schwartz, using also that $\ell^1(\Z_{\rm ph}) \subset \ell^2(\Z_{\rm ph})$, we get
\begin{align}
    ({\rm I}) & \lesssim_{s}  \sum_{j\in \Z_{\rm ph}}\bigg(  \sum_{j'\in \Z_{\rm ph}} \frac{1}{\braket{j-j'}^{2s_0}} \bigg) \bigg(  \sum_{j'\in \Z_{\rm ph}} \braket{j-j'}^{2s_0}|\cR_{j}^{j'}|^2 \braket{j'}^{2s} |u_{j'}|^{2} \bigg) \\
    & \lesssim_{s} \sum_{j\in \Z_{\rm ph}} \sup_{j'\in\Z_{\rm ph}} \braket{j-j'}^{2s_0} |\cR_{j}^{j'}|^2 \sum_{j'\in\Z_{\rm ph}} \braket{j'}^{2s} |u_j|^2 \\
    & \lesssim_{s} \bigg( \sup_{j'\in \Z_{\rm ph}} \sum_{j\in\Z_{\rm ph}} \braket{j-j'}^{s_0} |\cR_{j}^{j'}| \bigg)^2 \| u \|_{s}^2 \lesssim_{s} |\cR|_{0,s_0}^{2}\| u \|_{s}^2 \,, \label{I}
\end{align}
where we used $\sum_{j\in \Z_{\rm ph}}\braket{j-j'}^{-2s_0}<+\infty$ for $s_0$ as in \eqref{s0-embed}. To estimate $({\rm II})$, we note that
\begin{equation}
    \braket{j} > 2 \braket{j'} \quad \Rightarrow \quad \braket{j}\leq \braket{j-j'} + \braket{j'} \leq \braket{j-j'} + \tfrac12 \braket{j} \quad \Rightarrow \quad  \tfrac12 \braket{j} \leq \braket{j-j'} \,.
\end{equation}
Then, using again Cauchy-Schwartz and $\ell^1(\Z_{\rm ph}) \subset \ell^2(\Z_{\rm ph})$, we get
\begin{align}
    ({\rm II}) & \lesssim_{s}  \sum_{j\in \Z_{\rm ph}}\bigg(  \sum_{j'\in \Z_{\rm ph}} \frac{1}{\braket{j'}^{2s_0}} \bigg) \bigg(  \sum_{j'\in \Z_{\rm ph}} \braket{j-j'}^{2s}|\cR_{j}^{j'}|^2 \braket{j'}^{2s_0} |u_{j'}|^{2} \bigg) \\
    & \lesssim_{s} \sum_{j\in \Z_{\rm ph}} \sup_{j'\in\Z_{\rm ph}} \braket{j-j'}^{2s} |\cR_{j}^{j'}|^2 \sum_{j'\in\Z_{\rm ph}} \braket{j'}^{2s_0} |u_j|^2 \\
    & \lesssim_{s} \bigg( \sup_{j'\in \Z_{\rm ph}} \sum_{j\in\Z_{\rm ph}} \braket{j-j'}^{s} |\cR_{j}^{j'}| \bigg)^2 \| u \|_{s_0}^2 \lesssim_{s} |\cR|_{0,s}^{2}\| u \|_{s_0}^2 \,, \label{II}
\end{align}
where we used $\sum_{j\in \Z_{\rm ph}}\braket{j'}^{-2s_0}<+\infty$. Collecting \eqref{I} and \eqref{II} into \eqref{bas}, we conclude that
\begin{equation}\label{action.easy}
    \| \cR u \|_{s} \lesssim_s |\cR|_{0,s_0} \|u \|_{s}+ |\cR|_{0,s}\| u \|_{s_0} \lesssim |\cR|_{0,s} \|u \|_{s}\,,
\end{equation}
by monotonicity of the Sobolev topologies.

Then, we deduce the claimed estimate in \eqref{action.matrix.decay.est} when $m=0$ and $\cR=\cR(\lambda)$ depends on parameters by Definition \ref{block norm},{\eqref{norm.weighted.def}}, Leibniz rule,  \eqref{action.easy} and an induction argument.
 To prove the general claim for $\cR \in \OpM_{s+m}^{-m}$ with $m>0$, we note that
\begin{align}
    \| |{\rm D}|^{m}\cR u \|_{s}^{5,\gamma} \lesssim_{s} \big| |{\rm D}|^m \cR \big|_{0,s}^{5,\gamma}  \| u \|_{s}^{5,\gamma} \lesssim_{s,m} |\cR|_{-m,s+m}^{5,\gamma}  \| u \|_{s}^{5,\gamma} \,,
\end{align}
using that $|{\rm D}|^m \cR \in \OpM_{s}^{0}$ when  $\cR \in \OpM_{s+m}^{-m}$. This concludes the proof of item $(i)$.
\\[1mm]
We now prove item $(ii)$.
 We start by assuming that both $\cR$ and $\cQ$ do not depend on the parameter $\lambda$. The matrix elements for the composition operator $\cR\cQ$ follow the rule
\begin{equation}
	(\cR\cQ)_j^{j'} = \sum_{i \in \Z_{\rm Ph}} \cR_{j}^{i} \cQ_{i}^{j'}\,.
\end{equation}
We have
\begin{align}
\langle j' \rangle^{ - (m + m')} \sum_{j\in \Z_{\rm Ph}} \langle  j - j' \rangle^{s} |	(\cR\cQ)_j^{j'}| \leq (A) + (B) \,, \label{paris 0} 
\end{align}
where 
\begin{align}
&
(A) :=  \langle j' \rangle^{- (m + m')} \sum_{j,i \in \Z_{\rm Ph} \atop \braket{j-j'}\leq 2^{1/s}\braket{j-i}} \langle j - j' \rangle^{s} |\cR_j^{i} | | \cQ_{i}^{j'} | \,,  \nonumber	
\\
&
(B) := \langle j' \rangle^{- (m + m')}\sum_{j,i \in \Z_{\rm Ph} \atop  \braket{j-j'}> 2^{1/s}\braket{j-i}} \langle  j - j'  \rangle^{s} |\cR_j^{i}| | \cQ_{i}^{j'} |  \,.\nonumber 
\end{align}
We start with estimating $(A)$. 
By the Peetre inequality 
$
\langle i \rangle^m\langle j' \rangle^{- m} \lesssim_m  \langle j' - i \rangle^{|m|}
$ and having the series $\sum_{i \in \Z_{\rm Ph}} \langle  i - j'\rangle^{-s_{0}}  = C(s_{0})<\infty$ with $s_0>1$ as in \eqref{s0-embed},
one has
\begin{align}
    (A) & \lesssim_{m}  \sum_{j,i \in \Z_{\rm Ph} \atop \braket{j-j'}\leq 2^{1/s}\braket{j-i}} \frac{\braket{j-j'}^{s}}{\braket{j-i}^{s}\braket{i-j'}^{s_{0}}} \braket{i}^{-m}\braket{j-i}^{s} |\cR_{j}^{i}| \braket{j'}^{-m}\braket{i-j'}^{s_{0}+|m|} |\cQ_{i}^{j'}| \\
    & \lesssim_{m} \braket{i}^{-m} \sup_{i\in\Z_{\rm ph}}\sum_{j\in\Z_{\rm}}\braket{j-i}^{s} |\cR_{j}^{i}| \sum_{i\in\Z_{\rm ph}}\frac{1}{\braket{i-j'}^{s_{0}}}  \braket{j'}^{-m}\braket{i-j'}^{s_{0}+|m|} |\cQ_{i}^{j'}| \\
    & \stackrel{\eqref{def decay norm}}{\lesssim_{m}}C(s_{0}) |\cR|_{-m,s} |\cQ|_{-m',s_{0}+|m|}\,.
\end{align}
By similar arguments, having that,
\begin{equation}
     \braket{j-j'}> 2^{1/s}\braket{j-i} \quad \Rightarrow \quad \braket{j-j'} < \frac{2^{1/s}}{2^{1/s}-1} \braket{i-j'}\,,
\end{equation}
one gets $(B) \lesssim_m C(s)|\cR|_{s_{0}, m}|\cQ|_{s + |m|, m'} $, with $C(s):= \frac{2}{(2^{1/s}-1)^s}C(s_{0})$, and hence the claimed estimate follows by taking the supremum over $j' \in \Z^2$ in \eqref{paris 0}. If we reintroduce the dependence on the parameter $\lambda$, the  estimate for the seminorms $\sup_{\lambda\in\Lambda}|\pa_{\lambda}^k(\cP\cQ)(\lambda)|_{m,s-k}$, $k=0,1,...,5$ follows as usual by Leibniz rule with an induction argument.
\\[1mm]
\noindent
The claim in item $(iii)$ follows by an induction argument and item $(ii)$.
\\[1mm]
\noindent  The claim in item $(iv)$ follows by a Neumann series argument, together with item $(iii)$.
\\[1mm]
\noindent The claim in item $(v)$ is a direct consequence of Definition \ref{block norm}.
\end{proof}

There are special classes of linear operators that belong to $\OpM_{s}^{-\infty}$.
\begin{lem}\label{lemma.kernelHS}
	Let $s\in\N$, with $s\geq s_0\geqslant7$. The following holds:
	\\[1mm]
	\noindent $(i)$ If a linear operator $\cR=\cR(\lambda)$ acting on  $L_{\rm ph}^2(\T)$ is an integral operator
	\begin{equation}\label{kernel.gen}
		\cR(\lambda) h(\theta) = \int_{\T} K(\theta,\eta;\gamma) h(\eta) \wrt\eta
	\end{equation}
	for some kernel $K(\,\cdot\,,\,\cdot\,;\gamma)\in \cC^\infty(\T^2)$, which  is $5$-times differentiable with respect to $\lambda\in\Lambda$,  then $\cR\in \OpM_{s}^{-\infty}$, with estimate, for any $N\in\N_{0}$,
	\begin{align*}
		|\cR|_{-N,s}^{5,\gamma} &\lesssim_{N,s} \| K \|_{\cC^{N+s+2}(\T^2)}^{5,\gamma} \,;
	\end{align*}
	\noindent $(ii)$ If a linear operator $\cR=\cR(\lambda)$ in $L_{\rm ph}^2(\T)$ is a finite rank operator of the form
	\begin{equation}\label{finite.rank.gen}
		\cR(\lambda) h(\theta) = \sum_{j=1}^{m} \braket{h,g_j(\,\cdot\,;\lambda)}_{L_\theta^2(\T)} \phi_{j}(\theta;\lambda)\,,
	\end{equation}
	for some $g_j(\,\cdot\,;\lambda),\phi_{j}(\,\cdot\,;\lambda) \in \cC^{\infty}(\T)$, $j=1,...,m$, , which  are $5$-times differentiable with respect to $\lambda\in\Lambda$,  then $\cR\in \OpM_{s}^{-\infty}$, with estimates, for any $N\in\N_{0}$,
	\begin{equation}
		|\cR|_{-N,s}^{5,\gamma} \lesssim_{N,s} \sum_{j=1}^{m} \| g_j \|_{\cC^{N+s+1}(\T)}^{5,\gamma}  \| \phi_{j} \|_{\cC^{s+2}(\T)}^{5,\gamma}  \,.
	\end{equation}
\end{lem}
\begin{proof}
	We note that it is enough to prove item $(i)$, since for item $(ii)$ we have that \eqref{finite.rank.gen} satisfies \eqref{kernel.gen} with $K(\theta,\eta;\lambda):= g_j(\eta;\lambda)\phi_{j}(\theta;\lambda)$. Moreover, we prove the claim only for kernels $K(\theta,\eta)$ independent of $\lambda\in\Lambda$, since the proof for parameter dependent kernels follows by similar arguments.
	For $s\in\N$ and $N\in\N_{0}$, we get through integration by parts 
	\begin{align}
		|j'|^{N} \braket{j-j'}^{s}\int_{\T^2} & K(\theta,\eta)e^{\im(j'\eta-j\theta)} \wrt\eta\wrt\theta  =  \int_{\T^2} K(\theta,\eta) |{\rm D}_\eta|^{N} \braket{{\rm D}_\theta-{\rm D}_\eta}^{s} e^{\im(j'\eta-j\theta)} \wrt\eta \wrt \theta \\
		& =  \int_{\T^2}  |{\rm D}_\eta|^{N} \braket{{\rm D}_\theta-{\rm D}_\eta}^{s} K(\theta,\eta)e^{\im(j'\eta-j\theta)} \wrt\eta\wrt\theta \\
		& = \frac{1}{|j|^2} \int_{\T^2}  |{\rm D}_\eta|^{N} \braket{{\rm D}_\theta-{\rm D}_\eta}^{s} | {\rm D}_\theta|^2  K(\theta,\eta)e^{\im(j'\eta-j\theta)} \wrt\eta\wrt\theta \,.
	\end{align}
	Therefore, we deduce that
	\begin{align}
		 |\cR|_{-N,s} & = \sup_{j'\in\Z_{\rm ph}}\sum_{j\in\Z_{\rm ph}}  |j'|^{N} \braket{j-j'}^{s}|\cR_{j}^{j'}| \lesssim {\sum_{j\in\Z_{\rm ph}} \frac{\| K \|_{\cC^{N+s+2}(\T^2)}}{|j|^2}} \\
         &\lesssim \| K \|_{\cC^{N+s+2}(\T^2)} \,,
	\end{align}
	and the claimed estimate follows.
\end{proof}

\subsection{Reversibility and homological-type equation}

With the spaces introduced in \eqref{FS-even} and \eqref{FS-odd} in mind, we can now introduce the following definitions.
\begin{defi}\label{defi.rev}
	$(i)$ 
    We say that a linear operator $\Phi$ is {\it reversible} if
\begin{equation}
    \Phi:H^s_{\mathrm{even}}(\T^2)\to H^{s'}_{\mathrm{odd}}(\T^2)
\qquad \text{and} \qquad
\Phi:H^s_{\mathrm{odd}}(\T^2)\to H^{s'}_{\rm even}(\T^2) \,,
\end{equation}
and that it is {\it reversibility preserving} if
\begin{equation}
    \Phi:H^s_{\mathrm{even}}(\T^2)\to H^{s'}_{\mathrm{even}}(\T^2)
\qquad \text{and} \qquad
\Phi:H^s_{\mathrm{odd}}(\T^2)\to H^{s'}_{\mathrm{odd}}(\T^2)\,,
\end{equation}
  for some $s,s'\in\R$;    
	\\[1mm]
	\noindent $(ii)$ We say that an operator $\Phi:L^2(\T^2)\to L^2(\T^2)$ is {\it real} if $\Phi(u)$ is real valued for any $u$ real valued.
\end{defi}

For time independent linear operators acting on $L_{\rm ph}^2(\T)$, we characterize these properties in terms of the matrix elements.
\begin{lem}\label{prop.matrix.revreal}
	A time independent linear operator $\cA$ acting on  $L_{\rm ph}^2(\T)\subset H_{\circ}^{0}(\T^2)$ satisfies the following properties:
	\\[1mm]
	\noindent $(i)$ It is real if and only if
    \[
\cA_{1}^{1}\in\mathbb{R},\qquad
\cA_{1}^{j'}=\overline{\cA_{1}^{-j'}},\qquad
\cA_{j}^{1}=\overline{\cA_{-j}^{1}},\qquad
\cA_{j}^{j'}=\overline{\cA_{-j}^{-j'}},\quad \forall\,|j|,|j'|\ge2 \,;
\]
    \noindent $(ii)$  It is reversible if and only if
\[
\cA_{1}^{1}=0, \qquad
\cA_{1}^{-j'}= -\,\cA_{1}^{j'}, \qquad
\cA_{-j}^{1}= -\,\cA_{j}^{1}, \qquad
\cA_{-j}^{-j'}= -\,\cA_{j}^{j'},
\quad \forall\, |j|,|j'|\geq 2 \,;
\]
	\noindent $(iii)$ It is reversibility preserving if and only if
\[
\cA_{1}^{-j'}=\cA_{1}^{j'}, \qquad
\cA_{-j}^{1}=\cA_{j}^{1}, \qquad
\cA_{-j}^{-j'}=\cA_{j}^{j'},
\quad \forall\,|j|,|j'|\ge 2 \,.
\]
\end{lem}
\begin{proof}
    The claim for item $(i)$ (resp. for items $(ii),(iii)$) follows by Definition \ref{defi.rev} and direct computations with \eqref{matrix.repr.ph} and \eqref{real.h} (resp. \eqref{even.h}, \eqref{odd.h}).
\end{proof}
At last, we analyze a homological equation within the class of operators introduced above, which will play a central role in the full reducibility procedure developed in Section \ref{section-Full-redu}. More precisely, given a diagonal operator $\cD\in \OpM_{s}^{1}$ and a perturbation operator $\cE\in \OpM_{s}^{m}$, acting on $L_{\rm ph}^2(\T),$ we seek an operator $\cY$ solving the commutator equation
\begin{equation}\label{homol.equa}
	[\cD,\cY]=\cE\,.
\end{equation}
This equation naturally arises in the normal form analysis as the mechanism allowing one to eliminate off-diagonal contributions through successive conjugations of the linearized operator.
We assume that the diagonal operator $\cD$ has spectrum
\[
\big(\lambda_j(\alpha)\big)_{j\in\Z_{\rm ph}},
\]
depending smoothly on the parameter $\alpha\in[\alpha_1,\alpha_2]$. We further assume that the eigenvalues are simple and satisfy the following uniform gap and regularity conditions:
\begin{equation}\label{gap.D}
\begin{cases}
\displaystyle 
\inf_{\alpha\in[\alpha_{1},\alpha_{2}]}
|\lambda_{j}(\alpha)-\lambda_{j'}(\alpha)|
\geq
\nu_{1}|j-j'| \,,
\\[0.3cm]
\displaystyle
\sup_{k=1,\dots,5}\;
\sup_{\alpha\in[\alpha_{1},\alpha_{2}]}
\big|
\pa_{\alpha}^{k}
\big(
\lambda_{j}(\alpha)-\lambda_{j'}(\alpha)
\big)
\big|
\leq
\nu_{2}|j-j'| \,,
\\[0.3cm]
\end{cases}
\qquad
\forall\, j\neq j'\in\Z_{\rm ph}\,,
\end{equation}
for some given   $\nu_{1},\nu_{2}>0$. 
These assumptions quantify both the spectral separation of the diagonal frequencies and their smooth dependence with respect to the parameter $\alpha$. We also assume that
\begin{equation}\label{pure.D}
    \lambda_{-j}(\alpha)=-\lambda_{j}(\alpha) = \overline{\lambda_{j}(\alpha)} \,.
\end{equation}

Then, our main result reads as follows.
\begin{lem}\label{lemma.sol.hom}
Let $s,m\in\R$ and let $\cD\in \OpM_{s}^{1}$ be a diagonal operator satisfying the spectral assumptions \eqref{gap.D}, \eqref{pure.D}. 
Assume that $\cE\in \OpM_{s}^{m}$ satisfies
\[
\cE_{j}^{j}=0,
\qquad \forall\, j\in\Z_{\rm ph} \,.
\]
Then there exists an operator $\cY\in \OpM_{s+1}^{m}$ solving the homological equation \eqref{homol.equa}, whose matrix elements are given by
\begin{equation}\label{solu.Y}
\cY_{j}^{j'}
:=
\begin{cases}
\displaystyle
\frac{\cE_{j}^{j'}}{\lambda_{j}-\lambda_{j'}},
& j\neq j',
\\[0.3cm]
0,
& j=j',
\end{cases}
\end{equation}
and satisfying the estimate
\begin{equation}
|\cY|_{m,s+1}^{5,\gamma}
\leq
\mu\, |\cE|_{m,s}^{5,\gamma},
\qquad
\mu=\mu(\nu_{1},\nu_{2})>0.
\end{equation}
Moreover, if  $\cE$ is a reversible operator, then $\cY$ is reversibility preserving.
\end{lem}

\begin{proof}
	A direct computation shows that \eqref{solu.Y} solves the equation \eqref{homol.equa} when projecting it on the basis elements.
	Therefore, the claimed  estimate follows easily from \eqref{solu.Y}, \eqref{gap.D} and Definition \ref{block norm}, using Leibniz rule and Fa\`a di Bruno formula 
    when differentiating with respect to the parameter. The symmetry property of $\cY$ follows from Lemma \ref{prop.matrix.revreal} and \eqref{solu.Y}.
\end{proof}

\section{Regularity estimates for the nonlocal system}\label{Ref-VP}
As discussed previously in Section \ref{model.derivation}, the symmetry assumptions allow us to reduce the dynamics of the vortex pair to the evolution of a single vortex patch.
In this section, we investigate the structure of the linearized operator associated with the functional \(\cG\) introduced in \eqref{Func.G}, in a neighborhood of a small state, and analyze its functional regularity properties. In particular, we identify the principal contributions arising from the self-induced dynamics and the interaction terms. These results will play a fundamental role in the analysis developed later and will be used, in particular, in the proof of the existence of traveling vortex pairs in Section \ref{Traveling pairs}.

\subsection{Linearization}

 In polar coordinates, the corresponding contour dynamics equation takes the form \eqref{Func.G}-\eqref{F-compact}. The purpose of this section is to investigate the structure of the associated linearized operator around a small state.
More precisely, we linearize the nonlinear functional $\cG$ with respect to the boundary variable and identify the different contributions arising from the self-induced motion and the interaction terms. 
Our main result in this direction is the following proposition.

 \begin{pro}\label{linea-prop}
     Let $\cG$ be the functional introduced in \eqref{Func.G}. 
Let $\overline r$ be a smooth and sufficiently small function, and let $h$ be a smooth test function. Then the linearization of $\cG$ with respect to the boundary variable $r$ around a small  state $\overline r$ takes the form

     \begin{align}\label{nonlinear.cG}
\di_{r}\cG(\tU,\overline{r};\alpha)h= \omega \pa_{\vf} h &+ \pa_{\theta} \bigg(  \tfrac{\tU(\vf)}{\sqrt{1+2 \overline{r}}} \cos(\theta)\, h +\tV[\overline r;\alpha] h + \cR(\overline r;\alpha)[h] \bigg)\,, 
         \end{align}
         where, setting $\overline R (\vf,\theta):=\sqrt{1+2\overline r(\vf,\theta)}$,
         \begin{align}\label{V.at.boundary}
			\nonumber \tV[\overline r;\alpha](\vf,\theta) :=  - \frac{1}{ \overline R(\vf,\theta) }\Bigg( &\fint_{\T} \log\big|\overline R(\vf,\theta) e^{\im\theta}- \overline R(\vf,\eta)e^{\im\eta} \big| \pa_{\eta} \big( \overline R(\vf,\eta) \sin(\eta-\theta) \big) \wrt\eta \\
			+ & \fint_{\T} K_\alpha\big(\overline R(\vf,\theta)e^{\im\theta}, \overline R(\vf,\eta)e^{\im\eta} \big) \pa_{\eta} \big( \overline R(\vf,\eta) \sin(\theta+\eta) \big) \wrt \eta \Bigg)\,,
		\end{align}
  and 
\begin{align}\label{cR.lin.app}
	\nonumber 	\cR(\overline{r};\alpha)[h](\vf,\theta)& :=  \fint_{\T} \log \big| \overline R(\vf,\theta)e^{\im\theta} - \overline R(\vf,\eta)e^{\im\eta} \big| h(\vf,\eta) \wrt \eta \\
	&\  - \fint_{\T} K_{\alpha}\big( \overline R(\vf,\theta)e^{\im\theta}, \overline R(\vf,\eta)e^{\im\eta} \big) h(\vf,\eta) \wrt \eta\,.
\end{align}
 \end{pro} 
 \begin{proof}
 The main difficulty in the linearization of \eqref{Func.G} arises from the nonlinear nonlocal contribution \(F\) introduced in \eqref{F.app}. In order to analyze this term, we shall use a polar coordinates parametrization of the patch together with suitable boundary integral representations, allowing us to rewrite the volume integrals in terms of contour quantities. By using a polar coordinates change of variables, we obtain
\begin{align}
		F[r;\alpha]\vf,\theta) 
		& =    \fint_{\T} \int_{0}^{R(\vf,\eta)} \log |  \sqrt{1+2r(\vf,\theta)}e^{\im\theta} - \rho e^{\im\eta} | \rho \wrt \rho \wrt \eta 
		\\
		& \quad -   \fint_{\T} \int_{0}^{R(\vf,\eta)} K_\alpha(\sqrt{1+2 r(\vf,\theta)}e^{\im\theta},\rho e^{\im\eta})  \rho \wrt \rho \wrt \eta \,. \label{polar.integrals}
	  \end{align}
     Using Gâteaux derivatives, straightforward computations yield
\begin{align}\label{nonlinear.cG}
\di_{r}\cG(\tU,\overline{r};\alpha)h= \omega \pa_{\vf} h &+ \pa_{\theta} \bigg(  \tfrac{\tU(\vf)}{\sqrt{1+2 \overline{r}}} \cos(\theta)\, h +\di_{r}  F[\overline r; \alpha] h \bigg)\,, 
         \end{align}
         with  $\di_{r} F[\overline r; \alpha]$  a self-adjoint operator that, by \eqref{polar.integrals}, is decomposed into two distinct components: a multiplication operator and a compact integral operator, namely
\begin{equation}
	\di_{r}  F[\overline r; \alpha] h =   \tV[\overline r;\alpha] h + \cR(\overline r;\alpha)[h]  \,,
\end{equation}
where the real-valued function $ \tV[\overline r;\alpha]$ is given by
\begin{align}
	 \tV[\overline r;\alpha](\vf,\theta) &  := \frac{1}{2\pi\overline{R}} \bigg(\int_{O_{1}(\varphi)} \nabla_z\log \big| z -\xi \big|\cdot e^{\im \theta} \wrta(\xi)\bigg)\bigg|_{z=\overline{R} e^{\im \theta}} \\
	&  -\frac{1}{2\pi \overline{R}}\bigg( \int_{O_{1}(\varphi)} \nabla_z K_\alpha \big(z, \xi \big)\cdot e^{\im \theta}   \wrta(\xi)\bigg)\bigg|_{z=\overline{R} e^{\im \theta}}  \,,\label{V.lin.app}
\end{align}
where $\nabla_{z}=(\pa_{z},\pa_{\overline{z}})$, and the nonlocal operator $\cR(\overline r;\alpha)$ takes the form
\begin{align*}
	\nonumber 	\cR(\overline{r};\alpha)[h](\vf,\theta)& :=  \fint_{\T} \log \big| \overline R(\vf,\theta)e^{\im\theta} - \overline R(\vf,\eta)e^{\im\eta} \big| h(\vf,\eta) \wrt \eta \\
	&\  - \fint_{\T} K_{\alpha}\big( \overline R(\vf,\theta)e^{\im\theta}, \overline R(\vf,\eta)e^{\im\eta} \big) h(\vf,\eta) \wrt \eta\,.
\end{align*}
To conclude the claimed formula \eqref{V.at.boundary}, we need to rewrite the quantity \( \tV[\overline r;\alpha]\) in terms of boundary integral operators. 
 Applying \eqref{diff-gener}, we infer that 
\begin{align}
	\tV[\overline r;\alpha](\vf,\theta) &  = \frac{1}{\overline{R}} {\rm Re}\bigg(\frac{1}{\pi}\int_{O_{1}(\varphi)} (\partial_{\overline{z}}\log \big| z -\xi \big|) e^{-\im \theta} \wrta(\xi)\bigg)\bigg|_{z=\overline{R} e^{\im \theta}} \label{V.lin.app.1}\\
	&  -\frac{1}{ \overline{R}} {\rm Re}\bigg(\frac{1}{\pi} \int_{O_{1}(\varphi)} \partial_z K_\alpha \big(z, \xi \big) e^{\im \theta}   \wrta(\xi)\bigg)\bigg|_{z=\overline{R} e^{\im \theta}}  \,.\label{V.lin.app.2}
\end{align}
We start by computing the term in \eqref{V.lin.app.1}. By \eqref{partiaz-log} and \eqref{stokes.complex.app} we deduce that
\begin{align*}
    \bigg(\frac{1}{\pi}\int_{O_{1}} (\partial_{\overline{z}}\log \big| z -\xi \big|) e^{-\im \theta} \wrta(\xi)\bigg)\bigg|_{z=\overline{R} e^{\im \theta}}&=-\frac{1}{2\im \pi}\int_{\partial O_{1}} \log \big| \overline{R}(\theta) e^{\im \theta} -\xi \big| e^{-\im \theta} \wrt \xi\\
    &=-\frac{1}{2\im \pi}\int_{\mathbb{T}} \log \big| \overline{R}(\theta) e^{\im \theta} -\overline{R}(\eta) e^{\im \eta} \big| \partial_\eta\big[\,\overline{R}(\eta) e^{\im (\eta-\theta)}\big] \wrt \eta\,.
\end{align*}
Therefore
\begin{equation} \label{heat1}
    \begin{footnotesize}
        \begin{aligned}
              {\rm Re}\bigg(\frac{1}{\pi}\int_{O_{1}} (\partial_{\overline{z}}\log \big| z -\xi \big|) e^{-\im \theta} \wrta(\xi)\bigg)\bigg|_{z=\overline{R} e^{\im \theta}}
    &=\fint_{\mathbb{T}} \log \big| \overline{R}(\theta) e^{\im \theta} -\overline{R}(\eta) e^{\im \eta} \big| \partial_\eta\big[\,\overline{R}(\eta) \sin (\theta-\eta)\big] \wrt \eta \,. 
        \end{aligned}
    \end{footnotesize}
\end{equation}
For the second term \eqref{V.lin.app.2}, we use the identity \eqref{specific-K} together with \eqref{stokes.complex.app}, leading to
\begin{align*}
\bigg(\frac{1}{\pi} \int_{O_{1}(\varphi)} \partial_z K_\alpha \big(z, \xi \big) e^{\im \theta}   \wrta(\xi)\bigg)\bigg|_{z=\overline{R}(\theta) e^{\im \theta}}&=\frac{1}{2\im \pi} \int_{\partial O_{1}}  K_\alpha \big(\overline{R}(\theta) e^{\im \theta}, \xi \big) e^{\im \theta}   d\xi  \\
&=\frac{1}{2\im \pi} \int_{\mathbb{T}}  K_\alpha \big(\overline{R}(\theta) e^{\im \theta}, \overline{R}(\eta) e^{\im \eta} \big) \partial_\eta\big[\, \overline{R}(\eta)e^{\im (\theta+\eta)}\big]   \wrt \eta\,.  
\end{align*}
Consequently
\begin{equation}\label{heat2}
    \begin{footnotesize}
        \begin{aligned}
            {\rm Re}\bigg(\frac{1}{\pi} \int_{O_{1}(\varphi)} \partial_z K_\alpha \big(z, \xi \big) e^{\im \theta}   \wrta(\xi)\bigg)\bigg|_{z=\overline{R}(\theta) e^{\im \theta}}&
= \fint_{\mathbb{T}}  K_\alpha \big(\overline{R}(\theta) e^{\im \theta}, \overline{R}(\eta) e^{\im \eta} \big) \partial_\eta\big[\, \overline{R}(\eta)\sin\big(\theta+\eta\big)\big]   \wrt \eta \,. 
        \end{aligned}
    \end{footnotesize}
\end{equation}
Collecting \eqref{heat1}, \eqref{heat2} into \eqref{V.lin.app.1}, \eqref{V.lin.app.2}, it follows that
\begin{align*}
			\nonumber \tV[\overline r;\alpha](\vf,\theta) =  - \frac{1}{ \overline R(\vf,\theta) }\bigg( &\fint_{\T} \log\big|\overline R(\vf,\theta) e^{\im\theta}- \overline R(\vf,\eta)e^{\im\eta} \big| \pa_{\eta} \big[ \,\overline R(\vf,\eta) \sin(\eta-\theta) \big] \wrt\eta \\
			+ & \fint_{\T} K_\alpha\big(\overline R(\vf,\theta)e^{\im\theta}, \overline R(\vf,\eta)e^{\im\eta} \big) \pa_{\eta} \big[ \,\overline R(\vf,\eta) \sin(\theta+\eta) \big] \wrt \eta \bigg).
		\end{align*}
This ends the proof of Proposition \ref{linea-prop}.
 \end{proof}

\subsection{Regularity estimates}

We will establish in  this section  regularity estimates for the nonlinear vector field associated with \eqref{Func.G} together with its linearization around small states. 
These estimates provide the functional framework required for the bifurcation analysis developed in Section \ref{Traveling pairs}. More precisely, they will be used to prove the existence of traveling vortex pairs bifurcating from infinity, corresponding to the regime in which the two patches are sufficiently far apart and the interaction between them remains weak. In this limit, the dynamics may be viewed as a perturbation of the decoupled stationary configuration, allowing us to implement an implicit function theorem argument in suitable functional spaces.
The first result is the following. The Sobolev-analytic spaces $X^{\rho,s}$ used below are defined in \eqref{X.space} and we introduce the ball
\begin{equation}
	B^{\rho}(\delta) := \{ f \in X^{\rho,1} \, : \, \| f\|_{\rho,1} \leqslant \delta \}\,.
\end{equation}
\begin{pro}\label{prop-Functional-est}
	Let $k\in\N^*$, $\rho\geq 0$ and $s\geq s_0$. Let $\alpha_{0}>0$ sufficiently small and $\omega \in \R\setminus \{0\}$. There exists $\delta\in(0,1)$ such that the following assertions hold: 
   %
	\\[1mm]
	\noindent $(i)$
	The functional
	\begin{equation*}
		\begin{aligned}
			&\cG: X^{0,s-1} \times (X^{\rho,s}\cap B^{\rho}(\delta)) \times [-\alpha_{0},\alpha_{0}] \to X^{\rho,s-1}\,, 
		\end{aligned}
	\end{equation*}
	given by  \eqref{Func.G} and \eqref{F-compact}, is well-defined and satisfies the estimate
	\begin{equation}
		\| \cG(\tU,r;\alpha) \|_{\rho,s-1} \lesssim_{s} \| \tU \|_{0,s-1} + \| r \|_{\rho,s} + |\alpha|  \,; 
	\end{equation}
	\noindent $(ii)$ For any $\alpha\in[-\alpha_{0},\alpha_{0}]$, the linearized  operator 
    $$
    \di_{(\tU,r)} \cG( \overline\tU, \overline{r}; \alpha) : X^{0,s-1} \times X^{\rho,s} \to X^{\rho,s-1}
    $$
    takes the form
	\begin{equation}
		\begin{aligned}
			&\di_{(\tU,r)} \cG( \overline\tU, \overline{r}; \alpha) [\wh\tU,\whr] =  \omega \pa_{\vf}\whr +  \pa_{\theta}\big( \,\whU\sqrt{1+2 r} \cos(\theta) \big) + \pa_{\theta}\bigg( \Big( \tfrac{ \overline{U} \cos(\theta)}{\sqrt{1+2 \overline{r}}} + \tV[ \overline{r};\alpha] \Big) \whr +  \cR(\overline{r};\alpha) \whr \bigg) \,,
		\end{aligned}
	\end{equation}
	with $\tV[\overline r]$ and $\cR(\overline r)$ as in \eqref{V.at.boundary} and \eqref{cR.lin.app} of Proposition \ref{linea-prop}, In addition, 
	\begin{align*}
		\| \di_{(\tU,r)} \cG( \overline\tU, \overline{r}; \alpha) [\wh\tU,\whr] \|_{\rho,s-1} \lesssim_{s} &(1+\| \overline r\|_{\rho,1}) \| \wh\tU \|_{0,s-1} + \big(1+ \| \overline \tU \|_{0,0} + \| \overline r \|_{\rho,1} \big)\| \whr \|_{\rho,s} \\
        &+\| \overline r\|_{\rho,s} \| \wh\tU \|_{0,0} + \big( \| \overline \tU \|_{0,s-1} + \| \overline r \|_{\rho,s} \big)\| \whr \|_{\rho,1}\,;
	\end{align*}

	\noindent $(iii)$ The map $(\tU,r,\alpha) \mapsto \cG(\tU,r;\alpha)$ is $\cC^k$ with respect to $(\tU,r)$.
\end{pro}
\begin{proof}
    We extensively prove item $(i)$ only, since items $(ii)$ and $(iii)$ follow by straightforward and classical arguments, which we therefore omit. \\
    Recalling the definition of the functional $\cG$ in \eqref{Func.G}-\eqref{F-compact}, by Lemma \ref{lem-productlaw} we get
    \begin{align*}
		\big\| \omega \,\pa_{\vf} r + \pa_{\theta} \big( \tU(\vf) \sqrt{1+2 r(\vf,\theta)} \cos(\theta)\big) \big\|_{\rho,s-1} &\lesssim_{s} \| r \|_{\rho,s}+\| \tU \|_{0,s-1} \| \sqrt{1+2 r} \|_{\rho,0}\\
        &+\| \tU \|_{0,0}\| \sqrt{1+2 r} \|_{\rho,s}\\
        &\lesssim_{s} \| r \|_{\rho,s}+\| \tU \|_{0,s-1}+\| U \|_{0,0}(1+\| r \|_{\rho,s}) \,.
    \end{align*}
    Consequently
     \begin{align*}
		\big\| \omega \pa_{\vf} r + \pa_{\theta} \big( \tU(\vf) \sqrt{1+2 r(\vf,\theta)} \cos(\theta)\big) \big\|_{\rho,s-1} &\lesssim_{s} \| r \|_{\rho,s}+\| \tU \|_{0,s-1} \,.
    \end{align*}
    Now, we proceed to estimate the term $\pa_\theta F[r(\vf,\theta);\alpha]$  in \eqref{F-compact}, which requires a more careful and detailed analysis. The main difficulty arises from the first term in \eqref{F-compact}, representing the induced effect
    $$
    \mathcal{I}_{\rm self}[r](\varphi,\theta):= \fint_{\T} \log | R(\vf,\theta) e^{\im\theta} - R(\vf,\eta) e^{\im\eta} | \pa_{\theta,\eta}^2 \big[ R(\vf,\theta) R(\vf,\eta) \sin(\theta-\eta) \big] \wrt \eta\,,
    $$
    where we recall that $R(\vf,\theta)=\sqrt{1+2r(\vf,\theta)}$. In contrast, the second term 
    \begin{equation}\label{I.int}
        \cI_{\rm int}[r] (\vf,\theta) :=  -  \fint_{\T} K_\alpha ( R(\vf,\theta) e^{\im\theta}, R(\vf,\eta) e^{\im\eta} )|\pa_{\theta,\eta}^2 \big[ R(\vf,\theta) R(\vf,\eta) \sin(\eta+\theta) \big] \wrt \eta
    \end{equation}
    will be handled later straightforwardly, since the associated kernel does not exhibit any singularity. We split $\mathcal{I}_{\rm self}[r]$ as follows:
\begin{equation}\label{I.beria}
	    \begin{aligned}
		\mathcal{I}_{\rm self}[r](\vf,\theta)&= \fint_{\T} \log |  e^{\im\theta} -  e^{\im\eta} | \pa_{\theta,\eta}^2 \big[ R(\vf,\theta) R(\vf,\eta) \sin(\theta-\eta) \big] \wrt \eta\\&+
		 \fint_{\T} \log \left| \tfrac{R(\vf,\theta) e^{\im\theta} - R(\vf,\eta) e^{\im\eta}}{e^{\im\theta}-e^{\im\eta}} \right| \pa_{\theta,\eta}^2 \big[ R(\vf,\theta) R(\vf,\eta) \sin(\theta-\eta) \big] \wrt \eta\\
		&:= \mathcal{I}_1[r](\vf,\theta)+ \mathcal{I}_2[r](\vf,\theta)\,.
	\end{aligned}
\end{equation}
    By differentiating and making the change of variables $ \eta\to \eta+\theta$ we get
\begin{align*}
    \mathcal{I}_1[r](\vf,\theta)&=- \fint_{\T} \log |  1 -  e^{\im\eta} |\sin(\eta) d\eta- \fint_{\T} \log |  1 -  e^{\im\eta} |\, \mathcal{K}[r](\vf,\theta,\eta) \wrt \eta\\
    &= -\fint_{\T} \log |  \sin(\tfrac\eta2) |\, \mathcal{K}[r](\vf,\theta,\eta) \wrt \eta \,,
    \end{align*}
    where, denoting $R'(\vf,\theta):= (\pa_{\theta}R)(\vf,\theta)$,
    \begin{align*}
     \mathcal{K}[r](\vf,\theta,\eta)&:=\Big(R^\prime(\vf,\theta)R^\prime(\vf,\theta+\eta)+\big[R(\vf,\theta)R(\vf,\theta+\eta)-1\big]\Big)\sin(\eta)\\
     &+\Big(R^\prime(\vf,\theta)R(\vf,\theta+\eta)-R(\vf,\theta)R^\prime(\vf,\theta+\eta)\Big)\cos(\eta)\,.
    \end{align*}
    Therefore
    \begin{align*}
    \|\mathcal{I}_1[r]\|_{\rho,s-1}&\leqslant  -\fint_{0}^{2\pi}  \log |  \sin(\tfrac\eta2) |\, \|\mathcal{K}[r](\cdot,\eta)\|_{\rho,s-1} \wrt \eta\\
    &\leqslant C \sup_{\eta\in\T}\|\mathcal{K}[r](\cdot,\eta)\|_{\rho,s-1}\,,
    \end{align*}
    for some constant $C>0$.
    Applying Lemma \ref{lem-productlaw} we infer under the smallness condition that
     \begin{align}\label{est-K}
     \sup_{\eta\in\T}\|\mathcal{K}[r](\cdot,\eta)\|_{\rho,s-1}
     & \lesssim_{s} \|r\|_{\rho,s}\,,
    \end{align}
    which implies that
    \begin{align}\label{I1-est}
    \|\mathcal{I}_1[r]\|_{\rho,s-1}&\lesssim_{s} \|r\|_{\rho,s}\,.
    \end{align}
    We move now to the estimate of $\mathcal{I}_2[r]$. We first write by making a change of variables as before
    \begin{align*}
\mathcal{I}_2[r](\vf,\theta)=& -\fint_{\T} \log \left| 1+\tfrac{R(\vf,\theta)  - R(\vf,\theta+\eta) e^{\im\eta}}{1-e^{\im\eta}} \right|  \sin(\eta)\wrt \eta\\
&-\fint_{\T} \log \left| 1+\tfrac{R(\vf,\theta)  - R(\vf,\theta+\eta) e^{\im\eta}}{1-e^{\im\eta}} \right|  \mathcal{K}[r](\vf,\theta,\eta)\wrt \eta \,.
    \end{align*}
    By writing
    \begin{align*}
 \frac{R(\vf,\theta)  - R(\vf,\theta+\eta) e^{\im\eta}}{1-e^{\im\eta}}&=\sum_{\ell,j} R_{\ell,j}e^{\im(\ell\vf+j\theta)}  \frac{1-e^{\im (1+j)\eta}}{1-e^{\im\eta}} 
\end{align*}
and using 
$$
\sup_{\eta\in\T}\left|\frac{1-e^{\im (1+j)\eta}}{1-e^{\im\eta}}\right|\lesssim |1+j| \,,
$$
together with Lemma \ref{lem-productlaw} and the definition of the function spaces, we deduce that
\begin{align*}
 \sup_{\eta\in\T}\Big\|\tfrac{R(*,\cdot)  - R(*,\cdot+\eta) e^{\im\eta}}{1-e^{\im\eta}}\Big\|_{\rho,s-1}\lesssim_s \|r\|_{\rho,s} \,.  
\end{align*}
 Applying  the analytic expansion of $\log(1+x)$ around zero and the product law in Lemma \ref{lem-productlaw} with the previous estimate, we find, under the assumption $ \|r\|_{\rho,1}\leqslant \delta$,
 that
 \begin{align*}
 \Big\|\log\Big(1+\tfrac{R(*,\cdot)  - R(*,\cdot+\eta) e^{\im\eta}}{1-e^{\im\eta}}\Big)\Big\|_{\rho,s-1}&\lesssim_s \|r\|_{\rho,s}\sum_{n\in\N}\|r\|_{\rho,1}^n   \\
 &\lesssim_s \|r\|_{\rho,s} \,.
\end{align*}
 Putting this estimate together with \eqref{est-K} and the product estimate in Lemma \ref{lem-productlaw} we infer that
  \begin{align*}
    \|\mathcal{I}_2[r]\|_{\rho,s-1}& \lesssim_{s} \|r\|_{\rho,s} \,.
    \end{align*}
    Collecting this estimate with \eqref{I1-est} we deduce that $\cI_{\rm self}[r]$ in \eqref{I.beria} satisfies
    \begin{align*}
    \|\mathcal{I}_{\rm self}[r]\|_{\rho,s-1}&\lesssim_s \|r\|_{\rho,s} \,.
    \end{align*}
    Finally, we estimate $\cI_{\rm int}[r]$ in \eqref{I.int}. Noting that the kernel $K_{\alpha}$ in \eqref{Kalpha.kernel} is not singular, with similar arguments as above we get the estimate
    \begin{align}
    \| \cI_{\rm int}[r] \|_{\rho,s}|\lesssim_s |\alpha|+\|r\|_{\rho,s}\,.
   \end{align}
   This concludes the proof of item $(i)$ and of the proposition.
\end{proof}

\section{Traveling pairs  and asymptotics}\label{Traveling pairs}
The aim of this section is to prove Theorem \ref{Theo-1} by  constructing a branch of analytic traveling vortex pair solutions to the nonlinear equation \eqref{Func.G} in a neighborhood of the equilibrium state. More precisely, we shall establish the existence of a family of solutions parameterized by the small parameter $\alpha$, which measures the strength of the perturbation away from the trivial configuration. Our goal is to derive an explicit expansion of the solution with respect to $\alpha$, as this information will play a decisive role later when investigating the emergence of non–rigid periodic motions bifurcating from the branch of traveling pairs. The strategy is to view the governing equation \eqref{Func.G} as a nonlinear functional equation defined on a suitable Banach space of analytic boundaries, and then invoke the Implicit Function Theorem to locally parametrize the solution set.
To formulate our result precisely, we begin by introducing the appropriate functional framework.
For  $\rho>0$ and $s\geq 0$, we define the spaces
\begin{equation}\label{spaces.evenodd}
	\begin{aligned}
		X_{\geq 2, \rm even}^{\rho,s} &:= \big\{ f \in X^{\rho,s} \, : \, f= \Pi_{\geq 2} f \,, \ f(\theta)=f(-\theta) \big\} \,, \\
		X_{\geq 2, \rm odd}^{\rho,s} & := \big\{ f \in X^{\rho,s} \, : \, f= \Pi_{\geq 2} f \,, \ f(-\theta)=-f(\theta) \big\} \,,
	\end{aligned}
\end{equation}
where the operator $\Pi_{\geq 2}$ is introduced  in \eqref{proj.geq2}, and the functional  space $X^{\rho,s}$ is defined in Section \ref{Section-Func-sett}.
For $\delta>0$, we introduce  the closed ball
\begin{equation}\label{ball.even.small}
	B_\delta^{\rho,s} := \{ f \in X_{\geq 2, \rm even}^s \, : \, \| f\|_{\rho,s} \leq \delta \}\,.
\end{equation}
The main result of this section reads as follows.
\begin{pro}\label{prop.ralpha}
   Let $k\in\N^*$, $\rho>0$ and $s\geq 1$. There exist $\delta>0, \alpha_0>0$ and two curves 
   $$\alpha\in[-\alpha_0,\alpha_0]\mapsto r_\alpha\in B_\delta^{\rho,s}\,, \quad \alpha\in[-\alpha_0,\alpha_0]\mapsto \tU_\alpha\in\R \,,
   $$
   of  class $C^k$ such that
   $$
   \forall \,\alpha\in [-\alpha_0,\alpha_0]\,,\quad \cG(\tU_\alpha,r_\alpha;\alpha)=0 \,.
   $$
   The expansion of $r_\alpha$ is detailed in Proposition $\ref{prop.asympt.ralpha}.$
\end{pro}
The proof proceeds in several steps. We first expand the kernels involved in the definition of the functional $\cG$ for small $\alpha$, which provides the leading-order behavior of the problem. We then establish the regularity properties of $\cG$ and examine the structure of its linearization around the equilibrium configuration. With these ingredients in place, the Implicit Function Theorem can be applied, leading to the construction of the desired analytic branch of traveling vortex pair solutions, as detailed in Subsection \ref{Sec-proof-TP}.

\subsection{ Kernels  expansions}
For notational simplicity, we suppress the dependence on $\vf$ in the functions and we set
\begin{equation}
	R(\theta) = \sqrt{1+2 r(\theta)} = 1 + \tg(\theta) \,, \quad \textnormal{with} \quad \tg(\theta):= \sqrt{1+2 r(\theta)} - 1\,,
\end{equation}
for $r$ small enough. We begin by stating the following standard identity and include its full proof for the sake of completeness.
\begin{lem}\label{lemma-1-Id}
	Given a real function $\theta\in\R\mapsto \tg(\theta)$, we have that
	\begin{equation}
		\forall\, \theta\neq\eta,\quad \frac{\tg(\theta)e^{\im\theta}-\tg(\eta)e^{\im\eta}}{e^{\im\theta}-e^{\im\eta}} = \frac{\tg(\theta)+\tg(\eta)}{2} - \im\, \frac{\tg(\theta)-\tg(\eta)}{2} {\rm cotg}\Big(\frac{\theta-\eta}{2}\Big) \,. 
	\end{equation}
\end{lem}
\begin{proof}
	First, we compute
    \begin{align*}
		\frac{e^{\im\theta}}{e^{\im\theta}-e^{\im\eta}}& = \frac{e^{\im\theta}}{e^{\im\frac{\theta+\eta}{2}}(e^{\im\frac{\theta-\eta}{2}}-e^{-\im\frac{\theta-\eta}{2}})}=-\frac{\im}{2}\,\frac{e^{\im\frac{\theta-\eta}{2}}}{\sin(\tfrac{\theta-\eta}{2})}\\
        &=\frac12-\frac{\im}{2}{\rm cotg}(\tfrac{\theta-\eta}{2}) \,.
	\end{align*}
	By a permutation argument, we also get
    \begin{align*}
		\frac{e^{\im\eta}}{e^{\im\theta}-e^{\im\eta}}
        &=-\frac12-\frac{\im}{2}{\rm cotg}(\tfrac{\theta-\eta}{2}) \,.
	\end{align*}
    Therefore, the desired claim follows easily by summation of the above contributions.
\end{proof}
The next result establishes an expansion that will play a central role in describing the asymptotic structure of the traveling pairs.
\begin{lem}\label{lemma.kernels}
	Given a smooth real function $\theta\in\R\mapsto \tg(\theta)$ with a small amplitude, then we have
	\begin{equation}\label{BS.kernel.expand}
		\begin{aligned}
			\log \big| (1+\tg(\theta))e^{\im\theta} & - (1+\tg(\eta))e^{\im\eta} \big|  = \log| 1-e^{\im(\eta-\theta)} |  + \tfrac{\tg(\theta)+\tg(\eta)}{2} \\
			& -\tfrac12 \Big( \tfrac{\tg(\theta)+\tg(\eta)}{2} \Big)^2 +\tfrac12 \Big( \tfrac{\tg(\theta)-\tg(\eta)}{2} {\rm cotg} \Big(\tfrac{\theta-\eta}{2} \Big) \Big)^2 + O({|\tg|^3}+|\tg^\prime|^3) \,. 
				\end{aligned}
		\end{equation}
		Moreover, for $\alpha$ small enough, we have the expansion
		\begin{equation}\label{inter.kernel}
			\begin{aligned}
			K_\alpha\big( & (1+\tg(\theta))e^{\im\theta} ,  (1+\tg(\eta))e^{\im\eta} \big)  = - \tfrac{\alpha^2}{2}\big[  \cos(2\theta) +2 \cos(\theta-\eta) + \cos(2\eta) \big] \\
			& \ \ +\tfrac{\alpha^3}{3} \big[ \cos(3\theta) + 3\cos(2\theta-\eta) + 3 \cos(\theta-2\eta) + \cos(3\eta) \big] \\
			& \ \ - \tfrac{\alpha^4}{4} \big[ \cos(4\theta) + 4 \cos(3\theta-\eta) + 6\cos(2(\theta-\eta))  + 4\cos(\theta-3\eta) + \cos(4\eta)\big] \\
			& \ \ -\alpha^2\big[ \tg(\theta)\cos(2\theta)  + (\tg(\theta)+\tg(\eta)) \cos(\theta-\eta) + \tg(\eta)\cos(2\eta) \big] \\
			& \ \ +O\big(|\alpha|^5+|\alpha|^3 |\tg| +|\alpha|^2 | \tg|^2 \big) \,.
		\end{aligned}
	\end{equation}
\end{lem}
\begin{proof}
   We begin with the proof of \eqref{BS.kernel.expand}. Observe first the elementary decomposition
    \begin{align*}
\log \big| (1+\tg(\theta))e^{\im\theta}  - (1+\tg(\eta))e^{\im\eta} \big|  &= \log| 1-e^{\im(\eta-\theta)} | +\log \Big|1+\tfrac{\tg(\theta)e^{\im\theta}  - \tg(\eta)e^{\im\eta}}{e^{\im\theta}-e^{\im\eta}} \Big|	\,.
    \end{align*}
    For the second factor, we use the expansion near the origin, as in Lemma \ref{loga.expa},
    $$
    \log|1+z|={\rm Re}(z)-\tfrac12{\rm Re}(z^2)+O(z^3)\,,$$
    with $z=\frac{\tg(\theta)e^{\im\theta}  - \tg(\eta)e^{\im\eta}}{e^{\im\theta}-e^{\im\eta}}$, and invoking Lemma \ref{lemma-1-Id}, it yields
   \begin{align*}
 \log \Big|1+&\tfrac{\tg(\theta)e^{\im\theta}  - \tg(\eta))e^{\im\eta}}{e^{\im\theta}-e^{\im\eta}} \Big| =\tfrac{\tg(\theta)+\tg(\eta)}{2} \\
 &-\tfrac12{\rm Re}\Big(\tfrac{\tg(\theta)+\tg(\eta)}{2} - \im\, \tfrac{\tg(\theta)-\tg(\eta)}{2} {\rm cotg}\Big(\tfrac{\theta-\eta}{2}\Big)\Big)^2+O\bigg(\Big(\tfrac{\tg(\theta)e^{\im\theta}  - \tg(\eta)e^{\im\eta}}{e^{\im\theta}-e^{\im\eta}}\Big)^3\bigg)\\
 &=\tfrac{\tg(\theta)+\tg(\eta)}{2}  -\tfrac12 \Big( \tfrac{\tg(\theta)+\tg(\eta)}{2} \Big)^2 +\tfrac12 \Big( \tfrac{\tg(\theta)-\tg(\eta)}{2} {\rm cotg} \Big(\tfrac{\theta-\eta}{2} \Big) \Big)^2+ O(|\tg|^3+|\tg^\prime|^3)\,,
    \end{align*} 
    where the last term has been estimated by the classical Mean Value Theorem.
    This proves the identity in \eqref{BS.kernel.expand}.
  We now prove the expansion in \eqref{inter.kernel}. In particular, by \eqref{Kalpha.kernel}, we may expand the non-singular kernel as follows
  \begin{equation}
  	\begin{aligned}
  		K_{\alpha}\big((1+\tg(\theta))e^{\im\theta},(1+\tg(\eta))e^{\im\eta}\big)  = \sum_{n\geq 2} \tfrac{(-1)^{n+1}\alpha^n}{n} {\rm Re} \Big( \big( (1+\tg(\theta))e^{\im\theta} + (1+\tg(\eta))e^{-\im\eta} \big)^n \Big) \,.
  	\end{aligned}
  \end{equation}
  From the binomial identity we infer
  \begin{equation}
  	\begin{aligned}
  		\big( (1+\tg(\theta))e^{\im\theta} & + (1+\tg(\eta))e^{-\im\eta} \big)^n \\ & = (e^{\im\theta}+e^{-\im\eta})^n +  n (e^{\im\theta}+e^{-\im\eta})^{n-1} \big( \tg(\theta)e^{\im\theta} +\tg(\eta)^{-\im\eta}  \big) + O(|\tg|^2)\,,
  	\end{aligned}
  \end{equation}
  from which we get
    \begin{equation}
  	\begin{aligned}
  		K_{\alpha}&\big((1+\tg(\theta))e^{\im\theta} ,(1+\tg(\eta))e^{\im\eta}\big) \\
  		&   = - \tfrac{\alpha^2}{2}{\rm Re}\big( (e^{\im\theta}+e^{-\im\eta})^2 \big)  + \tfrac{\alpha^3}{3}{\rm Re}\big( (e^{\im\theta}+e^{-\im\eta})^3 \big) - \tfrac{\alpha^4}{4}{\rm Re}\big( (e^{\im\theta}+e^{-\im\eta})^4 \big)  \\
  		& \quad - \alpha^2 {\rm Re} \Big( (e^{\im\theta}+e^{-\im\eta}) \big( \tg(\theta)e^{\im\theta} +\tg(\eta)^{-\im\eta}  \big) \Big) + O(|\alpha|^5) + O(|\alpha|^3 |\tg|) + O(|\alpha|^2 |\tg|^2) \,.
  	\end{aligned}
  \end{equation}
  Explicit computations of these terms lead to \eqref{inter.kernel} and the proof is concluded.
\end{proof}
The following lemma is classical, see for instance \cite[Lemma 3.3]{Casrro}. It provides the Fourier representation of a convolution operator that arises in the linearization of the system around the equilibrium state.
\begin{lem}\label{lemma.cR.00}
For any $j\in \mathbb{Z}^\star,$ we have the following identity
$$
 \fint_{\T} \log\big| 1-  e^{\im\eta} \big| \cos( j\eta)  \wrt\eta=-\frac{1}{2|j|} \quad \text{and} \quad  \fint_{\T} \log\big| 1-  e^{\im\eta} \big| \sin( j\eta)  \wrt\eta=0\,.
$$
Moreover, for any smooth periodic function $h:\mathbb{T}\to \C$ with zero average, we define the operator
    $$
    \cR(0;0)h(\theta):=  \fint_{\T} \log|1-e^{\im(\eta-\theta)}| h(\eta) \wrt \eta \,.
    $$
    Then we have that
	\begin{equation}
		h(\theta)= \sum_{j\in\Z^\star} \whh(j) e^{\im j\theta} \quad  \Longrightarrow \quad  \cR(0;0)h(\theta)=- \sum_{j\in\Z^\star} \tfrac{\whh(j)}{2|j|} e^{\im j \theta} \,.
	\end{equation}
\end{lem}

\subsection{Regularity and linearization}
 We shall look for stationary  solutions to \eqref{Func.G} which solves the nonlinear equation
\begin{equation}\label{G0_station}
	\begin{aligned}
		\cG_{0}(\tU,r(\theta);\alpha) & := \pa_{\theta}\Big(   \tU  \sqrt{1+2 r(\theta)}\cos(\theta) + F[ r(\theta);\alpha ] \Big) = 0 \,,
	\end{aligned}
\end{equation}
with   $F[ r ;\alpha]$ as in  \eqref{F.app} and \eqref{F-compact}. Applying \eqref{F-compact} and Lemma \ref{lemma.cR.00}, we infer
\begin{align}\label{F1.onehalf}
\nonumber \pa_\theta F[r;\alpha] (\vf,\theta)
		& =    \fint_{\T} \log | e^{\im\theta} -  e^{\im\eta} |  \sin(\theta-\eta)  \wrt \eta\\
        &=- \fint_{\T} \log | 1 -  e^{\im\eta} |  \sin(\eta)  \wrt \eta\\
        &=0 \,.
        \end{align}
This implies that
\begin{equation}\label{G0.onehalf}
	\cG_{0}(0, 0 ; 0)=0 \,.
\end{equation}
The goal is to identify solutions to \eqref{G0_station} that bifurcate smoothly from the trivial solution $(\tU_{0}, r_0(\theta)) = (0,0)$ at $\alpha = 0$. To achieve this, we will apply the Implicit Function Theorem within the small parameter regime $|\alpha| \ll 1$. This approach allows us to rigorously establish the existence and uniqueness of nontrivial solutions near the trivial equilibrium.

We begin by analyzing the regularity properties of the nonlinear functional $\cG_{0}$ defined in \eqref{G0_station}. Recalling \eqref{spaces.evenodd}-\eqref{ball.even.small}, we have to ensure that $\cG_{0}$ is a map between the spaces $X_{\geq 2,{\rm even}}^{\rho,s}$ and $X_{\geq 2,{\rm odd}}^{\rho,s}$ with respect to $r$: we achieve this by constraining the speed $\tU$ as a function of $r\in B_{\delta}^{\rho,s}$.

\begin{pro}\label{lemma.constrain}
	Let $k\in\N^\star, \rho>0$ and $s\geq 1$. There exist $\delta>0$ and  $\alpha_{0}>0$ such that the following results  hold:
	\\[1mm]
	\noindent ${(i)}$
	The nonlinear functional
	\begin{equation}\label{fixed.constrain}
		B_\delta^{\rho,s} \times [-\alpha_{0},\alpha_{0}] \ni (r,\alpha) \mapsto \tU[r,\alpha] := - \frac{\int_{\T}F[r;\alpha](\theta)\cos(\theta) \wrt \theta}{\int_{\T}\sqrt{1+2 r(\theta)}\cos^2(\theta) \wrt\theta}\in \R
	\end{equation}
	is of class  {$\cC^k$.} In addition, the linear operator $\di_{r} \tU[\overline r;\alpha]: X_{\geq 2,\rm even}^{\rho,s} \to \R$  acts as
	\begin{align*}
		\di_{r} \tU[\overline r;\alpha] h &=  \frac{ \, \int_{\T}F[\overline r;\alpha](\theta) \cos(\theta) \wrt \theta  }{\big( \int_{\T} \sqrt{1+2\overline r} \cos^2(\theta)\wrt \theta \big)^2}\int_{\T} \frac{h(\theta)}{\sqrt{1+2\overline r}(\theta) }\cos^2(\theta)\wrt \theta \\&\qquad- \frac{\int_{\T}\big( \tV[\overline r;\alpha]h + \cR(\overline r;\alpha)[h](\theta) \big) \cos(\theta) \wrt \theta}{\int_{\T}\sqrt{1+2 \overline{r}(\theta)}\cos^2(\theta) \wrt\theta} \,;
	\end{align*}
	\noindent $(ii)$
 The associated composition map
	\begin{equation}\label{G0.reduced.def}
		\wt\cG_{0}: B_\delta^{\rho,s} \times [-\alpha_{0},\alpha_{0}] \to X_{\geq 2,\rm odd}^{\rho,s-1} \,, \quad \wt\cG_{0}(r;\alpha) := \cG_{0}\big( \tU[r;\alpha], r;\alpha \big)\,,
	\end{equation}
	is  of class $\cC^k$.
	
\end{pro}
\begin{rem}
The choice of the function $\tU$ imposed in the preceding proposition will become transparent in Proposition $\ref{lemma.constrain1}$. It plays the role of a Lagrange multiplier ensuring that the first sine mode is eliminated from the nonlinear functional. This orthogonality constraint is crucial for removing the degeneracy generated by translation invariance and for implementing the Implicit Function Theorem in a suitable functional setting.
\end{rem}
\begin{proof}
    ${{(i)}}$ The proof follows directly by the explicit definition of $\tU$ in \eqref{fixed.constrain} together with Proposition \ref{prop-Functional-est}.
    \\[1mm]
   ${{(ii)}}$ The regularity properties of the functional $\wt\cG_{0}$ follow directly from Proposition \ref{prop-Functional-est}-$(i)$. The remaining task is to verify the symmetry condition and to ensure that the range of this functional does not contain the mode-one component. The latter is guaranteed by the appropriate choice of the speed $\tU$ in \eqref{fixed.constrain}.

   By \eqref{ball.even.small} and \eqref{spaces.evenodd}, any function $r\in B_\delta^{\rho,s}$ is even in $\theta$. Therefore, from the structure of $F$ detailed in \eqref{polar.integrals} we deduce that
   $$
   F[r,\alpha](-\theta)=F[r,\alpha](\theta)\quad\hbox{and}\quad \tU[r,\alpha](-\theta)=\tU[r,\alpha](\theta) \,.
   $$
   Therefore, we obtain  that
   $$
   \cG_0(\tU,r,\alpha)(-\theta)=-\cG_0(\tU,r,\alpha)(\theta)\quad\hbox{and}\quad \wt\cG_0(r,\alpha)(-\theta)=-\wt\cG_0(r,\alpha)(\theta) \,.$$
   Regarding the absence of the mode-one component, it is enough to show that
   $$
   \int_{\mathbb{T}}\wt\cG_0(r,\alpha)(\theta)\sin(\theta)\wrt\theta=0 \,.
   $$
   According to \eqref{G0.reduced.def} and \eqref{G0_station}, this reduces  to checking that
   \begin{equation*}
	\begin{aligned}
	\int_{\mathbb{T}} \pa_{\theta}\Big(   \tU  \sqrt{1+2 r(\theta)}\cos(\theta) + F[ r(\theta);\alpha ] \Big) \sin(\theta) \wrt\theta = 0 \,.
	\end{aligned}
\end{equation*}
 By carrying out an integration by parts, this is equivalent to the identity
 $$
 \tU\int_{\T}\sqrt{1+2 r(\theta)}\cos^2(\theta) \wrt\theta =-\int_{\T}F[r(\theta);\alpha]\cos(\theta) \wrt \theta \,,$$
which is satisfied by the choice of $\tU$ in \eqref{fixed.constrain}. This concludes the proof.
\end{proof}
The next lemma establishes the asymptotic expansions in the parameter $\alpha$, which are needed for the later asymptotic expansions.
\begin{lem}\label{Id-lem-int}
    For $|\alpha|<1$, the following identities hold:
   \\[1mm]
   \noindent $(i)$ $ \wt\cG_0(0,\alpha)=\pa_\theta F[0;\alpha] 
   =\tfrac{1}{2}\sum_{n\geqslant 2} \alpha^{n}(-1)^{n+1}  \sin(n \theta)$;
   \\[1mm]
   \noindent $(ii)$ $ \tU[0,\alpha]=0$;
   \\[1mm]
   \noindent $(iii)$ $ \tV[0;\alpha](\theta) = \tfrac12+\tfrac{1}{2}\sum_{n\geq 2} (-\alpha)^{n}  \cos(n \theta)$;
\\[1mm]
\noindent $(iv)$ For any smooth  real function  $h(\theta)=\sum_{m\geqslant2}h_m \cos(m\theta)$,
\begin{align*}
		\cR(0;\alpha) h(\theta)& = -\sum_{m\geqslant2}\frac{1}{2m}h_m \cos(m\theta) + \sum_{n\geq m\atop
        m\geqslant2} (-\alpha)^{n}\frac{1}{2n}\binom{n}{m} h_m \cos\big((n-m) \theta\big)\,.
	\end{align*}
\end{lem}
\begin{proof}
We start with proving the identity in item $(i)$. Using \eqref{F-compact} and \eqref{Kalpha.kernel}  we infer
    \begin{align}\label{F-compactDD}
\nonumber \pa_\theta F[0;\alpha] 
		& =    \fint_{\T} \log |  e^{\im\theta} -  e^{\im\eta} |  \sin(\theta-\eta)  \wrt \eta +  \fint_{\T} K_\alpha (  e^{\im\theta},  e^{\im\eta} ) \sin(\eta+\theta)  \wrt \eta\\
        &=
\sum_{n\geq 2} \alpha^{n}\tfrac{(-1)^{n+1}}{n} {\rm Re}\fint_{\T} \big( e^{\im\theta}+e^{-\im\eta} \big)^n \sin(\eta+\theta)  \wrt \eta\,, 
\end{align}
where we have used, according to Lemma \ref{lemma.cR.00},
$$
\int_{\T} \log |  e^{\im\theta} -  e^{\im\eta} |  \sin(\theta-\eta)  \wrt \eta =0 \,.
$$
Using the binomial formula and making a change of variables yield
\begin{align}\label{binom.RE}
 {\rm Re}\fint_{\T}\big( e^{\im\theta}+e^{-\im\eta} \big)^n \sin(\eta+\theta)  \wrt \eta&= n{\rm Re}\Big(e^{\im n\theta}\fint_{\T}e^{-\im  \eta}  \sin(\eta)  \wrt \eta\Big) =\frac{n}{2}\sin(n\theta) \,. 
\end{align}
Plugging this formula into  \eqref{F-compactDD} we deduce that
\begin{align*}
 \pa_\theta F[0;\alpha] 
        &=\frac{1}{2}\sum_{n\geq 2} \alpha^{n}(-1)^{n+1}  \sin(n \theta) \,.
\end{align*}
To prove  the identity in item $(ii)$, we apply Proposition \ref{lemma.constrain} and item $(i)$, deducing that
\begin{align*}
    \tU[0,\alpha] &=\frac{1}{\pi} \int_{\T}\partial_\theta F[0;\alpha]\sin(\theta) \wrt \theta =0 \,.
\end{align*}
We now prove the identity in item $(iii)$. From \eqref{V.at.boundary},  we deduce that
\begin{align*}
			\nonumber \tV[0;\alpha](\theta) =& -  \fint_{\T} \log\big| e^{\im\theta}-  e^{\im\eta} \big| \cos(\eta-\theta)  \wrt\eta 
			-   \fint_{\T} K_\alpha\big( e^{\im\theta},  e^{\im\eta} \big)   \cos(\theta+\eta)  \wrt \eta  \,.
		\end{align*}
        By a change of variables and Lemma \ref{lemma.cR.00}, we get
        \begin{align*}
        - \fint_{\T} \log\big| e^{\im\theta}-  e^{\im\eta} \big| \cos(\eta-\theta)  \wrt\eta&=- \fint_{\T} \log\big| 1-  e^{\im\eta} \big| \cos(\eta)  \wrt\eta =\frac12 \,.
        \end{align*}
        Proceeding as for \eqref{F-compactDD} and \eqref{binom.RE}, we arrive at
         \begin{align*}
\nonumber   \fint_{\T} K_\alpha (  e^{\im\theta},  e^{\im\eta} ) \cos(\eta+\theta)  \wrt \eta
        &=\sum_{n\geq 2} \alpha^{n}\tfrac{(-1)^{n+1}}{n} {\rm Re}\fint_{\T} \big( e^{\im\theta}+e^{-\im\eta} \big)^n \cos(\eta+\theta)  \wrt \eta\\
        &=\frac{1}{2}\sum_{n\geq 2} \alpha^{n}(-1)^{n+1}  \cos(n \theta)\,. 
\end{align*}
Putting together the preceding identities allows us to get
\begin{equation}
			\nonumber \tV[0;\alpha](\theta) = \frac12+\frac{1}{2}\sum_{n\geq 2} (-\alpha)^{n}  \cos(n \theta)\,.
		\end{equation}
Finally, to establish the identity in item $(iv)$, we recall the definition of $\cR(\overline{r};\alpha)$ given \mbox{in \eqref{cR.lin.app},} from which we obtain
\begin{align}\label{cabo1}
		\cR(0;\alpha) h(\vf,\theta)& = \cR(0;0) h(\vf,\theta)   -\fint_{\T} K_{\alpha}\big( e^{\im\theta}, e^{\im\eta} \big) h(\vf,\eta) \wrt \eta\,.
	\end{align}
    For each $m\geqslant 2$, we have
    \begin{align*}
\nonumber   \fint_{\T} & K_\alpha (  e^{\im\theta},  e^{\im\eta} ) \cos(m\eta)  \wrt \eta
        =\sum_{n\geq 2} \alpha^{n}\tfrac{(-1)^{n+1}}{n} {\rm Re}\fint_{\T} \big( e^{\im\theta}+e^{-\im\eta} \big)^n \cos(m\eta)  \wrt \eta\\
        &=\frac{1}{2}\sum_{n\geq 2} \alpha^{n}\tfrac{(-1)^{n+1}}{n} {\rm Re}\fint_{\T} \big( e^{\im\theta}+e^{-\im\eta} \big)^n \big(e^{\im m \eta}+e^{-\im m \eta}\big)  \wrt \eta\,. 
\end{align*}
At this stage, the application of the binomial formula yields
\begin{align*}
\nonumber   \fint_{\T}  K_\alpha (  e^{\im\theta},  e^{\im\eta} ) \cos(m\eta)  \wrt \eta
        & =\frac{1}{2}\sum_{n\geq 2} \alpha^{n}\tfrac{(-1)^{n+1}}{n} \sum_{j=0}^{n} \binom{n}{j} {\rm Re}\Big( e^{\im j\theta} \fint_{\T} e^{\im(j-n)\eta} \big(e^{\im m \eta}+e^{-\im m \eta}\big)  \wrt \eta \Big) \\
        &=\sum_{n\geq m} \alpha^{n}(-1)^{n+1}\frac{1}{2n} \binom{n}{m}  \cos\big((n-m) \theta\big)\,. 
\end{align*}
Inserting it into \eqref{cabo1} and applying Lemma \ref{lemma.cR.00}, we get the claimed identity in item $(iv)$ for any smooth real function $h(\theta)=\sum_{m\geqslant2}h_m \cos(m\theta)$. This concludes the proof.
\end{proof}
Our next step is to examine the linearization of the functional $\wt\cG_{0}$ in the regime of small perturbations.
\begin{pro}\label{lemma.constrain1}
	Let $\rho>0$ and $s\geq 1$, and consider the functionals $\wt\cG_{0}$ and $\tU[r;\alpha]$ as defined in Proposition $\ref{lemma.constrain}.$ Then, there exist $\delta>0$ and  $\alpha_{0}>0$ such that the following holds:
	\\[1mm]
	$(i)$ 
	 For any $(\overline r,\alpha) \in B_\delta^{\rho,s} \times [-\alpha_{0},\alpha_{0}]$, the linearized operator $$\di_{r} \wt\cG_{0}(\overline r;\alpha):X_{\geq 2,\rm even}^{\rho,s}\to X_{\geq 2,\rm odd}^{\rho,s-1} 
     $$
     takes the form\begin{equation}\label{linear.with.constrain}
		\di_{r} \wt\cG_{0}(\overline r;\alpha) = \pa_{\theta}\circ \bigg(  \sqrt{1+2\overline r(\theta)} \cos(\theta) \, \di_{r}\tU[\overline r;\alpha] + \frac{\tU[\overline r;\alpha] \, \cos(\theta)}{\sqrt{1+2\overline r(\theta)}} + \tV[\overline r;\alpha] + \cR(\overline r;\alpha) \bigg) \,,
	\end{equation}
	 where $\tU[\overline r;\alpha]$ and the operator $\di_{r}\tU[\overline r;\alpha]$ are defined in Proposition $\ref{lemma.constrain},$ the function $\tV[\overline r;\alpha]$ is defined by \eqref{V.at.boundary} and the operator $\cR(\overline r;\alpha)$ is given by \eqref{cR.lin.app};
	\\[1mm]
	$(ii)$ The linear operator $\di_{r} \tU[0;\alpha]: X_{\geq 2,\rm even}^{\rho,s} \to \R$  acts as follows, 
	\begin{align*}
	h(\theta)=\sum_{m\geqslant2}h_m \cos(m\theta)\Longrightarrow 	\di_{r} \tU[0;\alpha] h 
        &{=-\tfrac14{\alpha^2}h_3  -\tfrac{1}{4}\sum_{m\geq 3} (-\alpha)^{m}\big(h_{m+1}+3h_{m-1}\big)\,;}
	\end{align*}
    $(iii)$ The speed $\tU[\overline{r};\alpha]$ satisfies  $$	\tU[\overline{r};\alpha]=O\big(\alpha^2\|\overline{r}\|_{L^\infty}+\|\overline{r}\|_{L^\infty}^2\big)\,;$$
	  $(iv)$ At the  equilibrium state, the linearized operator is an isomorphism   taking the form,
	\begin{align}
		\di_{r} \wt\cG_{0}(0;0) h & = \pa_{\theta}\circ \Pi_{\geq 2}\Big(  \tfrac12 h  + \cR(0;0)h \Big) \,, \label{bristol.0}
	\end{align}  
    where the operator $\cR(0;0)$ is described in Lemma $\ref{Id-lem-int}.$.  In particular, for any $j\in\Z$ with $|j|\geq 2$, we have that
   {
	\begin{equation}\label{dir.wtcG0}
		\di_{r} \wt\cG_{0}(0;0) [e^{\im j\theta}] =  \tfrac{\im\, j}{2|j|} (|j|-1) e^{\im j\theta} \,, \quad \big( \di_{r} \wt\cG_{0}(0;0) \big)^{-1} [e^{\im j \theta}] = - \tfrac{2 \im\, |j|}{j} \tfrac{1}{|j|-1} e^{\im j\theta}\,.
	\end{equation}}
\end{pro}
\begin{proof}
The linearization formula for $\wt\cG_{0}$ stated in item $(i)$ follows directly from Proposition \ref{prop-Functional-est}, combined with the chain rule, the linearization of the modulation parameter $\tU$ given in Proposition \ref{lemma.constrain}, and the linear expansion of the nonlinear term $F$ described in Proposition \ref{prop-Functional-est}.
\\
We now prove item $(ii)$.
    By Proposition \ref{lemma.constrain}, integration by parts and Lemma \ref{Id-lem-int}-$(i)$, we deduce that
    \begin{align}\label{alme1}
		\di_{r} \tU[0;\alpha] h
        & = - \frac{1}{\pi}{\int_{\T}\big( \tV[0;\alpha]h + \cR(0;\alpha)[h] \big) \cos(\theta) \wrt \theta} \,.
	\end{align}
For any smooth  real function 
        $\displaystyle{h(\theta)=\sum_{m\geqslant2}h_m \cos(m\theta)}$, we have, by Lemma \ref{Id-lem-int}-$(iii)$,
    \begin{align*}
 - \frac{1}{\pi}{\int_{\T}\big( \tV[0;\alpha]h \big) \cos(\theta) \wrt \theta}&=- \frac{1}{2\pi} \int_{\T}h(\theta)\cos(\theta)\wrt\theta -\frac{1}{2\pi}\sum_{n\geq 2} (-\alpha)^{n}  \int_{\T}h(\theta)\cos(n \theta)\cos(\theta)\wrt \theta \\
 &=-\frac{1}{4\pi}\sum_{n\geq 2} (-\alpha)^{n}  \int_{\T}h(\theta)\big(\cos((n+1) \theta)+\cos((n-1)\theta)) \big)\wrt \theta\\
 & =-\frac{\alpha^2}{4}h_3  -\frac{1}{4}\sum_{m\geq 3} (-\alpha)^{m}\big(h_{m+1}+h_{m-1}\big)  \,,
\end{align*}
and, by Lemma \ref{Id-lem-int}-$(iv)$,
\begin{align*}
		- \frac{1}{\pi}{\int_{\T}(\cR(0;\alpha)[h] ) \cos(\theta) \wrt \theta}& = -\frac{1}{\pi}\sum_{n\geq m\atop
        m\geq2} (-\alpha)^{n}\frac{1}{2n}\binom{n}{m} h_m \int_{\T}\cos\big((n-m) \theta\big)\cos(\theta)\wrt \theta\\
       &= -\frac{1}{2}\sum_{
        m\geqslant2} (-\alpha)^{m+1} h_m = -\frac{1}{2}\sum_{
        m\geqslant3} (-\alpha)^{m} h_{m-1} \,.
	\end{align*}
    Putting together the previous identities into \ref{alme1} yields
    \begin{align}\label{linera-drU}
		\nonumber\di_{r} \tU[0;\alpha] h &=-\frac{\alpha^2}{4}h_3  -\frac{1}{4}\sum_{m\geq 3} (-\alpha)^{m}\big(h_{m+1}+h_{m-1}\big)-\frac{1}{2}\sum_{
        m\geqslant3} (-\alpha)^{m} h_{m-1}\\
        &=-\frac{\alpha^2}{4}h_3  -\frac{1}{4}\sum_{m\geq 3} (-\alpha)^{m}\big(h_{m+1}+3h_{m-1}\big)\,.
	\end{align}
For the asymptotic in item $(iii)$, we use   Taylor expansion with respect to $r$, leading to
$$
\tU[\overline{r};\alpha]=\tU[0;\alpha]+\di_{r}\tU[0;\alpha][\overline{r}]+\int_0^1(1-s)\di_{r}^2\tU[s\overline{r};\alpha][\overline{r},\overline{r}]\wrt s \,.
$$
Then, by applying  Lemma \ref{Id-lem-int}-$(ii)$ and Proposition \ref{lemma.constrain}, we find that
$$
\tU[\overline{r};\alpha]=O\big(\alpha^2\|\overline{r}\|_{L^\infty}+\|\overline{r}\|_{L^\infty}^2\big) \,.
$$
Finally, we will move to the proof of  item $(iv)$. To obtain the identity  \eqref{bristol.0}, we first note that, by
\eqref{linear.with.constrain} and Lemma \ref{Id-lem-int}-$(ii)$, we are led to
\begin{align*}
	\di_{r} \wt\cG_{0}(0;\alpha) h = \pa_{\theta} \Big(   \cos(\theta) \, \di_{r}\tU[0;\alpha]h  + \tV[0;\alpha]\,h  + \cR(0;\alpha)h \Big) \,.
\end{align*}   
	By Lemma \ref{Id-lem-int}-$(iii)$,$(iv)$ and  Lemma \ref{lemma.cR.00},  we have that
	\begin{equation}
		\tV[0;0] = \tfrac12 \,, \quad \cR(0;0)[e^{\im j \theta}] =- \tfrac{1}{2|j|} \quad \forall\,j\in\Z\setminus\{0\}\,,
	\end{equation}
and, by item $(ii)$, we have $\di_{r}\tU[0;0]=0$. Now, it is straightforward to check that 
$$\di_{r} \wt\cG_{0}(0;0):X_{\geq 2,\rm even}^{\rho,s}\to X_{\geq 2,\rm odd}^{\rho,s-1}$$
defines an isomorphism. This concludes the proof. 
\end{proof}
We also need  to discuss the second variation of the functional $\wt\cG_{0}$.
{
\begin{pro}\label{lemma.constrain22}
	Under the assumptions of Proposition $\ref{lemma.constrain},$ we have
	\begin{align}\label{pula2}
		\di_{r}^2 \wt\cG_{0}(0;0)[h,h]& = -\tfrac12\Pi_{\geqslant2}\partial_\theta(h^2)\,.
	\end{align}
\end{pro}
}
\begin{proof}
We recall from Proposition \ref{lemma.constrain}-$(i)$ that
\begin{align*}
		\di_{r} \wt\cG_{0}(\overline r;0)h =& \pa_{\theta} \bigg(  \sqrt{1+2\overline r(\theta)} \cos(\theta) \, \di_{r}\tU[\overline r;0]h + \frac{\tU[\overline r;0]\, \cos(\theta)}{\sqrt{1+2\overline r(\theta)}}h(\theta)  \\
        &+ \tV[\overline r;0]h(\theta) + \cR(\overline r;0)h \bigg) \,.
	\end{align*}
    Differentiating one more time with respect to $r$ and using Lemma \ref{Id-lem-int}-$(ii)$, we infer
    
\begin{align*}
	\di_{r}^2 \wt\cG_{0}(0,0)[h,h]= \pa_{\theta} &\Big(  2 h(\theta)\cos(\theta) \, \di_{r}\tU[0;0]h +\cos(\theta)\di_{r}^2\tU[0;0][h,h]  \\
	&+ \big(\di_{r}\tV[0;0]h\big) h(\theta) + \di_{r}\cR(0;0)[h,h] \Big) \,.
\end{align*}
  Using items $(ii)$ and $(iii)$ of Proposition  \ref{lemma.constrain1}, we deduce that
    \begin{align}\label{1102}
		\di_{r}^2 \wt\cG_{0}(0;0)[h,h]& =-\sin(\theta)\di_{r}^2\tU[0;0][h,h]\\
        &+\pa_{\theta} \Big(  \big(\di_{r}\tV[0;0]h\big) h(\theta) + \di_{r}\cR(0;0)[h,h]  \Big)\,.
	\end{align}
    By construction, we have
    $$
    \di_{r}^2 \wt\cG_{0}(0;0)[h,h]=\Pi_{\geqslant2}\left(\di_{r}^2 \wt\cG_{0}[0;\alpha][h,h]\right).
    $$
    Therefore, we get from \eqref{1102}
    \begin{align}\label{112}
		\di_{r}^2 \wt\cG_{0}(0;0)[h,h]
        &=\Pi_{\geqslant2}\pa_{\theta} \Big(  \big(\di_{r}\tV[0;0]h\big) h(\theta) + \di_{r}\cR(0;0)[h,h]  \Big)\,.
	\end{align}
   Invoking \eqref{V.at.boundary} and Lemma \ref{lemma.kernels}, and noting  that $\overline R(\theta)= \sqrt{1+2\overline r(\theta)}$, we find
	\begin{align*}
			\nonumber \tV[\overline r;\alpha](\theta) = - \frac{1}{ \overline R(\theta) } &\fint_{\T} \log\big|\overline R(\theta) e^{\im\theta}- \overline R(\eta)e^{\im\eta} \big| \pa_{\eta} \big( \overline R(\eta) \sin(\eta-\theta) \big) \wrt\eta + O(\alpha^2) \,.
		\end{align*}
    By linearization, Lemma \ref{lemma.cR.00} and using that $ \tV[0;0]=\frac12$ by Lemma \ref{Id-lem-int}-$(ii)$,  we get
		\begin{equation}\label{pula1}
			\begin{aligned}
					\di_{r}\tV[0;0]h(\theta) &=- \tfrac12 h(\theta) -  \fint_{\T} \log\big| e^{\im\theta}- e^{\im\eta} \big| \pa_{\eta} \big( h(\eta) \sin(\eta-\theta) \big) \wrt\eta  \\
				& - \fint_{\T} \big(\di_{r}\log\big| R(\theta) e^{\im\theta}-  R(\eta)e^{\im\eta} \big| h\big)\big|_{ r=0}  \cos(\eta-\theta)  \wrt\eta \,.
			\end{aligned}
		\end{equation}
        On the other hand, using Lemma \ref{lemma.kernels}, we deduce that
        \begin{align}\label{log-1}
\big(\di_{r}\log\big| R(\theta) e^{\im\theta}- R(\eta)e^{\im\eta} \big| h \big) \big|_{ r=0}=\tfrac12(h(\theta)+h(\eta)) \,.
            \end{align}
            Therefore, we get, for any $h\in X_{\geq 2,\rm even}^{\rho,s}$,
        \begin{align*}
 \fint_{\T} \big(\di_{r}\log\big|\overline R(\theta) e^{\im\theta}- \overline R(\eta)e^{\im\eta} \big| h\big)  \cos(\eta-\theta)  \wrt\eta&= \tfrac12\fint_{\T} h(\eta) \cos(\eta-\theta)  \wrt\eta=0 \,.
            \end{align*}
           Plugging this identity into \eqref{pula1} and  integrating  by parts, we obtain
            \begin{align*}
			\di_{r}\tV[0;0]h(\theta) &= -\tfrac12 h(\theta)-  \fint_{\T} \log\big| e^{\im\theta}- e^{\im\eta} \big| \pa_{\eta} \big( h(\eta) \sin(\eta-\theta) \big) \wrt\eta\\
           &= - \tfrac12 h(\theta)+ \fint_{\T} {\rm Im}\Big(\tfrac{e^{\im \eta}}{ e^{\im\theta}- e^{\im\eta}} \Big)  h(\eta) \sin(\eta-\theta)  \wrt\eta \,.
		\end{align*}
     From the proof of Lemma \ref{lemma-1-Id} we have
    $$
    {\rm Im}\Big(\tfrac{e^{\im \eta}}{ e^{\im\theta}- e^{\im\eta}} \Big) =\tfrac12{\rm cotg}\Big(\tfrac{\eta-\theta}{2}\Big) \,.$$
    By classical trigonometric identities we deduce that
    \begin{align}
{\rm Im}\Big(\tfrac{e^{\im \eta}}{ e^{\im\theta}- e^{\im\eta}} \Big)\sin(\eta-\theta) =\tfrac12 \big(1+\cos(\eta-\theta)\big) \,.
    \end{align}
    Hence, we get, for any $h\in X_{\geq 2,\rm even}^{\rho,s}$
    \begin{align*}
		 \fint_{\T} {\rm Im}\Big(\tfrac{e^{\im \eta}}{ e^{\im\theta}- e^{\im\eta}} \Big)  h(\eta) \sin(\eta-\theta)  \wrt\eta&=\frac{1}{2  } \fint_{\T}\big(1+\cos(\eta-\theta)\big)  h(\eta)  \wrt\eta =0 \,.
		\end{align*}
        Therefore, we conclude that
        \begin{align*}
			\di_{r}\tV[0;0]h(\theta) &= - \tfrac12 h(\theta) \,.
            \end{align*}
            Inserting this into \eqref{112} allows to get
             \begin{align*}
		\di_{r}^2 \wt\cG_{0}(0;0)[h,h]& = \Pi_{\geqslant2}\pa_{\theta}  \big( - \tfrac12 h^2(\theta)+\di_{r}\cR(0;0)[h,h]\big) \,.
	\end{align*}
   Linearizing \eqref{cR.lin.app}, and using  Lemma \ref{lemma.kernels} together with \eqref{log-1}, we infer that  for $h\in X_{\geq 2,\rm even}^{\rho,s}$,
    \begin{align*}
	\partial_\theta\di_{r}\cR(0;0)[h,h](\theta)& = \frac{1}{4\pi} \partial_\theta\int_{\T} \big(h(\theta)+h(\eta)\big) h(\eta) \wrt \eta  =0 \,,
	\end{align*}
which implies in turn that 
$$\di_{r}^2 \wt\cG_{0}(0;0)[h,h] = \Pi_{\geqslant2}\pa_{\theta}  \big(-\tfrac12 h^2(\theta)\big) \,. $$
This ends the proof of the  formula in \eqref{pula2}.
    \end{proof}
\subsection{Construction of the traveling pairs}\label{Sec-proof-TP}
We are now in a position to detail the proof of  Proposition \ref{prop.ralpha}, to construct a curve of solutions $\alpha \in [-\alpha_0, \alpha_0] \mapsto r_{\alpha}(\theta)$ to \eqref{G0_station}, close to the trivial equilibrium given by $r_0(\theta) = 0$ at $\alpha = 0$.
\begin{proof}[Proof of Proposition $\ref{prop.ralpha}$ ]
	The claim follows from the standard Implicit Function Theorem applied with  the nonlinear map $\wt\cG_{0}$ introduced in Proposition \ref{lemma.constrain}, based on Proposition \ref{lemma.constrain1}.  As a consequence, we get a unique branch of functions $\alpha\in[-\alpha_0,\alpha_0]\mapsto r_\alpha\in X_{\geq 2,\rm even}^{\rho,s},$ satisfying
    \begin{equation}\label{IFT.equa}
		\begin{aligned}
			& \wt\cG_{0}(r_{\alpha};\alpha)=0, \quad \forall\,\alpha\in[-\alpha_{0},\alpha_{0}]\,,
	 	\\
 	& r_0(\theta) = 0 \,, \quad \tU[0,0] = 0\,.
		\end{aligned}
	\end{equation}
    Hence, we get
    $$
   \forall \, \alpha\in [-\alpha_0,\alpha_0]\,,\quad \cG(\tU_\alpha,r_\alpha;\alpha)=0 \,,\quad\hbox{with}\quad \tU_\alpha:=\tU(r_{\alpha},\alpha) \,.
   $$
   The $C^k$-regularity  of  the solution curve is guaranteed by the Implicit Function Theorem, as  the map $(\tU_\alpha,r_\alpha,\alpha)\mapsto \cG(\tU_\alpha,r_\alpha;\alpha)$ is analytic.   
\end{proof}

\subsection{Small amplitude expansion of the traveling waves}
The main objective here is to derive the expansion of the traveling wave profile $ r_\alpha $, as introduced in \mbox{Proposition \ref{prop.ralpha}}, with respect to the parameter $\alpha$. A first-order approximation, valid at $ O(\alpha) $, can be formally obtained by differentiating the implicit equation \eqref{IFT.equa} with respect to $\alpha$. This yields
\begin{equation}
	\pa_{\alpha} r_\alpha(\theta)|_{\alpha=0} = - \big(\di_{r} \wt\cG_{0}(0;0) \big)^{-1}\big[ \pa_{\alpha} \wt\cG_{0}(0;0) \big]\,.
\end{equation}
However, since the only explicit dependence on $\alpha$ in the equation comes from the term $ K_\alpha = O(\alpha^2)$, we find that $\partial_{\alpha} \wt\cG_{0}(0;0) = 0 $, and thus
$\pa_{\alpha} r_\alpha(\theta)|_{\alpha=0} =0$.
This implies that the leading order contribution to $ r_\alpha $ begins at $ O(\alpha^2) $. Although higher-order terms could in principle be computed by repeated differentiation of the implicit equation, the process quickly becomes algebraically intensive. Since we will only require the expansion of $ r_\alpha $ up to order $ O(\alpha^5) $ for subsequent analysis, we opt for a different—but still equivalent—approach rooted in the spirit of the Implicit Function Theorem. Our main result reads as follows.
{
\begin{pro}\label{prop.asympt.ralpha}
	The  solution $r_\alpha$ constructed in Proposition $\ref{prop.ralpha}$ for $|\alpha|\leq \alpha_{0}$ satisfies the asymptotic expansion
	\begin{equation}\label{ralpha.asymp}
		r_\alpha(\theta) = - \alpha^2 \cos(2\theta)  +\tfrac12 \alpha^3 \cos(3\theta)+\tfrac23 \alpha^4 \cos(4\theta)  + O(\alpha^5)\,.
	\end{equation}
    Moreover, the speed introduced in Proposition $\ref{lemma.constrain1}$ satisfies
   $$
   \tU[{r}_\alpha;\alpha]:=\tU_\alpha=O(\alpha^4) \,.
   $$
\end{pro}}
\begin{proof}
    The analysis presented below is formal in nature, but can be made rigorously by carefully estimating the associated error terms. To keep the exposition clear and concise, we adopt a purely formal approach to the expansion. We seek an asymptotic expansion of \( r_\alpha \) of the form,
	\begin{equation}\label{ansatz.ralpha}
		r_\alpha(\theta) = \alpha^2 r_2(\theta)  + \alpha^3 r_3(\theta) + \alpha^4 r_4(\theta) + O(\alpha^5) \,.
	\end{equation}
	By expanding \eqref{IFT.equa} in Taylor series with respect to $r$, we get
	\begin{equation}
		0 = \wt\cG_{0}(r_\alpha;\alpha) = \wt\cG_{0}(0;\alpha) + \di_{r} \wt\cG_{0}(0;\alpha) [r_{\alpha}] + \tfrac12 \di_{r}^2 \wt\cG_{0}(0;\alpha) [r_\alpha,r_\alpha] + O(|r_\alpha|^3) \,.
	\end{equation}
    We first write $r_\alpha:=\alpha^2\rho$, so that $\rho$ satisfies the equation
    \begin{equation}\label{equa.for.rho1}
		 \wt\cG_{0}(0;\alpha) +\alpha^2 \di_{r} \wt\cG_{0}(0;\alpha) [\rho] +\alpha^4 \tfrac12 \di_{r}^2 \wt\cG_{0}(0;0) [\rho,\rho] + O(\alpha^6)=0 \,.
	\end{equation}
    By Lemma \ref{Id-lem-int}-$(i)$, we have
    	\begin{equation}
    		\wt\cG_{0}(0;\alpha) =- \tfrac12 \alpha^2 \sin(2\theta) + \tfrac12 \alpha^3 \sin(3\theta) - \tfrac12\alpha^4 \sin(4\theta) + O(\alpha^5) \,.
    	\end{equation}
    	Applying   Proposition \ref{lemma.constrain1}, Lemma \ref{Id-lem-int}-$(iv)$ and \eqref{linera-drU}, allow to get  for a $h(\theta)=\sum_{m\geq 2} h_m \cos(m\theta)$,
    		\begin{align}
    			\di_{r} &\wt\cG_{0}(0;\alpha) h  = \pa_{\theta} \circ \Pi_{\geq 2}\Big( \cos(\theta) \di_{r}\tU[0;\alpha]h+  \tV[0;\alpha]h +  \cR(0;\alpha) h  \Big) \\
    			& = \pa_{\theta} \circ \Pi_{\geq 2} \Big(    \tfrac12 \big(1+\alpha^2\cos(2\theta)\big)h - \sum_{m\geq 2} \tfrac{h_m}{2m}  \cos(m\theta) + \tfrac{\alpha^2}{4} h_2 + O(\alpha^3) \Big) \\
    			& = \di_{r}\wt\cG_{0}(0;0)h { + \tfrac{\alpha^2}{2}\pa_{\theta}\circ \Pi_{\geq 2}\big(\cos(2\theta)h\big) } + O(\alpha^3) \,.
    		\end{align}
     Using Proposition \ref{lemma.constrain22}, we get
    	\begin{equation}
    	\di_{r}^2\wt\cG_{0}(0;0) [h,h] = - \tfrac12 \Pi_{\geq 2} \pa_{\theta}(h^2) \,.
    	\end{equation}
    Recalling \eqref{ansatz.ralpha}, and adopting the ansatz
    \begin{equation}\label{ansatz.rho}
    	\rho(\theta) = \rho_{0}(\theta) + \alpha \rho_{1}(\theta) + \alpha^2 \rho_{2}(\theta) + O(\alpha^3)\,,
    \end{equation}
    we find that equation  \eqref{equa.for.rho1} yields, order by order, 
    	\begin{align}
    		\alpha^2\big( &-\tfrac12 \sin(2\theta) + \di_{r}\wt\cG_{0}(0;0) \rho_{0} \big) + \alpha^3 \big( \tfrac12 \sin(3\theta) + \di_{r}\wt\cG_{0}(0;0) \rho_{1} \big) \\
    		& + \alpha^4 \big( -\tfrac12 \sin(4\theta) + \di_{r}\wt\cG_{0}(0;0) \rho_{2} +{ \tfrac12 \pa_{\theta}\Pi_{\geq 2}\big(\cos(2\theta)\rho_{0}\big) }  - \tfrac14 \Pi_{\geq 2} \pa_{\theta} (\rho_{0}^2) \big) + O(\alpha^5) = 0\,.
    	\end{align}
    Solving at orders $\alpha^2$ and $\alpha^3$, and invoking Proposition \ref{lemma.constrain1}-$(iv)$, we  obtain 
    	\begin{align}
    		\rho_{0}(\theta) &:=  \tfrac12 \big( \di_{r}\wt\cG_{0}(0;0) \big)^{-1} [\sin(2\theta)]     	= \tfrac{1}{4\im}  \big( \di_{r}\wt\cG_{0}(0;0) \big)^{-1}[e^{2\im\theta}-e^{-2\im\theta}]\nonumber \\ 
    		& = - \frac{1}{4\im} \frac{2\im}{2-1} \big( e^{2\im\theta} + e^{-2\im\theta} \big) = - \cos(2\theta) \,, \label{rho0.2theta} 
            \end{align}
            and
            \begin{align}
    			\rho_{1}(\theta) &:= - \tfrac12 \big( \di_{r}\wt\cG_{0}(0;0) \big)^{-1} [\sin(3\theta)]     	= -\tfrac{1}{4\im}  \big( \di_{r}\wt\cG_{0}(0;0) \big)^{-1}[e^{3\im\theta}-e^{-3\im\theta}] \nonumber \\
    		& = \frac{1}{4\im} \frac{2\im}{3-1} \big( e^{3\im\theta} + e^{-3\im\theta} \big) =\tfrac12 \cos(3\theta) \,.\label{rho1.3theta}
    	\end{align}
    To compute $\rho_{2}$, we first note that
    \begin{align}
    { \tfrac12 \pa_{\theta}\Pi_{\geq 2}\big(\cos(2\theta)\rho_{0}\big) } & = {- \tfrac14 \pa_{\theta}\Pi_{\geq 2}\big( 1+ \cos(4\theta) \big)  = \sin(4\theta)} \,, \\
       \tfrac14	\Pi_{\geq 2}\pa_{\theta}\big( \rho_{0}^2(\theta) \big) & =  \tfrac18  \Pi_{\geq 2}\pa_{\theta}\big(  1+\cos(4\theta)  \big)= - \tfrac12 \sin(4\theta)\,.
    \end{align}
    Therefore, solving at order $\alpha^4$, we obtain 
 	\begin{align}
	\rho_{2}(\theta)  &:=   \big( \di_{r}\wt\cG_{0}(0;0) \big)^{-1} \big[ \tfrac12 \sin(4\theta)- { \tfrac12 \pa_{\theta}\Pi_{\geq 2}\big(\cos(2\theta)\rho_{0}\big) }  + \tfrac14 \Pi_{\geq 2} \pa_{\theta} (\rho_{0}^2)   \big]  \\
	&=   \big( \di_{r}\wt\cG_{0}(0;0) \big)^{-1} \big[ \tfrac12 \sin(4\theta)-{\sin(4\theta)}  -\tfrac12 \sin(4\theta) \big]  \\
    & =  -\big( \di_{r}\wt\cG_{0}(0;0) \big)^{-1} \big[ {\sin(4\theta)}  \big]  = - \frac{1}{2\im}  \big( \di_{r}\wt\cG_{0}(0;0) \big)^{-1}[e^{4\im\theta}-e^{-4\im\theta}] \nonumber \\
    		& = \frac{1}{2\im} \frac{2\im}{4-1} \big( e^{4\im\theta} + e^{-4\im\theta} \big) = \tfrac23 \cos(4\theta)  \,.   \label{rho2.4theta} 
\end{align}
 Inserting \eqref{rho0.2theta}, \eqref{rho1.3theta}, \eqref{rho2.4theta} into \eqref{ansatz.rho} and recalling that $r_\alpha(\theta)=\alpha^2 \rho(\theta)$, we conclude the claimed expansion in \eqref{ralpha.asymp}. The estimate for $\tU_\alpha$ follows directly from {Proposition \ref{lemma.constrain1}-(iii)} together with the bound $r_\alpha = O(\alpha^2)$. 
\end{proof}

\subsection{Expansion of the linearized operator at the traveling wave}

In Section \ref{sect.linear.eq}, our main objective is to analyze the dynamics in a neighborhood of the traveling waves constructed in Proposition \ref{prop.ralpha}. To this end, we investigate the small-amplitude expansion of the associated linearized operator, whose general structure was derived in \eqref{nonlinear.cG}. 
At this stage, we focus only on the partial linearization with respect to the boundary variable \(r\), while the propagation speed \(U\) is frozen at the traveling-wave value \(\tU_{\alpha}\). The full linearization, including the modulation of the speed parameter, will be addressed later in Section \ref{sect.linear.eq}. 
So from \eqref{nonlinear.cG}, the partial linearization with respect to the shape variable $r$ takes the form  
\begin{align}\label{nonlinear.cGL}
\partial_{r}\cG(\tU_\alpha,{r}_\alpha;\alpha)h
   &= \omega \pa_{\vf} h 
   + \pa_{\theta} \bigg(  \tfrac{\tU_{\alpha}\,\cos(\theta)}{\sqrt{1+2 {r}_\alpha}} \, h 
   + \di_{r}  F[ r_\alpha; \alpha] h \bigg) \\
   &:= \omega \pa_{\vf} h + \pa_{\theta}\,\cT(r_\alpha;\alpha)[h] \,.
\end{align}
The operator \(\cT(r_\alpha;\alpha)\) admits the expansion  
\begin{align}
	\cT(r_\alpha;\alpha)[h] 
		&= \Big( \tfrac{\tU_{\alpha} \cos(\theta)}{\sqrt{1+2r_\alpha(\theta)}} 
               + \tV[r_\alpha;\alpha]  \Big)h 
               +  \cR(r_\alpha;\alpha)[h] \\
        &= \tV_{0,\alpha}(\theta)\,h + \cR(r_\alpha;\alpha)[h], 
        \qquad \tU_{\alpha}:=\tU[r_\alpha;\alpha],	\label{T.alpha.operator}
\end{align}
where \(\tV[r_\alpha;\alpha]\) and \(\cR(r_\alpha;\alpha)\) are defined in \eqref{V.at.boundary} and \eqref{cR.lin.app}, respectively.  For later convenience, we introduce the notation
\begin{equation}\label{wttV}
	\tV_{0,\alpha}(\theta) := \frac{\tU_{\alpha} \cos(\theta)}{\sqrt{1+2r_\alpha(\theta)}} + \tV[r_\alpha;\alpha] \,.
\end{equation}
This decomposition highlights two distinct components:  a leading-order multiplication operator given by the smooth function \(\tV_{0,\alpha}\),  
and a remainder term \(\cR(r_\alpha;\alpha)\), which, as we shall see, acts as a singular yet smoothing perturbation of the Hilbert transform.  

In what follows, we provide explicit asympotics for the terms in \eqref{T.alpha.operator} with respect to $\alpha\to 0$ up to order $\alpha^4$. In particular, the function $\tV_{0,\alpha}(\theta)$ is expanded in Lemma \ref{speed.V.expan} and the action of the operator $\cR(r_\alpha;\alpha)$ in Lemma \ref{lemma.R.expand}. These asymptotics will be used in Section \ref{sect.linear.eq} to compute the linear dynamics at the equilibrium. Furthermore, in Corollary \ref{cor.T.at.cos} we shall derive the expansion of the action of \(\cT(r_\alpha;\alpha)\) on the first Fourier mode, which is crucial to address technical difficulties related to the resonance of this mode at leading order.  

For later use, it is convenient to introduce the function $\tg_{\alpha}$ defined by
\begin{equation}\label{tg.alpha}
	\tg_{\alpha}(\theta) := \sqrt{1+2r_\alpha(\theta)} -1 \,,
\end{equation}
whose asymptotic behavior as $\alpha \to 0$ is described in the following result.
	\begin{lem}\label{sqrt.1+ralpha}
	For small values of $\alpha$ we have the expansion,
		\begin{align*}
			1 + \tg_{\alpha}(\theta) 
			& = 1 - \alpha^2 \cos(2\theta) + \tfrac12 \alpha^3 \cos(3\theta) + \alpha^4 \big( \tfrac{5}{12}\cos(4\theta) -  \tfrac14  \big) + O(\alpha^5) \,.
		\end{align*}
\end{lem}
\begin{proof}
	From \eqref{tg.alpha} and Proposition \ref{prop.asympt.ralpha}, we expand
		\begin{align}
			1+ \tg_{\alpha}(\theta) & = \sqrt{1+2r_\alpha(\theta)} = 1+ r_\alpha(\theta) - \tfrac12 \big( r_\alpha(\theta) \big)^2 + O(|r_\alpha|^3) \\
			& = 1 - \alpha^2 \cos(2\theta) + \tfrac12 \alpha^3 \cos(3\theta) + \alpha^4 \big( \tfrac23 \cos(4\theta)  - \tfrac12 \cos^2(2\theta) \big) + O(\alpha^5)  \\
			& = 1 - \alpha^2 \cos(2\theta) + \tfrac12 \alpha^3 \cos(3\theta) +\alpha^4 \big(  \tfrac{5}{12}\cos(4\theta) - \tfrac14   \big) + O(\alpha^5) \,,
		\end{align}
	and the claim is proved.
\end{proof}
The next goal is to expand the function $\tV_{0,\alpha}$ described through  \eqref{wttV}.
\begin{lem}\label{speed.V.expan}
	The function $\tV_{0,\alpha}$ 
	satisfies for small $\alpha$ the asymptotic
	\begin{equation}
	\tV_{0,\alpha}(\theta)= \tfrac12 + \tU_{\alpha}\cos(\theta) +\alpha^2\cos(2\theta) -\tfrac34 \alpha^3 \cos(3\theta)  + \tfrac{2}{3}\alpha^4\cos(4\theta)  + O(\alpha^5) \,. \label{speed.expand}
	\end{equation}
    In addition,
	\begin{align}
	 (\tV_{0,\alpha}(\theta))^{-1} 
    = 2 -4 \tU_{\alpha} \cos(\theta) -4\alpha^2 \cos(2\theta)  + 3 \alpha^3 \cos(3\theta) +\alpha^4 \big( 4 +\tfrac{4}{3} \cos(4\theta)  \big) + O(\alpha^5) \,.\label{inv.speed.expan}
	\end{align}
\end{lem}
\begin{proof}
    By \eqref{wttV}, Proposition \ref{prop.asympt.ralpha} and Lemma \ref{sqrt.1+ralpha} we get
    \begin{align}
	\tV_{0,\alpha}(\theta)&= \frac{\tU_{\alpha} \cos(\theta)}{1+\tg_{\alpha}(\theta)} + \tV[r_\alpha;\alpha](\theta)\\
    &=\tU_{\alpha} \cos(\theta)+\tV[r_\alpha;\alpha](\theta)+O(\alpha^5)\,.\label{form-V-0al}
    \end{align}
On the other hand, the formula  \eqref{V.at.boundary} allows to get

		\begin{align}
			\tV[r_\alpha;\alpha](\theta)=& \frac{-1}{1+\tg_{\alpha}(\theta)} \Bigg\{   \fint_{\T} \log\big|\mathtt{K}(\theta,\eta)  \big|\, \pa_{\eta} \big( (1+\tg_{\alpha}(\eta)) \sin(\eta-\theta) \big) \wrt\eta \\
			 &+ \fint_{\T} K_\alpha\big((1+\tg_{\alpha}(\theta)) e^{\im\theta}, (1+\tg_{\alpha}(\eta))e^{\im\eta} \big) \pa_{\eta} \big( (1+\tg_{\alpha}(\eta)) \sin(\theta+\eta) \big) \wrt \eta \Bigg\}\,,
		\end{align}
	where 
    $$
    \mathtt{K}(\theta,\eta):=(1+\tg_{\alpha}(\theta))  e^{\im\theta}- (1+\tg_{\alpha}(\eta))e^{\im\eta}
    $$
    and $\tg_{\alpha}(\theta)$ is as in \eqref{tg.alpha} and satisfies the expansion in Lemma \ref{sqrt.1+ralpha}. To get the expansion of $\tV[r_\alpha;\alpha]$ we first split it as follows
   \begin{align}\label{V-ralpha}
    \tV[r_\alpha;\alpha](\theta) = {=\frac{\tV_0(\theta)+\tV_{1}(\theta)}{1+\tg_{\alpha}(\theta)} } \,,
    \end{align}
with
\begin{align}\label{V-0-form}
\tV_0(\theta):=-\fint_{\T} \log\big|\mathtt{K}(\theta,\eta)  \big|\, \pa_{\eta} \big( (1+\tg_{\alpha}(\eta)) \sin(\eta-\theta) \big) \wrt\eta
\end{align}
and
\begin{align}\label{V-1-form}
\tV_1(\theta):=\fint_{\T} K_\alpha\big((1+\tg_{\alpha}(\theta)) e^{\im\theta}, (1+\tg_{\alpha}(\eta))e^{\im\eta} \big) \pa_{\eta} \big( (1+\tg_{\alpha}(\eta)) \sin(\theta+\eta) \big) \wrt \eta\,.
\end{align}
Our next goal is to determine the asymptotic structure of each component in this decomposition and collect them in the end.
\\[1mm]
\noindent $\blacktriangleright$ {\sc Expansion for $(1+\tg_{\alpha}(\theta))^{-1}$.}
	 By Taylor expansion, we obtain
	\begin{align*}
		(1 +& \tg_{\alpha}(\theta))^{-1}  = 1- \tg_{\alpha}(\theta) + \tg_{\alpha}^2(\theta) + O(|r_\alpha|^3).
	\end{align*}
Then, applying Lemma \ref{sqrt.1+ralpha} yields
\begin{align}\label{calcetto1}
		\nonumber (1 + \tg_{\alpha}(\theta))^{-1}  &= 1 + \alpha^2 \cos(2\theta) - \tfrac12 \alpha^3 \cos(3\theta) + \alpha^4 \big( -\tfrac{5}{12} \cos(4\theta)+\tfrac14 +\cos^2(2\theta) \big)  + O(\alpha^5) \\
		& = 1 + \alpha^2 \cos(2\theta) - \tfrac12 \alpha^3 \cos(3\theta) + \alpha^4 \big( \tfrac{1}{12}\cos(4\theta) + \tfrac{3}{4} \big) + O(\alpha^5) \,.
	\end{align}
\noindent $\blacktriangleright$ {\sc Expansion for $\tV_{0}(\theta)$.}
	 First, using \eqref{BS.kernel.expand} in Lemma \ref{lemma.kernels}, together with a change of variables, we expand $\tV_0(\theta)$ in \eqref{V-0-form} as
	 \begin{align}
	 		\tV_{0}(\theta) 
	 		& = -\fint_{\T} \log\big| 1- e^{\im\eta} \big| \pa_{\eta} \big( (1+\tg_{\alpha}(\eta+\theta)) \sin(\eta) \big) \wrt\eta \\
	 		&\quad - \frac{1}{2} \fint_{\T} (\tg_{\alpha}(\theta)+\tg_{\alpha}(\eta)) \pa_{\eta} \big( (1+\tg_{\alpha}(\eta)) \sin(\eta-\theta) \big) \wrt\eta \\
	 		& \quad + \frac{1}{8} \fint_{\T} (\tg_{\alpha}(\theta)+\tg_{\alpha}(\eta))^2 \cos(\eta-\theta) \wrt\eta \\
	 		& \quad -\frac{1}{8} \fint_{\T}  (\tg_{\alpha}(\theta)-\tg_{\alpha}(\eta))^2 {\rm cotg}^2\big(\tfrac{\theta-\eta}{2}\big)  \cos(\eta-\theta) \wrt\eta + O( |\tg_{\alpha}|^3 + |\tg_{\alpha}'|^3) \\
	 		& =: \sum_{j=1}^4\tV_{0}^{j} + O(\alpha^6) \,.
	 \end{align}
	 We shall analyze these terms one  by one. First, we note that, for any integer  $k\geq 2$, by using Lemma \ref{lemma.cR.00} and item $1.$ in Appendix \ref{app.trigo},
	 	\begin{align}
	 		&\fint_{\T} \log\big| 1- e^{\im\eta} \big| \pa_{\eta} \big( \cos(k(\eta+\theta)) \sin(\eta) \big) \wrt\eta \\
	 		& = \frac{1}{2 }\fint_{\T} \log\big| 1- e^{\im\eta} \big|   \big( (k+1)\cos(k\theta+(k+1)\eta) - (k-1)\cos(k\theta+(k-1)\eta) \big) \wrt\eta \\
	 		& = \frac{1}{4} \big( \cos(k\theta)-\cos(k\theta) \big) = 0\,.
	 	\end{align}
	 Combining this identity with  Lemma \ref{sqrt.1+ralpha} and Lemma \ref{lemma.cR.00}, we find
	 	\begin{align}
	 		\tV_{0}^{1} &= - \fint_{\T} \log\big| 1- e^{\im\eta} \big| \pa_{\eta} \big( (1+\tg_{\alpha}(\eta+\theta)) \sin(\eta) \big) \wrt\eta \\
	 		&= -  (1-\tfrac14 \alpha^4)\fint_{\T} \log\big| 1- e^{\im\eta} \big|  \cos(\eta)  \wrt\eta +O(\alpha^5) = \tfrac12 \big(1- \tfrac14 \alpha^4\big) + O(\alpha^5) \,.
	 	\end{align}
	 Next, we move to the expansion of $\tV_{0}^{2}$. First, we observe that
    \begin{align}
	 	\tV_{0}^{2}(\theta)	&= - \frac{1}{2}  \fint_{\T} (\tg_{\alpha}(\theta)+\tg_{\alpha}(\eta)) \pa_{\eta} \big( (1+\tg_{\alpha}(\eta)) \sin(\eta-\theta) \big) \wrt\eta \\
	 		& = - \frac{1}{2} \fint_{\T} \tg_{\alpha}(\eta) \pa_{\eta} \big( (1+\tg_{\alpha}(\eta)) \sin(\eta-\theta) \big) \wrt\eta  \,.
	 	\end{align}     
     Thus,   Lemma \ref{sqrt.1+ralpha} yields
	 	\begin{align}
	 	\tV_{0}^{2}(\theta)	
	 		& = - \frac{1}{2} \fint_{\T} \big(- \alpha^2 \cos(2\eta) +\tfrac12 \alpha^3 \cos(3\eta) + \alpha^4\big( \tfrac{5}{12} \cos(4\eta) -\tfrac14 \big) \big) \\ & \quad \quad \quad \quad\pa_{\eta} \big( (1-\alpha^2 \cos(2\eta)) \sin(\eta-\theta) \big) \wrt\eta + O(\alpha^5) \\
	 		& = - \frac{\alpha^4}{2} \fint_{\T}\cos(2\eta) \pa_{\eta} \big( \cos(2\eta) \sin(\eta-\theta) \big) \wrt\eta + O(\alpha^5) \\
            & = O(\alpha^5) \,,
	 	\end{align}
since the last integral identically vanishes.
	Concerning the term $\tV_{0}^{3}$, we expand,  by Lemma \ref{sqrt.1+ralpha},
	 	\begin{align}
	 		\tV_{0}^{3}& = \frac{1}{8} \fint_{\T} \big( \tg_{\alpha}(\theta) + \tg_{\alpha}(\eta) \big)^2 \cos(\eta-\theta) \wrt\eta \\
	 		& =  \frac{\alpha^4}{8} \fint_{\T} \big( \cos(2\theta) + \cos(2\eta) \big)^2 \cos(\eta-\theta) \wrt\eta + O(\alpha^5) \\
           & = O(\alpha^5) \,,
	 	\end{align}
since the last integral identically vanishes.
	 We are left to analyze the term $\tV_{0}^{4}$. By item $3.$ in Appendix \ref{app.trigo} and checking the trigonometric reductions, we find
	 	\begin{align}
	 			(	\cos(2\theta)&-\cos(2\eta))^2  {\rm cotg}^2\big(\tfrac{\theta-\eta}{2}\big) \cos(\eta-\theta)  =\\
	 		& = 4\sin^2(\theta+\eta)(\cos(\theta-\eta)+1)^2 \cos(\eta-\theta) \\
	 		& = 2- \cos(4\theta) -\cos(4\eta)+ \tfrac12 \cos(3\theta-3\eta) + 2\big( \cos(2\theta-2\eta) - \cos(2\theta+2\eta)  \big) \\
	 		& \quad + \tfrac72 \cos(\theta-\eta) -\tfrac14 \big( \cos(\theta-5\eta) + \cos(5\theta-\eta) \big) - \tfrac74 \big( \sin(\theta+3\eta)+\sin(3\theta+\eta)\big) \,.
	 	\end{align}
	 It follows that
	 	\begin{align}
	 		\tV_{0}^{4}(\theta)& = - \frac{1}{8}  \fint_{\T}  (\tg_{\alpha}(\theta)-\tg_{\alpha}(\eta))^2 {\rm cotg}^2\big(\tfrac{\theta-\eta}{2}\big)  \cos(\eta-\theta) \wrt\eta \\
	 		& = - \frac{\alpha^4}{8} \fint_{\T}  (\cos(2\theta)-\cos(2\eta))^2 {\rm cotg}^2\big(\tfrac{\theta-\eta}{2}\big)  \cos(\eta-\theta) \wrt\eta + O(\alpha^5)\\
	 		& =  - \frac{\alpha^4}{8} \fint_{\T}(2-\cos(4\theta)) \wrt \eta + O(\alpha^5)  = \alpha^4\big( -\tfrac14 + \tfrac18 \cos(4\theta)\big) + O(\alpha^5)\,.
	 	\end{align}
	 Summing up all the contributions, we conclude that
	 	\begin{align}
	 	\tV_{0} (\theta)&  = \tfrac12 - \tfrac18 \alpha^4+ \alpha^4 (-\tfrac14 + \tfrac18\cos(4\theta)) + O(\alpha^5) \\
	 		&= \tfrac12 +\alpha^4(-\tfrac38 + \tfrac18\cos(4\theta)) + O(\alpha^5) \,.\label{calcetto2}
	 	\end{align}
\noindent $\blacktriangleright$ {\sc Expansion for $\tV_{1}(\theta)$.}
 Using \eqref{inter.kernel} in Lemma \ref{lemma.kernels} and recalling Lemma \ref{sqrt.1+ralpha}, we expand $\tV_{1}(\theta)$ in \eqref{V-1-form} as
\begin{align*}
		\tV_{1}(\theta) =&	- \fint_{\T} K_\alpha\big((1+\tg_{\alpha}(\theta)) e^{\im\theta}, (1+\tg_{\alpha}(\eta))e^{\im\eta} \big) \cos(\theta+\eta)  \wrt \eta  \\
       & +\alpha^2 \fint_{\T} K_\alpha\big((1+\tg_{\alpha}(\theta)) e^{\im\theta}, (1+\tg_{\alpha}(\eta))e^{\im\eta} \big) \pa_{\eta} \big( \cos(2\eta) \sin(\theta+\eta) \big) \wrt \eta +O(\alpha^5)\\
   :=&\tV_{1}^1(\theta)+\tV_{1}^2(\theta)+O(\alpha^5)\,.
	\end{align*}
    Applying \eqref{inter.kernel} once again leads to
    \begin{align*}
		\tV_{1}^2(\theta) =&-\tfrac12\alpha^4 \fint_{\T} \big[\cos(2\theta)+2\cos(\theta-\eta)+\cos(2\eta)] \pa_{\eta}\big( \cos(2\eta) \sin(\theta+\eta) \big) \wrt \eta +O(\alpha^5)\\
       =&-\alpha^4 \fint_{\T} \cos(\theta-\eta) \pa_{\eta}\big( \cos(2\eta) \sin(\theta+\eta) \big) \wrt \eta +O(\alpha^5)\,.
	\end{align*}
    Integration by parts and straightforward computations yield
    \begin{align*}
		\tV_{1}^2(\theta) 
       &= \alpha^4 \fint_{\T} \sin(\theta-\eta)  \cos(2\eta) \sin(\theta+\eta)  \wrt \eta +O(\alpha^5)\\
      &=\tfrac14 \alpha^4 +O(\alpha^5)\,. 
	\end{align*}
    Concerning the expansion of $\tV_{1}^1$, we use \eqref{inter.kernel} together with Lemma \ref{sqrt.1+ralpha}, which allows us to obtain
    \begin{equation*}
			\begin{aligned}
			K_\alpha\big( & (1+\tg(\theta))e^{\im\theta} ,  (1+\tg(\eta))e^{\im\eta} \big)  = - \tfrac{\alpha^2}{2}\big[  \cos(2\theta) +2 \cos(\theta-\eta) + \cos(2\eta) \big] \\
			& \ \ +\tfrac{\alpha^3}{3} \big[ \cos(3\theta) + 3\cos(2\theta-\eta) + 3 \cos(\theta-2\eta) + \cos(3\eta) \big] \\
			& \ \ - \tfrac{\alpha^4}{4} \big[ \cos(4\theta) + 4 \cos(3\theta-\eta) + 6\cos(2(\theta-\eta))  + 4\cos(\theta-3\eta) + \cos(4\eta)\big] \\
			& \ \ +\alpha^4\big[ \cos^2(2\theta)  + (\cos(2\theta)+\cos(2\eta)) \cos(\theta-\eta) + \cos^2(2\eta) \big] +O\big(\alpha^5 \big) \,.
		\end{aligned}
	\end{equation*}
    Therefore
    \begin{align*}
		\tV_{1}^1(\theta) =&	\alpha^2 \fint_{\T} \cos(\theta-\eta) \cos(\theta+\eta)  \wrt \eta-\alpha^3  \fint_{\T} \cos(2\theta-\eta) \cos(\theta+\eta)  \wrt \eta \\
       &+\alpha^4  \fint_{\T} \cos(3\theta-\eta) \cos(\theta+\eta)  \wrt \eta\\
       &-\alpha^4\fint_{\T} (\cos(2\theta)+\cos(2\eta)) \cos(\theta-\eta) \cos(\theta+\eta)  \wrt \eta\,.
	\end{align*}
   Using classical  trigonometric identities we infer
      \begin{align*} \tV_{1}^1(\theta) &=	\tfrac12\alpha^2 \cos(2\theta)-\tfrac12\alpha^3  \cos(3\theta) +\tfrac12\alpha^4  \cos(4\theta)-\alpha^4\big(\tfrac12+\tfrac14\cos(4\theta)\big)\\
        &=\tfrac12\alpha^2 \cos(2\theta)-\tfrac12\alpha^3  \cos(3\theta) +\alpha^4\big(-\tfrac12+\tfrac14\cos(4\theta)\big)\,.
	\end{align*}
    Putting together the preceding results yields
   \begin{align}\label{calcetto3}
		\nonumber \tV_{1}(\theta) =&	\tfrac14\alpha^4+\tfrac12\alpha^2 \cos(2\theta)-\tfrac12\alpha^3  \cos(3\theta) +\alpha^4\big(-\tfrac12+\tfrac14\cos(4\theta)\big)+O(\alpha^5)\\
       =&	\tfrac12\alpha^2 \cos(2\theta)-\tfrac12\alpha^3  \cos(3\theta) +\alpha^4\big(-\tfrac14+\tfrac14\cos(4\theta)\big)+O(\alpha^5) \,.
	\end{align} 
\noindent $\blacktriangleright$ {\sc Conclusion.}
Plugging \eqref{calcetto1}, \eqref{calcetto2} and \eqref{calcetto3} into \eqref{V-ralpha} gives
\begin{align*}
    \tV[r_\alpha;\alpha]=&\Big(\tfrac12 +\tfrac12\alpha^2 \cos(2\theta)-\tfrac12\alpha^3  \cos(3\theta) +\alpha^4\big(-\tfrac58+\tfrac38\cos(4\theta)\big)\Big)\\
    &\times\Big(1 + \alpha^2 \cos(2\theta) - \tfrac12 \alpha^3 \cos(3\theta) + \alpha^4 \big( \tfrac{1}{12}\cos(4\theta) + \tfrac{3}{4} \big)\Big)+O(\alpha^5)\\
    =&\tfrac12 + \alpha^2\cos(2\theta) - \tfrac34\,\alpha^3\cos(3\theta) + \tfrac{2}{3}\,\alpha^4\cos(4\theta) + O(\alpha^5)\,.
    \end{align*}
    Inserting this expansion into \eqref{form-V-0al} allows to get
    \begin{align}\label{VV-0}
	{\tV_{0,\alpha}(\theta)}
    &=\tfrac12+\tU_{\alpha} \cos(\theta)+ \alpha^2\cos(2\theta) - \tfrac34\,\alpha^3\cos(3\theta) + \tfrac{2}{3}\,\alpha^4\cos(4\theta)+O(\alpha^5)\,.
    \end{align} 
	 This concludes the proof of the expansion in \eqref{speed.expand}. We now prove the last expansion \eqref{inv.speed.expan}.
By \eqref{VV-0}, \eqref{prop.asympt.ralpha} and Taylor expansion, we get
	\begin{align*}
		(\tV_{0,\alpha}(\theta))^{-1} & = 2\big[  1 + 2\tU_{\alpha} \cos(\theta) + 2\alpha^2 \cos(2\theta) - \tfrac32 \alpha^3 \cos(3\theta) +\tfrac{4}{3} \alpha^4 \cos(4\theta) + O(\alpha^5) \big]^{-1} \\
		& = 2 -4\tU_{\alpha} \cos(\theta) -4\alpha^2 \cos(2\theta) + 3 \alpha^3 \cos(3\theta)  + \alpha^4\big( 4 +\tfrac{4}{3}  \cos(4\theta)\big) + O(\alpha^5) \,,
	\end{align*}
	and the claim is proved.
\end{proof}
In the next lemma we derive the asymptotic expansion of the nonlocal operator $\cR(r_\alpha;\alpha)$ defined in \eqref{cR.lin.app}. In particular, we identify its principal contribution together with the lower-order interaction terms generated by the weak coupling between the two vortex patches. We remark that, for sake of convenience, we shall use here the notation \eqref{notation.tR} that we will extensively employ in Section \ref{sect.linear.eq}.
\begin{lem}\label{lemma.R.expand}
The operator $\cR(r_\alpha;\alpha)$ in \eqref{cR.lin.app}, with $\bar r = r_\alpha$, admits the following expansion when acting on zero-average functions $h:\theta\in\mathbb{T}\mapsto h(\theta)$, 	
    \begin{align}
			&\cR(r_\alpha;\alpha) [h]  = \cR(0;0)[h] +	\tA[h]+ {\alpha^2}\fint_{\T}\cos(\theta-\eta) h(\eta) \wrt \eta \\
			&\quad  - {\alpha^3} \fint_{\T}\big( \cos(2\theta-\eta) + \cos(\theta-2\eta) \big)h(\eta) \wrt \eta +{\alpha^4}\fint_{\T} \mathtt{K}_4(\theta,\eta)h(\eta) \wrt\eta  +\alpha^5 \tR[h]\,,
			\label{exp.R.total}
	\end{align}
  for some $\tR \in \OpM_{s}^{-\infty}$, $s\geq s_0$,  with
    \begin{align}
        \mathtt{K}_4(\theta,\eta)&= \tfrac12 \cos(\theta-\eta) - \cos(\theta+\eta)+\tfrac32\cos(2\theta-2\eta)  -\tfrac12\cos(2\theta+2\eta)\\
        &+\tfrac12\cos(3\theta-\eta) -\tfrac14\cos(3\theta+\eta)  +\tfrac12 \cos(\theta-3\eta)-\tfrac14\cos(\theta+3\eta) \,, \label{K4.kernel}
    \end{align}
and    $\tA[h]\in \R$ is a constant in $\theta$. 
	In particular, when $h= \Pi_{\geq 2}h$, we have
			\begin{align}
			\cR(r_\alpha;\alpha) [h] & = \cR(0;0)  [h] +	\tA[ h]  - {\alpha^3} \fint_{\T} \cos(\theta-2\eta)h(\eta) \wrt \eta    \\
			&\quad +{\alpha^4}\fint_{\T}\big(\tfrac32\cos(2\theta-2\eta)  -\tfrac12\cos(2\theta+2\eta)\big)  h(\eta) \wrt\eta \\
			&  \quad +{\alpha^4} \fint_{\T} \big(\tfrac12 \cos(\theta-3\eta)-\tfrac14\cos(\theta+3\eta) \big)  h(\eta) \wrt\eta +\alpha^5 \tR[h] \,.
			\label{exp.R.geq2}
			\end{align}
	Moreover, when acting on $h(\theta)=\cos(\theta)$ and $h(\theta)=\sin(\theta)$, it takes the form
		\begin{align}
			\cR(r_\alpha;\alpha) [\cos](\theta) & = \tA[\cos] + {\tfrac12(-1+\alpha^2)\cos(\theta)} \\
			&\quad - \tfrac12\alpha^3 \cos(2\theta)+ \alpha^4\big( -\tfrac14 \cos(\theta) +\tfrac18 \cos(3\theta) \big) + O(\alpha^5) \,, \label{exp.R.modecos}  \\
				\cR(r_\alpha;\alpha) [\sin](\theta) & = \tA[\sin] +{ \tfrac12(-1+\alpha^2)\sin(\theta)} \\
				& \quad - \tfrac12\alpha^3\sin(2\theta)+\alpha^4\big( \tfrac34\sin(\theta) +\tfrac38 \sin(3\theta)\big) + O(\alpha^5) \,. \label{exp.R.modesin}
		\end{align}
\end{lem}

\begin{proof}
    By \eqref{cR.lin.app} and \eqref{tg.alpha}, we obtain the decomposition
		\begin{align}\label{Split-R0-1}
				\nonumber\cR(r_\alpha;\alpha) [h](\theta)  & =  \fint_{\T} \log\big| (1+\tg_{\alpha}(\theta))e^{\im\theta} - (1+\tg_{\alpha}(\eta))e^{\im\eta}  \big| h(\eta) \wrt\eta \\
			\nonumber	& \quad -  \fint_{\T} K_\alpha \big( (1+\tg_{\alpha}(\theta))e^{\im\theta} , (1+\tg_{\alpha}(\eta)) e^{\im\eta} \big) h(\eta) \wrt\eta \\
				& =: \cR_{0} [h](\theta) + \cR_{1} [h](\theta) \,.
  		\end{align}
    Let us start with the expansion of the first operator \(\mathcal{R}_0\). Using the kernel expansion \eqref{BS.kernel.expand} from Lemma \ref{lemma.kernels}, together with a suitable change of variables, we infer
		\begin{align*}
			\cR_{0} [h] 
			& =  \fint_{\T} \log\big|1 - e^{\im\eta}  \big| h(\eta+\theta) \wrt\eta   + \frac{1}{2} \fint_{\T} \big( \tg_{\alpha}(\theta) + \tg_{\alpha}(\eta) \big) h(\eta) \wrt\eta \\
			& \ \ + \frac{1}{8} \fint_{\T} \big( \tg_{\alpha}(\theta)- \tg_{\alpha}(\eta) \big)^2 {\rm cotg}^2\big(\tfrac{\theta-\eta}{2}\big) h(\eta) \wrt\eta\\
            &\ \ - \frac{1}{8} \fint_{\T} \big( \tg_{\alpha}(\theta)+ \tg_{\alpha}(\eta) \big)^2 h(\eta) \wrt\eta     + \alpha^5 \tR_{0}[h] \,,
		\end{align*}
        for some $\tR_{0}\in\OpM_{s}^{-\infty}$. It can also be rewritten in the following form,
        \begin{align}
			\cR_{0} [h]
			& = \cR(0;0) [h]  + \tfrac{1}{8} \fint_{\T} \big( \tg_{\alpha}(\theta)- \tg_{\alpha}(\eta) \big)^2 {\rm cotg}^2\big(\tfrac{\theta-\eta}{2}\big) h(\eta) \wrt\eta \\
            & \ \ - \tfrac{1}{4} \tg_{\alpha}(\theta)\fint_{\T} \tg_{\alpha}(\eta)h(\eta) \wrt\eta+ \fint_{\T}\big( \tfrac12 \tg_{\alpha}(\eta)- \tfrac18 \tg_{\alpha}^2(\eta) \big) h(\eta) \wrt\eta + \alpha^5 \tR_{0}[h]\,.\label{polase.R0}
		\end{align}
	Recall that the action of the Fourier multiplier $\cR(0;0)$ is described in Lemma \ref{lemma.cR.00}. 
    Moreover, the term
\[
\fint_{\T}\Big( \tfrac12 \tg_{\alpha}(\eta)- \tfrac18 \tg_{\alpha}^2(\eta) \Big) h(\eta)\,\wrt\eta
\]
contributes only through a real constant. Therefore, there is no need to expand it further, since it disappears after applying the operator \(\partial_\theta\).
    On the other hand, applying  Lemma \ref{sqrt.1+ralpha} leads to
	\begin{align}
			- \tfrac{1}{4} \tg_{\alpha}(\theta&)\fint_{\T} \tg_{\alpha}(\eta)h(\eta) \wrt\eta  = -\tfrac{1}{4}\alpha^4 \cos(2\theta) \fint_{\T} \cos(2\eta) h(\eta) \wrt\eta + O(\alpha^5)[h] \\
			&=-\tfrac{1}{8} \alpha^4\fint_{\T}\big( \cos(2\theta-2\eta) + \cos(2\theta+2\eta) \big)h(\eta) \wrt\eta + O(\alpha^5)[h]  \,.\label{polase1}
	\end{align}
	By item $3.$ in Appendix \ref{app.trigo}, we have that
	\begin{equation}
		\begin{aligned}
			(	\cos(2\theta)&-\cos(2\eta))^2  {\rm cotg}^2\big(\tfrac{\theta-\eta}{2}\big)  
			 = 3 - \tfrac12\big( \cos(4\theta)+\cos(4\eta)\big) + 4 \cos(\theta-\eta)  \\
			&\quad + \cos(2\theta-2\eta) -3\cos(2\theta+2\eta) -2\big( \cos(3\theta+\eta)+\cos(\eta+3\eta) \big)\,.
		\end{aligned}
	\end{equation}
Then, exploiting the fact that \(h\) has vanishing spatial average, we obtain
	\begin{align}
		\tfrac{1}{8} & \fint_{\T} \big( \tg_{\alpha}(\theta)- \tg_{\alpha}(\eta) \big)^2 {\rm cotg}^2\big(\tfrac{\theta-\eta}{2}\big) h(\eta) \wrt\eta \\
		&= \tfrac{1}{8}\alpha^4 \fint_{\T} \big( \cos(2\theta)- \cos(2\eta) \big)^2 {\rm cotg}^2\big(\tfrac{\theta-\eta}{2}\big) h(\eta) \wrt\eta + O(\alpha^5)[h] \\
		& =-\tfrac{1}{16}\alpha^4 \fint_{\T} \cos(4\eta) h(\eta) \wrt\eta 
		+\tfrac{1}{2}\alpha^4\fint_{\T}\cos(\theta-\eta)h(\eta) \wrt\eta \\
		& \quad +\tfrac{1}{8}\alpha^4\fint_{\T}\big( \cos(2\theta-2\eta) -3\cos(2\theta+2\eta) \big)h(\eta) \wrt\eta \\
		& \quad -\tfrac{1}{4}\alpha^4 \fint_{\T} \big( \cos(3\theta+\eta) +\cos(\theta+3\eta) \big)h(\eta) \wrt\eta + \alpha^5 \tR_{0}[h]\,.\label{polase2}
	\end{align}
Collecting all the terms \eqref{polase1}, \eqref{polase2} into \eqref{polase.R0}, 
we conclude that
	\begin{align}\label{Expand-R-0}
		\cR_{0} [h] & = \cR(0;0)[h] +	 \fint_{\T}\big( \tfrac12 \tg_{\alpha}(\eta)- \tfrac18 \tg_{\alpha}^2(\eta) -\tfrac{1}{16} \alpha^4 \cos(4\eta) \big) h(\eta) \wrt\eta 
		\\
		& \quad+{\alpha^4}\fint_{\T}\tfrac12 \cos(\theta-\eta)h(\eta) \wrt\eta  -{\alpha^4}\fint_{\T}\tfrac12 \cos(2\theta+2\eta) h(\eta) \wrt\eta \\
		&  \quad -{\alpha^4} \fint_{\T} \tfrac14\big( \cos(3\theta+\eta) +\cos(\theta+3\eta) \big)h(\eta) \wrt\eta + \alpha^5 \tR_{0}[h] \,.
 	\end{align}
Next, we move to the asymptotic expansion of the operator $\cR_{1}$ defined in \eqref{Split-R0-1}. For this aim, we use \eqref{inter.kernel} from Lemma \ref{lemma.kernels} and recalling that $h(\theta)$ has zero average,
	\begin{align}
		\cR_{1} [h]    
		&=  \fint_{\T} \big( \tfrac12 \alpha^2 \cos(2\eta) -\tfrac13 \alpha^3 \cos(3\eta) + \tfrac14 \alpha^4 \cos(4\eta) +\alpha^2 \tg_{\alpha}(\eta) \cos(2\eta)  \big) h(\eta) \wrt\eta \\
		&\ \  + {\alpha^2}\fint_{\T}\cos(\theta-\eta) h(\eta) \wrt \eta - {\alpha^3} \fint_{\T}\big( \cos(2\theta-\eta) + \cos(\theta-2\eta) \big)h(\eta) \wrt \eta    \\
		& \ \   + \tfrac{1}{2} \alpha^4 \fint_{\T} \big( 2\cos(3\theta-\eta) + 3\cos(2\theta-2\eta) +2\cos(\theta-3\eta) \big) h(\eta) \wrt\eta   \\
		& \ \ + {\alpha^2} \fint_{\T} \big( \tg_{\alpha}(\theta) + \tg_{\alpha}(\eta)\big) \cos(\theta-\eta) h(\eta) \wrt \eta + \alpha^5 \tR_{1}[h]\,,
	\end{align}
for some $\tR_{1}\in \OpM_{s}^{-\infty}$. According to Lemma \ref{sqrt.1+ralpha} and classical trigonometric identities, one can expand
	\begin{align}
		{\alpha^2} \fint_{\T}& \big( \tg_{\alpha}(\theta) + \tg_{\alpha}(\eta) \big)\cos(\theta-\eta) h(\eta) \wrt \eta  = -{\alpha^4} \fint_{\T} \cos(2\theta) \cos(\theta-\eta) h(\eta) \wrt \eta  \\
        &\ \ -{\alpha^4} \fint_{\T} \cos(2\eta) \cos(\theta-\eta) h(\eta) \wrt \eta + O(\alpha^5)[h] \\
		&=-{\alpha^4} \fint_{\T} \big( \cos(\theta+\eta) + \tfrac12 \cos(\theta-3\eta) +\tfrac12 \cos(3\theta-\eta)\big) h(\eta) \wrt \eta  + \alpha^5 \tR_{1}[h] \,.
	\end{align}
Summing up, we obtain that
	\begin{align*}
	&	\cR_{1}[h]  =   \fint_{\T} \big( \tfrac12 \alpha^2 \cos(2\eta) -\tfrac13 \alpha^3 \cos(3\eta) + \tfrac14 \alpha^4 \cos(4\eta) +\alpha^2 \tg_{\alpha}(\eta) \cos(2\eta)  \big) h(\eta) \wrt\eta \\
		& \ \  + {\alpha^2}\fint_{\T}\cos(\theta-\eta) h(\eta) \wrt \eta - {\alpha^3} \fint_{\T}\big( \cos(2\theta-\eta) + \cos(\theta-2\eta) \big)h(\eta) \wrt \eta  \\
	& \ \    + {\alpha^4} \fint_{\T} \big( \tfrac12\cos(3\theta-\eta) + \tfrac32\cos(2\theta-2\eta) +\tfrac12 \cos(\theta-3\eta)  -\cos(\theta+\eta) \big) h(\eta) \wrt\eta \\
	& \ \  + \alpha^5 \tR_{1}[h]\,.
	\end{align*}
Therefore, substituting this expansion together with \eqref{Expand-R-0} into \eqref{Split-R0-1} yields the desired expansion \eqref{exp.R.total}.
As for the expressions in \eqref{exp.R.geq2}, \eqref{exp.R.modecos}, and \eqref{exp.R.modesin}, they are obtained by straightforward computations using \eqref{exp.R.total} and Lemma \ref{lemma.cR.00}.
\end{proof}
We conclude this section with the following expansion, obtained through direct computations based on Lemma \ref{speed.V.expan} and Lemma \ref{lemma.R.expand}.
\begin{cor}\label{cor.T.at.cos}
	The operator $\cT(r_\alpha;\alpha)$ in \eqref{T.alpha.operator} restricted on the mode one expands as follows,
	\begin{align}
		 	\cT(r_\alpha;\alpha)[\cos](\theta) &=   \tfrac12\tU_{\alpha}\big( 1+\cos(2\theta) \big) + \alpha^2 \big( \cos(\theta)+\tfrac12\cos(3\theta)  \big)   -  \alpha^3 \big(   \tfrac78 \cos(2\theta) +\tfrac38 \cos(4\theta)\big) \\
		 	& \ \ + \alpha^4\big(-\tfrac14\cos(\theta)+ \tfrac{31}{48} \cos(3\theta) + \tfrac{25}{48} \cos(5\theta) \big) + O(\alpha^5)  \\
		 	& = \alpha^2 \big( \cos(\theta)+ \tfrac12\cos(3\theta)  \big) + O(\alpha^3) \,,
	\end{align}
	and,
		\begin{align}
		\cT(r_\alpha;\alpha)[\sin](\theta) &= \tfrac12\tU_{\alpha} \sin(2\theta)+ \tfrac12 \alpha^2 \sin(3\theta) -  \alpha^3 \big(  \tfrac18 \sin(2\theta) + \tfrac38 \sin(4\theta)  \big)    \\
		& \ \ + \alpha^4\big( \tfrac34\sin(\theta) - \tfrac{7}{48} \sin(3\theta) + \tfrac{25}{48} \sin(5\theta)  \big) + O(\alpha^5) \\
		& = \tfrac12 \alpha^2 \sin(3\theta) + O(\alpha^3) \,.
	\end{align}
\end{cor}

\section{Spectral theory around traveling pairs}\label{sect.linear.eq}

Having established in the previous section the existence of a family of stationary solutions (traveling pairs) together with their asymptotic behavior, we now turn to the problem of identifying periodic-in-time solutions bifurcating from these stationary profiles. The first step in this construction is to analyze the structure of the linearized operator around the traveling pairs. As we shall see, the spectral properties of this operator already reveal the presence of time-periodic solutions at the linear level.
The main challenge, however, lies in proving that these solutions persist when nonlinear effects are taken into account. This is precisely where KAM techniques play a central role, allowing us to control the delicate interplay between small divisors and nonlinear interactions. 
Before exploring the spectral analysis of the linearized operator, we first need to reformulate the nonlinear problem in a framework centered at the traveling pairs. In this setting, the speed naturally plays the role of a Lagrangian multiplier, ensuring the elimination of any  component along the first sine mode. Once this reformulation is in place, we will proceed to describe the structure of the linearized operator at the traveling pairs and show how it can be diagonalized by means of specific, time-independent transformations. This step provides the essential groundwork for the subsequent construction of time-periodic solutions in a neighborhood of the traveling pairs.
\subsection{Nonlinear equation for the perturbation}
Proposition \ref{prop.ralpha} provides the existence of a unique analytic branch of traveling pairs for small values of $\alpha$. More precisely, it ensures the construction of a couple $(\tU_\alpha, r_\alpha)$ such that  
\begin{equation*}
   \forall \, \alpha \in [-\alpha_0,\alpha_0] \,, \qquad \cG(\tU_\alpha,r_\alpha;\alpha)=0 \,,
\end{equation*}
The next goal is to reformulate the problem in a neighborhood of these stationary solutions within a suitable functional framework. In particular, we introduce symmetry-adapted phase spaces and suitable orthogonality conditions allowing us to isolate the relevant dynamical directions and remove the degeneracies generated by the invariance of the system. This reformulation will provide the appropriate setting for the spectral analysis and the nonlinear perturbative arguments developed in the subsequent sections.
 Throughout the remainder of the paper, we shall use \(\alpha\) as the main bifurcation parameter and assume that it belongs to a sufficiently small interval of the form $[\alpha_{1},\alpha_{2}]$ such that
\begin{equation}\label{alpha1alpha2}
	  \alpha_1:=\tfrac14 \alpha_{0} ,\quad  \alpha_{2} := \tfrac34 \alpha_{0} \,\quad\hbox{and}\quad  0<\alpha_{0} \ll 1.
\end{equation}
This choice is made for convenience, and the results remain valid for any interval of the form $[\alpha_1,\alpha_2]$ with $0<\alpha_1<\alpha_2\ll1.$ 
Moving $\alpha$ as a parameter in $[\alpha_{1},\alpha_{2}]$ will allow us later to impose non-resonance conditions in the construction of the periodic solutions.
The choice of the interval $[\alpha_{1},\alpha_{2}]$ in \eqref{alpha1alpha2} ensures that, for any $0<\alpha\leq \alpha_{2}<\alpha_{0}$, there exists a nontrivial traveling pair $(\tU_{\alpha},r_{\alpha})$, whereas $\alpha\geq \alpha_{1}=\tfrac14 \alpha_{0}$ is to fix its size of order $\alpha_{0}$ and to search for time-periodic bifurcations around this equilibrium of size $\varepsilon>0$ much smaller than $\alpha_{0}$.

\medskip

\noindent
For \(0<\varepsilon\ll \alpha_{0}\) sufficiently small, we seek reversible solutions \((\tU,r)\) to \eqref{Func.G} of the form
\begin{equation}\label{perturbed.solution}
	\tU(\vf) = \tU_{\alpha} + \varepsilon \tW(\vf) \,, \quad r(\vf,\theta) = r_\alpha(\theta) + \varepsilon \rho(\vf,\theta)\,,
\end{equation}
where $\rho\in H^{s}_{\circ,\textnormal{even}}$ can be decomposed, recalling the definition of $\Pi_{\geq 2}$ in \eqref{proj.geq2}, as
\begin{equation}
	\rho(\vf,\theta) =\rho_{1}(\vf) \cos(\theta) + \rho_{\geq 2}(\vf,\theta) \,, \quad \rho_{\geq 2}(\vf,\theta)= \Pi_{\geq 2} \rho_{\geq 2}(\vf,\theta)\,,
\end{equation}
 under the following symmetry conditions
\begin{equation}
	\tW(\vf) = \tW(-\vf)\,, \quad \rho_{1}(\vf)=\rho_{1}(-\vf)\,, \quad \rho_{\geq 2}(\vf,\theta) = \rho_{\geq 2}(-\vf,-\theta)\,.
\end{equation}
Then, solving equation \eqref{Func.G} in a neighborhood of the traveling pairs $(\tU_{\alpha},r_\alpha)$ is equivalent to solving the following problem in terms of $(\tW, \rho)$,
\begin{equation}\label{full.nonlin.eq.pert}
	\wt\cG (\tW,\rho;\alpha,\varepsilon) := \varepsilon^{-1} \cG(\tU_{\alpha}+ \varepsilon \tW, r_\alpha+ \varepsilon \rho;\alpha)=0\,.
\end{equation}
To make sure that
\begin{equation}\label{phase.space.cond}
	\wt\cG:  H_{\rm even}^{s}(\T)\times H^{s}_{\circ,\textnormal{even}}\to H^{s-1}_{\circ,\textnormal{odd}}\,,
\end{equation}
where the sine mode one is absent from the range, we need to use $\tW(\vf)$ as a Lagrange multiplier 
by imposing that
\begin{equation}
	\int_{\T} \wt\cG(\tW,\rho;\alpha,\varepsilon)(\vf,\eta) \sin(\eta) \wrt\eta = 0 \,,
\end{equation}
that is
\begin{equation*}
	\int_{\T} \Big( \big(\tU_{\alpha}+\varepsilon\tW(\vf)\big) \sqrt{1+2(r_\alpha(\eta)+\varepsilon\rho(\vf,\eta))} \cos(\eta) + F(r_\alpha+\varepsilon\rho;\alpha)  \Big) \cos(\eta) \wrt\eta = 0 \,.
\end{equation*}
By virtue of \eqref{fixed.constrain}, this equation is equivalent to 
\begin{equation}\label{tU1.def}
		    \begin{aligned}
		         \tW(\varphi)=\tW[\rho](\vf) = - \frac{1}{\varepsilon}\Bigg( \frac{\int_{\T} F(r_\alpha+\varepsilon\rho;\alpha)\cos(\eta)\wrt\eta}{\int_{\T}\sqrt{1+2(r_\alpha(\eta)+\varepsilon\rho(\vf,\eta))}\cos^2(\eta)\wrt\eta} -  \frac{\int_{\T} F(r_\alpha;\alpha)\cos(\eta)\wrt\eta}{\int_{\T}\sqrt{1+2r_\alpha(\eta)}\cos^2(\eta)\wrt\eta}  \Bigg)\,.
		    \end{aligned}
	\end{equation}
We intend to establish  the following result.
\begin{pro}\label{prop.U1}
	Let $k\in\mathbb{N}$ and $s\geqslant s_0.$  Then the map
	\begin{equation}
		\tW:  H_{\circ,\textnormal{even}}^s(\T^2) \mapsto H^s_\textnormal{even}(\T)\,, \quad  \rho \mapsto \tW[\rho]\,,
	\end{equation}
	defined in \eqref{tU1.def} is of class $C^k$ and satisfies
	\begin{equation}\label{U10}
		\tW[0] = 0,\quad
		\di_{\rho} \tW[\overline \rho][h] = - \frac{\int_{\T} \cT_{\varepsilon,\alpha}[h](\vf,\eta) \cos(\eta) \wrt\eta }{ \int_{\T} \sqrt{1+ 2\overline{r}(\vf,\eta)} \cos^2(\eta) \wrt\eta} \,,
	\end{equation}
	where
	\begin{equation}
		\cT_{\varepsilon,\alpha} := \tfrac{\tU_{\alpha}+\varepsilon \tW[\overline\rho]}{\sqrt{1+ 2 \overline{r}}} \cos(\theta)+ \di_{r} F(\overline{r};\alpha) \,, \quad  \overline{r}(\vf,\theta):= r_\alpha(\theta) + \varepsilon \overline{\rho}(\vf,\theta) \,.
	\end{equation}
\end{pro}
\begin{proof}
   The regularity of $\tW$ is deduced in a standard manner from \eqref{tU1.def}, \eqref{full.nonlin.eq.pert} and \mbox{Proposition \ref{prop-Functional-est}}. Its symmetry follows from the fact that
   $$
   F(r_\alpha+\varepsilon\rho;\alpha)(-\varphi,-\theta)=F(r_\alpha+\varepsilon\rho;\alpha)(\varphi,\theta)\,,
   $$
   together with change of variables.
   The corresponding linear operator $\di_{\rho}\tW[\overline \rho]$ is obtained by a straightforward linearization
   from \eqref{tU1.def}, noting that
    \begin{equation}
        - \frac{\int_{\T} F(r_\alpha+\varepsilon \overline \rho;\alpha) \cos(\eta) \wrt\eta}{\int_{\T}\sqrt{1+2 (r_\alpha(\eta) +\varepsilon \overline \rho(\vf,\eta))} \cos^2(\eta) \wrt\eta} = \tU_{\alpha} +\varepsilon \tW[\overline \rho](\varphi)\,,
    \end{equation}
and using Proposition \ref{prop-Functional-est}. This ends the proof of the claimed result.
 \end{proof}
After this reduction, and clearly having $\partial_\varphi r_\alpha=0$, we reformulate equation \eqref{full.nonlin.eq.pert} as
\begin{equation}\label{equation.rho.full}
			\begin{aligned}
			&\tG(\rho;\alpha,\varepsilon)  := \wt\cG(\tW[\rho],\rho;\alpha,\varepsilon)  \\
			&	= \omega\,\pa_{\vf} \rho + \varepsilon^{-1} \pa_{\theta} \Big( \big(\tU_{\alpha}+\varepsilon\tW[\rho]\big) \sqrt{1+2(r_\alpha+\varepsilon\rho)} \cos(\theta) + F(r_\alpha+\varepsilon\rho;\alpha) \Big) = 0 \,.
		\end{aligned}
\end{equation}
\subsection{Linearized operator at the equilibrium}
The main goal of this section is to analyze the functional \(\tG\) defined in \eqref{equation.rho.full}. In particular, we identify its linear part around the equilibrium configuration and describe the structure of the corresponding linearized operator.\\ 
In view  of Proposition \ref{prop.ralpha} and \eqref{equation.rho.full}, we obtain
  $$\tG(0;\alpha,\varepsilon)=0 \,.
$$
For $s\in\R$ and $\delta>0$ we define the closed ball
\begin{equation*}
	B_{\delta,\textnormal{even}}^{s} := \big\{ f \in H^s_{\circ,\textnormal{even}}(\T^2) \, : \, \| f\|_{s} \leqslant \delta \big\}\,.
\end{equation*}

\begin{pro}\label{lemma.constrain-Lun1}
	Let $k\in\N$ and $s\geq 1$. There exist $\delta>0,\, \alpha_{0}>0$ and  $\varepsilon_{0}>0$ such that the following results  hold:
    \\[1mm]
    \noindent $(i)$ The map
    \begin{equation}
		\tG: B_{\delta,\textnormal{even}}^{s} \times [\alpha_{1},\alpha_{2}]\times [-\varepsilon_{0},\varepsilon_{0}] \to H^{s-1}_{\circ,\textnormal{odd}}(\T^2)
	\end{equation}
	is  of class $\cC^k$;
    \\[1mm]
    \noindent $(ii)$  For any $\overline{\rho} \in B_{\delta,{\rm even}}^{s}$, the linearized operator $\partial_\rho\tG(\overline{\rho};\alpha,\varepsilon)$ takes the form
    \begin{align}
	& \partial_\rho\tG(\overline{\rho};\alpha,\varepsilon) :  H_{\circ,{\rm even}}^{s}(\T^2)\to H_{\circ,{\rm odd}}^{s-1}(\T^2)\,, \\ &h \mapsto \partial_\rho\tG(\overline{\rho};\alpha,\varepsilon)[h]=\omega\,\pa_{\vf} h + \pa_{\theta} \Big( \tF_{\varepsilon,\alpha}(\vf,\theta) \cK_{\varepsilon,\alpha}[h] + \cT_{\varepsilon,\alpha}[h] \Big)\,,
\end{align}
where $\cT_{\varepsilon,\alpha}$ is defined in Proposition \ref{prop.U1} and 
$$
		\tF_{\varepsilon,\alpha}(\vf,\theta) := - \tfrac{\sqrt{1+2(r_\alpha+\varepsilon \overline{\rho}(\vf,\,\cdot\,))(\theta)}\cos(\theta)}{\int_{\T}\sqrt{1+2(r_\alpha+\varepsilon \overline{\rho}(\vf,\,\cdot\,))(\eta)}\cos^2(\eta)\wrt\eta}  \,,\quad
		\cK_{\varepsilon,\alpha}[h](\vf)  := \int_{\T} \cos(\eta) \cT_{\varepsilon,\alpha}[h](\vf,\eta) \wrt\eta \,.
	$$
    \end{pro}
\begin{proof}
The proof of items $(i)$ and $(ii)$ follows directly from \eqref{equation.rho.full}, \eqref{full.nonlin.eq.pert}, Proposition \ref{prop.U1}, and Proposition \ref{prop-Functional-est}. In particular, we note that 
    \begin{equation}
        \tF_{\varepsilon,\alpha}(\theta) \cK_{\varepsilon,\alpha}[h] = \sqrt{1+2(r_{\alpha}+\varepsilon\overline{\rho})} \,\cos(\theta) \di_{\rho} \tW[\overline{\rho}][h] \,,
    \end{equation}
    with $\di_{\rho} \tW[\overline{\rho}][h]$ as in \eqref{U10}.
\end{proof}
From Proposition \ref{lemma.constrain-Lun1}, we infer that, upon linearizing equation \eqref{equation.rho.full} around the equilibrium state \(\rho=0\), we obtain the linear equation
\begin{equation}\label{linear.rho.full}
	\partial_\rho\tG(0;\alpha,\varepsilon)[h]=\omega\,\pa_{\vf} h + \pa_{\theta} \Big( \tF_{0,\alpha}(\theta) \cK_{0,\alpha}[h] + \cT_{0,\alpha}[h] \Big) = 0 \,,
\end{equation}
where, see also \eqref{T.alpha.operator}, \eqref{wttV}, \eqref{V.at.boundary}, \eqref{cR.lin.app},
\begin{equation}\label{cT.0.alpha}
	\cT_{0,\alpha}  =  \tfrac{\tU_{\alpha} \cos(\theta)}{\sqrt{1+2r_\alpha(\theta)}} + \tV[r_\alpha;\alpha] + \cR(r_\alpha;\alpha) =: \tV_{0,\alpha}(\theta) + \cR(r_\alpha;\alpha)\,,
\end{equation}
and
	\begin{equation}\label{smoothF0}
	\begin{aligned}
		\tF_{0,\alpha}(\theta) & := - \frac{\sqrt{1+2r_\alpha(\theta)}\cos(\theta)}{\int_{\T}\sqrt{1+2r_\alpha(\eta)}\cos^2(\eta)\wrt\eta}  \,,
	\end{aligned}
\end{equation}
\begin{equation}\label{smoothK0}
	\begin{aligned}
		\cK_{0,\alpha}[h](\vf) & := \int_{\T} \cos(\eta) \cT_{0,\alpha}[h](\vf,\eta) \wrt\eta \,.
	\end{aligned}
\end{equation}
We note that the function $\tF_{0,\alpha}$ is even in $\theta$, while $\cK_{0,\alpha}[h]$ is even in $\varphi$. Furthermore, the spatial operator
\begin{equation}
\pa_{\theta}\circ \Big( \tF_{0,\alpha}(\theta)\cK_{0,\alpha} + \cT_{0,\alpha} \Big)
= \pa_{\theta}\circ \Big( \tF_{0,\alpha}(\theta)\cK_{0,\alpha}+\tV_{0,\alpha}(\theta) + \cR(r_\alpha;\alpha) \Big)
\end{equation}
has variable coefficients and is not diagonal. The spectral analysis of this operator is essential for determining whether its kernel elements, which provide the first approximation to the nonlinear equation \eqref{equation.rho.full}, can generate time-periodic solutions. Because this analysis is technically demanding, obtaining a complete description of the spectrum requires constructing a sequence of transformations that diagonalize the operator when acting on the phase space \mbox{space $H^s_{\circ,\textnormal{even}}(\T^2)$} introduced in \eqref{Hs-circ}.

\medskip

\noindent
For later use, we have the following intermediate result.
\begin{lem}\label{sa.pp.T.eps}
Let $\varepsilon$ and $\overline{r}\in H^s_{\circ,\textnormal{even}}(\T^2)$ be sufficiently small. Then the operator $\cT_{\varepsilon,\alpha}$, defined in Proposition $\ref{prop.U1}$, is self-adjoint on $H^0_{\circ,\textnormal{even}}(\T^2)$ and parity preserving.  
\end{lem}
\begin{proof}
	The claim of self-adjointness follows from the fact that the function $\tV_{\varepsilon,\alpha}$ is real and that the kernel of the operator $\cR(\overline r;\alpha)$ in \eqref{cR.lin.app} is real and symmetric. 	
	We now prove that $\cT_{\varepsilon,\alpha}$ is parity preserving.
	Since $\overline r = {\rm even}(\vf,\theta)$, we clearly have $\sqrt{1+2\overline r} ={\rm even}(\vf,\theta)$. It follows immediately that $ \frac{\overline \tU(\vf) \cos(\theta)}{\sqrt{1+2\overline r(\vf,\theta)}}={\rm even}(\vf,\theta)$. The claim that $\tV_{\varepsilon,\alpha}={\rm even}(\vf,\theta)$ follows by directly checking the parity condition on $\tV[\overline r;\alpha](\vf,\theta)$. Indeed, we have $\tV[\overline r;\alpha]={\rm even}(\vf,\theta)$ by \eqref{V.at.boundary},  the assumption $\overline r = {\rm even}(\vf,\theta)$ and using that
	\begin{equation}
		\begin{aligned}
			&\log \big| \sqrt{1+2\overline r(\vf,\theta)}e^{-\im\theta} - \sqrt{1+2\overline r(\vf,\eta)}e^{-\im\eta} \big| = \log \big| \sqrt{1+2\overline r(\vf,\theta)}e^{\im\theta} - \sqrt{1+2\overline r(\vf,\eta)}e^{\im\eta} \big| \,, \\
			& K_{\alpha} \big( \sqrt{1+2\overline r(\vf,\theta)}e^{-\im\theta}, \sqrt{1+2\overline r(\vf,\eta)}e^{-\im\eta} \big) = K_{\alpha} \big( \sqrt{1+2\overline r(\vf,\theta)}e^{\im\theta}, \sqrt{1+2\overline r(\vf,\eta)}e^{\im\eta} \big) \,,
		\end{aligned}
	\end{equation}
	with $K_\alpha$ as in \eqref{Kalpha.kernel}. The claim that $\cR(\overline r;\alpha)$ is parity preserving follows by the same arguments. This concludes the proof.
\end{proof}

\subsection{Reduction of the transport operator}\label{Redu-transpo}
We now construct the first transformation in the diagonalization procedure for the linearized operator \eqref{linear.rho.full}.  
To this end, for each $|\alpha|\leq \alpha_{0}$, we introduce a family of $\alpha$-dependent diffeomorphisms of the torus $\T$, defined by
\begin{equation}
	 \theta\in \T \mapsto \vartheta := \theta + \beta(\theta;\alpha) \ \leftrightsquigarrow \ \theta = \vartheta + \breve{\beta}(\vartheta;\alpha)\,,
\end{equation}
where $\vartheta\mapsto  \theta:=\vartheta+\breve{\beta}(\vartheta;\alpha)$ is the inverse diffeomorphism and the real function $\theta\in\T\mapsto \beta(\theta;\alpha)\in\R$ is defined by the relation
\begin{equation*}
	\tV_{0,\alpha}(\theta)(1+ \beta'(\theta;\alpha)) = \tm_{1}(\alpha)\,,
\end{equation*}
with $\tV_{0,\alpha}(\theta)$ the function defined in \eqref{cT.0.alpha}, see also Lemma \ref{speed.V.expan}. The constant $\tm_{1}(\alpha)$ is to be determined. 
Then this equation can be solved  in the periodic setting as follows
\begin{equation}\label{defn.beta}
	 \beta(\theta)= \beta(\theta;\alpha) := \pa_{\theta}^{-1} \bigg(  \frac{\tm_{1}(\alpha)}{\tV_{0,\alpha}(\theta)}- 1 \bigg)\,,
\end{equation}
 provided that the constant $\tm_{1,\alpha}$ is given by
\begin{equation}\label{m1.alpha}
	\tm_{1}(\alpha):=\big( \langle 
	\tV_{0,\alpha}^{-1} \rangle_{\theta}\big)^{-1} \in \R\,. 
\end{equation}
We remark that the invertibility of the map
$\theta \mapsto \theta + \beta(\theta;\alpha)$
requires ensuring that
\begin{equation}
	\inf_{\theta \in \T} (1 + \beta'(\theta;\alpha)) = \inf_{\theta\in\T}\bigg(  \frac{\tm_{1}(\alpha)}{\tV_{0,\alpha}(\theta)} \bigg) := c_0 > 0\,.
\end{equation}
This latter property occurs provided that  $|\alpha|\leq \alpha_{0} $ small enough, according to  Lemma \ref{speed.V.expan} and using the fact that $ \tV_{0,\alpha}(\theta) = \tfrac12 + O(\alpha^2) >0$.
Therefore, we deduce that the mapping $\T\ni \theta\mapsto \theta+\beta(\theta)$ is invertible for any $|\alpha| \leq  \alpha_{0}$.
Now, we  define the symplectic of change of coordinates  operator $\cS$  and  its inverse $\cS^{-1}$
\begin{equation}\label{S.diffeo}
	\cS := (1+ \beta'(\theta)) \circ \cB \,, \quad \cS^{-1} := \cB^{-1} \circ (1 + \beta'(\theta))^{-1},
\end{equation}
where $\cB$ is the
composition map and $\cB^{-1}$ its inverse, given by
\begin{equation}\label{B.diffeo}
	\cB u(\theta):= u(\theta + \beta(\theta)) \,, \quad \cB^{-1} u(\vartheta) := u(\vartheta+ \breve{\beta}(\vartheta))\,,
\end{equation}
The main goal of this section is to derive asymptotic expansions of the quantities introduced above in the regime of small $\alpha.$

\begin{lem}\label{lemma.reparam.0}
	Let $s\geq s_0$. The following hold:
	\\[1mm]
	\noindent $(i)$ We have the following estimates
	\begin{align}
		&	\big| \tm_{1} - \big( \tfrac12 -  \alpha^4 \big) \big|^{5,\alpha_{0}} \lesssim \alpha_{0}^5 \,, \label{m1alpha0.est} \\
		&	\big\| \beta -\big( -2\tU_{\alpha} \sin(\theta) -\alpha^2 \sin(2\theta) +\tfrac12 \alpha^3 \sin(3\theta) + \tfrac{1}{6} \alpha^4 \sin(4\theta) \big)\ \big\|_{s}^{5,\alpha_{0}}  \lesssim_{s} \alpha_{0}^5 \,, \label{beta0.est} \\
		&	\big\| \breve\beta - \big(2\tU_{\alpha} \sin(\theta) + \alpha^2 \sin(2\theta) - \tfrac12 \alpha^3 \sin(3\theta) + \tfrac{5}{6} \alpha^4 \sin(4\theta) \big)\ \big\|_{s}^{5,\alpha_{0}}  \lesssim_{s} \alpha_{0}^5 \,. \label{brevebeta0.est} 
	\end{align}
	Moreover, $\beta(\theta)={\rm odd}(\theta)$ and $\breve{\beta}(\vartheta)={\rm odd}(\vartheta)$;
	\\[1mm]
	\noindent $(ii)$
    The maps $\cS$ and $\cB$, defined by \eqref{S.diffeo}, \eqref{B.diffeo}, and \eqref{defn.beta}, are bounded and invertible operators. Moreover, for any $h\in H^s(\T^2)$, they satisfy the estimates
	\begin{align}
		&	\| \cS^{\pm 1} h \|_{s}^{5,\alpha_{0}}  \lesssim_{s} \| h \|_{s}^{5,\alpha_{0}} \,, \quad 	\| (\cS^{\pm 1} -{\rm Id}) h \|_{s}^{5,\alpha_{0}}  \lesssim_{s} \alpha_{0}^2  \| h \|_{s+1}^{5,\alpha_{0}} \label{bounded.S.est}\,, \\
        &	\| \cB^{\pm 1} h \|_{s}^{5,\alpha_{0}}  \lesssim_{s} \| h \|_{s}^{5,\alpha_{0}} \,, \quad 	\| (\cB^{\pm 1} -{\rm Id}) h \|_{s}^{5,\alpha_{0}}  \lesssim_{s} \alpha_{0}^2  \| h \|_{s+1}^{5,\alpha_{0}} \label{bounded.B.est}\,.
	\end{align}
	In particular, we write
	\begin{equation}\label{cS.minus.Id}
		\cS^{\pm 1} = {\rm Id} + \cP_{\pm} \,, \quad \cB^{\pm 1} = {\rm Id} + \cQ_{\pm} \,,
	\end{equation}
	where the operators $\cP_{\pm}$ satisfy the estimates, for any $h\in H^s(\T^2)$,
	\begin{align}
	& 	\big\| \big(\cP_{+} -\big( -\alpha^2\pa_{\theta}\circ \sin(2\theta) \big) \big)h \big\|_{s}^{5,\alpha_{0}}  \lesssim_{s} \alpha_{0}^3 \| h \|_{s+2}^{5,\alpha_{0}} \,, \label{S.expand.est}  \\
	& 	\big\| \big( \cP_{-} -\alpha^2 \pa_{\vartheta}\circ \sin(2\vartheta) \big)h \big\|_{s}^{5,\alpha_{0}}  \lesssim_{s} \alpha_{0}^3 \| h \|_{s+2}^{5,\alpha_{0}} \,, \label{invS.expand.est} \\
    & 	\big\| \big(\cQ_{+} -\big( -\alpha^2\sin(2\theta)\pa_{\theta} \big) \big)h \big\|_{s}^{5,\alpha_{0}}  \lesssim_{s} \alpha_{0}^3 \| h \|_{s+2}^{5,\alpha_{0}} \,, \label{B.expand.est}  \\
	& 	\big\| \big( \cQ_{-} -\alpha^2 \sin(2\vartheta)\pa_{\vartheta} \big)h \big\|_{s}^{5,\alpha_{0}}  \lesssim_{s} \alpha_{0}^3 \| h \|_{s+2}^{5,\alpha_{0}} \,. \label{invB.expand.est}
\end{align}
	In addition, $\cS^{\pm 1}$ and $\cB^{\pm 1}$ are real and  reversibility preserving;
	\\[1mm]
	\noindent $(iii)$ We have that
	\begin{equation}
		\cS^{-1} \pa_{\theta}\circ \tV_{0,\alpha}(\theta) \,\cS = \pa_{\vartheta}\circ \tm_{1,\alpha} \,.
	\end{equation}
\end{lem}
\begin{proof}
	We start with item $(i)$.   
    The estimate \eqref{m1alpha0.est} for $\tm_{1,\alpha}$ and the estimate \eqref{beta0.est} for $\beta$ follow by \eqref{m1.alpha}, \eqref{defn.beta} and \eqref{inv.speed.expan} in Lemma \ref{speed.V.expan}. 
    Alternatively, we search for their expansion from equation $\tV_{0,\alpha}(\theta)(1+\beta'(\theta;\alpha))=\tm_{1}(\alpha)$ in \eqref{defn.beta}. That is, we search for an expansion of $\beta(\theta;\alpha)$ of the form $\beta(\theta;\alpha)= \tb_{2}(\theta) \alpha^2 + \tb_{3}(\theta)\alpha^3 + \tb_{4}(\theta;\alpha)\alpha^4 + O(\alpha^5)$ and we insert it in \eqref{defn.beta}, with $\tm_{1,\alpha}$ obtained by fixing
    \begin{equation}\label{average.m1}
        \tm_{1}(\alpha) = \braket{\tV_{0,\alpha}(\,\cdot\,)(1+\beta'(\,\cdot\,;\alpha))}_{\theta} \,.
    \end{equation}
    By Lemma \ref{speed.V.expan}, we compute
    \begin{align}
        \tV_{0,\alpha}(\theta)(1+\beta'(\theta;\alpha)) & = \tfrac12 + \alpha^2 \big( \tfrac12 \tb_{2}'(\theta) + \cos(2\theta) \big)+ \alpha^3\big( \tfrac12 \tb_{3}'(\theta) - \tfrac34 \cos(3\theta) \big) \\
        & +\alpha^4 \big( \tfrac12 \tb_{4}'(\theta;\alpha) +\alpha^{-4}\tU_{\alpha}\cos(\theta) +\tfrac{2}{3}\cos(4\theta) + \tb_{2}'(\theta) \cos(2\theta) \big) + O(\alpha^5)\,.
    \end{align}
    We immediately obtain that
    \begin{equation}
        \tb_{2}(\theta) = - \sin(2\theta) \,, \quad \tb_{3}(\theta) = \tfrac12 \sin(3\theta) \,,
    \end{equation}
    which also implies that 
    $$\tb_{2}'(\theta)\cos(2\theta) = - 2 \cos^2(2\theta)=-1-\cos(4\theta)\,.
    $$
    From \eqref{average.m1}, it follows that $\tm_{1}(\alpha) = \tfrac12 - \alpha^4 + O(\alpha^5) $ and we determine
    \begin{equation}
        \tb_{4}'(\theta)= -2\alpha^{-4}\tU_{\alpha} \cos(\theta) - \tfrac{4}{3}\cos(4\theta) + 2\cos(4\theta) \ \Rightarrow \tb_{4}(\theta) = -2\alpha^{-4}\tU_{\alpha} \sin(\theta) + \tfrac{1}{6} \sin(4\theta) \,.
    \end{equation}
    This concludes the proof of \eqref{m1alpha0.est} and \eqref{beta0.est}.
    \\
    To prove the estimate \eqref{brevebeta0.est} for $\breve{\beta}(\vartheta)=\breve{\beta}(\vartheta;\alpha)$, we search for expansion of the form $$\breve{\beta}(\vartheta)= \breve\tb_{2}(\vartheta) \alpha^2 + \breve\tb_{3}(\vartheta) \alpha^3 + \breve\tb_{4}(\vartheta;\alpha)\alpha^4 + O(\alpha^5)
    $$ by the implicit definition of $\breve{\beta}(\vartheta)$, which is
	\begin{equation}\label{implicit.breve}
		\breve{\beta}(\vartheta) = - \beta\big(\vartheta+ \breve{\beta}(\vartheta)\big) \,.
	\end{equation}
	By expanding \eqref{implicit.breve} in Taylor series around $\vartheta$ and in powers of $\alpha$, together with \eqref{beta0.est}, we find that
	\begin{equation}
		\begin{aligned}
			\breve\tb_{2}(\vartheta) & =- \big(- \sin(2\vartheta)\big)\,, \quad \breve\tb_{3}(\vartheta) = - \tfrac12 \sin(3\vartheta) \,, \\
			 \breve\tb_{4}(\vartheta) & = - \big(-2\alpha^{-4}\tU_{\alpha} \sin(\vartheta) + \tfrac{1}{6} \sin(4\vartheta) \big)  - \breve{\tb_{2}}(\vartheta)  \pa_{\vartheta}(-\sin(2\vartheta))  \\
		&	= 2\alpha^{-4}\tU_{\alpha} \sin(\vartheta) +  \big(-\tfrac{1}{6}+1) \sin(4\vartheta) \,,
		\end{aligned}
	\end{equation}
and we conclude the desired expansion.\\	In particular, the claimed estimate in \eqref{brevebeta0.est} follows by an implicit function argument on the $O(\alpha^5)$-remainder for the expansion of $\breve{\beta}(\vartheta)$.
	\\
	We now move to item $(ii)$. The estimates \eqref{bounded.S.est} and \eqref{bounded.B.est} follow by \eqref{S.diffeo}, \eqref{B.diffeo}, Lemma \ref{lemma.diffeo.gen} and estimates \eqref{beta0.est}, \eqref{brevebeta0.est}.  To prove \eqref{S.expand.est} and \eqref{B.expand.est}, first we note that, by Taylor expansion, we have that
	\begin{equation}\label{B.Taylor}
	\cB = {\rm Id} + \beta(\theta) \pa_{\theta} + \beta^2(\theta)  \int_{0}^{1} (1-\tau) \cB_{\tau}\circ \pa_{\theta}^2 \wrt\tau \,, \quad \cB_{\tau}h(\theta):= h(\theta+\tau\beta(\theta)) \,.
\end{equation}
	Consequently, we obtain that
	\begin{align}
		\cS & = (1+\beta ') \Big( {\rm Id} + \beta \pa_{\theta} + \beta ^2 \int_{0}^{1} (1-\tau) \cB_{\tau}\circ \pa_{\theta}^2 \wrt\tau \Big) \\
		& = {\rm Id } + \beta'+\beta \pa_{\theta} + \beta \beta' \pa_{\theta}  + (1+\beta') \beta ^2 \int_{0}^{1} (1-\tau) \cB_{\tau}\circ \pa_{\theta}^2 \wrt\tau \\
		& = {\rm Id} + \pa_{\theta} \circ \beta +\beta\beta'\pa_{\theta} + (1+\beta') \beta ^2 \int_{0}^{1} (1-\tau) \cB_{\tau}\circ \pa_{\theta}^2 \wrt\tau\,.\label{S.Taylor}
	\end{align}
	Therefore, \eqref{S.expand.est} and \eqref{B.expand.est} follow by \eqref{B.Taylor}, \eqref{S.Taylor}, Lemma \ref{lem-productlaw}-$(ii)$ and estimate \eqref{beta0.est}. 
	Finally, to establish \eqref{invS.expand.est} and \eqref{invB.expand.est}, we observe that, by \eqref{S.diffeo} and the fact that
    	\begin{equation}\label{deri.brevebeta}
		1+ \breve{\beta}'(\vartheta)= \frac{1}{1+\beta'(\vartheta +\breve{\beta}(\vartheta))} = \cB^{-1} \Big\{ \frac{1}{1+\beta'(\theta)} \Big\}\,,
	\end{equation}
     we have that
	\begin{equation}
		\cS^{-1} = \cB^{-1} \circ (1+\beta')^{-1} = (1+\breve{\beta}') \circ \cB^{-1}.
	\end{equation}
	Hence, by arguing as \eqref{B.Taylor} and \eqref{S.Taylor}, we obtain that
	\begin{align}
    \cB^{-1} &= {\rm Id} + \breve\beta(\vartheta) \pa_{\vartheta} + \breve\beta^2(\vartheta)  \int_{0}^{1} (1-\tau) \cB_{\tau}^{-1}\circ \pa_{\vartheta}^2 \wrt\tau \,, \\
		\cS^{-1}  & = {\rm Id} + \pa_{\vartheta} \circ \breve\beta + \breve{\beta}\breve{\beta}'\pa_{\vartheta}+ (1+\breve\beta') \breve\beta ^2 \int_{0}^{1} (1-\tau) \cB_{\tau}^{-1}\circ \pa_{\vartheta}^2 \wrt\tau\,,
	\end{align}
	from which we deduce \eqref{invS.expand.est} and \eqref{invB.expand.est}, using also Lemma \ref{lem-productlaw}-$(ii)$ and estimate \eqref{brevebeta0.est}.	
	\\
	As for item $(iii)$, it follows straightforwardly from \eqref{S.diffeo}, \eqref{B.diffeo}, and \eqref{defn.beta}. This ends  the proof of the claimed results.
\end{proof}

The next goal is to give an expansion of the conjugation of the remainder $\cS^{-1} \pa_{\theta}\cR(0;0)\cS$, that will be used later.
The following notation will be adopted throughout this section, in order to streamline the presentation.

\smallskip 

\noindent {{\bf Notation.} Given two operators  $A=A(\alpha),B=B(\alpha) \in \OpM_{s}^{m}$, $\alpha\in[\alpha_{1},\alpha_{2}]$,  we will write
\begin{equation}\label{notation.tR}
	A(\alpha) = B(\alpha) + \alpha^c \tR\,, \quad c \in \N\,,
\end{equation}
where $\tR=\tR(\alpha)$ is a smoothing operator in $\OpM_{s}^{-\infty}$ that we do not explicitly define from step to step, with estimates $|\tR|_{-N,s}^{5,\alpha_0}\leq C$ for any $N\in\N_{0}$ and  for some constant $C=C(N,s)>0$ independent of $\alpha_{2}$. For the definition of this topology, see Definition \ref{block norm}.}

\begin{lem}\label{lemma.conjR00}
	For $s\geq s_0$, we have that
	\begin{equation}
 \cE_{\cR(0;0)}:= \cS^{-1} \pa_{\theta} \cR(0;0) \cS - \pa_{\vartheta}\cR(0;0) \in \OpM_{s}^{-\infty}\,,
	\end{equation}
	with estimates, for any $N\in\N_{0}$,
	\begin{equation}
		\begin{aligned}
			&	|\cE_{\cR(0;0)}  |_{-N,s}^{5,\alpha_{0}} \lesssim_{N,s} \alpha_{0}^2 \,, \quad  |   \cE_{\cR(0;0)}- \cI_{\cR(0;0)}  |_{-N,s}^{5,\alpha_{0}} \lesssim_{N,s} \alpha_{0}^5\,,
		\end{aligned}
	\end{equation}
	where $\cI_{\cR(0;0)}$ is defined, for any smooth  $h$ with zero average, by
	\begin{align}
			\cI_{\cR(0;0)}&[h](\vartheta)  = - \alpha^2 \fint_{\T}\sin(\vartheta+\eta)  h(\eta) \wrt \eta  +\alpha^3 \fint_{\T}\Big( \tfrac12\sin(\vartheta+2\eta) + \sin(2\vartheta+\eta) \Big)h(\eta) \wrt \eta    \\
			& -\alpha^4\fint_{\T} \Big( \tfrac{7}{12}\sin(\vartheta+3\eta) + \tfrac{2}{3} \sin(2\vartheta+2\eta) + \tfrac{7}{4} \sin(3\vartheta+\eta)  -\tfrac12\sin(\vartheta-\eta) \Big)h(\eta) \wrt\eta		\,.
	\end{align}

\end{lem}

\begin{proof}
	First, using  \eqref{S.diffeo} and \eqref{B.diffeo} we get the identity
	$$
		\cS^{-1} \pa_{\theta} =\pa_{\vartheta} \cB^{-1}.
	$$
	Consequently,
	\begin{equation}\label{lunch1}
		\cS^{-1} \pa_{\theta} \cR(0;0) \cS = \pa_{\vartheta} \cB^{-1} \cR(0;0) \cS \,.
	\end{equation}
We now expand the action of  $\cB^{-1}\cR(0;0)[ \cS  h  ]$. 
By  \eqref{S.diffeo}, \eqref{B.diffeo} and Lemma \ref{lemma.cR.00}, we compute
	\begin{align}
		\cB^{-1}  \cR(0;0) [\cS h] & = \cB^{-1} \Big[ \fint_{\T}  \log\big| 1-e^{\im(\eta-\theta)}  \big| (1+\beta'(\eta)) h(\eta+\beta(\eta))  \wrt\eta  \Big]   \\
		& = \cB^{-1} \Big[  \fint_{\T}  \log\big| 1-e^{\im(\eta +\breve{\beta}(\eta)-\theta)}  \big| h(\eta)  \wrt\eta \Big]  \\
		& = \fint_{\T}  \log\big| 1-e^{\im(\eta +\breve{\beta}(\eta)-\vartheta-\breve{\beta}(\vartheta))}  \big| h(\eta)  \wrt\eta \,.
	\end{align}
At this stage, we shall  expand the kernel using Taylor expansion, for $|\alpha|\leq \alpha_{0} $ small enough, see also \eqref{brevebeta0.est},
\begin{equation}\label{rain1}
	\begin{aligned}
		\log\big| 1&-e^{\im(\eta +\breve{\beta}(\eta)-\vartheta-\breve{\beta}(\vartheta))}  \big| = \log \big| 1- e^{\im(\eta-\vartheta)} - \tH_{\breve{\beta}}(\eta,\vartheta) \big| \\
		&= \log\big| 1-e^{\im(\eta -\vartheta)}\big|  + {\rm Re} \Big( \tfrac{\tH_{\breve{\beta}}(\eta,\vartheta)}{1-e^{\im(\eta-\vartheta)}}  \Big) - \tfrac12  {\rm Re}  \Big( \tfrac{\tH_{\breve{\beta}}(\eta,\vartheta)}{1-e^{\im(\eta-\vartheta)}} \Big)^2   + \tW_{\breve{\beta}}(\eta,\vartheta) \,,
	\end{aligned}
\end{equation}
where
$$
		\tH_{\breve{\beta}} (\eta,\vartheta)  := -e^{\im(\eta -\vartheta +\breve{\beta}(\eta)-\breve{\beta}(\vartheta))} + e^{\im(\eta-\vartheta)}
$$
and
\begin{equation}
	\tW_{\breve{\beta}}(\eta,\vartheta) := - \tfrac12 \Big( \tfrac{\tH_{\breve{\beta}}(\eta,\vartheta)}{1-e^{\im(\eta-\vartheta)}} \Big)^3\int_{0}^{1} (1-\tau)^2 \big( \pa_{x}^3 \log|1-\tau x| \big)|_{x=\frac{\tH_{\breve{\beta}}(\eta,\vartheta)}{1-e^{\im(\eta-\vartheta)}}} \wrt \tau \,  .
\end{equation}
Applying once again  Taylor expansion, we infer
\begin{equation}\label{tH.beta}
	\begin{aligned}
		\tH_{\breve{\beta}} (\eta,\vartheta)  		=&  - \im\,(\breve{\beta}(\eta)-\breve{\beta}(\vartheta)) e^{\im(\eta-\vartheta)} + \tfrac12 (\breve{\beta}(\eta)-\breve{\beta}(\vartheta))^2 e^{\im(\eta-\vartheta)}  \\
		& + \tfrac{\im}{2} \big( \breve{\beta}(\eta)-\breve{\beta}(\vartheta) \big)^3\, \int_{0}^{1} (1-\tau)^2 e^{\im ( \eta-\vartheta  + \tau ( \breve{\beta}(\eta)-\breve{\beta}(\vartheta) ) )} \wrt \tau  \,.
	\end{aligned}
\end{equation}
On the other hand,  by standard  trigonometric identities, we get
\begin{equation}\label{powerade4}
	\frac{e^{\im(\eta-\vartheta)}}{1-e^{\im(\eta-\vartheta)}} = -\frac12 +\frac{\im}{2} \frac{\sin(\eta-\vartheta)}{1-\cos(\eta-\vartheta)} = -\frac12 +\frac{\im}{2} {\rm cotg} \Big( \tfrac{\eta-\vartheta}{2}\Big) \,.
\end{equation}
Combining   \eqref{tH.beta} and \eqref{powerade4}, we obtain
	\begin{align*}
		{\rm Re} \Big( \frac{\tH_{\breve{\beta}}(\eta,\vartheta)}{1-e^{\im(\eta-\vartheta)}}  \Big)  &= (\breve{\beta}(\eta)-\breve{\beta}(\vartheta)) {\rm Im} \Big( \frac{e^{\im(\eta-\vartheta)}} {1-e^{\im(\eta-\vartheta)}}  \Big)  \\
		& \quad + \frac12 (\breve{\beta}(\eta)-\breve{\beta}(\vartheta))^2 {\rm Re} \Big( \frac{e^{\im(\eta-\vartheta)}}{1-e^{\im(\eta-\vartheta)}}  \Big)  + O(\alpha^5)  \\
		& =\frac12 (\breve{\beta}(\eta)-\breve{\beta}(\vartheta)) {\rm cotg} \big( \tfrac{\eta-\vartheta}{2} \big) - \tfrac14  (\breve{\beta}(\eta)-\breve{\beta}(\vartheta))^2 + O(\alpha^5) \,,	\end{align*}
	and 
	\begin{align}
				{\rm Im} \Big( \frac{\tH_{\breve{\beta}}(\eta,\vartheta)}{1-e^{\im(\eta-\vartheta)}}  \Big)  &= - (\breve{\beta}(\eta)-\breve{\beta}(\vartheta)) {\rm Re} \Big( \frac{e^{\im(\eta-\vartheta)}} {1-e^{\im(\eta-\vartheta)}}  \Big) + O(\alpha^3)  \\ 
		& = \frac12 (\breve{\beta}(\eta)-\breve{\beta}(\vartheta)) + O(\alpha^3) \,, \\
		\frac12 {\rm Re} \Big( \Big( \frac{\tH_{\breve{\beta}}(\eta,\vartheta)}{1-e^{\im(\eta-\vartheta)}} \Big)^2  \Big) & =\frac12 \Big(   {\rm Re} \Big( \frac{\tH_{\breve{\beta}}(\eta,\vartheta)}{1-e^{\im(\eta-\vartheta)}} \Big)  \Big)^2 - \frac12 \Big(  {\rm Im} \Big( \frac{\tH_{\breve{\beta}}(\eta,\vartheta)}{1-e^{\im(\eta-\vartheta)}} \Big)  \Big)^2 \\
		& = \frac18 \Big(  (\breve{\beta}(\eta)-\breve{\beta}(\vartheta)) {\rm cotg} \big( \tfrac{\eta-\vartheta}{2} \big)  \Big)^2 - \frac18 \big( \breve{\beta}(\eta)-\breve{\beta}(\vartheta) \big)^2 + O(\alpha^5) \,.\label{rain2}
	\end{align}	
To get the suitable expansion, we have to study terms of the form
\begin{equation}
	(\sin(k\eta)-\sin(k\vartheta)) {\rm cotg} \big( \tfrac{\eta-\vartheta}{2} \big) \,, \quad k =1,2,3,4 \,.
\end{equation}
Using item $6.$ in Appendix \ref{app.trigo}, we obtain the following standard trigonometric identities
\begin{equation*}
	\begin{aligned}
		(\sin(\eta)-\sin(\vartheta)) {\rm cotg} \big( \tfrac{\eta-\vartheta}{2} \big) & =   \cos(\eta)+ \cos(\vartheta) \,, \\
		(\sin(2\eta)-\sin(2\vartheta)) {\rm cotg} \big( \tfrac{\eta-\vartheta}{2} \big) & = \cos(2\eta) + 2\cos(\eta+\vartheta) + \cos(2\vartheta) \,, \\
		(\sin(3\eta)-\sin(3\vartheta)) {\rm cotg} \big( \tfrac{\eta-\vartheta}{2} \big) & = \cos(3\eta) + 2 \cos(2\eta+\vartheta) + 2 \cos(\eta+2\vartheta)  + \cos(3\vartheta) \,,  \\
		(\sin(4\eta)-\sin(4\vartheta)) {\rm cotg} \big( \tfrac{\eta-\vartheta}{2} \big) & = \cos(4\eta) + 2 \cos(3\eta+\vartheta) + 2 \cos(2\eta+2\vartheta) \\
		& \quad + 2 \cos(\eta+3\vartheta) + \cos(4\vartheta) \,, 
	\end{aligned}
\end{equation*}
and 
\begin{equation*}
	\begin{aligned}
		\Big(	(\sin(2\eta)&-\sin(2\vartheta)) {\rm cotg} \big( \tfrac{\eta-\vartheta}{2} \big)  \Big)^2 =  \tfrac12(\cos(4\eta)+\cos(4\vartheta)) +2  \cos(3\eta+\vartheta)
			+2\cos(\eta+3\vartheta)  \\
		&+ \cos(2\eta-2\vartheta) + 3\cos(2\eta+2\vartheta) + 4 \cos(\eta-\vartheta) + 3 \,, \\
		(\sin(2\eta)&-\sin(2\vartheta))^2 = -\tfrac12 (\cos(4\eta)+\cos(4\vartheta)) + \cos(2\eta+2\vartheta) - \cos(2\eta-2\vartheta) +1 \,.
	\end{aligned}
\end{equation*}
It follows from the above identities together with \eqref{rain1}, \eqref{rain2} and \eqref{brevebeta0.est}, that
	\begin{align}
		\log\big| 1&-e^{\im(\eta +\breve{\beta}(\eta)-\vartheta-\breve{\beta}(\vartheta))} \big|  =  \log \big| 1- e^{\im(\eta-\vartheta)} \big| +\tfrac12 (\breve{\beta}(\eta)-\breve{\beta}(\vartheta)) {\rm cotg} \big( \tfrac{\eta-\vartheta}{2}\big) \\
		&	\quad -\tfrac18 (\breve{\beta}(\eta)-\breve{\beta}(\vartheta))^2 - \tfrac18 \Big( (\breve{\beta}(\eta)-\breve{\beta}(\vartheta)) {\rm cotg} \big( \tfrac{\eta-\vartheta}{2}\big)  \Big)^2 + O(|\breve{\beta}|^3) \,, \end{align}
			and
		\begin{align}
		\log\big| 1&-e^{\im(\eta +\breve{\beta}(\eta)-\vartheta-\breve{\beta}(\vartheta))} \big| 
		= \log \big| 1- e^{\im(\eta-\vartheta)} \big|  + \tU_{\alpha} \big( \cos(\eta)+\cos(\vartheta) \big) \\
		& + \tfrac12\alpha^2 \big(  \cos(2\eta)+2 \cos(\eta+\vartheta)+\cos(2\vartheta) \big) \\
		& -\tfrac14 \alpha^3 \big( \cos(3\eta) + 2 \cos(2\eta+\vartheta) + 2 \cos(\eta+2\vartheta)  + \cos(3\vartheta)  \big)  \\
		&+\alpha^4 \Big( \tfrac{5}{12}\big( \cos(4\eta) + 2 \cos(\vartheta+3\eta) + 2 \cos(2\vartheta+2\eta) + 2 \cos(3\vartheta + \eta) + \cos(4\vartheta) \big) \\
		& \quad -\tfrac18 \big( - \tfrac12(\cos(4\eta)+\cos(4\theta))  + \cos(2\eta+2\vartheta) - \cos(2\eta-2\vartheta) +1 \big) \\
		& \quad -\tfrac18 \big(   \tfrac12(\cos(4\eta)+\cos(4\vartheta)) +2  \cos(3\eta+\vartheta)		+2\cos(\eta+3\vartheta)  \\
		&\quad \quad + \cos(2\eta-2\vartheta) + 3\cos(2\eta+2\vartheta) + 4 \cos(\eta-\vartheta) + 3   \big) \Big) + O(\alpha^5) \,.	\end{align}
	By gathering the terms, we get
	\begin{align}
		\log\big| 1&-e^{\im(\eta +\breve{\beta}(\eta)-\vartheta-\breve{\beta}(\vartheta))} \big|  = \log \big| 1- e^{\im(\eta-\vartheta)} \big|+\tU_{\alpha} \big( \cos(\eta)+\cos(\vartheta) \big)\\
		& + \tfrac12\alpha^2 \big(  \cos(2\eta)+2 \cos(\eta+\vartheta)+\cos(2\vartheta) \big) \\
& -\tfrac14 \alpha^3 \big( \cos(3\eta) + 2 \cos(2\eta+\vartheta) + 2 \cos(\eta+2\vartheta)  + \cos(3\vartheta)  \big)  \\
&+\alpha^4 \Big( \tfrac{5}{12} \big(\cos(4\eta)+\cos(4\vartheta)\big)  + \tfrac{7}{12}\cos(3\eta+\vartheta)+\tfrac{1}{3}\cos(2\eta+2\vartheta)+\tfrac{7}{12}\cos(\eta+3\vartheta) \\
& \quad- \tfrac12 \cos(\eta-\vartheta) -\tfrac12 \Big)   + O(\alpha^5) \,.
	\end{align}
	We point out that $O(\alpha^5)$ is a smooth function of order $\alpha^5.$ 
It follows that  for a zero average function $h$ and recalling the notation in \ref{notation.tR}), we obtain
	\begin{align}
		\cB^{-1} & \cR(0;0) [\cS h]  = \cR(0;0)[h]  + \tU_{\alpha} \fint_{\T} \cos(\eta) h(\eta) \wrt\eta \\
		&+ \alpha^2 \fint_{\T}\big( \tfrac12\cos(2\eta) + \cos(\vartheta+\eta)\big)h(\eta)  \wrt\eta \\
		& -\alpha^3\fint_{\T} \big( \tfrac14\cos(3\eta) + \tfrac12 \cos(\vartheta+2\eta)+ \tfrac12 \cos(2\vartheta + \eta ) \big) h(\eta) \wrt\eta  \\
		&  +\alpha^4 \fint_{\T} \big(  \tfrac{5}{12}\cos(4\eta) + \tfrac{7}{12} \cos(\vartheta+3\eta)+ \tfrac{1}{3} \cos(2\vartheta+2\eta) \\
		& \quad \quad \quad + \tfrac{7}{12}\cos(3\vartheta+\eta) - \tfrac12  \cos(\vartheta-\eta)\big) h(\eta) \wrt\eta  + \alpha^5 \tR[h] \,,\label{powerade5}
	\end{align}
	where $\tR \in \OpM_{s}^{-\infty}$ is a smoothing kernel operator.
Inserting \eqref{powerade5} into \eqref{lunch1}, we find
\begin{equation}
	\cS^{-1} \pa_{\theta} \cR(0;0) \cS = \pa_{\vartheta} \cR(0;0) + \cE_{\cR(0;0)} \,,
\end{equation}
with
	\begin{align}
		&	\cE_{\cR(0;0)}[h](\vartheta)  =- \alpha^2\fint_{\T}\sin(\vartheta+\eta)  h(\eta) \wrt \eta   +\alpha^3 \fint_{\T}\big( \tfrac12\sin(\vartheta+2\eta) + \sin(2\vartheta+\eta) \big)h(\eta) \wrt \eta    \\
		& \quad-\alpha^4\fint_{\T} \Big( \tfrac{7}{12}\sin(\vartheta+3\eta) + \tfrac{2}{3} \sin(2\vartheta+2\eta) + \tfrac{7}{4} \sin(3\vartheta+\eta)  -\tfrac12\sin(\vartheta-\eta) \Big)h(\eta) \wrt\eta  \\
      &\qquad\qquad + \alpha^5 \tR[h]\,,
	\end{align}
	from which we deduce the claims, using also Lemma \ref{lemma.kernelHS}.
\end{proof}
The next step is to study the localized transformation induced by $\cS$ on the phase space, obtained by projecting it onto $H_{\circ}^{s}(\T^2)$. The latter space is defined in \eqref{Hcirc.def}. To this end, we define the operator
$$
\cS_{\rm ph} := \Pi_{\rm ph} \cS \Pi_{\rm ph},
$$
with $\Pi_{\rm ph}$ as in \eqref{proj.ph}. The   following result was partially proved in \cite{HR21}.
	\begin{lem}\label{Lema-decompMon}
			 Let $s\geq 0$. The following hold:
    \\[1mm]
    \noindent $(i)$ The operator
	\begin{equation}\label{S.proj}
		\cS_{\rm ph}  : H_\circ^s(\T^2) \to H_\circ^s(\T^2) 
	\end{equation}
	is well defined and invertible, satisfying the following estimates,
    \begin{align}
		&	\| \cS_{\rm ph}^{\pm 1} h \|_{s}^{5,\alpha_{0}}\lesssim_{s} \| h \|_{s}^{5,\alpha_{0}}\,, \quad  	\| \big(\cS_{\rm ph}^{\pm 1} -\Pi_{\rm ph}\big)h \|_{s}^{5,\alpha_{0}}\lesssim_{s} \alpha_{0}^2 \| h \|_{s+1}^{5,\alpha_{0}} \,.
	\end{align}
     Furthermore, the operators $\cS_{\rm ph}^{\pm 1}$ are real and reversibility preserving;
     \\[1mm]
    \noindent $(ii)$
		We have the representations
        \begin{align}
            \cS_{\rm ph} h = \cS h - \cS_{0}[h] \bs_{1} \,, \quad \cS_{\rm ph}^{-1} h = \cS^{-1} h - \cS_{1}[h] \cS^{-1}\bs_{1}\,, \label{representations}
        \end{align}
        with
        \begin{align}
            \cS_{0}[h]&:= \braket{h,(\cB^{-1}-{\rm Id})\bs_{1}}_{L_\theta^2(\T)} \,, \\
            \cS_{1}[h]&:= \braket{h,(\cB-{\rm Id})g_{1}}_{L_\theta^2(\T)} \,, \quad g_1(\theta):=\frac{\mathbf{s}_1(\theta)}{\big\langle \mathbf{s}_1,\mathcal{B}\mathbf{s}_1\big\rangle_{L^2_\theta(\T)}}\,, \label{g1.S}
        \end{align}
        such that the following estimates hold, for any $m\geq 0$,
        \begin{equation}
            |(\cS_{0}[\,\cdot\,]) \bs_{1}|_{-m,s}^{5,\alpha_{0}} \lesssim_{m,s} \alpha_{0}^2 \,, \quad |(\cS_{1}[\,\cdot\,]) \cS^{-1}\bs_{1}|_{-m,s}^{5,\alpha_{0}} \lesssim_{m,s} \alpha_{0}^2 \,. \label{est.repre}
        \end{equation}
			\end{lem}
We postpone the proof of this lemma to Appendix \ref{app.diff.U}.
Now,  recalling \eqref{linear.rho.full}, \eqref{cT.0.alpha}, \eqref{smoothF0}, \eqref{smoothK0}, we denote by $\cL_{0,\alpha}$ the operator
	\begin{align}\label{cL0.alpha}
		\cL_{0,\alpha} & := \omega\, \pa_{\vf} + \pa_{\theta}\circ \big( \tV_{0,\alpha}(\theta) + \cR_{0,\alpha} \big) \,, \quad \cR_{0,\alpha} := \cR(r_\alpha;\alpha) + \tF_{0,\alpha}(\theta) \cK_{0,\alpha} \,.
	\end{align}
We have the following conjugation result.
\begin{pro}\label{diffeo.conj.no.proj}
	{\bf (Reduction of the transport term - unperturbed).}
	Consider the operator
	$$
	\cL_{0,\alpha}^{(1)}  := 	\cS_{\rm ph}^{-1} \cL_{0,\alpha} \cS_{\rm ph} \,,
	$$
    with $\cL_{0,\alpha}$ as in \eqref{cL0.alpha} and $\cS_{\rm ph}$ as in \eqref{S.proj}, \eqref{S.diffeo}, \eqref{B.diffeo}.
	Then, for any $s\geqslant s_0$, we have that  $\cL_{0,\alpha}^{(1)}:H_{\circ,{\rm even}}^s\to H_{\circ,{\rm odd}}^{s-1}$ is  continuous, real and reversible. In particular, we have that 	
	\begin{align}
			\cL_{0,\alpha}^{(1)}  
			&= \omega\,\pa_{\vf} + \cD_{0,\alpha}+ \cE_{0,\alpha} \,,\label{conj.no.proj}
	\end{align}
	where:
    \\[1mm]
    \noindent $\bullet$ The diagonal operator $\cD_{0,\alpha}$ takes the form
	\begin{equation}\label{diag.0.alpha}
		\begin{aligned}
			\cD_{0,\alpha}& :=  \Pi_{\rm ph}\pa_{\vartheta} \big( \tm_{1}(\alpha) + \cR(0;0) \big)\Pi_{\rm ph} =  {\rm diag}\big\{ \lambda_{j}(\alpha) \,: \,\ j \in \Z_{\rm ph} \big\} \,,  \\
			\lambda_{j}(\alpha) & := \begin{cases}
            0 & j=1\,, \\
			    \im \big(  \tm_{1}(\alpha) j - \tfrac{1}{2} {\rm sgn}(j) \big)\,, &  |j|\geq 2 \,,
			\end{cases} 
		\end{aligned}
	\end{equation}
	with   $\tm_{1}(\alpha)\in\R$ being a constant  satisfying the estimate \eqref{m1alpha0.est};
	\\[1mm]
    \noindent $\bullet$ The  smoothing operator $\cE_{0,\alpha}\in\OpM_{s}^{-\infty}$
	satisfies the estimates, for any $m\in \N_{0}$,
	\begin{equation}\label{cE.0.alpha.est}
		\begin{aligned}
			&	|  \cE_{0,\alpha}  |_{-m,s}^{5,\alpha_{0}} \lesssim_{m} \alpha_{0}^3  \,,  \quad	| \cE_{0,\alpha}-  \cI_{0,\alpha}  |_{-m,s}^{5,\alpha_{0}} \lesssim_{m} \alpha_{0}^5  \,,
		\end{aligned}
	\end{equation}
	where the leading part $\cI_{0,\alpha}$ is given, for any $h\in H_{\circ,{\rm even}}^s $, by
   \begin{align}
    \cI_{0,\alpha}[h]  
     =&    \alpha^3 \fint_{\T} \big(-\tfrac12 \cos(\vartheta) \sin(2\eta)+3\sin(2\vartheta)\cos(\eta)\big) h(\eta) \wrt\eta \\
     &
     + \alpha^4 \fint_{\T}\big( \tfrac{1}{3} \sin(2\vartheta+2\eta)  - 3\sin(2\vartheta-2\eta)\big) h(\eta) \wrt\eta \\
     &  + \alpha^4 \fint_{\T} \big(-4 \sin(3\vartheta)\cos(\eta) + \tfrac{2}{3} \cos(\vartheta) \sin(3\eta)  \big)h(\eta) \wrt\eta  \\
     & - \tfrac{3}{2}\alpha^4\fint_{\T}  \sin(3\vartheta+3\eta) h(\eta)\wrt \eta\,.\label{I.0.alpha}
        \end{align}
\end{pro}
	 \begin{proof}
	 In the  proof, we extensively use the notation introduced in \eqref{notation.tR} without specifying the expressions of the several operators in $\OpM_{s}^{-\infty}$. First, the continuity, the reality and the reversibility properties of the operator
\begin{equation}\label{map.cL1}
	\cL_{0,\alpha}^{(1)} : H_{\circ,{\rm even}}^s \to H_{\circ,{\rm odd}}^{s-1}
\end{equation}
follow directly from the corresponding properties of $\cL_{0,\alpha}$ in Lemma \ref{sa.pp.T.eps} (with $\overline{r}=r_\alpha$), together with the mapping properties of the transformation $\cS_{\rm ph}$ established in Lemma \ref{Lema-decompMon}-$(i)$.
For notational  convenience, we introduce the operator 
     \begin{align}\label{LL11}
     \mathbb{L}[h]:=\pa_{\theta}\big[(\tV_{0,\alpha}(\theta)+\cR_{0,\alpha})h\big]\,.
     \end{align}
  Using Lemma \ref{Lema-decompMon}-$(ii)$ and the fact that the operators $\cS_{\rm ph}^{\pm}$ commute with $\omega\,\pa_{\vf},$ we compute
\begin{align}
    \cL_{0,\alpha}^{(1)} & := \omega\,\pa_{\vf} + \cS_{\rm ph}^{-1} \LL \cS_{\rm ph} \\
    & = \omega \,\pa_{\vf} + \cS_{\rm ph}^{-1}\LL \cS - (\cS_{0}[\,\cdot\,])\cS_{\rm ph}^{-1} [\LL[\bs_{1}]] \\
    & = \omega \,\pa_{\vf} + \cS^{-1}\LL \cS - \big( \cS_{1}[\LL \cS[\,\cdot\,]] \big) \cS^{-1}[\bs_{1}] - (\cS_{0}[\,\cdot\,])\cS_{\rm ph}^{-1} [\LL[\bs_{1}]] \,. \label{conj.cL1}
\end{align}
Applying Lemma \ref{lemma.reparam.0}-$(iii)$, we find
\begin{align}
	\cS^{-1} \LL \cS & = \pa_{\vartheta} \circ \tm_{1}(\alpha) + \cS^{-1} \pa_{\theta} \,\cR_{0,\alpha} \cS \\
	& = \pa_{\vartheta} \big( \tm_{1}(\alpha) + \cR(0;0) \big) +  \cS^{-1} \pa_{\theta} \cR_{0,\alpha} \cS - \pa_{\vartheta} \cR(0;0) \,. \label{conj.LL}
\end{align}
Inserting \eqref{conj.LL} into \eqref{conj.cL1}, we deduce that $\cL_{0,\alpha}^{(1)}$ takes the form in \eqref{conj.no.proj} and \eqref{diag.0.alpha},
where the new remainder $\cE_{0,\alpha}$, recalling \eqref{map.cL1}, is given by
\begin{equation}
	\cE_{0,\alpha}:= \cE_{0,\alpha}^{\rm I} + \cE_{0,\alpha}^{\rm II} + \cE_{0,\alpha}^{\rm III} \,,\label{cE0.split}
\end{equation}
with
\begin{align}
	\cE_{0,\alpha}^{\rm I} & := \Pi_{\rm ph} \big( \cS^{-1} \pa_{\theta} \cR_{0,\alpha} \cS - \pa_{\vartheta} \cR(0;0) \big) \Pi_{\rm ph}\,, \label{cE0.I} \\
	\cE_{0,\alpha}^{\rm II} & := - \big( \cS_{1}\big[\LL \cS[\Pi_{\rm ph}[\,\cdot\,]]\big] \big)\Pi_{\rm ph} \cS^{-1}[\bs_{1}] \,, \label{cE0.II} \\
	\cE_{0,\alpha}^{\rm III} & := - (\cS_{0}[\Pi_{\rm ph}[\,\cdot\,]])\Pi_{\rm ph}\cS_{\rm ph}^{-1} [\LL[\bs_{1}]]  \label{cE0.III} \,.
\end{align}
In the following, we shall expand each of these contributions to the remainder $\cE_{0,\alpha}$.
\\[1mm]
\noindent  $\blacktriangleright$ {\bf Expansion of $\cE_{0,\alpha}^{\rm I}$.} 
First, by \eqref{cE0.I}, \eqref{S.diffeo} and  \eqref{B.diffeo}, we get
	\begin{align}
		\cE_{0,\alpha}^{\rm I}  
		& = \Pi_{\rm ph} \big( \cS^{-1} \pa_{\theta}(\cR_{0,\alpha}-\cR(0;0))\cS + \cS^{-1}\pa_{\theta}\,\cR(0;0) \cS -  \pa_{\vartheta} \,\cR(0;0)  \big) \Pi_{\rm ph} \\
		& = \Pi_{\rm ph}  \pa_{\vartheta} \, \cB^{-1}\big(\cR_{0,\alpha}-\cR(0;0)\big)\cS  \Pi_{\rm ph} + \Pi_{\rm ph} \cE_{\cR(0;0)} \Pi_{\rm ph} \\
		&  := \Pi_{\rm ph} \pa_{\vartheta}\, \big( \cR_{0,\alpha}^{(1)} + \cR_{0,\alpha}^{(2)}  \big)\Pi_{\rm ph} + \Pi_{\rm ph} \cE_{\cR(0;0)} \Pi_{\rm ph}  \,, \label{leone1}  
	\end{align}
	where $\cE_{\cR(0;0)}$ has already been introduced and estimated  in Lemma \ref{lemma.conjR00}. Moreover the operators $\cR_{0,\alpha}^{(1)}$ and $ \cR_{0,\alpha}^{(2)} $ are defined through   \eqref{cL0.alpha}, \eqref{cT.0.alpha}, \eqref{smoothF0}, \eqref{smoothK0} as follows
	\begin{equation}\label{leone5}
		\begin{aligned}
			\cR_{0,\alpha}^{(1)} & :=  \cB^{-1} \big(\cR(r_\alpha;\alpha) -\cR(0;0) \big)\cS \,, \quad 	\cR_{0,\alpha}^{(2)}  := (\cB^{-1} \tF_{0,\alpha})  \cK_{0,\alpha} \cS \,.
		\end{aligned}
	\end{equation}	
	We start with the study of $\cR_{0,\alpha}^{(1)}.$
	By Lemma \ref{lemma.R.expand} and \eqref{B.diffeo}, we have that
	\begin{align}
		&	\cB^{-1}	\big(\cR(r_\alpha;\alpha)- \cR(0;0)\big [\cS h]  = 	\tA[\cS h]+ \alpha^2\fint_{\T}\cos(\vartheta+\breve{\beta}(\vartheta)-\eta-\breve{\beta}(\eta)) h(\eta) \wrt \eta  \\
		&\quad  - \alpha^3 \fint_{\T}\big( \cos(2\vartheta+2\breve{\beta}(\vartheta)-\eta-\breve{\beta}(\eta)) + \cos(\vartheta+\breve{\beta}(\vartheta)-2\eta-2\breve{\beta}(\eta)) \big)h(\eta) \wrt \eta    \\
		& \quad+\alpha^4\fint_{\T} \mathtt{K}_4\big(\vartheta+\breve{\beta}(\vartheta), \eta+\breve{\beta}(\eta) \big)h(\eta) \wrt\eta+\alpha^5 \tR_{1}[h]\,,
	\end{align}  
	for some $\tR_{1}\in\OpM_{s}^{-\infty}$, with the kernel $\tK_{4}$ defined as in \eqref{K4.kernel}.
	We note that, for $m,k\in\Z$, by \eqref{brevebeta0.est} and Taylor expansion,
	\begin{align}
		\cos &\big(m(\vartheta+\breve{\beta}(\vartheta)) +k(\eta+\breve{\beta}(\eta) \big)   = \cos(m\vartheta+k\eta) -\big(m \breve{\beta}(\vartheta) + k\breve{\beta}(\eta)\big) \sin(m\vartheta+k\eta) + O(|\breve{\beta}|^2) \\
		\label{powerade6} & = \cos(m\vartheta+k\eta) -  \alpha^2 \big(m\sin(2\vartheta)+k\sin(2\eta)\big)\sin(m\vartheta+k\eta) + O(\alpha^3) \\
		& = \cos(m\vartheta+k\eta) + \tfrac{m}{2}\alpha^2\big( \cos((m+2)\vartheta+k\eta) - \cos((m-2)\vartheta+k\eta) \big) \\
		 &\quad + \tfrac{k}{2}\alpha^2\big( \cos(m\vartheta+(k+2)\eta) - \cos(m\vartheta+(k-2)\eta) \big) + O(\alpha^3) = \cos(m\vartheta+k\eta) + O(\alpha^2)\,.
	\end{align}
	In particular, when $m=1$ and $k=-1$, we have
	\begin{equation}
		\begin{aligned}
			\cos\big(\vartheta-\eta  + \breve{\beta}(\vartheta)-\breve{\beta}(\eta)\big)   = \cos(\vartheta-\eta) &+ \alpha^2 \big( \tfrac12 \cos(3\vartheta-\eta)  + \tfrac12 \cos(\vartheta-3\eta) - \cos(\vartheta+\eta) \big)+ O(\alpha^3 ) \,.
		\end{aligned}
	\end{equation}
	This implies that
	\begin{align*}
		&\alpha^2 \fint_{\T}   \cos(\vartheta+\breve{\beta}(\vartheta)-\eta - \breve{\beta}(\eta))   h(\eta) \wrt \eta    = \alpha^2 \fint_{\T}   \cos(\vartheta-\eta)   h(\eta) \wrt \eta \\
		& + \alpha^4 \fint_{\T} \big( \tfrac12 \cos(3\vartheta-\eta)  + \tfrac12 \cos(\vartheta-3\eta) - \cos(\vartheta+\eta) \big)  h(\eta) \wrt\eta  + O(\alpha^5)[h]\,.
	\end{align*}
	Moreover, we deduce from \eqref{powerade6} that
	$$
	\mathtt{K}_4\big(\vartheta+\breve{\beta}(\vartheta), \eta+\breve{\beta}(\eta) \big)=\mathtt{K}_4\big(\vartheta, \eta \big)+O(\alpha^2)\,.
	$$
	Summing up, recalling \eqref{K4.kernel}, we obtain that 
	\begin{align}
		\cR_{0,\alpha}^{(1)} [h]&  = \tA[\cS h] +	 \alpha^2\fint_{\T}\cos(\vartheta-\eta) h(\eta) \wrt \eta    - \alpha^3 \fint_{\T}\big( \cos(2\vartheta-\eta) + \cos(\vartheta-2\eta) \big)h(\eta) \wrt \eta    \\
		&+\alpha^4\fint_{\T} \big( \tfrac12 \cos(\vartheta-\eta) - 2\cos(\vartheta+\eta)\big)h(\eta) \wrt\eta   +\alpha^4\fint_{\T} \big(\tfrac32\cos(2\vartheta-2\eta)  -\tfrac12\cos(2\vartheta+2\eta)\big) h(\eta) \wrt\eta \\
		&   +\alpha^4\fint_{\T} \big(\cos(3\vartheta-\eta) -\tfrac14\cos(3\vartheta+\eta) \big)h(\eta) \wrt\eta +\alpha^4\fint_{\T} \big( \cos(\vartheta-3\eta)-\tfrac14\cos(\vartheta+3\eta) \big)h(\eta) \wrt\eta \\
		& + \alpha^5 \tR_{1}[h] \,.
	\end{align} 
	It follows that
	\begin{align}
		 \pa_{\vartheta} \cR_{0,\alpha}^{(1)} [h]  &= 	- \alpha^2\fint_{\T}\sin(\vartheta-\eta) h(\eta) \wrt \eta    + \alpha^3 \fint_{\T}\big( 2\sin(2\vartheta-\eta) + \sin(\vartheta-2\eta) \big)h(\eta) \wrt \eta    \\
			 &+\alpha^4\fint_{\T} \big( 2 \sin(\vartheta+\eta) -\tfrac12 \sin(\vartheta-\eta)\big)h(\eta) \wrt\eta   +\alpha^4\fint_{\T} \big(\sin(2\vartheta+2\eta)  -3\sin(2\vartheta-2\eta)\big) h(\eta) \wrt\eta \\
			  &+\alpha^4\fint_{\T} \big(\tfrac34 \sin(3\vartheta+\eta) -3\sin(3\vartheta-\eta) \big)h(\eta) \wrt\eta +\alpha^4\fint_{\T} \big( \tfrac14 \sin(\vartheta+3\eta)-\sin(\vartheta-3\eta) \big)h(\eta) \wrt\eta\\
	 &+\alpha^5 \tR_{1}[h] \,. \label{leone02}
	\end{align} 
	Using the projection $\Pi_{\rm ph}$, defined in \eqref{proj.ph} it is easy to check that
	\begin{align} 	
		&	\Pi_{\rm ph}\pa_{\vartheta} \cR_{0,\alpha}^{(1)} [\Pi_{\rm ph} h]  = \alpha^3 \fint_{\T}\big( 2\sin(2\vartheta)\cos(\eta) -\cos(\vartheta) \sin(2\eta) \big) (\Pi_{\rm ph} h) (\eta) \wrt \eta    \\
		&\quad +\alpha^4 \fint_{\T}\big(\sin(2\vartheta+2\eta)  -3\sin(2\vartheta-2\eta)\big) (\Pi_{\rm ph} h)(\eta) \wrt\eta \\
		&  \quad +\alpha^4 \fint_{\T} \big(-\tfrac94 \sin(3\vartheta)\cos(\eta) + \tfrac54 \cos(\vartheta) \sin(3\eta)  \big)(\Pi_{\rm ph} h)(\eta) \wrt\eta +\alpha^5\Pi_{\rm ph}\tR_{1}\Pi_{\rm ph}[h]  \,.\label{leone14} 
	\end{align} 
	We now estimate  $\cR_{0,\alpha}^{(2)}$ in \eqref{leone5}.
	First, we compute, by \eqref{smoothF0}, \eqref{B.diffeo}, \eqref{powerade6} (with $m=1,2$ and $k=0$) and Lemma \eqref{sqrt.1+ralpha},
		\begin{align}
				\{\cB^{-1}\tF_{0,\alpha} \}(\vartheta)  &= \tF_{0,\alpha}(\vartheta+\breve{\beta}(\vartheta))  = - \frac{\sqrt{1+2r_\alpha(\vartheta+\breve{\beta}(\vartheta))} \cos(\vartheta+\breve{\beta}(\vartheta))}{\int_{\T}\sqrt{1+2r_\alpha(\eta)}\cos^2(\eta)\wrt\eta}  \\
			& = - \frac{\big( 1- \alpha^2 \cos(2\vartheta) + O(\alpha^3) \big) \big( \cos(\vartheta) +\tfrac{\alpha^2}{2} (\cos(3\vartheta)-\cos(\vartheta)) + O(\alpha^3) \big)}{\frac12\int_{\T}\big(1-\alpha^2\cos(2\eta)+ O(\alpha^3)\big)(1+\cos(2\eta))\wrt\eta } \\
			& = - \frac{1}{\pi\big( 1 - \frac{\alpha^2}{2} + O(\alpha^3) \big)} \big( \cos(\vartheta) - \alpha^2 \cos(\vartheta) + O(\alpha^3) \big)\\
			& =- \tfrac{1}{2\pi}(2  - \alpha^2 )\cos(\vartheta)  + O(\alpha^3)\,.\label{cicap1}
		\end{align}
		To compute $\cK_{0,\alpha}[ \cS h]$,  we combine  \eqref{smoothK0}, \eqref{S.diffeo}, \eqref{B.diffeo}, \eqref{powerade6} (with $m=0$ and $k=1,3$) and Corollary \ref{cor.T.at.cos}, in order to get
			\begin{align}
				\cK_{0,\alpha}[ \cS h] & = \int_{\T} (\cS h)(\eta) \cT_{0,\alpha} [\cos(\eta)] \wrt \eta = \int_{\T} h (\eta) \cB^{-1} \{\cT_{0,\alpha}[\cos(\eta)] \}\wrt\eta \\
				& = \int_{\T} h (\eta) \Big( \alpha^2 \big(\cos(\eta+\breve{\beta}(\eta)) + \tfrac12 \cos(3(\eta+\breve{\beta}(\eta)))   \big) + O(\alpha^3)  \Big) \wrt\eta\\
				& =  \int_{\T} h(\eta) \big( \alpha^2 \big( \cos(\eta) +  \tfrac{1}{2}\cos(3\eta) ) + O(\alpha^3) \big) \wrt \eta \,.
			\end{align}
			Therefore,
			\begin{align}
				\pa_{\vartheta}\cR_{0,\alpha}^{(2)}[h] & = \pa_{\vartheta}\big( - \tfrac{1}{2\pi}(2-\alpha^2) \cos(\vartheta) + O(\alpha^3) \big)     \int_{\T} \big( \alpha^2 \big( \cos(\eta) +  \tfrac{1}{2}\cos(3\eta) ) + O(\alpha^3) \big) h(\eta) \wrt \eta \\
				& =  (2-\alpha^2)\alpha^2 \sin(\vartheta)  \,  \fint_{\T} \big( \cos(\eta)+ \tfrac12\cos(3\eta)  \big) h(\eta) \wrt\eta  +\alpha^5 \tR_{2}[h] \,,\label{leone3}
			\end{align}
			for some $\tR_{2}\in\OpM_{s}^{-\infty}$. It yields that
			\begin{align}\label{proj-hh}
				\Pi_{\rm ph}\pa_{\vartheta} \cR_{0,\alpha}^{2}[\Pi_{\rm ph}h] = \alpha^5 \Pi_{\rm ph}\tR_{2}\Pi_{\rm ph}[h]:=\alpha^5\tR_{3}[h] \,,
			\end{align}
			for some $\tR_{3}\in\OpM_{s}^{-\infty}$. 
		By Lemma \ref{lemma.conjR00}, we have that
		\begin{align}
			\Pi_{\rm ph} & \cE_{\cR(0;0)} [\Pi_{\rm ph} h ] = \alpha^3\fint_{\T} \big( \tfrac12\cos(\vartheta) \sin(2\eta) +  \sin(2\vartheta) \cos(\eta)  \big)(\Pi_{\rm ph} h)(\eta) \wrt\eta \\
			& - \alpha^4 \fint_{\T} \big( \tfrac{7}{12} \cos(\vartheta) \sin(3\eta) + \tfrac{2}{3} \sin(2\vartheta+2\eta) + \tfrac{7}{4} \sin(3\vartheta) \cos(\eta) \big) (\Pi_{\rm ph} h)(\eta) \wrt\eta + \alpha^5 \tR_{4}[h]\,, \label{proj.ER00}
		\end{align}
		for some $\tR_{4}\in\OpM_{s}^{-\infty}$. Finally, collecting \eqref{leone14}, \eqref{proj-hh} and \eqref{proj.ER00} into \eqref{leone1}, we find
		\begin{align}
			\cE_{0,\alpha}^{\rm I}[h]& = \Pi_{\rm ph} \big( \pa_{\vartheta}\, \big( \cR_{0,\alpha}^{1} + \cR_{0,\alpha}^{2} \big)+\cE_{\cR(0;0)}   \big)[\Pi_{\rm ph} [h] ] \\
			& = \alpha^3 \fint_{\T} \big(-\tfrac12 \cos(\vartheta) \sin(2\eta)+3\sin(2\vartheta)\cos(\eta)\big) (\Pi_{\rm ph} h)(\eta) \wrt\eta \\
			& + \alpha^4 \fint_{\T}\big( \tfrac{1}{3} \sin(2\vartheta+2\eta)  -3\sin(2\vartheta-2\eta)\big) (\Pi_{\rm ph} h)(\eta) \wrt\eta \\
			\label{stepI}	&  + \alpha^4 \fint_{\T} \big(-4 \sin(3\vartheta)\cos(\eta) + \tfrac{2}{3} \cos(\vartheta) \sin(3\eta)  \big)(\Pi_{\rm ph} h)(\eta) \wrt\eta +\alpha^5 \tR_5[h]  \,,
		\end{align}
		for some $\tR_5\in\OpM_{s}^{-\infty}$. 
\\[1mm]
\noindent  $\blacktriangleright$ {\bf Expansion of $\cE_{0,\alpha}^{\rm II}$.} We now expand \eqref{cE0.II}. First, by Lemma \ref{lemma.reparam.0}-(ii), Lemma \ref{Lema-decompMon}-$(ii)$ and $\Pi_{\rm ph}\bs_{1}=0$, we note that
\begin{align}
   \cE_{0,\alpha}^{\rm II} & =- (\mathcal{S}_1[\mathbb{L}\, \cS \Pi_{\rm ph}h])\,\Pi_{\rm ph}\mathcal{S}^{-1}\mathbf{s}_1\\
   & = -(\mathcal{S}_1[\mathbb{L}\,  \Pi_{\rm ph}h])\,\Pi_{\rm ph}(\mathcal{S}^{-1}-\hbox{Id})\mathbf{s}_1+\alpha^6 \tR[h] \,, \label{II.cE}
\end{align}
for some $\tR \in\OpM_{s}^{-\infty}$. By Lemma \ref{lemma.reparam.0}-$(ii)$, we have that
    \begin{align}
      \Pi_{\rm ph}(\cS^{-1}-{\rm Id})[\mathbf{s}_1] & = \sqrt 2\Pi_{\rm ph} \big( \pa_{\vartheta} (\sin(2\vartheta)\sin(\vartheta) \big)\\
      &=\tfrac{3\sqrt{2}}{2}\alpha^2\sin(3\vartheta) +O(\alpha^3)\,. \label{SSLL1}
    \end{align}
Recalling \eqref{g1.S}, Lemma \ref{Lema-decompMon}-$(ii)$ and \ref{lemma.reparam.0}-$(ii)$, we now look at
\begin{align}
    \cS_{1}[\Pi_{\rm ph}[\LL[h]]] & =  \frac{\braket{\Pi_{\rm ph}[\LL[h]], (\cB-{\rm Id})\bs_{1}}_{L_\theta^2(\T)}}{\braket{\bs_{1},\cB \bs_{1}}_{L_\theta^2(\T)}}  \\
    & =  \braket{\LL[h], \Pi_{\rm ph}(\cB-{\rm Id})\bs_{1}}_{L_\theta^2(\T)} + \alpha^4 \tR[h] \,, \label{cS1.expand}
\end{align}
for some $\tR \in\OpM_{s}^{-\infty}$. By Lemma \ref{lemma.reparam.0}-$(ii)$, we have that
\begin{align}
    \Pi_{\rm ph}(\cB-{\rm Id})\bs_{1} & = -\sqrt 2 \alpha^{2} \Pi_{\rm ph}\big(\sin(2\theta) \cos(\theta) \big)+ O(\alpha^3) = - \tfrac{\sqrt 2}{2}\alpha^{2} \sin(3\theta)  + O(\alpha^3) \,. \label{B-Id.s1}
\end{align}
We now compute, recalling \eqref{LL11}, \eqref{linear.rho.full}, \eqref{cT.0.alpha}, \eqref{smoothF0}, \eqref{smoothK0},
\begin{equation}
    \LL[h] = \pa_{\theta}\big(\tV_{0,\alpha}h + \cR_{0,\alpha}[h]\big) = \pa_{\theta}\big( \cT_{0,\alpha}[h] + \tF_{0,\alpha}(\theta) \cK_{0,\alpha}[h] \big) \,.
\end{equation}
By Lemma \ref{lemma.R.expand} and Lemma \ref{speed.V.expan}, we have
\begin{equation}
    \pa_{\theta}\cT_{0,\alpha} = \pa_{\theta}\cT_{0,0} + \pa_{\theta}\big(\cT_{0,\alpha} -\cT_{0,0}\big)\,, \quad \cT_{0,0}:=\tfrac12 + \cR(0;0)\,,
\end{equation}
with
\begin{equation}
    \big\| \pa_{\theta}\big(\cT_{0,\alpha} - \cT_{0,0}\big)[h]\big\|_{s}^{5,\alpha_{0}} \lesssim_{s} \alpha_{0}^2 \| h \|_{s+1}^{5,\alpha_{0}} \,.
\end{equation}
    From Lemma \ref{sa.pp.T.eps}, we deduce that $\cT_{0,0}=\tfrac12 +\cR(0;0)$ is self-adjoint. Therefore, by \eqref{smoothK0}, it follows that 
	\begin{align}
		\cK_{0,\alpha}[h] & = \int_{\T} \cos(\eta) \cT_{0,\alpha}[h](\eta) \wrt\theta \\
        &= \int_{\T} \cos(\eta) \big(\tfrac12 + \cR(0;0)\big)[h](\eta) \wrt\eta + \alpha^2 \tR \\ 
		& = \int_{\T}  \big(\tfrac12 + \cR(0;0)\big)[ \cos(\eta)] h(\eta) \wrt\eta + \alpha^2 \tR = \alpha^2 \tR \,,
	\end{align}
    for some $\tR \in\OpM_{s}^{-\infty}$. We obtain that
    \begin{align}
		\tF_{0,\alpha}(\theta)\cK_{0,\alpha}[h]
        &= \alpha^2 \tR \,,\label{smoothKM0L}
	\end{align}
    and, consequently,
    \begin{align}
    \mathbb{L} [h]&=\partial_\theta\Big( \cT_{0,\alpha}[h]+\tF_{0,\alpha}(\theta)\cK_{0,\alpha}[h]\Big)  = \pa_{\theta}\cT_{0,0}[h]+ \pa_{\theta}\big( \cT_{0,\alpha}-\cT_{0,0} \big) [h]+ \alpha^2 \tR[h] \,.\label{LL.h.expand}
    \end{align}
    Inserting \eqref{LL.h.expand} and \eqref{B-Id.s1} into \eqref{cS1.expand}, we compute, together with Proposition \ref{lemma.constrain1}-$(iv)$,
    \begin{align}
        \cS_{1}[\LL[h]] & = \braket{\LL[h],\Pi_{\rm ph}(\cB-{\rm Id})\bs_{1}}_{L_\theta^2(\T)} + \alpha^3 \tR \\
        & = - \tfrac{\sqrt 2}{2} \alpha^2 \fint_{\T} \big(\pa_{\theta}\cT_{0,0} \big)[h](\eta) \sin(3\eta) \wrt \eta + \alpha^3 \tR \\
        & =  \tfrac{\sqrt 2}{2} \alpha^2 \fint_{\T} h(\eta) \big(\pa_{\theta}\cT_{0,0} \big)[\sin(3\,\cdot\,)](\eta) \wrt \eta + \alpha^3 \tR \\
        & =  \tfrac{\sqrt 2}{2} \alpha^2 \fint_{\T} h(\eta) \cos(3\eta) \wrt\eta + \alpha^3 \tR \,, \label{cS1.expand2}
     \end{align}
     for some $\tR \in\OpM_{s}^{-\infty}$. Finally, collecting \eqref{cS1.expand2} and \eqref{SSLL1} back into \eqref{II.cE}, we obtain that
     \begin{align}
         \cE_{0,\alpha}^{\rm II} & = -\tfrac32 \alpha^4 \fint_{\T} \sin(3\vartheta) \cos(3\theta) h(\eta) \wrt\eta+ \alpha^5 \tR \,, \label{stepII}
     \end{align}
     for some $\tR \in\OpM_{s}^{-\infty}$.
\\[1mm]
\noindent $\blacktriangleright$ {\bf Expansion of $\cE_{0,\alpha}^{\rm III}$.} We finally expand \eqref{cE0.III}. First, we claim that
	\begin{align}\label{LL[s1]}
	\mathbb{L}[\mathbf{s}_1]= \tfrac{3\sqrt{2}}{2}\alpha^2\cos(3\theta)+O(\alpha^3)  \,.
\end{align}
In view of \eqref{T.alpha.operator} and \eqref{cL0.alpha}, we write
     \begin{align}
	\mathbb{L}[\mathbf{s}_1]= \pa_{\theta} \big(\mathcal{T}(r_\alpha;\alpha)[\mathbf{s}_1] + \tF_{0,\alpha}(\theta) \cK_{0,\alpha}[\mathbf{s}_1]\big)  \,. \label{LL.at.s1}
\end{align}
By Corollary \ref{cor.T.at.cos}, recalling the notation \eqref{cT.0.alpha}, we infer that	
\begin{align}\label{form-tham}
	\cT_{0,\alpha}[\mathbf{s}_1]=\mathcal{T}(r_\alpha;\alpha)[\mathbf{s}_1] &=   \tfrac{\sqrt{2}}{2} \alpha^2\sin(3\theta)  +O( \alpha^3)\,.
\end{align}	
By \eqref{smoothF0} and \eqref{ralpha.asymp}, with similar computation as in \eqref{cicap1}, we find
	\begin{align}\label{f0-forma}
		\tF_{0,\alpha}(\theta) & := - \frac{\sqrt{1+2r_\alpha(\theta)}\cos(\theta)}{\int_{\T}\sqrt{1+2r_\alpha(\eta)}\cos^2(\eta)\wrt\eta} =-\tfrac{1}{\pi}\cos(\theta)+O(\alpha^2)  \,.
	\end{align}	
From  \eqref{form-tham} and \eqref{smoothK0}, we deduce that
	\begin{align}
		\cK_{0,\alpha}[\bs_{1}](\vf) & = \int_{\T} \cos(\eta) \cT_{0,\alpha}[\mathbf{s}_1](\eta) \wrt\eta \\
		& =\alpha^2\tfrac{\sqrt{2}}{2} \int_{\T} \cos(\eta)\sin(3\eta)  \wrt\eta+O(\alpha^3) =O(\alpha^3)\,.\label{smoothKM0}
	\end{align}
Therefore, the claimed expansion in \eqref{LL[s1]} follows by \eqref{LL.at.s1}, \eqref{form-tham} and \eqref{smoothKM0}.  By Lemma \ref{Lema-decompMon}-$(ii)$, Lemma \ref{lemma.reparam.0}-$(ii)$ and \eqref{scalar.prod.L2}, we find that
        \begin{align}
        \nonumber\mathcal{S}_0[h]&=\big\langle h,  \big(\mathcal{B}^{-1}-\textnormal{Id}\big)\mathbf{s}_1  \big\rangle_{L^2_\theta(\T)}\,\\
        &=\alpha^2\sqrt{2}\fint_{\T} h(\eta)  \sin(2\eta)\cos(\eta)\wrt \eta+ \alpha^3 \tR_{6} \,,
        \end{align}
        for some $\tR_{6} \in \OpM_{s}^{-\infty}$. 
        Consequently, we get
        \begin{align}\label{form-So}
            \mathcal{S}_0[\Pi_{\rm ph}h]=\alpha^2\tfrac{\sqrt{2}}{2}\fint_{\T} (\Pi_{\rm ph}h)(\eta)  \sin(3\eta)\wrt\eta+ \alpha^3 \tR_{6}\,.
        \end{align}
Collecting \eqref{LL[s1]} and \eqref{form-So} into \eqref{cE0.III}, using also Lemma \ref{Lema-decompMon}-$(i)$, we obtain
\begin{align}
    \cE_{0,\alpha}^{\rm III} & = - (\cS_{0}[\Pi_{\rm ph}[\,\cdot\,]])\Pi_{\rm ph}[\LL[\bs_{1}]]  + \alpha^5 \tR_{7} \\
    & = - \tfrac32  \alpha^4\fint_{\T} \cos(3\vartheta) \sin(3\eta) (\Pi_{\rm ph}h)(\eta) \wrt\eta + \alpha^5 \tR_{7}\,, \label{stepIII}
\end{align}
for some $\tR_{7} \in \OpM_{s}^{-\infty}$. 
\\[1mm]
\noindent  $\blacktriangleright$ {\bf Conclusion.} We finally collect the contributions from  \eqref{stepI}, \eqref{stepII} and \eqref{stepIII} into \eqref{cE0.split} to obtain that 
\begin{equation}\label{final.cE.exp}
    \cE_{0,\alpha} = \cI_{0,\alpha} + \alpha^5 \,\tR
\end{equation}
for some $\tR_{7} \in \OpM_{s}^{-\infty}$, where $\cI_{0,\alpha}$ is the integral operator as in \eqref{I.0.alpha}. The estimates in \eqref{cE.0.alpha.est} then follow by the explicit expansion in \eqref{final.cE.exp}, \eqref{I.0.alpha} and by Lemma \ref{lemma.kernelHS}.
This concludes the proof.             
\end{proof}

In the following lemma, we compute the action of the operator $\cL_{0,\alpha}^{(1)}$ in \eqref{conj.no.proj} on the first three modes, which are the ones that give nontrivial contributions in the $\alpha^4$-expansion  $\cI_{0,\alpha}$, defined in \eqref{I.0.alpha}, of the smoothing remainder operator $\cE_{0,\alpha}$. 

\begin{lem} \label{lemma.modes123.E}
	$(i)$ We have that
	\begin{align}
	\cD_{0,\alpha}[\cos(\vartheta)] &= 0\,,\\
		 \cD_{0,\alpha}[\sin(j\vartheta)] & = \im \lambda_{j}^{(0)}(\alpha) \cos(j\vartheta) \,, \quad j\geq 2 \,,  \\
	 \cD_{0,\alpha}[\cos(j\vartheta)] & = -\im \lambda_{j}^{(0)}(\alpha) \sin(j\vartheta) \,, \quad j\geq 2 \,;
	\end{align}
\noindent $(ii)$	We have that 
\begin{align}
			\cE_{0,\alpha}[\cos(\vartheta)] & = \tfrac32\alpha^3\sin(2\vartheta)- 2 \alpha^4\sin(3\vartheta)+ O(\alpha^5)\,,\quad \Pi_{1,\bc}\cE_{0,\alpha}[\cos(\,\cdot\,)]=0 \,, \\
		\cE_{0,\alpha}[e^{\pm \im 2\vartheta}] & = \mp \tfrac14 \im \alpha^3 \cos(\vartheta) \pm \im \alpha^4 \big( \tfrac{3}{2}e^{\pm \im 2 \vartheta} + \tfrac{1}{6}e^{\mp\im 2 \vartheta} \big) + O(\alpha^5) \,, \\
	\cE_{0,\alpha}[e^{\pm \im 3\vartheta}] & =  \pm \im \tfrac{1}{3} \alpha^4 \cos(\vartheta) \mp \im \tfrac34 \alpha^4 e^{\mp\im 3 \vartheta}  + O(\alpha^5) \,.
	\end{align}
\end{lem}
\begin{proof}
	The claims follow by  straightforward computations using  \eqref{conj.no.proj}, \eqref{diag.0.alpha}, \eqref{I.0.alpha}
    and \eqref{dir.wtcG0}.
\end{proof}

\subsection{Full reduction}\label{section-Full-redu}
We are now in a position to diagonalize the operator $\cL_{0,\alpha}^{(1)}$ in \eqref{conj.no.proj}
through a time–independent transformation.  
It is worth emphasizing that no small–divisor problems arise at this stage,
and therefore a linear Nash--Moser scheme is unnecessary.  
In fact, the reduction is carried out in a single step by solving a
nonlinear homological equation.  
Its linear part is explicitly solved at the level of the matrix elements of the
reversible operator $\cE_{0,\alpha}$ with respect to the Fourier basis of
$L_{\rm ph}^2(\T)$, see \eqref{matrix.rep.compact}–\eqref{matrix.elem.compact}.  
The main result of this section is stated in the following proposition.
\begin{pro}{\bf (Diagonal reducibility - unperturbed).}\label{prop.diag.red.equi}
	Let $   S\geq s_0$ and $s\in [s_0,S]$. There exists $\alpha_{0}=\alpha_{0}(S)\in (0,1)$ small enough such that, for any $m\geq 0$  and for any $\alpha\in [\alpha_{1},\alpha_{2}]$ as in \eqref{alpha1alpha2}, the following holds. There exists an invertible, reversibility preserving map $\cU={\rm Id} + \cY:H_{\circ}^s(\T^2)\to H_{\circ}^s(\T^2)$ with $\cY \in \OpM_{s}^{-m}$ for any $m\in\mathbb{N}$, and such that the operator $\cL_{0,\alpha}^{(1)}$ in \eqref{conj.no.proj} is conjugated to the real and reversible operator
	\begin{align}
		\cL_{0,\alpha}^{(2)} & := \cU^{-1} \cL_{0,\alpha}^{(1)} \cU = \omega\,\pa_{\vf} + \cD_{0,\alpha}^{(\infty)} :H_{\circ,{\rm even}}^s(\T^2) \to H_{\circ,{\rm odd}}^{s-1}(\T^2) \,,
	\end{align}
	where $\cD_{0,\alpha}^{(\infty)} $ is the diagonal operator
	\begin{align}
		\cD_{0,\alpha}^{(\infty)} & := \cD_{0,\alpha} + \cZ_{0,\alpha} =  {\rm diag}\big\{ \lambda_{j}^{(\infty)}(\alpha) \, : \, j\in\Z_{\rm ph} \big\} \,, \\
		\lambda_{j}^{(\infty)}(\alpha) &:= \begin{cases}
			0\,, & j= 1\,, \\
			\lambda_{j}(\alpha) + \tz_{j}(\alpha) \,, & |j|\geq 2\,,
		\end{cases} \label{lambda.linear}
	\end{align}
	with $\cZ_{0,\alpha}\in\OpM_{s}^{-m}$ and $\lambda_{j}(\alpha)$  as in \eqref{diag.0.alpha}, 
	satisfying the following property, \begin{equation}\label{pure.im.lambda}
		\begin{aligned}
			&	\lambda_{-j}^{(\infty)}(\alpha)=	- \lambda_{j}^{(\infty)}(\alpha) = \overline{\lambda_{j}^{(\infty)}}(\alpha) \,, 
			\quad |j|\geq 2\,.
		\end{aligned}
	\end{equation}
	Furthermore, we have the following estimates
	\begin{align}
		& |\tz_{\pm 2} \mp \im \tfrac32 \alpha^4|^{5,\alpha_{0}} \lesssim \alpha_{0}^5 \,, \quad  |\tz_{j}|^{5,\alpha_{0}} \lesssim \alpha_{0}^5 \,,  \label{tz.est}\\
		& |\cY|_{-m,s}^{5,\alpha_{0}} \lesssim_{s,m} \alpha_{0}^3 \,, \quad |\cZ_{0,\alpha}|_{-m,s}^{5,\alpha_{0}}\lesssim_{s,m}\alpha_{0}^{4}\,. \label{Y.gen.est}
	\end{align}
\end{pro}

\begin{proof}
	Without further relabeling, we rescale
	\begin{equation}\label{rescale}
		\cE_{0,\alpha} \rightsquigarrow \e \, \cE_{0,\alpha} \,, \quad  \cY \rightsquigarrow \e \, \cY \,, \quad \e := \alpha_{0}^3 \,,
	\end{equation}
	so that $|\cE_{0,\alpha}|_{0,{s_0}}^{5,\alpha_{0}}\lesssim 1$. For simplicity in the notation, we also drop the subscripts. We compute
	\begin{align}
		\cL^{(2)} & = \cU^{-1} \cL^{(1)} \cU \\
		& = \omega\,\pa_{\vf} + \cU^{-1}\big( \cD+ \e\, \cD \cY + \e \, \cE+ {\e^2} \, \cE\cY \big) \\
		& = \omega \,\pa_{\vf} + \cD+ \e \,\cU^{-1} \big(  [\cD,\cY]  +\cE+ \e \, \cE \cY \big)\,.
	\end{align}
	We seek  operators $\cY$ and $\cZ$ that solve the quadratic homological equation
	\begin{equation}\label{nonlin.hom.eq}
	[\cD,\cY]  +\cE+ \e \, \cE \cY = \cZ+ \e \, \cY \cZ=\cU\cZ\,.
	\end{equation}
	It is convenient to look for $\cY$ and $\cZ$ in the form of formal power series expansions in $\varepsilon$
	\begin{equation}\label{YZ.series}
		\cY = \sum_{n\geq 0} \e^n \cY_{n} \,, \quad \cZ= \sum_{n\geq 0} \e^n \cZ_{n} \,.
	\end{equation}
	Inserting \eqref{YZ.series} into \eqref{nonlin.hom.eq}, we get
	\begin{equation}
		\sum_{n\geq 0} [\cD,\cY_{n}] + \cE +  \sum_{n\geq 0} \e^{n+1} \cE \cY_{n} = \sum_{n\geq 0} \e^n \cZ_n +  \sum_{n\geq 0} \e^{n+1}\cX_{n} \,, \quad \cX_{n}:= \sum_{k=0}^{n} \cY_{k} \cZ_{n-k} \,.
	\end{equation}
	By equating the powers of $\e$, we obtain the recursive system
	\begin{align}
		[\cD,\cY_{0}] &= -\cE + \cZ_{0} \,, \label{Y0z0.eq} \\
		[\cD,\cY_{n}] & = - \cE \cY_{n-1} +\cX_{n-1} +\cZ_{n} \,, \quad \forall\, n\geq 1 \,. \label{YnZn.eq}
	\end{align}
	We now split the proof into three steps.
	\\[1mm]
	\noindent {\sc Step I. }{\bf Verification of the gap condition.} To solve \eqref{Y0z0.eq}, we want to apply Lemma \ref{lemma.sol.hom}. To this end, we have to check the gap condition \eqref{gap.D} for the operator $\cD$ in \eqref{diag.0.alpha}, that is
    \begin{equation}\label{gap.claim}
           \begin{cases}
        \inf_{\alpha\in[\alpha_{1},\alpha_{2}]}|\lambda_{j}(\alpha)-\lambda_{j'}(\alpha)| \geq  \nu_{1}|j-j'|\,, \\
         \sup_{k=1,...,k_0} \sup_{\alpha\in [\alpha_{1},\alpha_{2}]}\big|\pa_{\alpha}^k(\lambda_{j}(\alpha)-\lambda_{j'}(\alpha))\big|  \leq \nu_{2}|j-j'|\,,
    \end{cases} \quad      \forall\,j,j'\in \Z_{\rm ph}\,, \ j\neq j'\,,
    \end{equation}
    for some constants $\nu_1,\nu_2>0.$ To this end, we distinguish between two different cases:
    \\[1mm]
   \noindent $\blacktriangleright$ If $|j|\geq 2$ and $j'=1$, then
\begin{align}
\lambda_{j}(\alpha)-\lambda_{j'}(\alpha)
= \lambda_{j}(\alpha)
= \tfrac12 \big(j-\mathrm{sgn}(j)\big)
+ \big(\tm_{1}(\alpha)-\tfrac12\big)j \,.
\end{align}
Using estimate \eqref{m1alpha0.est} in Lemma~\ref{lemma.reparam.0}-(i), we obtain
\begin{align}
\big| \lambda_{j}(\alpha) \big|
&\geq \tfrac12 \big||j|-1\big|
     - \big|\tm_{1}(\alpha)-\tfrac12\big|\,|j| \\
&\geq \tfrac14 |j| - C \alpha_0^4 |j|
 \geq \tfrac18 |j|
 \geq \tfrac{1}{12} |j-1| \,,
\end{align}
for $\alpha_0 \ll 1$ sufficiently small. Here we used that
\[
\big||j|-1\big| \geq \tfrac12 |j| \geq \tfrac13 |j-1| \qquad \text{for } |j|\geq 2 \,,
\]
Moreover, again by Lemma~\ref{lemma.reparam.0}-(i), for any $k=1,\dots,5$ we have
\begin{align}
\big|\partial_\alpha^k \lambda_{j}(\alpha)\big|
= \big|\partial_\alpha^k \tm_{1}\big|\,|j|
\leq C \alpha^{5-k} |j|
\leq C |j|
\leq \tfrac23 C |j-1| \,,
\end{align}
for some constant $C>0$, uniformly for all $\alpha\in[\alpha_1,\alpha_2]$ as defined in
\eqref{alpha1alpha2}.
\\[1mm]
		\noindent $\blacktriangleright$ If $|j|, |j'|\geq 2$ and $j\neq j'$, then
		\begin{align}
				\lambda_{j}(\alpha)-\lambda_{j'}(\alpha) & = \tfrac12 (j-j') - \tfrac12\big({\rm sgn}(j) -{\rm sgn}(j') \big) + (\tm_{1}(\alpha)-\tfrac12)(j-j') \\
				& = \tfrac12(j-j')\Big( 1 - \tfrac{{\rm sgn}(j) -{\rm sgn}(j')}{j-j'} + 2 (\tm_{1}(\alpha)-\tfrac12) \Big) \,.
		\end{align}
		Note that ${\rm sgn}(j) = - {\rm sgn}(j') $ and $|j|, |j'|\geq 2$ imply $|j-j'|\geq 4$. Therefore, by \eqref{m1alpha0.est} in Lemma \eqref{lemma.reparam.0}-$(i)$, we get
		\begin{align}
				|\lambda_{j}(\alpha)-\lambda_{j'}(\alpha)| & \geq \tfrac12 |j-j'|\Big( 1 -  \tfrac{|{\rm sgn}(j) -{\rm sgn}(j')|}{|j-j'|} - 2 |\tm_{1}(\alpha)-\tfrac12| \Big)  \\
				& \geq \tfrac12 |j-j'| \big( 1 - \tfrac12 - C \alpha_{0}^4 \big) \geq \tfrac18 |j-j'|  \,,
		\end{align}
			for $\alpha_{0}\ll 1$ sufficiently small.  Moreover, by Lemma \ref{lemma.reparam.0}-$(i)$, we deduce, for any $k=1,...5$,
            \begin{align}
                \big|\pa_{\alpha}^k(\lambda_{j}(\alpha)-\lambda_{j'}(\alpha))\big| = |\pa_{\alpha}^k \tm_{1}| |j-j'| \leq C \alpha^{5-k} |j-j'| \leq C |j-j'| \,,
            \end{align}
            for some $C>0$ and uniformly for any $\alpha\in [\alpha_{1},\alpha_{2}]$ as in \eqref{alpha1alpha2}.
 	Therefore, we obtain that \eqref{gap.claim} holds with $\nu_{1}=\tfrac{1}{12}$ and $\nu_{2}=C$. 
 		\\[1mm]
\noindent {\sc Step II. }{\bf Estimates in the low regularity norm.} 
 	Let $m\geq 0$ be arbitrary and define
 	\begin{equation}\label{Z0Zn.def}
 		\cZ_{0}:= {\rm diag} \,\cE :={\rm diag}\{ \cE_{j}^{j} \, : \, j \in \Z_{\rm ph} \} \,, \quad \cZ_{n} := {\rm diag} (\cE \cY_{n-1} -\cX_{n-1} )\,.
 	\end{equation}
    Then, we apply Lemma \ref{lemma.sol.hom},
 	to solve \eqref{Y0z0.eq} for $\cY_{0}$ and, recursively, \eqref{YnZn.eq} for $\cY_{n}$. The resulting operators are reversibility preserving and  satisfy the required estimates at the regularity level $s=s_0$, namely,
 	 	\begin{align}\label{low.Y0Z0}
 		|\cY_{0}|_{-m,s_0}^{5,\alpha_{0}} \leq \mu\, |\cE|_{-m,s_0}^{5,\alpha_{0}} \,, \quad 	|\cZ_{0}|_{-m,s_0}^{5,\alpha_{0}} \leq  |\cE|_{-m,s_0}^{5,\alpha_{0}} \,.
 	\end{align}
 	Moreover,  for $n\geq 1$,  Lemma \ref{standard prop decay norm}-$(ii)$ yields
 	\begin{align}
	|\cY_{n}|_{-m,s_0}^{5,\alpha_{0}} & \leq \mu\, C(s_0)  |\cE|_{0,s_0}^{5,\alpha_{0}}   |\cY_{n-1}|_{-m,s_0}^{5,\alpha_{0}}  + \mu\, C(s_0) \sum_{k=0}^{n-1}  |\cY_{k}|_{0,s_0}^{5,\alpha_{0}}   |\cZ_{n-1-k}|_{-m,s_0}^{5,\alpha_{0}} \,, \\
	|\cZ_{n}|_{-m,s_0}^{5,\alpha_{0}} & \leq C(s_0)  |\cE|_{0,s_0}^{5,\alpha_{0}}   |\cY_{n-1}|_{-m,s_0}^{5,\alpha_{0}}  + C(s_0) \sum_{k=0}^{n-1}  |\cY_{k}|_{0,s_0}^{5,\alpha_{0}}   |\cZ_{n-1-k}|_{-m,s_0}^{5,\alpha_{0}}  \,. \label{low.YnZn} 
   \end{align}
	We remark that the constants $C(s_0)>0$ and $\mu=\mu(\nu_{1},\nu_{2})>0$ are independent of $m\geq 0$. We define the positive quantities
	\begin{equation}\label{an.def.s0}
		a_{n,m}(s_0) := 	|\cY_{n}|_{-m,s_0}^{5,\alpha_{0}}  + 	|\cZ_{n}|_{-m,s_0}^{5,\alpha_{0}} \,, \quad \gamma= 1 + \mu \,,
	\end{equation}
	which, by \eqref{low.Y0Z0}, \eqref{low.YnZn}, satisfy the estimates
	\begin{align}
		a_{0,m}(s_0) & \leq \gamma |\cE|_{-m,s_0}^{5,\alpha_{0}} \,, \\
		a_{n,m}(s_0) &  \leq \gamma C(s_0)  |\cE|_{0,s_0}^{5,\alpha_{0}} a_{n-1,m}(s_0)  + \gamma C(s_0) \sum_{k=0}^{n-1}  a_{k,0}(s_0) a_{n-1-k,m}(s_0) \,, \quad n\geq 1 \,.\label{an-s0}
	\end{align}
	Moreover, in terms of the quantities $(b_{n}(s_0))_{n\geq 0}$ defined by
	\begin{equation}\label{anbn}
		a_{n,m}(s_0) = \big( C(s_0) |\cE|_{0,s_0}^{5,\alpha_{0}} \big)^{n} |\cE|_{-m,s_0}^{5,\alpha_{0}}b_{n}(s_0) \,,
	\end{equation}
	the estimates in \eqref{an-s0} can be rewritten in the form
		\begin{align}
	0\leq 	b_0(s_0) & \leq \gamma  \,, \\
	0\leq 	b_n(s_0) &  \leq \gamma  b_{n-1}(s_0)  + \gamma  \sum_{k=0}^{n-1}  b_{k}(s_0) b_{n-1-k}(s_0) \,, \quad n\geq 1 \,.\label{bn-s0}
	\end{align}
	Note that the system \eqref{bn-s0} governing the sequence $(b_{n}(s_0))$ is independent of the parameter $m\geq 0$.
	We now introduce the auxilary discrete dynamical system
	\begin{align}
		c_0(s_0) &= \gamma \,, \\
		c_n(s_0) & = \gamma  c_{n-1}(s_0)  + \gamma  \sum_{k=0}^{n-1}  c_{k}(s_0) c_{n-1-k}(s_0) \,, \quad n\geq 1 \,. \label{cn-s0}
	\end{align}
	A straightforward induction argument shows that \begin{align}\label{bn-cn}0\leq b_n(s_0) \leq c_n(s_0)\,, \quad  \forall \, n\in\N_{0}\,.
    \end{align}
    To analyze the growth of $(c_n(s_0))_{n\in\N_0}$, we introduce the associated generating function 
	\begin{equation}\label{F0-gener}
		F_0(x) := \sum_{n\geq 0} c_n(s_0) x^n \,.
	\end{equation}
	By \eqref{cn-s0}, we get
	\begin{align}
		F_0(x) & =  \gamma + \sum_{n\geq 1} c_{n}(s_0) x^n \\
		& = \gamma + \gamma \sum_{n\geq 1} c_{n-1}(s_0) x^n + \gamma \sum_{n\geq 1} \Big( \sum_{k=0}^{n-1} c_{k}(s_0) c_{n-1-k}(s_0) \Big) x^{n} \\
		& = \gamma + \gamma x F_0(x) + \gamma x F_0^2(x) \,.
	\end{align}
	That is, $F_0(x)$ satisfies the quadratic  equation
	\begin{equation}
		\gamma x F_0^2(x) + (\gamma x -1)F_0(x) + \gamma = 0 \,,
	\end{equation}
	which, together with the condition $F_0(0)=\gamma$, admits the explicit solution, 
	\begin{align}
		F_0(x) & = \frac{1}{2\gamma x} \Big( 1 - \gamma x - \sqrt{ (1-\gamma x)^2 - 4\gamma^2 x } \Big) \\
		& =  2\Big( \gamma^{-1} -  x + \sqrt{ (\gamma^{-1}- x)^2 - 4 x } \Big)^{-1} \,.\label{F0-form}
	\end{align}
The domain of definition of $F_0$ containing $x=0$, which corresponds to the region where the power series in \eqref{F0-gener} converges, determined by requiring that $F_0(x)$ be real-valued. This amounts to imposing the nonnegativity of the discriminant of the quadratic equation, namely
	\begin{equation}
	\Delta = 	(\gamma^{-1}- x)^2 - 4 x  = x^2 - 2(2+\gamma^{-1}) +\gamma^{-2} \geq 0\,.
	\end{equation}
	The smallest positive zero of $\Delta= 0$ is given by
	\begin{equation}\label{x.star}
		x_* \equiv x_*(\gamma)= 2 + \gamma^{-1} - 2 \sqrt{1+ \gamma^{-1}} \in (0, \gamma^{-1}) \,,
	\end{equation}
	which ensures that, for any $x\in [0,x_*]$ the series in \eqref{F0-gener} converges, with 
	\begin{align}\label{F0xstar}
    F_0(x_*)=\big( \sqrt{1+\gamma^{-1}} -1 \big)^{-1} = \gamma \big( \sqrt{1+\gamma^{-1}} +1 \big)> \gamma >0\,.
    \end{align} 
	Therefore, by \eqref{YZ.series}, \eqref{an.def.s0},\eqref{anbn} and \eqref{bn-cn}, we conclude that
	\begin{align}
		|\cY|_{-m,s_0}^{5,\alpha_{0}} + |\cZ|_{-m,s_0}^{5,\alpha_{0}} &\leq \sum_{n\geq 0}  \e^n \big( |\cY_{n}|_{-m,s_0}^{5,\alpha_{0}} + |\cZ_{n}|_{-m,s_0}^{5,\alpha_{0}}  \big) \\
        &\leq |\cE|_{-m,s_0}^{5,\alpha_{0}}\sum_{n\geq 0}   \big( \e\, C(s_0) |\cE|_{0,s_0}^{5,\alpha_{0}} \big)^{n} c_{n}(s_0) \\
		\label{lownorm.est}& \leq\, |\cE|_{-m,s_0}^{5,\alpha_{0}} F_0\big( \e\, C(s_0) |\cE|_{0,s_0}^{5,\alpha_{0}} \big) \,,
	\end{align}
	provided that the smallness condition \begin{align}\label{smalln-ep}{}\alpha_{0}^3 = \e \leq x_* \big({C(s_0)}|\cE|_{0,s_0}^{5,\alpha_{0}} \big)^{-1}
    \end{align}
    is satisfied. Remark that this assumption is   independent of $m\geq 0$.
     Since the coefficients $c_{n}(s_0)$ in \eqref{F0-gener} are positive, we deduce that
	 \begin{equation}
	 	c_{n}(s_0) \leq x_*^{-n} F_0(x_*)  \quad  \forall\, n\in \N_{0}  \,,
	 \end{equation}
	 which implies in view of  \eqref{anbn}  and \eqref{bn-cn} that
	 \begin{align}\label{anm-est}
	 	a_{n,m}(s_0) \leq  F_0(x_*) |\cE|_{-m,s_0}^{5,\alpha_{0}} \Big( {C(s_0)}\,x_*^{-1} |\cE|_{0,s_0}^{5,\alpha_{0}} \Big)^n \quad \forall \, n\in\N_{0} \,.
	 \end{align}
	\noindent {\sc Step III. }{\bf Estimates in high regularity norms.}
    We now extend the previous estimates to the regime $s\in[s_0,S]$ and establish the convergence of the series under smallness assumptions that depend only on the upper  regularity threshold  $S$. By Lemma~\ref{lemma.sol.hom} together with  the tame estimates provided by  Lemma \ref{standard prop decay norm}-$(ii)$, we deduce  that the sequences of operators $(\cY_{n})_{n\geq 0}$ and $(\cZ_{n})_{n\geq 0}$ solving \eqref{Y0z0.eq} and \eqref{YnZn.eq} satisfy the following estimates,
	 	\begin{align}
		 		|\cY_{0}|_{-m,s}^{5,\alpha_{0}} \leq \mu\, |\cE|_{-m,s}^{5,\alpha_{0}} \,, \quad 	|\cZ_{0}|_{-m,s}^{5,\alpha_{0}} \leq  |\cE|_{-m,s}^{5,\alpha_{0}} \,, \label{high.Y0Z0}
		 	\end{align}
	and,  for all $n\geq 1$,  
	 	\begin{align}
		 		|\cY_{n}|_{-m,s}^{5,\alpha_{0}} & \leq \mu \, \big( C(s_0)  |\cE|_{0,s_0}^{5,\alpha_{0}}   |\cY_{n-1}|_{-m,s}^{5,\alpha_{0}} + C(s)  |\cE|_{0,s}^{5,\alpha_{0}}   |\cY_{n-1}|_{-m,s_0}^{5,\alpha_{0}} \big) \\
		 		&  + \mu\,  \sum_{k=0}^{n-1} \big( C(s_0) |\cY_{k}|_{0,s_0}^{5,\alpha_{0}}   |\cZ_{n-1-k}|_{-m,s}^{5,\alpha_{0}} + C(s) |\cY_{k}|_{0,s}^{5,\alpha_{0}}   |\cZ_{n-1-k}|_{-m,s_0}^{5,\alpha_{0}} \big) \,, \\
		 			|\cZ_{n}|_{-m,s}^{5,\alpha_{0}} & \leq C(s_0)  |\cE|_{0,s_0}^{5,\alpha_{0}}   |\cY_{n-1}|_{-m,s}^{5,\alpha_{0}} + C(s)  |\cE|_{0,s}^{5,\alpha_{0}}   |\cY_{n-1}|_{-m,s_0}^{5,\alpha_{0}} \\
		 		&  +  \sum_{k=0}^{n-1} \big( C(s_0) |\cY_{k}|_{0,s_0}^{5,\alpha_{0}}   |\cZ_{n-1-k}|_{-m,s}^{5,\alpha_{0}} + C(s) |\cY_{k}|_{0,s}^{5,\alpha_{0}}   |\cZ_{n-1-k}|_{-m,s_0}^{5,\alpha_{0}} \big)  \,. \label{high.YnZn} 
		 	\end{align}
		 	We remark that the constants $C(s_0), C(s)>0$ and $\mu=\mu(\nu_{1},\nu_{2})>0$ are independent of $m\geq 0$. 
		We define the positive quantities
		 	\begin{equation}\label{an.def.s}
		 		a_{n,m}(s) := 	|\cY_{n}|_{-m,s}^{5,\alpha_{0}}  + 	|\cZ_{n}|_{-m,s}^{5,\alpha_{0}} \,, \quad \gamma= 1 + \mu \,,
		 	\end{equation}
		 	which, by \eqref{high.Y0Z0}, \eqref{high.YnZn}, satisfy the estimates
		 	\begin{align}
		 	a_{0,m}(s) & \leq \gamma |\cE|_{-m,s}^{5,\alpha_{0}} \,, \\
		 	a_{n,m}(s) & \leq \gamma \big( C(s_0)  |\cE|_{0,s_0}^{5,\alpha_{0}}  a_{n-1,m}(s) + C(s)  |\cE|_{0,s}^{5,\alpha_{0}}  a_{n-1,m}(s_0) \big) \\
		 		&  +\gamma  \sum_{k=0}^{n-1} \big( C(s_0) a_{k,0}(s_0)  a_{n-1-k,m}(s) + C(s)a_{n-1-k,0}(s) a_{k,m}(s_0) \big) \,. \label{an-s}
		 	\end{align}
            Applying \eqref{anm-est} from Step II,  together with the fact  that $ \gamma x_* \in (0,1)$ by \eqref{x.star}, and $F_0(x_*)>\gamma>1$ by \eqref{F0xstar} and \eqref{an.def.s},
            yields, for $n\geq 1$,
            \begin{align}\label{an-m,s}
		 	a_{n,m}(s) & \leq   C(\gamma, s_0)\,  \delta\,  a_{n-1,m}(s) + C(\gamma,s) \gamma C(s_0)  |\cE|_{-m,s_0}^{5,\alpha_{0}}   |\cE|_{0,s}^{5,\alpha_{0}} \delta^{n-1} \\
		 	\nonumber 	&  +  \sum_{k=0}^{n-1}  C(\gamma,s_0)   \,\delta^{k+1} \, a_{n-1-k,m}(s) +\sum_{k=0}^{n-1} C(\gamma,s) \gamma C(s_0)|\cE|_{-m,s_0}^{5,\alpha_{0}} \delta^{k}\,a_{n-1-k,0}(s)  \,,
		 	\end{align}
            where the constants 
            \begin{equation}\label{C.gamma.s}
              \forall\, s \in [s_0,S]\,, \quad   C(\gamma,s):= C(s)C(s_0)^{-1} F_0(x_*(\gamma
            )) \geq F_0(x_*(\gamma)) =: C(\gamma,s_0) >1 \,,
            \end{equation}
            and we have  used the notation 
            \begin{align}\label{small-delta}\delta:={C(s_0)}\,x_*^{-1} |\cE|_{0,s_0}^{5,\alpha_{0}} \,.
            \end{align} We introduce the generating function
		 	\begin{equation}\label{F-gener}
		 		\wtF_m(x) := \sum_{n\geq 0} a_{n,m}(s) x^n \,.
		 	\end{equation}
            When $m=0$, using also \eqref{C.gamma.s}, \eqref{small-delta}, we get
            \begin{align*}
		 	\forall \, n\geq1\,, \quad a_{n,0}(s) & \leq   C(\gamma, s_0)\,  \delta\,  a_{n-1,0}(s) + C(\gamma,s)  |\cE|_{0,s}^{5,\alpha_{0}}\,   \delta^{n}   +  2C(\gamma,s)\sum_{k=0}^{n-1}     \,\delta^{k+1} \, a_{n-1-k,0}(s)\,.
		 	\end{align*}
            Then
            \begin{align*}
		 	\wtF_0(x) & \leq  \gamma |\cE|_{0,s}^{5,\alpha_{0}}+ C(\gamma, s_0)\,  \delta\,  x\wtF_0(x) + C(\gamma,s)  |\cE|_{0,s}^{5,\alpha_{0}}\,   \sum_{n\geq 1}\delta^{n} x^n\\
            &+ 2 C(\gamma,s)\delta \sum_{n\geq1}x^n\sum_{k=0}^{n-1}     \,\delta^{k} \, a_{n-1-k,0}(s)\,.
		 	\end{align*}
            We note that, by the standard Cauchy product,
            \begin{align*}
            \sum_{n\geq1}x^n\sum_{k=0}^{n-1}     \,\delta^{k} \, a_{n-1-k,0}(s)&=x \wtF_0(x)\sum_{n\geq0}
        \delta^n x^n =\frac{x}{1-\delta x} \wtF_0(x)\,.
            \end{align*}
            Consequently
            \begin{align}\label{F1=gen}
		  	 \wtF_0(x) & \leq  \gamma |\cE|_{0,s}^{5,\alpha_{0}}+ C(\gamma, s_0)\,  \delta\,  x\wtF_0(x) + C(\gamma,s)  |\cE|_{0,s}^{5,\alpha_{0}}\,   \frac{\delta x}{1-\delta x} + 2 C(\gamma,s)\frac{\delta x}{1-\delta x} \wtF_0(x)\,.
		 	\end{align}
            From \eqref{YZ.series} and \eqref{an.def.s} we deduce that
            \begin{align*}
		|\cY|_{0,s}^{5,\alpha_{0}} + |\cZ|_{0,s}^{5,\alpha_{0}} &\leq \sum_{n\geq 0}  \e^n \big( |\cY_{n}|_{0,s}^{5,\alpha_{0}} + |\cZ_{n}|_{0,s}^{5,\alpha_{0}}  \big) \leq \wtF_0(\e)\,.
	\end{align*}
    It follows from \eqref{F1=gen} that
    \begin{align}
	 	\wtF_0(\e) & \leq  \gamma |\cE|_{0,s}^{5,\alpha_{0}}+ C(\gamma, s_0)\,  \delta\,  \e \wtF_0(\e) + C(\gamma,s)  |\cE|_{0,s}^{5,\alpha_{0}}\,   \frac{\delta \e}{1-\delta \e}+  2C(\gamma,s)\frac{\delta \e}{1-\delta \e} \wtF_0(\e)\,.
		 	\end{align}
            Now, by imposing the smallness  condition
            \begin{align}\label{assum-pet}
            2C(\gamma,s)\frac{\delta \e}{1-\delta\e}\leq\frac14  \quad \Leftrightarrow \quad  {\delta \e}\leq\frac{1}{8C(\gamma,s)+1} \,,
        \end{align}
            we get
            \begin{align}\label{F1=gen1}
		  |\cY|_{0,s}^{5,\alpha_{0}} + |\cZ|_{0,s}^{5,\alpha_{0}}\le 	\wtF_0(\e) & \leq 3 \gamma |\cE|_{0,s}^{5,\alpha_{0}}\,.
		 	\end{align}
            It remains to address the general case $m\geq 0$. Going  back  to \eqref{an-m,s}, it follows that
            \begin{align*}
		 	\wtF_m(x) & \leq \gamma |\cE|_{-m,s}^{5,\alpha_{0}}+  C(\gamma, s_0)\,  \delta\,  x\wtF_m(x) + C(\gamma,s)   \gamma C(s_0) |\cE|_{-m,s_0}^{5,\alpha_{0}} |\cE|_{0,s}^{5,\alpha_{0}} \frac{x}{1-\delta x} \\
		 	\nonumber 	&  +  \delta C(\gamma,s_0)\sum_{n\geq 1}x^n\sum_{k=0}^{n-1}     \,\delta^{k} \, a_{n-1-k,m}(s) + C(\gamma,s) \gamma C(s_0)|\cE|_{-m,s_0}^{5,\alpha_{0}}\sum_{n\geq 1}x^n\sum_{k=0}^{n-1}  \delta^{k}\,a_{n-1-k,0}(s)  \,.
		 	\end{align*}
          Using once again the Cauchy product we deduce that  
             \begin{align}
		 	 \wtF_m(x) & \leq \gamma |\cE|_{-m,s}^{5,\alpha_{0}}+  C(\gamma, s_0)\,  \delta\,  x\wtF_m(x) + C(\gamma,s)  \gamma C(s_0)  |\cE|_{-m,s_0}^{5,\alpha_{0}} |\cE|_{0,s}^{5,\alpha_{0}}\frac{x}{1-\delta x} \\
		 	&  +  \delta C(\gamma,s_0)\frac{x}{1-\delta x}\wtF_m(x) + C(\gamma,s) \gamma C(s_0)|\cE|_{-m,s_0}^{5,\alpha_{0}}\frac{x}{1-\delta x}\wtF_0(x)  \,.\label{F-m-x}
		 	\end{align}
            From \eqref{YZ.series} and \eqref{an.def.s} we deduce that
            \begin{align}
		|\cY|_{-m,s}^{5,\alpha_{0}} + |\cZ|_{-m,s}^{5,\alpha_{0}} &\leq \sum_{n\geq 0}  \e^n \big( |\cY_{n}|_{-m,s}^{5,\alpha_{0}} + |\cZ_{n}|_{-m,s}^{5,\alpha_{0}}  \big) \leq \wtF_m(\e) \,.
	\end{align}
    Thus, from \eqref{F-m-x} we get that
    \begin{align}\label{F-m-x0}
		 	\nonumber \wtF_m(\e) & \leq \gamma |\cE|_{-m,s}^{5,\alpha_{0}}+  C(\gamma, s_0)\,  \delta\,  \e \wtF_m(\e) + C(\gamma,s)  \gamma C(s_0) |\cE|_{-m,s_0}^{5,\alpha_{0}} |\cE|_{0,s}^{5,\alpha_{0}}  \frac{\e}{1-\delta \e} \\
		 	 	&  +   C(\gamma,s_0)\frac{\delta\e}{1-\delta \e}\wtF_m(\e) + C(\gamma,s) \gamma C(s_0) |\cE|_{-m,s_0}^{5,\alpha_{0}}\frac{\e}{1-\delta \e}\wtF_0(\e)  \,.
		 	\end{align}
            Therefore, combining  the condition \eqref{assum-pet}, \eqref{F-m-x0} and \eqref{F1=gen1}, and  adjusting the constant $C(\gamma,s)$, we obtain
            \begin{align*}
		  \wtF_m(\e) & \leq C(\gamma,s) |\cE|_{-m,s}^{5,\alpha_{0}}+ C(\gamma,s)  |\cE|_{0,s}^{5,\alpha_{0}}  |\cE|_{-m,s_0}^{5,\alpha_{0}}    \,.
		 	\end{align*}
            We remark that, in view of \eqref{small-delta}, the smallness condition \eqref{assum-pet} is ensured provided that
            \begin{align} \label{last.small}
            { \e} |\cE|_{0,s_0}^{5,\alpha_{0}}\leq\frac{x_*}{C(s_0)(4C(\gamma,S)+1)} \,.
        \end{align}
\noindent Since $ |\cE|_{0,s_0}^{5,\alpha_{0}}$ is of order $O(1)$ and $ s \in [s_0,S]$, \eqref{last.small} reduces to be a smallness condition on $ \alpha_0$ alone.

		 Finally, we invert the rescaling in \eqref{rescale} to deduce \eqref{Y.gen.est} from \eqref{lownorm.est}, \eqref{I.0.alpha} and \eqref{cE.0.alpha.est}.
        We also get   the estimates in \eqref{tz.est} by virtue of \eqref{Z0Zn.def}, \eqref{lownorm.est}, Lemma \ref{standard prop decay norm}-$(v)$, and Lemma \ref{lemma.modes123.E}-$(ii)$, We are left to conclude that $\cU={\rm Id} + \cY$ is real and reversibility preserving and that \eqref{lambda.linear}, \eqref{pure.im.lambda} hold. First, we note that $\cY_{0}$ (resp. $\cZ_{0}$) is real and reversibility preserving (resp. real and reversible), respectively, by \eqref{solu.Y} in Lemma \ref{lemma.sol.hom}, \eqref{Z0Zn.def}, Lemma \ref{prop.matrix.revreal} and the fact that $\cE$ is real and reversible. By induction with \eqref{Y0z0.eq}, \eqref{Z0Zn.def} and by \eqref{YZ.series}, we deduce with the same arguments that $\cY_{n}$ and $\cU={\rm Id} + \cY$ (resp. $\cZ_{n}$ and $\cZ$) are real and reversibility preserving (resp. real and reversible). We conclude that \eqref{lambda.linear}, \eqref{pure.im.lambda} hold by \eqref{diag.0.alpha} in Proposition \ref{diffeo.conj.no.proj}, Lemma \ref{lemma.modes123.E}-$(i)$ and Lemma \ref{prop.matrix.revreal}-$(i)$,$(ii)$. This concludes the proof.
\end{proof}
In the following corollary, we collect the expansion of the unperturbed eigenvalues that we obtained in Proposition \ref{diffeo.conj.no.proj} and Proposition \ref{prop.diag.red.equi}.
\begin{cor}\label{cor.eigen.asympt}
For $\alpha \in [\alpha_{1},\alpha_{2}]$ as in \eqref{alpha1alpha2},	we have the following expansions for small $\alpha$:
\\[1mm]
\noindent $\blacktriangleright$ for $|j|= 2$,
\begin{equation}
	\begin{aligned}
		-\im\,\lambda_{j}^{(\infty)}(\alpha) &= j \big(\tfrac12 - \alpha^4 +\tc_{1}(\alpha) \big) + {\rm sgn}(j) \big(-\tfrac12 + \tfrac{3}{2} \alpha^4\big) + \tr_{j}(\alpha) \in\R\,;
	\end{aligned}
\end{equation}
\noindent $\blacktriangleright$ for $|j|\geq 3$, 
\begin{equation}
	\begin{aligned}
		-\im\,\lambda_{j}^{(\infty)}(\alpha) &= j \big(\tfrac12 -  \alpha^4 +\tc_{1}(\alpha) \big) - \tfrac12 {\rm sgn}(j) + \tr_{j}(\alpha)  \in \R \,. 
	\end{aligned}
\end{equation}
The constants $\tc_{1}(\alpha)$ and $\tr_{j}(\alpha)$   satisfy  the property
\begin{equation}\label{pure.im.cond}
	\begin{aligned}
		\tc_{1}(\alpha)\in\mathbb{R},& \quad \tr_{-j}(\alpha) = -{\tr}_{j}(\alpha)\in\mathbb{R}\,,
	\end{aligned}
\end{equation}
together with the pointwise estimates, for any $\alpha \in [\alpha_{1},\alpha_{2}]$ and  for any $k=0,...,5$,
\begin{align}\label{small.const.alpha}
	\big| \pa_{\alpha}^k \tc_{1}(\alpha) \big|+\sup_{|j|\geq 2}\big| \pa_{\alpha}^k \tr_{j}(\alpha)\big| & \leq C \alpha^{5-k} \,, 
	\end{align}
for some constant $C>0$ independent of $\alpha$.
\end{cor}
\begin{proof}
	The symmetry property in \eqref{pure.im.cond} follows by \eqref{diag.0.alpha}  in Proposition \ref{diffeo.conj.no.proj} and \eqref{lambda.linear}, \eqref{pure.im.lambda} in Proposition \ref{prop.diag.red.equi}. The estimates in \eqref{small.const.alpha} follow by \eqref{lambda.linear}, \eqref{tz.est} in Proposition \ref{prop.diag.red.equi},  \eqref{diag.0.alpha}  in Proposition \ref{diffeo.conj.no.proj},  and \eqref{m1alpha0.est} in Lemma \ref{lemma.reparam.0}.
\end{proof}
\subsection{Linear dynamics around the equilibrium}

We are now in a position to provide a complete description of the linear dynamics around the equilibrium $\rho = 0$ for \eqref{equation.rho.full}. Our goal is to analyze the solutions $h(\varphi,\theta)$ of the linear system \eqref{linear.rho.full}. To this end, we collect the transformations constructed in the previous subsections, which allow us to diagonalize the time-independent variable-coefficient operator $\cL_{0,\alpha}$ appearing in \eqref{cL0.alpha} (see \eqref{cT.0.alpha}, \eqref{smoothF0}, and \eqref{smoothK0}).

\begin{pro}\label{prop.Malpha}
	Let $   S\geq s_0$ and $s\in [s_0,S]$. There exists $\alpha_{0}=\alpha_{0}(S)\in (0,1)$ small enough such that, for any $\alpha\in [\alpha_{1},\alpha_{2}]$ as in \eqref{alpha1alpha2}, the following hold:
	\\[1mm]
	\noindent $(i)$ There exists an invertible, reversibility preserving operator $\cM_{\alpha}:H_{\circ}^{s}(\T^2)\to H_{\circ}^{s}(\T^2)$  given by
	\begin{equation}\label{Malpha}
		\cM_{\alpha} := \cS_{\rm ph} \cU \,,
	\end{equation}
	with $\cS_{\rm ph}$ as in Lemma $\ref{Lema-decompMon}$ and $\cU$ as in Proposition $\ref{prop.diag.red.equi}$, such that $\cL_{0,\alpha}$ in \eqref{cL0.alpha} is conjugated to
	\begin{equation}
		\cL_{0,\alpha}^{(2)} := \cM_{\alpha}^{-1} \cL_{0,\alpha} \cM_{\alpha} = \omega\,\pa_{\vf} + \cD_{0,\alpha}^{(\infty)} : H_{\circ,{\rm even}}^{s}(\T^2) \to H_{\circ,{\rm odd}}^{s-1} (\T^2)\,,
	\end{equation}
	where the diagonal operator $\cD_{0,\alpha}^{(\infty)}$ is as in Proposition $\ref{prop.diag.red.equi}.$ Moreover, the following estimates hold, for every $h\in H_{\circ,{\rm even}}^{s}(\T^2)$,
	\begin{equation}
		\| \cM_{\alpha}^{\pm 1} h \|_{s}^{5,\alpha_{0}} \lesssim \| h\|_{s}^{5,\alpha_{0}} \,, \quad \| \big(\cM_{\alpha}^{\pm 1} -{\rm Id}\big)h \|_{s}^{5,\alpha_{0}} \lesssim \alpha_{0}^2 \| h\|_{s+1}^{5,\alpha_{0}} \,;
	\end{equation}
    \noindent $(ii)$ We have the following identities
    \begin{align}
       &  \cU^{-1}\Pi_{1,{\bf s}} = \Pi_{1,{\bf s}} \cU = 0 \,, \\
       & \cM_{\alpha}= \cS \,\cU - \Pi_{1,{\bf s}}  \cS \Pi_{\rm ph} \cU \,, \\
       & \cM_{\alpha}^{-1} = \cU^{-1} \cS^{-1} - \cU^{-1} \Pi_{\rm ph} \cS^{-1} \Pi_{1,{\bf s}}  + \cU^{-1} \cR_{\cS,1,{\bf s}} \,,
    \end{align}
    implying that  $\cM_{\alpha}^{\pm 1} - (\cS\,\cU)^{\pm 1}$ is a finite rank operator.
\end{pro}
\begin{proof}
    Item $(i)$ follows by Lemma \ref{Lema-decompMon}, Lemma \ref{diffeo.conj.no.proj} and Proposition \ref{prop.diag.red.equi}. Item $(ii)$ follows by Lemma \ref{Lema-decompMon} and Proposition \ref{prop.diag.red.equi}, using that the map $\cU$ in Proposition \ref{prop.diag.red.equi} satisfies $\cU^{\pm 1}= \Pi_{\rm ph}\cU^{\pm 1}\Pi_{\rm ph}$ by construction.
\end{proof}

We are finally able to conclude the linear analysis at the equilibrium.
\begin{pro}
    For any choice of the frequency $\omega\in\R\setminus\{0\}$ given by
    \begin{equation}
        \omega = -\im \lambda_{n}^{(\infty)}(\alpha) \in \R \setminus\{0\}\,, \quad n \geq 2\,,
    \end{equation}
    with $\lambda_{n}^{(\infty)}(\alpha)$ as in  Corollary $\ref{cor.eigen.asympt}$, there exists a smooth, even, time-periodic solution $h$ of the linear equation (see also \eqref{linear.rho.full})
    \begin{equation}
        \cL_{0,\alpha}[h] = \omega\,\pa_{\vf}h + \pa_{\theta} \big(\tV_{0,\alpha}h + \cR_{0,\alpha}[h] \big) = 0 \,,
    \end{equation}
    taking the form
    \begin{equation}
	\begin{aligned}
		h(\vf,\theta) 
		& =\ta_{1}\cM_{\alpha}[\cos(\,\cdot\,)](\theta) +  \ta_{n}\cos(\vf) \cM_{\alpha}[\cos(n\,\cdot\,)](\theta)+ \ta_{n}\sin(\vf) \cM_{\alpha}[\sin(n\,\cdot\,)](\theta) \,,
	\end{aligned}
\end{equation}
with $\ta_{1},\ta_{n}\in\R$ two arbitrary constants.
\end{pro}
\begin{proof}
    By Proposition \ref{prop.Malpha}, we have that
    \begin{equation}
	\cL_{0,\alpha}[h] = 0 \quad \stackrel{h = \cM_{\alpha}[\tth]}{\Leftrightarrow} \quad \cL_{0,\alpha}^{(2)} [\tth]  = \omega\,\pa_{\vf}\tth + \cD_{0,\alpha}^{(\infty)}[h] = 0 \,.
\end{equation}
We note that
\begin{equation}
	\begin{aligned}
			\cL_{0,\alpha}^{(2)}[\ta_{1}(\vf)\cos(\vartheta)] & =  (\omega\,\pa_{\vf}\ta_{1}(\vf)) \cos(\vartheta) \,, 
	\end{aligned}
\end{equation}
and, for $n\geq 2$, $\ell \in\Z$, using also \eqref{pure.im.lambda} in Proposition \ref{prop.diag.red.equi},
\begin{equation}
	\begin{aligned}
		\cL_{0,\alpha}^{(2)}[e^{\im(n\vartheta+\ell\vf)}] & = \im \,\big( \omega\,\ell -\im\,\lambda_{n}^{(\infty)}(\alpha) \big) e^{\im(n\vartheta+\ell\vf)} \,, \\
		\cL_{0,\alpha}^{(\infty)}[e^{-\im(n\vartheta+\ell\vf)}] & = \big( -\im\,\omega\,\ell + \lambda_{-n}^{(\infty)}(\alpha) \big)e^{-\im(n\vartheta+\ell\vf)} \\
		& = -\im \, \big( \omega\,\ell -\im\,\lambda_{n}^{(\infty)}(\alpha) \big) e^{-\im(n\vartheta+\ell\vf)}\,,
	\end{aligned}
\end{equation}
from which we deduce that
\begin{equation}
	\begin{aligned}
		\cL_{0,\alpha}^{(\infty)}[\cos(n\vartheta+\ell\vf)] & = - \big( \omega\,\ell -\im\,\lambda_{n}^{(\infty)}(\alpha) \big) \sin(n\vartheta+\ell\vf) \,.
	\end{aligned}
\end{equation}
Therefore, any trigonometric function of the form
\begin{equation}
	\tth(\vf,\vartheta) = \ta_{1} \cos(\vartheta) + \ta_{n} \cos(n\vartheta - \vf) \,, \quad \ta_{1}\in\R\,, \ \ta_{n}\in\R\setminus\{0\}\,, \quad n\geq 2\,,
\end{equation}
is an element, ${\rm even}$ in $(\vf,\vartheta),$ in the kernel of $\cL_{0,\alpha}^{(2)}$ provided that 
$\omega= -\im \,\lambda_{n}^{(\infty)}(\alpha) \in \R$. It follows that the function
\begin{equation}
	\begin{aligned}
		h(\vf,\theta) &= \cM_{\alpha}[\ta_{1}\cos(\,\cdot\,)+ \ta_{n} \cos(-\vf+n\,\cdot\,) ](\theta) \\
		& =\ta_{1}\cM_{\alpha}[\cos(\,\cdot\,)](\theta) +  \ta_{n}\cos(\vf) \cM_{\alpha}[\cos(n\,\cdot\,)](\theta)+ \ta_{n}\sin(\vf) \cM_{\alpha}[\sin(n\,\cdot\,)](\theta)
	\end{aligned}
\end{equation}
is a time-periodic solution of the linear equation $\cL_{0,\alpha}[h]=0$ as soon as $\omega = -\im\, \lambda_{n}^{(\infty)}(\alpha)$. 
\end{proof}
The goal of the following sections is now to construct a small amplitude time-periodic solution of the nonlinear equation \ref{equation.rho.full} that bifurcates from the equilibrium configuration and is supported at the leading order of the same mode $\tJ \in \N$ of a given linear solution, with a frequency of oscillation $\omega_{\varepsilon}\in\R\setminus\{0\}$ close to $\omega=-\im\,\lambda_{\tJ}^{(\infty)}(\alpha)$. The mode $\tJ\in\N$ will be referred to as {\it tangential mode}, and, from now on, it is fixed satisfying the condition
\begin{equation}\label{tangential.mode}
	\tJ \geq 2 \,.
\end{equation}
Complementarily, we define as {\it normal modes} all the remaining modes of oscillations in the phase space different from the mode 1 and from the tangential mode $\tJ$. In particular, we denote
\begin{equation}\label{Z.perp}
	\Z_{\perp} := \Z_{\rm ph} \setminus\{1,\pm \tJ\} = \big\{ j \in \Z \,:\, |j|\geq 2 , \ j\neq \pm \tJ \big\} \,.
\end{equation}

\section{Nonlinear reformulation and Lyapunov-Schmidt reduction}\label{Luap-Schm}
In the previous section, we linearized the nonlinear equations in \eqref{equation.rho.full} around the stationary solution and analyzed the kernel of the corresponding linear operator \eqref{linear.rho.full}, thereby obtaining a family of linear time-periodic solutions near the equilibrium. After fixing the tangential site as in \eqref{tangential.mode}, our next objective is to construct, around  one of these small-amplitude periodic solutions  of the linearized system, a time-periodic solution of the full nonlinear problem. 
To carry out this bifurcation analysis, it is convenient to reformulate the nonlinear equations so as to achieve two goals:
\begin{itemize}
    \item  The linearization at any state is already diagonal at the leading order; 
     \item The system is expressed in coordinates that clearly separate the tangential and normal components of the solution.
\end{itemize}
\paragraph{Linear change of variable.} This preparatory step will allow us to apply the iterative scheme in a framework adapted to the spectral properties of the linearized operator.
We achieve the first point by changing the unknown {\it at the nonlinear level}. Specifically, we set
\begin{equation}
	\rho(\vf,\theta) = \cM_{\alpha}[\fu(\vf,\,\cdot\,)](\theta) \,,
\end{equation}
with $\cM_{\alpha}$ being the linear map given by  \eqref{Malpha}. Then the function $\rho$ is a solution of \eqref{equation.rho.full} if and only if $\fu$ solves the equation
\begin{equation}\label{G(u).equa}
		\bG(\fu):=\bG(\fu;\omega,\alpha,\varepsilon) = 0\,,
\end{equation}
where
\begin{equation}\label{G(u).def}
	\begin{aligned}
		\bG(\fu) & := \cM_{\alpha}^{-1}\big[ \tG(\cM_{\alpha}[\fu])\big] \\
		&= \omega\,\pa_{\vf} \fu +\varepsilon^{-1} \cM_{\alpha}^{-1} \circ \pa_{\theta}\Big( \fU[\fu](\vf) \sqrt{1+2 \fr(\vf,\theta)}\cos(\theta) + F(\fr;\alpha) \Big) \,,
	\end{aligned}
\end{equation}
with
\begin{equation}\label{tU.fr}
	\begin{aligned}
		\fU[\fu](\vf) := \tU_{\alpha} +\varepsilon \tW[\cM_{\alpha}[\fu]](\vf)\,, \quad \fr(\vf,\theta):= r_\alpha(\theta) + \varepsilon\cM_{\alpha}[\fu(\vf,\,\cdot\,)](\theta) \,,
	\end{aligned}
\end{equation}
see also Proposition \ref{prop.U1}.
By virtue of Proposition \ref{prop.ralpha} we have  
\begin{equation}\label{bG.equizero}
	\forall \, \omega, \varepsilon, |\alpha|\leqslant\alpha_0,\quad \bG(0;\omega,\alpha,\varepsilon)=0 \,.
\end{equation}
For $s\geq s_0$ and $\delta>0$ we define the closed ball,
\begin{equation*}
	B_{\delta,\textnormal{even}}^{s} := \{ f \in {W_\upsilon^5(\Lambda;H^s_{\circ,\textnormal{even}}(\T^2)} \, : \, \| f\|_{s}^{5,\upsilon}\leq \delta \}\,.
\end{equation*}
Let $\Omega\subset \R\setminus\{0\}$ be an arbitrary closed set.
Recall also the choice of $[\alpha_{1},\alpha_{2}]$ in \eqref{alpha1alpha2}. The first result deals with the regularity of the functional $\bG.$

\begin{pro}\label{lemma.map.bG}
	Let $k\in\N, S>s_0$ and $s\in[s_0,S]$. There exist $\delta>0$ and  $\alpha_{0},\varepsilon_{0}>0$ small enough  such that the following results  hold:
    \\[1mm]
    \noindent $(i)$ The map
    \begin{equation}
		\bG: B_{\delta,\textnormal{even}}^{s} \times \Omega \times [\alpha_{1},\alpha_{2}]\times [-\varepsilon_{0},\varepsilon_{0}] \to{W_\upsilon^5(\Lambda;H^{s-1}_{\circ,\textnormal{odd}}(\T^2))}
	\end{equation}
	is  of class $\cC^k$; 
    \\[1mm]
    \noindent $(ii)$  For any $\overline{\fu} \in B_{\delta,{\rm even}}^{s}$,  the linearized operator $\partial_\fu\bG(\overline{\fu};\omega,\alpha,\varepsilon)$ takes the form
    \begin{align} \partial_\fu\bG(\overline{\fu};\omega,\alpha,\varepsilon) 
   [h]  
    &=\omega\,\pa_{\vf} h + \cM_{\alpha}^{-1}\pa_{\theta} \Big( \tV_{\varepsilon,\alpha}(\vf,\theta) \cM_{\alpha}[h]  + \cR_{\varepsilon,\alpha}\big[ \cM_{\alpha}[h] \big] \Big)\,,
\end{align}
where
\begin{align}
	\tV_{\varepsilon,\alpha}(\vf,\theta) & := \tfrac{\tU_{\alpha} + \varepsilon \tW[\cM_{\alpha}[\overline{\fu}]](\vf) }{{\sqrt{1+2\overline{\fr}(\vf,\theta)}}} \cos(\theta) + \tV[\overline{\fr};\alpha](\vf,\theta) \,, \label{tV.vare} \\
	\cR_{\varepsilon,\alpha} & := \cR(\overline{\fu};\alpha) +\tF_{\varepsilon,\alpha}\,\cK_{\varepsilon,\alpha} \,, \quad \overline{\fr}(\vf,\theta) = r_{\alpha}(\theta) + \varepsilon \cM_{\alpha}[\overline{\fu}(\vf,\,\cdot\,)](\theta)\,, \label{cR.vare}
\end{align}
with $\tV[\overline{\fr};\alpha]$  defined in \eqref{V.at.boundary}, $\cR(\overline{\fr};\alpha)$ defined in \eqref{cR.lin.app}, $\tW[\,\cdot\,]$ in Proposition \ref{prop.U1}, and $\tF_{\varepsilon,\alpha}$, $\cK_{\varepsilon,\alpha}$ in Proposition $\ref{lemma.constrain-Lun1},$ with $\overline{\rho}(\vf,\theta)$ replaced by $\cM_{\alpha}[{\overline{\fu}}({\vf},\,\cdot\,)](\theta)$. Moreover, the following estimates hold,  for any $m\in\N_{0}$ and  $h\in {W_\upsilon^5(\Lambda;H_{\circ,{\rm even}}^s(\T^2))}$,
\begin{align}
	\| \tV_{\varepsilon,\alpha} - \tV_{0,\alpha} \|_{s}^{5,\upsilon} & \lesssim_{s} \varepsilon \big( 1 + \| \overline{\fu} \|_{s}^{5,\upsilon} \big) \,, \\
	\| \partial_\theta^m \big( \cR_{\varepsilon,\alpha}-\cR_{0,\alpha} \big) [h] \|_{s}^{5,\upsilon} & \lesssim_{m,s} \varepsilon \Big( {\big(1+\| \overline{\fu}\|_{s_0+m}^{5,\upsilon}\big)}\| h\|_{s}^{5,\upsilon} + \| \overline{\fu}\|_{s{+m}}^{5,\upsilon}  \| h\|_{s_0}^{5,\upsilon}  \Big) \,;
\end{align}
\noindent $(iii)$ For any $\overline{\fu} \in B_{\delta,{\rm even}}^{s}$ and $(\omega,\alpha,\varepsilon)\in \Omega\times [\alpha_1,\alpha_2]\times [-\varepsilon_0,\varepsilon_0]$,  the second variation 
$$\pa_{\fu}^2\bG(\overline{\fu};\omega,\alpha,\varepsilon):W_\upsilon^5(\Lambda;H_{\circ,{\rm even}}^s(\T^2)) \times W_\upsilon^5(\Lambda;H_{\circ,{\rm even}}^s(\T^2)) \to W_\upsilon^5(\Lambda;H_{\circ,{\rm even}}^{s-1}(\T^2)) $$
satisfies the estimates, for any  $h_1,h_2\in {W_\upsilon^5(\Lambda;H_{\circ,{\rm even}}^s(\T^2))}$,
\begin{equation}
    \|\pa_{\fu}^2\bG(\overline{\fu};\omega,\alpha,\varepsilon)[h_1,h_2]\|_{s}^{5,\upsilon} \lesssim_s \varepsilon^2 \|h_1 \|_{s+1}^{5,\upsilon} \|h_2 \|_{s+1}^{5,\upsilon} \,.
\end{equation}
    \end{pro}
\begin{proof}
    The proof of items $(i)$ and $(ii)$ follows by \eqref{G(u).def}, \eqref{equation.rho.full}, \eqref{full.nonlin.eq.pert}, \eqref{tV.vare}, \eqref{cR.vare}, \eqref{cT.0.alpha}, \eqref{smoothF0}, \eqref{smoothK0} and Proposition \ref{prop.U1}, Proposition \ref{prop-Functional-est}, Proposition \ref{lemma.constrain-Lun1}, Proposition \ref{prop.Malpha} and Lemma \ref{lem-productlaw}. The proof of item $(iii)$ follows by standard arguments, together with item $(i)$ and Lemma \ref{lem-productlaw}.
\end{proof}

 \paragraph{Lyapunov-Schmidt reduction.}
Next, we solve equation \eqref{G(u).equa} by separating the tangential mode and the degenerate first mode from the normal directions. More precisely, we seek solutions \(\fu(\vf,\vartheta)\) of the form
\begin{equation}\label{split.fu}
	\fu(\vf,\vartheta) = \fu_{\intercal}(\vf,\vartheta) + \fu_{\perp}(\vf,\vartheta) \,,
\end{equation}
where the tangential part $\fu_{\intercal}(\vf,\vartheta)$ (including also the mode $\cos(\vartheta)$) has the form
\begin{equation}\label{fu.tang}
	\fu_{\intercal}(\vf,\vartheta) = \fc_{1}(\vf) \cos(\vartheta) + \fa_{\tJ}(\vf)\cos(\tJ\vartheta)+\fb_{\tJ}(\vf)\sin(\tJ\vartheta) \,, \quad \fc_{1}, \fa_{\tJ},\fb_{\tJ}:\T\to\R\,,
\end{equation}
and the normal part $\fu_{\perp}(\vf,\vartheta)$ satisfies (recalling \eqref{Z.perp})
\begin{equation}\label{fu.perp}
	\fu_{\perp}(\vf,\vartheta) = \Pi_{\perp} \fu_{\perp} (\vf,\vartheta):= \sum_{j\in\Z_{\perp}} \wh\fu_{\perp}(\vf,j) \te_{j}(\vartheta) \,.
\end{equation}
We seek solutions in the space 
\({W_\upsilon^5(\Lambda; H_{\circ,{\rm even}}^s(\mathbb{T}^2))}\)
(see \eqref{Hs-circ}). In view of the decomposition introduced in
\eqref{split.fu}, together with \eqref{fu.tang} and \eqref{fu.perp}, we impose
the following parity conditions:
\begin{equation}\label{parities.1}
\fc_{1}(\vf), \, \fa_{\tJ}(\vf)\ \text{ are even in } \vf\,,
\qquad
\fb_{\tJ}(\vf)\ \text{ is odd in } \vf\,,
\qquad
\fu_{\perp}(\vf,\vartheta)\ \text{ is even in } (\vf,\vartheta) \,.
\end{equation}
We then define the tangential projection by
$$
{\Pi_{\intercal}} := \Pi_{\rm ph} - \Pi_{\perp},
$$
where \(\Pi_{\rm ph}\) is defined in \eqref{proj.ph} and \(\Pi_{\perp}\) in
\eqref{fu.perp}. With these definitions, we perform a Lyapunov--Schmidt
reduction: solving equation \eqref{G(u).equa} is equivalent to solving the
system
\begin{equation}
\begin{cases}\label{LyaSch}
	\bG_{\intercal}(\fu_{\intercal},\fu_{\perp};\omega,\alpha,\varepsilon) := \Pi_{\intercal} \bG(\fu_{\intercal}+\fu_{\perp}) = 0\,, \\
	\bG_{\perp}(\fu_{\intercal},\fu_{\perp};\omega,\alpha,\varepsilon) := \Pi_{\perp} \bG(\fu_{\intercal}+\fu_{\perp}) = 0 \,.
\end{cases}
\end{equation}
In the next sections, we proceed as follows:
\begin{enumerate}
	\item For a fixed normal component $\fu_{\perp}$, we solve the
	{\it{bifurcation equation}}
	\begin{equation}
		\bG_{\intercal}(\fu_{\intercal},\fu_{\perp};\omega,\alpha,\varepsilon)=0
	\end{equation}
	with respect to $\fu_{\intercal}=\fu_{\intercal}(\vf,\vartheta;\fu_{\perp})$.
	The tangential component of the solution is obtained by a fixed-point argument,
	after introducing suitable action--angle coordinates in order to invert the
	restricted linearized operator at equilibrium. A crucial point is to track the
	regularity dependence of $\fu_{\intercal}$ on the normal component $\fu_{\perp}$, which, at this stage, 
	plays  the role of a parameter. This analysis is carried out in
	Section~\ref{section.bifurcation};
	\item Next, we solve the {\it{range equation}}
	\begin{equation}\label{bifurcation.eq}
		\wt\bG_{\perp}(\fu_{\perp};\omega,\alpha,\varepsilon) := \bG_{\perp}(\fu_{\intercal}(\fu_{\perp}),\fu_{\perp};\omega,\alpha,\varepsilon) = 0 
	\end{equation}
	with the unknown $\fu_{\perp}$. The normal component is determined
	via a nonlinear Nash--Moser iteration. This requires constructing an approximate
	right inverse of the linearized operator at each approximate solution. Since the
	non-resonance conditions needed to invert these operators induce regularity losses only in
	spatial regularity, it suffices to reduce the transport part-that is, the
	leading-order term-of the linearized operator, and to show that the remaining
	terms are  smoothing in space. This is the content of Section~\ref{section.normal};
	\item While solving the bifurcation equation, we use a hypothetical conjugation argument to
	fix the final frequency of oscillation
	\begin{equation}
		\omega_{\varepsilon}
		=\omega_{\varepsilon}(\fu_{\perp})
		=\omega_{0}+O(\varepsilon),
		\qquad
		\omega_{0}=-\im\,\lambda_{\tJ}^{(\infty)}(\alpha),
	\end{equation}
	which allows us to treat $\omega\in\mathbb{R}$ as a free parameter when solving
	the range equation in Section~\ref{section.normal}. The measure estimates
	associated with the final non-resonance conditions (first Melnikov and transport
	conditions) are established in Section~\ref{section-measures}.
\end{enumerate}

\section{Solution of the bifurcation equation}\label{section.bifurcation}
Following the Lyapunov–Schmidt reduction performed in Section 6, the construction of time-periodic vortex patches is reduced to two coupled problems: a finite-dimensional bifurcation equation, governing the tangential directions associated with symmetries and phase dynamics, and an infinite-dimensional range equation, describing the normal perturbations of the vortex boundary.
In this section, we address the first of these components.
The bifurcation equation is obtained by projecting the nonlinear functional equation onto the kernel of the linearized operator at the approximate traveling dipole. This kernel is three-dimensional in space, reflecting in particular the trivial degeneracy of the first mode inherited from translation invariance. At leading order, the reduced system exhibits a nearly integrable structure. Indeed, in the regime of large separation between the two patches (that is, for small interaction parameter $\alpha$), the dynamics can be regarded as a perturbation of an integrable transport system.
As a consequence, the tangential modes admit a natural formulation in terms of action–angle variables, enabling us to exploit the Hamiltonian and reversible structure of the problem.\\
In this description, we treat the normal component \( \fu_{\perp} \) as a parameter taking values in a suitable function space. We assume that it satisfies the bound (see also \eqref{ansatz.approx})
\begin{equation}\label{ansatz.fu.bifurc}
| \fu_{\perp} |_{s_0+ \mu_0}^{5,\upsilon} {\leqslant  1}\,, \quad \textnormal{for some constant} \quad \mu_0 >0 \,,
\end{equation}
and we proceed to solve the  {\it bifurcation equation}
\begin{equation}
	\bG_{\intercal}(\fu_{\intercal},\fu_{\perp};\omega,\alpha,\varepsilon) = 0
\end{equation}
for $\fu_{\intercal} = \fu_{\intercal}(\vf,\vartheta;\fu_{\perp})$, with $\bG_{\intercal}$ as in \eqref{LyaSch}. In view of \eqref{fu.tang}, and for regularity  purposes, it is convenient to seek functions  $(\fc_{1}(\vf),\fa_{\tJ}(\vf),\fb_{\tJ}(\vf))$ in $H^s(\T)$ such that
	\begin{align}
	\fu_{\intercal}(\vf,\vartheta) = 	\tA(\fc_{1},\fa_{\tJ},\fb_{\tJ})(\vf,\vartheta)& :=  \fc_{1}(\vf) \cos(\vartheta) + \fa_{\tJ}(\vf)\cos(\tJ\vartheta)+\fb_{\tJ}(\vf)\sin(\tJ\vartheta)\,,\label{A.map.abc}
	\end{align}
 together with the parity conditions in \eqref{parities.1} to ensure that $\fu_{\intercal}\in H_{\circ,{\rm even}}^s(\T^2)$. Therefore, we are led to solve the $3\times 3$ nonlinear system 
				\begin{align}
			&\bG_{1}(\fc_{1},\fa_{\tJ},\fb_{\tJ},\fu_{\perp};\omega,\alpha,\varepsilon):= \tfrac{1}{\pi}\big( \bG_{\intercal}(\tA(\fc_{1},\fa_{\tJ},\fb_{\tJ}),\fu_{\perp};\omega,\alpha,\varepsilon),\cos(\,\cdot\,)\big)_{L^2} \label{bG1} \\
			=&  \omega\,\pa_{\vf}\fc_{1}(\vf)+ \frac{2}{\varepsilon} \fint_{\T} \cM_{\alpha}^{-1} \pa_{\theta}\Big[ (\tU_{\alpha}+ \varepsilon \tW(\vf))  \sqrt{1+2 \fr(\vf,\,\cdot\,)} \cos(\,\cdot\,) + F[\fr(\vf,\,\cdot\,);\alpha]  \Big](\eta) \cos(\eta) \wrt\eta = 0  \,, \\
			&\bG_{A}(\fc_{1},\fa_{\tJ},\fb_{\tJ},\fu_{\perp};\omega,\alpha,\varepsilon):= \tfrac{1}{\pi}\big( \bG_{\intercal}(\tA(\fc_{1},\fa_{\tJ},\fb_{\tJ}),\fu_{\perp};\omega,\alpha,\varepsilon),\tfrac{\te_{\tJ}+\te_{-\tJ}}{2}\big)_{L^2} \label{bGA} \\
			=& \omega\,\pa_{\vf}\fa_{\tJ}(\vf) + \frac{2}{\varepsilon} \fint_{\T} \cM_{\alpha}^{-1} \pa_{\theta}\Big[ (\tU_{\alpha}+ \varepsilon \tW(\vf))  \sqrt{1+2 \fr(\vf,\,\cdot\,)} \cos(\,\cdot\,) + F[\fr(\vf,\,\cdot\,);\alpha]  \Big](\eta) \cos(\tJ\eta) \wrt\eta = 0  \,,  \\ &\bG_{B}(\fc_{1},\fa_{\tJ},\fb_{\tJ},\fu_{\perp};\omega,\alpha,\varepsilon):= \tfrac{1}{\pi}\big( \bG_{\intercal}(\tA(\fc_{1},\fa_{\tJ},\fb_{\tJ}),\fu_{\perp};\omega,\alpha,\varepsilon),\tfrac{\te_{\tJ}-\te_{-\tJ}}{2\im}\big)_{L^2} \label{bGB} \\
			=& \omega\,\pa_{\vf}\fb_{\tJ}(\vf) + \frac{2}{\varepsilon} \fint_{\T} \cM_{\alpha}^{-1} \pa_{\theta}\Big[ (\tU_{\alpha}+ \varepsilon \tW(\vf))  \sqrt{1+2 \fr(\vf,\,\cdot\,)} \cos(\,\cdot\,) + F[\fr(\vf,\,\cdot\,);\alpha]  \Big](\eta) \sin(\tJ\eta) \wrt\eta = 0  \,, 
		\end{align}
where we remind from  \eqref{tU.fr} and \eqref{A.map.abc} that
\begin{equation}
\begin{aligned}
	\fr(\vf,\theta) & = r_\alpha(\theta) + \varepsilon \cM_{\alpha}\big[ \tA(\fc_{1},\fa_{\tJ},\fb_{\tJ})(\vf,\,\cdot\,) + \fu_{\perp}(\vf,\,\cdot\,)  \big](\theta)\,. 
\end{aligned}
\end{equation}
 In the following proposition, we rewrite the above system in a compact form in order to clearly separate the linear contributions around the equilibrium from the remaining quadratic nonlinear terms.
 \begin{pro}\label{prop.compact.bifurc}
The following hold:
\\[1mm]
\noindent $(i)$	The system
 	\begin{equation}\label{fX.def}
 			\begin{cases}
 			\bG_{1}(\fX,\fu_{\perp};\omega,\alpha,\varepsilon) = 0\,, \\
 			\bG_{A}(\fX,\fu_{\perp};\omega,\alpha,\varepsilon) = 0\,, \\ 
 			\bG_{B}(\fX,\fu_{\perp};\omega,\alpha,\varepsilon) = 0\,,
 		\end{cases} 
 		 \quad \fX(\vf)= \begin{pmatrix}
 		\fc_{1}(\vf) \\ \fa_{\tJ}(\vf) \\ \fb_{\tJ}(\vf)
 		\end{pmatrix} \,,
 	\end{equation}
 	reads as
 	\begin{equation}\label{compact.eq.1J}
 		\tL[\omega,\alpha] \fX(\vf)+ \varepsilon\, \tQ(\fX,\fu_{\perp};\alpha,\varepsilon) = 0 \,,
 	\end{equation}
 	where
 	the linear operator $\tL$ is given by
 	\begin{equation}\label{linear.tang.matrix}
 		\tL[\omega,\alpha] := \begin{pmatrix}
 			\omega\,\pa_{\vf} & 0 & 0 \\
 			0 & \omega\,\pa_{\vf} & -\im\,\lambda_{\tJ}^{(\infty)}(\alpha) \\
 			0 & \im\,\lambda_{\tJ}^{(\infty)}(\alpha) & \omega\,\pa_{\vf}
 		\end{pmatrix}\,,
 	\end{equation}
 	and the quadratic term
 		\begin{equation}\label{tQ.def}
 		\tQ(\fX,\fu_{\perp}):= \tQ(\fX,\fu_{\perp};\alpha,\varepsilon) := \begin{pmatrix}
 			\tQ_{1}(\fX,\fu_{\perp};\alpha,\varepsilon)  \\ 
 			\tQ_{A}(\fX,\fu_{\perp};\alpha,\varepsilon)  \\ 
 			\tQ_{B}(\fX,\fu_{\perp};\alpha,\varepsilon ) \\ 
 		\end{pmatrix}
 	\end{equation}
 	satisfies the estimates: under the smallness condition
    $$
    \| \fX\|_{s_0}^{5,\upsilon} \leqslant 1\,,\quad  \| \fu_{\perp}\|_{s_0}^{5,\upsilon} \leqslant1 \,,
    $$ 
    we get, for any $s\geq s_0$,
 	\begin{align}
 		\| \tQ(\fX,\fu_{\perp}) \|_{s}^{5,\upsilon}& \lesssim_{s} \,\big( \| \fX\|_{s_0}^{5,\upsilon} + \| \fu_{\perp}\|_{s_0}^{5,\upsilon}   \big)\big( \| \fX\|_{s}^{5,\upsilon} + \| \fu_{\perp}\|_{s}^{5,\upsilon}   \big)\,, \label{quad.est.Q} \\
 			\| \pa_{\fX}\tQ(\fX,\fu_{\perp})[h] \|_{s}^{5,\upsilon}  &\lesssim_{s} \| h \|_{s}^{5,\upsilon}+\big( \| \fX\|_{s}^{5,\upsilon} + \| \fu_{\perp}\|_{s}^{5,\upsilon} \big)\| h \|_{s_0}^{5,\upsilon} \,,  \label{dX.est.Q}  \\
 				\| \pa_{\fu_{\perp}}\tQ(\fX,\fu_{\perp})[h] \|_{s}^{5,\upsilon} &\lesssim_{s} \| h \|_{s}^{5,\upsilon} +  \big( \| \fX\|_{s}^{5,\upsilon} + \| \fu_{\perp}\|_{s}^{5,\upsilon} \big)\| h \|_{s_0}^{5,\upsilon} \,;  \label{du.est.Q} 
 	\end{align}
 	\noindent $(ii)$ Assume that  $\fX$ and  $\fu_{\perp}$ satisfy the symmetry properties \eqref{parities.1}, then we have
 	\begin{equation}
 		\begin{aligned}
 			\tQ_{1}(\fX,\fu_{\perp};\alpha,\varepsilon)= {\rm odd}(\vf) \,, \quad	\tQ_{A}(\fX,\fu_{\perp};\alpha,\varepsilon)= {\rm odd}(\vf)\,, \quad		\tQ_{B}(\fX,\fu_{\perp};\alpha,\varepsilon)=  {\rm even}(\vf) \,.
 		\end{aligned}
 	\end{equation}
 \end{pro}

\begin{proof}
	We start with item $(i)$. To obtain the form in \eqref{compact.eq.1J}, we expand, using also \eqref{bG.equizero},
	\begin{equation}
		\begin{aligned}
			\bG_{1}&(\fX,\fu_{\perp};\omega,\alpha,\varepsilon)= \di_{\fX}\bG_{1}(0,0;\omega,\alpha,\varepsilon)[\fX]+ \di_{\fu_{\perp}}\bG_{1}(0,0;\omega,\alpha,\varepsilon)[\fu_{\perp}]+ \varepsilon \, \tQ_{1} (\fX,\fu_{\perp};\alpha,\varepsilon)  \,,\\
		\bG_{A}&(\fX,\fu_{\perp};\omega,\alpha,\varepsilon)= \di_{\fX}\bG_{A}(0,0;\omega,\alpha,\varepsilon)[\fX]+ \di_{\fu_{\perp}}\bG_{A}(0,0;\omega,\alpha,\varepsilon)[\fu_{\perp}] + \varepsilon \, \tQ_{A} (\fX,\fu_{\perp};\alpha,\varepsilon)  \,,\\
		\bG_{B}&(\fX,\fu_{\perp};\omega,\alpha,\varepsilon)= \di_{\fX}\bG_{B}(0,0;\omega,\alpha,\varepsilon)[\fX]+ \di_{\fu_{\perp}}\bG_{B}(0,0;\omega,\alpha,\varepsilon)[\fu_{\perp}] + \varepsilon \, \tQ_{B} (\fX,\fu_{\perp};\alpha,\varepsilon)  \,,
		\end{aligned}
	\end{equation}
	where the differentials  with respect to $\fX=(\fc_{1},\fa_{\tJ},\fb_{\tJ})^{\top}$ are given by, recalling Proposition \ref{prop.Malpha},
	\begin{equation}
				\begin{aligned}
				\di_{\fX}\bG_{1}(0,0)[\fX]& := \omega\,\pa_{\vf}\fc_{1}(\vf) 
				+  \frac{1}{\pi}\int_{\T} \cM_{\alpha}^{-1}\circ \pa_{\theta} \circ (\tV_{0,\alpha}+\cR_{0,\alpha}) \circ \cM_{\alpha}[ \tA(\fX)(\,\cdot\,) ](\eta) \cos(\eta) \wrt \eta \\
				& =  \omega\,\pa_{\vf} \fc_{1}\,, \\
				\di_{\fX}\bG_{A}(0,0)[\fX]& :=  \omega\,\pa_{\vf}\fa_{\tJ}(\vf)+ \frac{1}{\pi} \int_{\T} \cM_{\alpha}^{-1}\circ \pa_{\theta} \circ (\tV_{0,\alpha}+\cR_{0,\alpha}) \circ  \cM_{\alpha}[ \tA(\fX)(\,\cdot\,) ] (\eta) \cos(\tJ \eta) \wrt \eta \\
				& = \omega\,\pa_{\vf}\fa_{\tJ}  - \im\, \lambda_{\tJ}^{(\infty)}(\alpha) \fb_{\tJ}  \,, \\
				\di_{\fX}\bG_{B}(0,0)[\fX] &:= \omega\,\pa_{\vf}\fb_{\tJ}(\vf)+  \frac{1}{\pi} \int_{\T} \cM_{\alpha}^{-1}\circ \pa_{\theta} \circ (\tV_{0,\alpha}+\cR_{0,\alpha}) \circ \cM_{\alpha}[ \tA(\fX)(\,\cdot\,) ] (\eta) \sin(\tJ \eta) \wrt \eta\\
				& =  \omega\,\pa_{\vf}\fb_{\tJ}  + \im\, \lambda_{\tJ}^{(\infty)}(\alpha) \fa_{\tJ}  \,.
			\end{aligned}
	\end{equation}
	On the other hand, by Proposition \ref{prop.Malpha}, the diagonal structure of the linearized operator implies that the derivatives with respect to \( \fu_{\perp} \) vanish,
	\begin{equation}
				\begin{aligned}
				\di_{\fu_{\perp}}\bG_{1}(0,0)[\fu_{\perp}] & ={\frac1\pi}\int_{\T} \cM_{\alpha}^{-1}\circ \pa_{\theta} \circ (\tV_{0,\alpha}+\cR_{0,\alpha}) \circ \cM_{\alpha}[ \fu_{\perp}(\vf,\,\cdot\,) ](\eta) \cos(\eta) \wrt \eta = 0\,,\\ 
				\di_{\fu_{\perp}}\bG_{A}(0,0)[\fu_{\perp}] & ={\frac1\pi}\int_{\T} \cM_{\alpha}^{-1}\circ \pa_{\theta} \circ (\tV_{0,\alpha}+\cR_{0,\alpha}) \circ \cM_{\alpha}[ \fu_{\perp}(\vf,\,\cdot\,) ](\eta) \cos(\tJ\eta) \wrt \eta = 0\,,\\ 
				\di_{\fu_{\perp}}\bG_{B}(0,0)[\fu_{\perp}]& ={\frac1\pi}\int_{\T} \cM_{\alpha}^{-1}\circ \pa_{\theta} \circ (\tV_{0,\alpha}+\cR_{0,\alpha}) \circ \cM_{\alpha}[ \fu_{\perp}(\vf,\,\cdot\,) ](\eta) \sin(\tJ\eta) \wrt \eta = 0 \,.
			\end{aligned}
	\end{equation}
The  quadratic nonlinearities (with respect to $(\fX,\fu_{\perp})$) are given by 
	\begin{align}
		&
			\begin{aligned}
				\tQ_{1} (\fX,\fu_{\perp}) &:= 	\bG_{1}(\fX,\fu_{\perp})- 	\di_{\fX}\bG_{1}(0,0)[\fX]\\
				&=\frac{1}{{\pi\varepsilon}} \int_{\T} \cM_{\alpha}^{-1}\circ \pa_{\theta}\Big[ (\tU_{\alpha}+ \varepsilon{\tW(\vf)})  \sqrt{1+2 \fr(\vf,\,\cdot\,)} \cos(\,\cdot\,) + F[\fr(\vf,\,\cdot\,);\alpha] \\
				&\quad \quad \quad  -\varepsilon \, (\tV_{0,\alpha}+\cR_{0,\alpha})\circ \cM_{\alpha}[ \tA(\fX)(\,\cdot\,)  ] \Big](\eta) \cos(\eta) \wrt\eta \,,  
			\end{aligned}
        \label{Q1} \\
		& 
			\begin{aligned}
				\tQ_{A} (\fX,\fu_{\perp}) &:= 	\bG_{A}(\fX,\fu_{\perp})- 	\di_{\fX}\bG_{A}(0,0)[\fX] \\
				&	=  \frac{1}{{\pi\varepsilon}} \int_{\T} \cM_{\alpha}^{-1}\circ \pa_{\theta}\Big[ (\tU_{\alpha}+ \varepsilon {\tW(\vf)})  \sqrt{1+2 \fr(\vf,\,\cdot\,)} \cos(\,\cdot\,) + F[\fr(\vf,\,\cdot\,);\alpha]  \\
				&\quad \quad \quad -\varepsilon \, (\tV_{0,\alpha}+\cR_{0,\alpha})\circ \cM_{\alpha}[ \tA(\fX)(\,\cdot\,)  ] \Big](\eta) \cos(\tJ\eta) \wrt\eta \,,
			\end{aligned}
        \label{QA}\\
		&
			\begin{aligned}
				\tQ_{B} (\fX,\fu_{\perp})  &:= 	\bG_{B}(\fX,\fu_{\perp})- 	\di_{\fX}\bG_{B}(0,0)[\fX]\\
				&	=  \frac{1}{{\pi\varepsilon}} \int_{\T} \cM_{\alpha}^{-1}\circ \pa_{\theta}\Big[ (\tU_{\alpha}+ \varepsilon{\tW(\vf)})  \sqrt{1+2 \fr(\vf,\,\cdot\,)} \cos(\,\cdot\,) + F[\fr(\vf,\,\cdot\,);\alpha]  \\
				&\quad \quad \quad -\varepsilon \, (\tV_{0,\alpha}+\cR_{0,\alpha})\circ \cM_{\alpha}[ \tA(\fX)(\,\cdot\,)  ] \Big](\eta) \sin(\tJ\eta) \wrt\eta \,.
			\end{aligned}
        \label{QB} 
	\end{align}
 In particular, we may rewrite $\tQ_{1}$, and similarly for the other quadratic functionals, as a second order remainder of the Taylor expansion
  $$ 
				\tQ_{1} (\fX,\fu_{\perp}) ={ \int_0^1 (1-t)\di^2_{\fX,\fu_{\perp}}\bG_{1}(t\fX,t\fu_{\perp})[(\fX,\fu_{\perp}),(\fX,\fu_{\perp})]\wrt  t} \,.
                $$
	 Therefore, the estimate in \eqref{quad.est.Q}  follows by \eqref{Q1}, \eqref{QA}, \eqref{QB}, \eqref{bG1}, \eqref{bGA}, \eqref{bGB}, \eqref{LyaSch} and Proposition \ref{lemma.map.bG}-$(iii)$. 
The estimates \eqref{dX.est.Q}, \eqref{du.est.Q} follow by similar arguments and we omit the details.
\smallskip
     
	\noindent We now prove item $(ii)$. By Proposition \ref{lemma.constrain} and \eqref{G0_station}, we deduce that the nonlinear functional
	\begin{equation}
		(\tU(\vf),r(\vf,\theta) )\mapsto  \pa_{\theta}\Big(\tU(\vf) \sqrt{1+2r(t,\theta)} \cos(\theta)  + F[r(t,\theta);\alpha] \Big)
	\end{equation}
	maps functions $\big( {\rm even}(\vf), {\rm even}(\vf,\theta)\big)$ into functions ${\rm odd}(\vf,\theta)$. Moreover, by \eqref{parities.1}, \eqref{tU.fr} and the fact that $\cM_{\alpha}$ is parity preserving by Proposition \ref{prop.Malpha}, we have that $\fr(\vf,\theta)={\rm even}(\vf,\theta)$. Therefore, using also that the operators $ \tV_{0,\alpha},\cR_{0,\alpha}$ are parity preserving by Lemma \ref{sa.pp.T.eps} and Proposition \ref{lemma.constrain-Lun1}-$(ii)$, and the explicit formulae in \eqref{Q1}, \eqref{QA}, \eqref{QB},  the claim follows.
\end{proof}

\subsection{Action-angle coordinates} We now aim to solve the nonlinear system \eqref{compact.eq.1J} by means of  a fixed point argument. At this stage, however, the linear operator $\tL[\omega,\alpha]$ appearing in \eqref{linear.tang.matrix} is not diagonal, and its invertibility requires suitable conditions  on the frequency $\omega$. To overcome this difficulty, it is convenient to introduce {\it action–angle coordinates} for the tangential subspace. This choice is motivated by the following considerations:
\begin{itemize}
	\item The linear dynamics as the equilibrium is integrable;
	\item Action–angle variables provide a natural framework for determining the oscillation frequency of the full nonlinear solution.

\end{itemize}
The action–angle variables are introduced through a nonlinear change of coordinates. In the planar setting, their construction naturally follows from polar coordinates.
Let $\xi > 0,$ and let $ \fa_{\tJ}, \fb_{\tJ} $ be two real-valued, $2\pi$-periodic functions. Assume that the associated action
\[
\tI(\vf)
=
\frac{1}{2}\Big( \fa_{\tJ}^2(\vf) + \fb_{\tJ}^2(\vf) - \xi \Big)
\]
is sufficiently small.
We then define the angle variable $\Theta$ as a $2\pi$-periodic function satisfying
\begin{align}\label{action.angle}
\fa_{\tJ}(\vf) + \im\,\fb_{\tJ}(\vf)
=
\sqrt{\xi + 2\tI(\vf)}\,
e^{\im(\vf + \Theta(\vf))}\,.
\end{align}
We intend to prove the following result.
\begin{pro}\label{prop.AA}
	Let the functions $\tI:\T\to\R$ and $\Theta:\T\to\mathbb{R}$ be as before. The following assertions hold:
	\\[1mm]
	\noindent $(i)$ 
    By introducing  the variable
	\begin{equation}
		\fY(\vf):= \big( \fc_{1}(\vf), \tI(\vf),\Theta(\vf) \big)^\top \,,
	\end{equation}
	the nonlinear system \eqref{compact.eq.1J} is transformed into
	\begin{equation}\label{compact.aa.sys}
		\omega\,\pa_{\vf}\fY(\vf) +\bigg(\begin{smallmatrix}
			0 \\ 0 \\ \omega+\im\,\lambda_{\tJ}^{(\infty)}(\alpha)
		\end{smallmatrix} \bigg)+ \varepsilon \, \tP(\fY,\fu_{\perp};\alpha,\varepsilon) = 0\,, 
	\end{equation}
	for some quadratic nonlinearity
	\begin{equation}
		\tP(\fY,\fu_{\perp};\alpha,\varepsilon) :=\begin{pmatrix}
			\tP_{1}(\fY,\fu_{\perp};\alpha,\varepsilon) \\ \tP_{\tI}(\fY,\fu_{\perp};\alpha,\varepsilon) \\ \tP_{\Theta}(\fY,\fu_{\perp};\alpha,\varepsilon)
		\end{pmatrix}\,.
	\end{equation}
	Furthermore, under the smallness condition
    $$
    \| \fY\|_{s_0}^{5,\upsilon} \leqslant 1\,, \quad  \| \fu_{\perp}\|_{s_0}^{5,\upsilon} \leqslant1 \,,
    $$
    we get the following estimates, for any $s\geq s_0$,
	\begin{equation}\label{quad.est.AA}
	    \begin{aligned}
        \| \fX \|_{s}^{5,\upsilon}&  \lesssim_{s}
        1+\|\fY \|_{s}^{5,\upsilon}  \,,   \\
			\| \tP(\fY,\fu_{\perp}) \|_{s}^{5,\upsilon} &  \lesssim_{s}   1 +	\| \fY \|_{s}^{5,\upsilon} + 	\| \fu_{\perp} \|_{s}^{5,\upsilon}  \,, \\
				\| \pa_{\fY}\tP(\fY,\fu_{\perp})[h] \|_{s}^{5,\upsilon} &\lesssim_{s} \| h \|_{s}^{5,\upsilon}+  \big(  \| \fY\|_{s}^{5,\upsilon} + \| \fu_{\perp}\|_{s}^{5,\upsilon} \big) \| h \|_{s_0}^{5,\upsilon} \,, \\
					\| \pa_{\fu_{\perp}}\tP(\fY,\fu_{\perp})[h] \|_{s}^{5,\upsilon} &\lesssim_{s} \| h \|_{s}^{5,\upsilon}+  \big(  \| \fY\|_{s}^{5,\upsilon} + \| \fu_{\perp}\|_{s}^{5,\upsilon} \big) \| h \|_{s_0}^{5,\upsilon} \,;
	\end{aligned}
	\end{equation}
	\noindent $(ii)$ Assuming that  $\fX$ and  $\fu_{\perp}$ satisfy the symmetry properties \eqref{parities.1}, then we have
	\begin{equation}\label{parity.AA}
		\tI(\vf)={\rm even}(\vf) \,, \quad \Theta(\vf)={\rm odd}(\vf) \,,
	\end{equation}
	and
	\begin{align}
		\tP_{1}(\fY,\fu_{\perp};\alpha,\varepsilon)= {\rm odd}(\vf) \,, \quad	\tP_{\tI}(\fY,\fu_{\perp};\alpha,\varepsilon)= {\rm odd}(\vf)\,, \quad		\tP_{\Theta}(\fY,\fu_{\perp};\alpha,\varepsilon)=  {\rm even}(\vf) \,.
	\end{align}
\end{pro}
\begin{proof}
	The proof of item $(i)$ follows by elementary computations. To prove item $(ii)$, first we differentiate \eqref{action.angle} with respect to $\vf$ and solve with respect to $\omega\,\pa_{\vf}\tI$ and $\omega\,\pa_{\vf}\Theta$, obtaining
	\begin{equation}\label{deri.aa}
		\begin{aligned}
			\omega\,\pa_{\vf} \tI &= \fa_{\tJ} \,\omega\,\pa_{\vf}\fa_{\tJ} + \fb_{\tJ}\,\omega\,\pa_{\vf}\fb_{\tJ} \,, \\
			\omega\,\pa_{\vf} \Theta& =\frac{1}{{\xi+2\tI}} \big(-\fb_{\tJ} \omega\,\pa_{\vf} \fa_{\tJ} + \fa_{\tJ} \omega\,\pa_{\vf}\fb_{\tJ} \big) - \omega\,.
		\end{aligned}
	\end{equation}
	By inserting \eqref{compact.eq.1J} into \eqref{deri.aa}, we obtain the system \eqref{compact.aa.sys}, with nonlinearities
	\begin{equation}\label{nonlin.P}
		\begin{aligned}
			\tP_{1}(\fY,\fu_{\perp};\alpha,\varepsilon)&:= \tQ_{1}(\fX,\fu_{\perp};\alpha,\varepsilon)|_{\fX=\fX(\fY)} \,, \\
			\tP_{\tI}(\fY,\fu_{\perp};\alpha,\varepsilon) &:= {\sqrt{\xi+2\tI}} \big( \cos(\vf+\Theta)  \tQ_{A}(\fX,\fu_{\perp};\alpha,\varepsilon)+\sin(\vf+\Theta)  \tQ_{B}(\fX,\fu_{\perp};\alpha,\varepsilon)\big)|_{\fX=\fX(\fY)} \,, \\
			\tP_{\Theta}(\fY,\fu_{\perp};\alpha,\varepsilon) & := \frac{1}{\sqrt{\xi+2\tI}} \big(-\sin(\vf+\Theta)  \tQ_{A}(\fX,\fu_{\perp};\alpha,\varepsilon)+ \cos(\vf+\Theta) \tQ_{B}(\fX,\fu_{\perp};\alpha,\varepsilon)  \big) |_{\fX=\fX(\fY)}\,.
		\end{aligned}
	\end{equation}
	The estimates in \eqref{quad.est.AA} follow by \eqref{nonlin.P}, Lemma \ref{lem-productlaw}-$(v)$ and Proposition \ref{prop.compact.bifurc}-$(i)$. The claimed parities in item $(iii)$ follow by \eqref{nonlin.P}, the parity conditions  \eqref{parities.1} and Proposition \ref{prop.compact.bifurc}-$(ii)$. This concludes the proof.
\end{proof}

\subsection{Modified bifurcation equation}\label{sec-MBE}

We now turn to the resolution of the system \eqref{compact.aa.sys}. In order to guarantee solvability, the frequency $\omega$ must be suitably chosen. Indeed, the third nonlinear component $\tP_{\Theta}(\fY,\fu_{\perp})$ is ${\rm even}(\vf)$ and may therefore possess a nonzero average that needs to be compensated.
At the same time, both the solution $\fY$ and the oscillation frequency depend on the normal component $\fu_{\perp}$, which is still to be determined. The latter requires suitable non-resonance conditions involving the parameter $\omega$.
For this reason, we employ a {\it hypothetical conjugation argument}: instead of directly solving system \eqref{compact.aa.sys}, we consider the modified problem
\begin{equation}\label{mod.aa.sys}
	\omega\,\pa_{\vf}\fY(\vf)
	+
	\begin{pmatrix}
		0 \\ 0 \\ \tw
	\end{pmatrix}
	+
	\varepsilon\, \tP(\fY,\fu_{\perp};\alpha,\varepsilon)
	= 0 \,,
\end{equation}
where we look for solutions with respect to the pair $(\fY(\vf),\tw)$, with $\tw \in \mathbb{R}$;
solutions of the original system \eqref{compact.aa.sys} are then recovered by imposing the compatibility condition
\begin{equation}\label{tw.back}
	\tw = \omega + \im\,\lambda_{\tJ}^{(\infty)}(\alpha) \,.
\end{equation}
In the following theorem, we solve the bifurcation equation in action–angle variables and derive quantitative estimates for the dependence of the solutions on the normal component $\fu_{\perp}$.
\begin{theo}	\label{theo.FP.tangential}
Let $S>s_0$, there exists $\varepsilon_{0} \in (0,1)$ sufficiently small such that, for any $\varepsilon \in (0,\varepsilon_{0})$, 
    the following holds:
\\[1mm]
\noindent $(i)$
	Assume  the smallness condition
    $$
    \| \fY\|_{s_0}^{5,\upsilon} \leqslant 1\,,\quad  \| \fu_{\perp}\|_{s_0}^{5,\upsilon} \leqslant1 \,,
    $$
    and let  $\tw$ be fixed as follows
	\begin{equation}\label{tw.constrain}
		\tw =-\varepsilon \braket{\tP_{\Theta}\big(\fY(\,\cdot\,;\fu_{\perp},\omega,\alpha,\varepsilon),\fu_{\perp};\alpha,\varepsilon\big)}_{\vf}.
	\end{equation}
	Then, for any $s\in[s_0,S],$ we have the estimates
	\begin{align}\label{tw.est}
		|\tw|^{5,\upsilon} &\lesssim \varepsilon,\\ 
			| \pa_{\fY}\tw(\fY,\fu_{\perp})[h] |^{5,\upsilon} &\lesssim  \| h \|_{s}^{5,\upsilon} +\big( \| \fY\|_{s}^{5,\upsilon} + \| \fu_{\perp}\|_{s}^{5,\upsilon} \big) \| h \|_{s_0}^{5,\upsilon} \,, \\
		| \pa_{\fu_{\perp}}\tw(\fY,\fu_{\perp})[h] |^{5,\upsilon} &\lesssim 
        \| h \|_{s}^{5,\upsilon}+  \big( \| \fY\|_{s}^{5,\upsilon} + \| \fu_{\perp}\|_{s}^{5,\upsilon} \big) \| h \|_{s_0}^{5,\upsilon} \,;
	\end{align}
   \noindent $(ii)$ Under the assumption \eqref{tw.constrain}, the nonlinear equation \eqref{mod.aa.sys} admits a  solution $$(\fY(\vf),\tw)=(\fY(\vf;\fu_{\perp},\omega,\alpha,\varepsilon),\tw(\fu_{\perp}, \omega,\alpha,\varepsilon))
   $$
   such that for any $s\in[s_0,S],$
   \begin{align*}
		\| \fY(\fu_{\perp}) \|_{s}^{5,\upsilon} &\leqslant\,C\varepsilon \big(1+ \| \fu_{\perp}\|_{s}^{5,\upsilon} \big)
	\end{align*}
    and
    \begin{align*}
   	\| \di_{\fu_{\perp}} \fY(\fu_{\perp}) [h] \|_{s}^{5,\upsilon} 
   	& \lesssim_{s}  \varepsilon\big(  \| h\|_{s}^{5,\upsilon}  +  \| \fu_{\perp}\|_{s}^{5,\upsilon}   \| h\|_{s_0}^{5,\upsilon}  \big) \,;
   \end{align*}
\noindent $(iii)$  Under the point $(ii)$, the components of $\fX(\vf;\fu_{\perp})$ defined in \eqref{fX.def} satisfy the estimates
\begin{equation}
	\| \di_{\fu_{\perp}}\fc_{1}(\fu_{\perp}) [h] \|_{s}^{5,\upsilon}, \| \di_{\fu_{\perp}}\fa_{\tJ}(\fu_{\perp}) [h] \|_{s}^{5,\upsilon} , \| \di_{\fu_{\perp}}\fb_{\tJ}(\fu_{\perp}) [h] \|_{s}^{5,\upsilon} \lesssim_{s} \varepsilon\big(  \| h\|_{s}^{5,\upsilon}  +  \| \fu_{\perp}\|_{s}^{5,\upsilon}   \| h\|_{s_0}^{5,\upsilon}  \big)\,.
\end{equation}
As a consequence, we obtain that
\begin{equation}
	\| \di_{\fu_{\perp}} \fu_{\intercal}(\fu_{\perp})[h] \|_{s}^{5,\upsilon} \lesssim_{s} \varepsilon \big(  \| h\|_{s}^{5,\upsilon}  +  \| \fu_{\perp}\|_{s}^{5,\upsilon}   \| h\|_{s_0}^{5,\upsilon}  \big)\,.
\end{equation}
\end{theo}

\begin{proof}
	The estimates of  item $(i)$ detailed in \eqref{tw.est} follow by \eqref{tw.constrain} and Proposition \ref{prop.AA}-$(ii)$. As to the second point $(ii)$, we rewrite system \eqref{mod.aa.sys} as
	\begin{equation}\label{tS.eps}
		\fY(\vf) = - (\omega\,\pa_{\vf})^{-1} \Big( \varepsilon \tP(\fY,\fu_{\perp};\alpha,\varepsilon) + \Big(\begin{smallmatrix}
			0  \\ 0 \\ \tw
		\end{smallmatrix}\Big) \Big) =: \tS_{\varepsilon} \big( \fY, \fu_{\perp}; \omega,\alpha,\varepsilon \big) \,,
	\end{equation}
	which is well defined by Proposition \ref{prop.AA}-$(iii)$ and the choice of $\tw$ in \eqref{tw.constrain}.
	We shall provide a solution $\fY(\vf)=\fY(\vf,\fu_{\perp};\omega,\alpha,\varepsilon)$ of \eqref{tS.eps} by using a fixed point argument. Actually, we will show that for $s\in [s_0,S]$, there exists  $\varepsilon_0=\varepsilon(s_0,S)$ small enough, such that for any $\varepsilon\in[0,\varepsilon_0],$   the map $\tS_{\varepsilon}$ is a contraction map in the closed ball 
	\begin{equation}
		B_{s_0}:= \Big\{ \fY \in {W_\upsilon^5(\Lambda; H^{s_0})} \, : \, \|\fY\|_{s_0}^{5,\upsilon} \leqslant 1 , \ \ \fY(\vf) \ \text{satisfies \eqref{parities.1}, \eqref{parity.AA}} \Big\} \,,
	\end{equation}
	for a given $\delta$, for any $\fu_{\perp}= \Pi_{\perp} \fu_{\perp} \in H_{\circ,{\rm even}}^{s}(\T^2)$ satisfying \eqref{ansatz.fu.bifurc}, and for any values of the parameters $(\omega,\alpha ) \in\Lambda\subset \R\setminus\{0\}\times [\alpha_{1},\alpha_{2}]$. By \eqref{tS.eps}, \eqref{tw.constrain} and  Proposition \ref{prop.AA}-$(ii)$, we have, assuming \eqref{ansatz.fu.bifurc}, that
	\begin{align}\label{Low-est-S}
		\nonumber\| \tS_{\varepsilon}(\fY,\fu_{\perp}) \|_{s_0}^{5,\upsilon} & \lesssim_{s} \varepsilon |\omega|^{-1} \| \tP(\fY,\fu_{\perp}) \|_{s_0}^{5,\upsilon} \\
		\nonumber& \lesssim_{s}\varepsilon  \big(1+ \| \fY \|_{s_0}^{5,\upsilon} + \| \fu_{\perp}\|_{s_0}^{5,\upsilon} \big) \\
		& \leq C_{1} \varepsilon\,,
	\end{align}
	for some constants $C_{1}>0$, which ensures that $\tS_{\varepsilon}(\,\cdot\,, \fu_{\perp})$ maps $B_{s_0}$ into itself, provided that $\varepsilon$ is small enough.
	The Frech\'et derivative of $\tS_{\varepsilon}(\fY,\fu_{\perp})$ in \eqref{tS.eps} with respect to $\fY$ is given, using also \eqref{tw.constrain}, by
	\begin{align}
		\pa_{\fY} \tS_{\varepsilon}(\fY,\fu_{\perp}) & =  - (\omega\,\pa_{\vf})^{-1} \Big( \varepsilon \pa_{\fY}\tP(\fY,\fu_{\perp}) + \Big(\begin{smallmatrix}
			0  \\ 0 \\ \pa_{\fY}\tw
		\end{smallmatrix}\Big) \Big) \\
		&  =  - (\omega\,\pa_{\vf})^{-1} \Big( \varepsilon \pa_{\fY}\tP(\fY,\fu_{\perp}) -\varepsilon \Big(\begin{smallmatrix}
			0  \\ 0 \\ \braket{ \pa_{\fY} \tP_{\Theta}(\fY,\fu_{\perp}) }_{\vf}
		\end{smallmatrix}\Big) \Big)  \,. \label{pa.fY.tS}
	\end{align}
	Therefore, by \eqref{tS.eps}, \eqref{pa.fY.tS} and Proposition \ref{prop.AA}-$(ii)$, we deduce that
	\begin{align}
		\| \pa_{\fY} \tS_{\varepsilon}(\fY,\fu_{\perp}) \|_{\cB((H^{s_0}(\T))^3)}^{5,\upsilon} & \lesssim_{s} \varepsilon |\omega|^{-1}  \| \pa_{\fY} \tP_{\varepsilon}(\fY,\fu_{\perp}) \|_{\cB((H^{s_0}(\T))^3)}^{5,\upsilon} \\
		& \lesssim \varepsilon  \big( 1+ \| \fY \|_{s_0}^{5,\upsilon} + \| \fu_{\perp}\|_{s_0}^{5,\upsilon} \big) \\
		& \leqslant C_{2} \varepsilon\leq \tfrac12\,, \label{contr.est}
	\end{align}
	provided that $\varepsilon $ is small enough. We conclude that $\tS_{\varepsilon}$ is a contraction and we deduce the existence of a solution $\fY(\vf)= \fY(\vf,\fu_{\perp};\omega,\alpha,\varepsilon)$ of \eqref{tS.eps} in the unit ball $B_{s_0}$. The regularity persistence of $\fY$ can be deduced from the equation \eqref{tS.eps} in a similar way to \eqref{Low-est-S} leading for any $s\in[s_0,S]$
    \begin{align*}
		\| \fY \|_{s}^{5,\upsilon} &\leqslant  \| \tS_{\varepsilon}(\fY,\fu_{\perp}) \|_{s}^{5,\upsilon}  \leqslant C\varepsilon  \big(1+ \| \fY \|_{s}^{5,\upsilon} + \| \fu_{\perp}\|_{s}^{5,\upsilon} \big) \,.
	\end{align*}
Thus,  for  $\varepsilon>0$ small enough, we get 
\begin{align*}
		\| \fY \|_{s}^{5,\upsilon} &\leqslant\,C\varepsilon \big(1+ \| \fu_{\perp}\|_{s}^{5,\upsilon} \big) \,.
	\end{align*}
   We now turn to the second estimate in $(ii)$ concerning the linearization.
     By differentiating \eqref{tS.eps} with respect to $\fu_{\perp}$ at the solution $(\fY(\vf),\tw)=\big( \fY(\vf,\fu_{\perp}), \tw(\fu_{\perp}) \big)$, we get
\begin{equation}
	\di_{\fu_{\perp}} \fY(\fu_{\perp}) =  \pa_{\fY}\tS_{\varepsilon}\big(\fY(\fu_{\perp}),\fu_{\perp}\big) \big[ \di_{\fu_{\perp}}\fY(\fu_{\perp}) \big] + \pa_{\fu_{\perp}} \tS_{\varepsilon}\big(\fY(\fu_{\perp}),\fu_{\perp}\big) \,,
\end{equation}
and, solving with respect to $\di_{\fu_{\perp}} \fY(\fu_{\perp})$,
\begin{equation}\label{night.eps}
	\big( {\rm Id} - \pa_{\fY}\tS_{\varepsilon}\big(\fY(\fu_{\perp}),\fu_{\perp}\big) \big)\big[ \di_{\fu_{\perp}}\fY(\fu_{\perp}) \big]  = \pa_{\fu_{\perp}} \tS_{\varepsilon}\big(\fY(\fu_{\perp}),\fu_{\perp}\big) \,.
\end{equation}
   First, by \eqref{contr.est}, we deduce that  $ {\rm Id} - \pa_{\fY}\tS_{\varepsilon}\big(\fY(\fu_{\perp}),\fu_{\perp}\big) $ is invertible by Neumann series, with estimate
   \begin{equation*}
   		\big\| \big( {\rm Id} -\pa_{\fY} \tS_{\varepsilon}(\fY,\fu_{\perp})\big)^{-1} \big\|_{\cB((H^{s_0}(\T))^3)}^{5,\upsilon} \leq \Big( 1 - 	\| \pa_{\fY} \tS_{\varepsilon}(\fY,\fu_{\perp}) \|_{\cB((H^{s_0}(\T))^3)}^{5,\upsilon} \Big)^{-1} \leq 2 \,,
   \end{equation*}
   under the same smallness condition for $\varepsilon>0$. In addition, one gets the tame estimates, for any $s\in[s_0,S],$
   \begin{equation}\label{nightL2.eps}
   		\big\| \big( {\rm Id} -\pa_{\fY} \tS_{\varepsilon}(\fY,\fu_{\perp})\big)^{-1}[h] \big\|_{s}^{5,\upsilon} \lesssim  \| h \|_{s}^{5,\upsilon}+ \| \fu_{\perp}\|_{s}^{5,\upsilon} \| h \|_{s_0}^{5,\upsilon} \,.
   \end{equation}
   The Frech\'et derivative of $\tS_{\varepsilon}(\fY,\fu_{\perp})$ in \eqref{tS.eps} with respect to $\fu_{\perp}$ is given, using also \eqref{tw.constrain}, by
   \begin{align}
   	\pa_{\fu_{\perp}} \tS_{\varepsilon}(\fY,\fu_{\perp}) & =  - (\omega\,\pa_{\vf})^{-1} \Big( \varepsilon \pa_{\fu_{\perp}}\tP(\fY,\fu_{\perp}) + \Big(\begin{smallmatrix}
   		0  \\ 0 \\ \pa_{\fu_{\perp}}\tw
   	\end{smallmatrix}\Big) \Big) \\
   	&  =  - (\omega\,\pa_{\vf})^{-1} \Big( \varepsilon \pa_{\fu_{\perp}}\tP(\fY,\fu_{\perp}) -\varepsilon \Big(\begin{smallmatrix}
   		0  \\ 0 \\ \braket{ \pa_{\fu_{\perp}} \tP_{\Theta}(\fY,\fu_{\perp}) }_{\vf}
   	\end{smallmatrix}\Big) \Big)  \,. \label{pa.fuperp.tS}
   \end{align}
   Therefore, by estimate \eqref{tw.est} in item $(i)$, \eqref{pa.fuperp.tS}, Proposition \ref{prop.AA}-$(ii)$ and assuming \eqref{ansatz.fu.bifurc}, we obtain 
   \begin{align}
   	\|  	\pa_{\fu_{\perp}} \tS_{\varepsilon}(\fY,\fu_{\perp}[h]\|_{s}^{5,\upsilon} & \lesssim_{s} \varepsilon \| 	\pa_{\fu_{\perp}} \tP(\fY,\fu_{\perp}[h] \|_{s}^{5,\upsilon} \\
   	& \lesssim_{s} \varepsilon \Big( \big( 1+\| \fY\|_{s_0}^{5,\upsilon} + \| \fu_{\perp}\|_{s_0}^{5,\upsilon} \big) \| h \|_{s}^{5,\upsilon}+  \big( 1+ \| \fY\|_{s}^{5,\upsilon} + \| \fu_{\perp}\|_{s}^{5,\upsilon} \big) \| h \|_{s_0}^{5,\upsilon} \Big) \\
   	& \lesssim_{s}  \varepsilon\big(  \| h\|_{s}^{5,\upsilon}  +  \| \fu_{\perp}\|_{s}^{5,\upsilon}   \| h\|_{s_0}^{5,\upsilon}  \big) \,.\label{night3.eps} 
   \end{align}
   By \eqref{night.eps}, \eqref{nightL2.eps} and \eqref{night3.eps}, we deduce that
   \begin{align}
   	\| \di_{\fu_{\perp}} \fY(\fu_{\perp}) [h] \|_{s}^{5,\upsilon} & = \big\|  \big( {\rm Id} -\pa_{\fY} \tS_{\varepsilon}(\fY,\fu_{\perp})\big)^{-1} \big[  	\pa_{\fu_{\perp}} \tS_{\varepsilon}(\fY,\fu_{\perp}[h] \big] \big\|_{s}^{5,\upsilon} \\
   	& \lesssim_{s}  \varepsilon\big(  \| h\|_{s}^{5,\upsilon}  +  \| \fu_{\perp}\|_{s}^{5,\upsilon}   \| h\|_{s_0}^{5,\upsilon}  \big) \,. \label{night4.eps}
   \end{align}
   This completes the proof of point $(ii)$.\\
Regarding the  point $(iii)$, it follows in a straightforward way  from \eqref{night4.eps}, Proposition \ref{prop.AA}-$(ii)$, \eqref{fX.def} and \eqref{fu.tang}. This concludes the proof.
\end{proof}


\section{Solution of the range equation}\label{section.normal}
In the previous section, we have  solved the bifurcation equation in the tangential directions. We now turn to the infinite-dimensional component of the problem: the range equation. While the bifurcation equation governs the neutral modes associated with symmetries and phase dynamics, the range equation captures the genuinely deformative part of the vortex boundary — namely, the normal oscillations around the translating dipole.
A central feature of the range equation is the structure of its linearization around the approximate solution. The linearized operator exhibits a leading transport term inherited from the underlying Euler dynamics, combined with lower-order nonlocal contributions generated by the interaction between the vortex patches. Although the dominant transport part suggests a simple dynamical behavior, the presence of small divisors and the loss of derivatives induced by the transport mechanism prevent a direct inversion.

To overcome this difficulty, we first analyze the linearized operator in the normal direction and identify its principal structure. Through a change of variables and a conjugation argument, we transform the operator into a form that is diagonal up to a smoothing remainder and admits an approximate inverse satisfying tame estimates. The key idea is to separate the leading transport dynamics from the smoothing remainders and to exploit the Hamiltonian and reversible structure of the system. This reduction is done in the spirit of \cite{BM20, BFM24, FMT25, HHM21, HR21}.

Once an approximate inverse is constructed, we implement a Nash–Moser iterative scheme to solve the nonlinear range equation. This iterative procedure compensates for the loss of derivatives and allows us to control the small divisor effects under suitable non-resonance conditions on the temporal frequency. 
This completes the infinite-dimensional step of the Lyapunov–Schmidt scheme and provides the full solution of the nonlinear problem, up to the measure estimates on the parameter set that will be addressed in Section \ref{section-measures}.

\medskip

We recall from \eqref{LyaSch} and \eqref{bifurcation.eq} that the {range equation} is given by
\begin{equation}\label{range.eq}
	\wt\bG_{\perp}(\fu_{\perp};\omega,\alpha,\varepsilon) = \bG_{\perp} (\fu_{\intercal}(\fu_{\perp}),\fu_{\perp};\omega,\alpha,\varepsilon)=0\,,
\end{equation}
where,  according to  \eqref{fu.tang},
	\begin{align}
		\fu_{\intercal}(\vf,\vartheta;\fu_{\perp}) & = \fu_{\intercal}(\vf,\vartheta;\fu_{\perp},\omega,\alpha,\varepsilon)\\
		& = \fc_{1}(\vf;\fu_{\perp}) \cos(\vartheta) + \fa_{\tJ}(\vf;\fu_{\perp})\cos(\tJ\vartheta)+\fb_{\tJ}(\vf;\fu_{\perp})\sin(\tJ\vartheta) \,, \label{fu.tang.sol.NM}
	\end{align}
Recall that the functions $\fu_{\intercal}(\vf,\vartheta;\fu_{\perp})$ and 
$(\fc_{1}(\vf;\fu_{\perp}), \fa_{\tJ}(\vf;\fu_{\perp}), \fb_{\tJ}(\vf;\fu_{\perp}))$ 
are constructed in Theorem \ref{theo.FP.tangential}-$(ii)$. 
In particular, the triple $(\fc_{1}(\vf), \fa_{\tJ}(\vf), \fb_{\tJ}(\vf))$ 
solves the bifurcation equation \eqref{compact.eq.1J} provided that it solves 
\eqref{compact.aa.sys} and that the parameter $\tw=\tw(\fu_{\perp})$, 
given by Theorem \ref{theo.FP.tangential}-$(i)$, satisfies the compatibility condition \eqref{tw.back}.\\
In this section, we solve \eqref{range.eq}, treating $\omega \in \mathbb{R} \setminus \{0\}$ as a free parameter. 
Only in Section \ref{section-measures} will we link back  $\omega$ to the final oscillation frequency, by imposing  the compatibility condition \eqref{tw.back}, in order to finally  produce a time-periodic solution.
%
From \eqref{G(u).def}, we recall that  the nonlinear functional $\wt\bG_{\perp}(\fu_{\perp})$ appearing in \eqref{range.eq} is given by 
\begin{equation}
	\begin{aligned}
		\wt\bG_{\perp}(\fu_{\perp};\omega,\alpha,\varepsilon) & = \omega\,\pa_{\vf} \fu_{\perp} + \varepsilon^{-1} \Pi_{\perp} \cM_{\alpha}^{-1} \circ \pa_{\theta} \Big(\fU_{\perp}[\fu_{\perp}](\vf) \sqrt{1+2\fr_{\perp}(\vf,\theta)}\cos(\theta) + F(\fr_{\perp};\alpha) \Big) \,,
	\end{aligned}
\end{equation}
where, in view of    \eqref{tU.fr} and \eqref{A.map.abc}, we have 
\begin{equation}
	\begin{aligned}
		\fU_{\perp}[\fu_{\perp}](\vf)& := \tU_{\alpha} + \varepsilon\tW\big[ \cM_{\alpha}\big[ \fu_{\intercal}(\fu_{\perp}) + \fu_{\perp} \big] \big](\vf)\,, \\
		\fr_{\perp}(\vf,\theta) & := r_\alpha(\theta) + \cM_{\alpha}\big[ \fu_{\intercal}(\vf,\,\cdot\,;\fu_{\perp}) + \fu_{\perp}(\vf,\,\cdot\,)  \big](\theta)\,.
	\end{aligned}
\end{equation}
To state our first result , we need to introduce the function spaces.
For $s\in\mathbb{R}$, recalling \eqref{FS-even}, \eqref{FS-odd} and \eqref{fu.perp}, we define
\begin{align}
	H_\perp^{s}&:=H_{\perp}^{s}(\mathbb{T}^{2},\mathbb{R}):=\big\{ h\in H^{s}(\mathbb{T}^{2};\mathbb{R}) \, : \, h = \Pi_{\perp} h \big\}\,, \label{Hnormal.def} \\
	H^{s}_{\perp,{\rm even}}& :=H^{s}_{{\rm even}}\cap H^s_{\perp}\quad \textnormal{ and }\quad H^{s}_{\perp,{\rm odd}}:=H^{s}_{{\rm odd}}\cap H^s_{\perp} \,. \label{Hs-normal}
\end{align}
For $\delta>0$, we define the closed ball
\begin{equation*}
	B_{\delta,{\rm even}}^{s,\perp} := \big\{ f \in W_\upsilon^5(\Lambda;H^s_{\perp,{\rm even}}(\T^2)) \, : \,  \| f\|_{s}^{5,\upsilon}\leqslant \delta \big\}\,,
\end{equation*}
with $(\omega,\alpha)\in\Lambda=\Omega\times[\alpha_1,\alpha_2]$ and  $\Omega\subset \R\setminus\{0\}$ be an arbitrary closed set. The choice of $[\alpha_{1},\alpha_{2}]$ is made in \eqref{alpha1alpha2}.
The first preliminary result is an immediate consequence of \eqref{range.eq}, \eqref{G(u).def}, Proposition \ref{lemma.map.bG} and Theorem \ref{theo.FP.tangential}.
\begin{pro}\label{lemma.map.bG.solo}
	Let $k\in\N, S>s_0$ and $s\in[s_0,S]$. There exist $\delta>0$ and  $\alpha_{0},\varepsilon_{0}>0$ small enough  such that the following results  hold.
	\\[1mm]
	\noindent $(i)$ The map
	\begin{equation}
		\wt\bG_{\perp}: B_{\delta,{\rm even}}^{s,\perp} \times \Lambda\times [-\varepsilon_{0},\varepsilon_{0}] \to W_\upsilon^5(\Lambda;H^{s-1}_{\perp,{\rm odd}}(\T^2))
	\end{equation}
	is  of class $\cC^k$;
	\\[1mm]
	\noindent $(ii)$  For any $\underline{\fu}_{\perp} \in B_{\delta,{\rm even},\perp}^{s}$,  the linearized operator $\di_{\fu_{\perp}}\wt\bG_{\perp}(\underline{\fu}_{\perp};\omega,\alpha,\varepsilon)$ takes the form 
	\begin{align}
		\di_{\fu_{\perp}}\wt\bG_{\perp}(\underline{\fu}_{\perp};\omega,\alpha,\varepsilon)
		[h] 
		&=\omega\,\pa_{\vf} h + \Pi_{\perp}\cM_{\alpha}^{-1}\pa_{\theta} \Big( \tV_{\varepsilon,\alpha}(\vf,\theta)   + \cR_{\varepsilon,\alpha}\Big)\big[ \cM_{\alpha}\circ \big( {\rm Id} + \di_{\fu_{\perp}} \fu_{\intercal}(\underline{\fu}_{\perp}) \big)[h] \big]  \,,
	\end{align}
	where $\tV_{\varepsilon,\alpha}(\vf,\theta)$ and $\cR_{\varepsilon,\alpha}$ are as in \eqref{tV.vare} and  \eqref{cR.vare}, respectively, with $\overline{\fu}$ and $\overline{\fr}$ replaced by
	\begin{equation}
		\overline{\fu} \rightsquigarrow   \underline{\fu}:= \fu_{\intercal}(\underline{\fu}_\perp) + \underline{\fu}_\perp \,, \quad \overline{\fr} \rightsquigarrow \underline{\fr}_\perp(\vf,\theta) := r_\alpha(\theta) + \cM_{\alpha}\big[ \underline{\fu}(\vf,\,\cdot\,) \big](\theta) \,,
	\end{equation}
	and where $\di_{\fu_{\perp}} \fu_{\intercal}(\underline{\fu}_{\perp})$ satisfies the estimates in Theorem $\ref{theo.FP.tangential}$-$(ii)$.
\end{pro}

The main objective of the next sections is to find zeroes of the nonlinear functional 
$\widetilde{\bG}_{\perp}(\fu_{\perp})$ defined in \eqref{range.eq} by means of a Nash–Moser iterative scheme, which will be carried out in Section \ref{section.NASH}. 
As is standard in such singular perturbative settings, the core difficulty lies in the inversion of the linearized operator around a suitable approximate solution. 
Therefore, we begin by analyzing the linearization of $\widetilde{\bG}_{\perp}$ at a generic approximate state in Section \ref{section.linear.eps}, and establishing its invertibility with tame estimates in Section \ref{sect.almost.inv}. 
This linear analysis provides the fundamental ingredient required to control the loss of derivatives and to close the Nash–Moser iteration.

\subsection{Linearization at the approximate solution}\label{section.linear.eps}
We impose the following smallness condition, which will be satisfied by each approximate solution $\underline{\fu}_{\perp}$ generated along the nonlinear Nash–Moser iteration,
\begin{equation}\label{ansatz.approx}
	\| \underline{\fu}_\perp \|_{s_0+\mu_{0}}^{5,\upsilon} \leq 1 \,, \quad \text{for some constant} \ \ \mu_{0}=\mu_{0}(\tau)>0 \,.
\end{equation}
The stability of the eigenvalues of the linearized operator with respect to the state $\fu_{\perp}$ plays a crucial role in the measure estimates of the final Cantor set. 
To quantify this dependence, we introduce the notation: for any functional $h$ depending on $\fu_{\perp}$, we set
\begin{equation}\label{Delta12}
    \Delta_{12} h := h(\fu_{\perp,2}) - h(\fu_{\perp,1})\quad\hbox{and}\quad \Delta_{12}\fu_{\perp}:=\fu_{\perp,2}-\fu_{\perp,1} \,.
\end{equation}
We estimate such variations in another low norm $\| \,\cdot\, \|_{s_1}$ for all Sobolev indexes $s_1$ such that
\begin{equation}\label{s1}
	s_1 + \sigma_{0} \leq s_0 + \mu_{0}\,, \quad \text{for some constant} \ \ \sigma_{0}=\sigma_{0}(\tau)>0 \,.
\end{equation}
The next result  deals with  the structure of the linearized operator $\di_{\fu_{\perp}}\wt\bG_{\perp}(\underline{\fu}_\perp;\omega,\alpha,\varepsilon)$ introduced in Proposition \ref{lemma.map.bG.solo}.


\begin{lem}\label{expand.sL.pert}
	Let $M\in\N, S>s_0$ and $s\in[s_0,S]$. The linearized operator $\di_{\fu_{\perp}}\wt\bG_{\perp}(\underline{\fu}_\perp;\omega,\alpha,\varepsilon)$ is given by
	\begin{align}
		\sL_{\perp} &:= \di_{\fu_{\perp}}\wt\bG_{\perp}(\underline{\fu}_\perp;\omega,\alpha,\varepsilon) 
		= \Pi_{\perp} \Big(  \omega\,\pa_{\vf}  + \cD_{0,\alpha}^{(\infty)}+ \pa_{\vartheta} \fv_{\varepsilon}(\vf,\vartheta) \Big) \Pi_{\perp} +\sE_{0}\,,\label{sL.perp}
	\end{align} 
	where the diagonal operator $\cD_{0,\alpha}^{(\infty)}$ is stated  in Proposition \ref{prop.diag.red.equi}, the real function $\fv_{\varepsilon}$ is even in $(\vf,\vartheta)$, satisfying the estimates 
	\begin{align}
		\| {\mathfrak v}_{\varepsilon} \|_{s}^{5,\upsilon} & \lesssim_{s} \varepsilon (1+ \|\underline{\fu}_{\perp} \|_s^{5,\upsilon}) \,,\label{fv.est} \\
		\| \Delta_{12}{\mathfrak v}_{\varepsilon} \|_{s_1} &\lesssim_{s} \varepsilon \| \Delta_{12}\fu_{\perp} \|_{s} \,, \label{fv.est12}
	\end{align}
	and the  operator $\sE_{0}$ is real and reversible, satisfying,   for any $h \in W_\upsilon^5(\Lambda;H^s_{\perp,\textnormal{even}}(\T^2))$ and for any $m=0,1,...,M$,
	\begin{align}
		\| \partial_\theta^{m}\sE_{0}[h] \|_{s}^{5,\upsilon} & \lesssim \varepsilon\big(\big(1+\| \underline{\fu}_{\perp}\|_{s_0+m}^{5,\upsilon}\big) \| h \|_{s}^{5,\upsilon} + \|\underline{\fu}_{\perp} \|_{s+m}^{5,\upsilon} \| h \|_{s_0}^{5,\upsilon} \big)\,,  \label{sE.smooth.est} \\
		\|\partial_\theta^{m} \Delta_{12}\sE_{0}[h] \|_{s}& \lesssim \varepsilon  \| h \|_{s}  \| \Delta_{12}\fu_{\perp} \|_{s+m}  \,.\label{sE.smooth.est12}
	\end{align}
\end{lem}
\begin{proof}
	Applying Proposition \ref{lemma.map.bG.solo}  yields to
    \begin{align}
		\sL_{\perp} & = \Pi_{\perp} \Big(  \omega\,\pa_{\vf} +\cM_{\alpha}^{-1} \pa_{\theta} \big( \tV_{\varepsilon,\alpha}  + \cR_{\varepsilon,\alpha} \big) \cM_{\alpha} \big( {\rm Id} + \di_{\fu_{\perp}} \fu_{\intercal}(\underline{\fu}_{\perp}) \big)  \Big) \Pi_{\perp} \\
		& = \Pi_{\perp} \Big(  \omega\,\pa_{\vf} +\cM_{\alpha}^{-1} \pa_{\theta} \big( \tV_{\varepsilon,\alpha}  + \cR_{\varepsilon,\alpha} \big) \cM_{\alpha}   \Big) \Pi_{\perp}\\
        &+\Pi_{\perp}   \cM_{\alpha}^{-1} \pa_{\theta} \big( \tV_{\varepsilon,\alpha}  + \cR_{\varepsilon,\alpha} \big) \cM_{\alpha}   \di_{\fu_{\perp}} \fu_{\intercal}(\underline{\fu}_{\perp})    \Pi_{\perp}\,. 
	\end{align}
	Thus, by Proposition~\ref{prop.Malpha}--$(ii)$ and straightforward computation we may expand
	\begin{align}
		\sL_{\perp} 
		& = \Pi_{\perp} \Big(  \omega\,\pa_{\vf} +\cM_{\alpha}^{-1} \pa_{\theta} \big( \tV_{0,\alpha}  + \cR_{0,\alpha} \big) \cM_{\alpha}   \Big) \Pi_{\perp} \\
		& + \Pi_{\perp} \Big(  \omega\,\pa_{\vf} + \cU^{-1} \cS^{-1} \pa_{\theta} \big( \tV_{\varepsilon,\alpha} - \tV_{0,\alpha}  + \cR_{\varepsilon,\alpha} -\cR_{0,\alpha}\big) \cS \cU\Pi_{\perp} +  \tR_{0} \\
		&  = \Pi_{\perp} \Big(  \omega\,\pa_{\vf}  + \cD_{0,\alpha}^{(\infty)}+ \pa_{\vartheta} \circ\,\fv_{\varepsilon}\Big) \Pi_{\perp} +\sE_{0}\,,
	\end{align}
	where
	\begin{align}
		\fv_{\varepsilon}(\vf,\vartheta) & := \cB^{-1} \{ \tV_{\varepsilon,\alpha}(\vf,\,\cdot\,) - \tV_{0,\alpha}(\,\cdot\,)\}(\vartheta) \,, \label{fv.def} 
	\end{align}
	and the remainder $\sE_{0}$ takes the form
	\begin{equation}\label{sE0.def}
		\sE_{0} := \Pi_{\perp} \cU^{-1} \cS^{-1} \pa_{\theta} \big(\cR_{\varepsilon,\alpha} - \cR_{0,\alpha} \big) \cS \cU \Pi_{\perp} + \tR_{0}\,, 
	\end{equation}
	with $\tR_{0}$ given by
	\begin{align}
		\tR_{0} & = \Pi_{\perp}  \cM_{\alpha}^{-1} \pa_{\theta} \big( \tV_{\varepsilon,\alpha}  + \cR_{\varepsilon,\alpha} \big) \cM_{\alpha} \di_{\fu_{\perp}} \fu_{\intercal}(\underline{\fu}_{\perp}) \Pi_{\perp} \\
		& + \Pi_{\perp}   \big( \cM_{\alpha}^{-1}  - \cU^{-1} \cS^{-1} \big)\pa_{\theta} \big( \tV_{\varepsilon,\alpha} - \tV_{0,\alpha}  + \cR_{\varepsilon,\alpha} -\cR_{0,\alpha}\big)\cM_{\alpha}^{-1}\Pi_{\perp} \\
		& + \Pi_{\perp} \cU^{-1} \cS^{-1} \pa_{\theta} \big( \tV_{\varepsilon,\alpha} - \tV_{0,\alpha}  + \cR_{\varepsilon,\alpha} -\cR_{0,\alpha}\big) \big(\cM_{\alpha}- \cS \cU\big)\Pi_{\perp}  \\
        & + \Pi_{\perp} \big( \pa_{\vartheta} \fv_{\varepsilon} (\cU-\Pi_{\rm ph}) + (\cU^{-1}-\Pi_{\rm ph})\pa_{\vartheta} \fv_{\varepsilon} \,\cU  \big)\Pi_{\perp} \,.
	\end{align}
	The estimate \eqref{fv.est}  follows by \eqref{fv.def}, Proposition \eqref{lemma.map.bG}-$(iii)$ and Lemma \ref{lemma.reparam.0}-$(ii)$. We note that $ \di_{\fu_{\perp}} \fu_{\intercal}(\underline{\fu}_{\perp})$ is a finite rank operator by \eqref{fu.tang}, \eqref{fX.def}, \eqref{action.angle} and Theorem \ref{theo.FP.tangential}. Therefore, together with Proposition \ref{prop.Malpha}-$(ii)$, we deduce that $\tR_{0}$ is a finite rank operator and, consequently by \eqref{sE0.def}, that $\sE_{0}$ is a smoothing operator. The estimate  \eqref{sE.smooth.est} then follows by Proposition \ref{prop.diag.red.equi}, Proposition \ref{prop.Malpha}, Proposition \ref{lemma.map.bG}-$(iii)$, Theorem \ref{theo.FP.tangential}-$(ii)$ and Lemma \ref{standard prop decay norm}-$(i),(iv)$. The estimates \eqref{fv.est12} and \eqref{sE.smooth.est12} follow by similar arguments, together with Lemma \ref{lem-productlaw}-$(v)$.
\end{proof}

\subsection{Reduction of the transport term}\label{sect.redu.eps}
Combining \eqref{lambda.linear}, \eqref{diag.0.alpha}, and \eqref{conj.no.proj}, we deduce that the operator $\sL_{\perp}$ defined in \eqref{sL.perp} can be decomposed as
\begin{align}\label{sL.perp.intro}
	\sL_{\perp} & =  \Pi_{\perp} \Big(  \omega\,\pa_{\vf}  + \cD_{0,\alpha}^{(\infty)}+\pa_{\vartheta} \fv_{\varepsilon}(\vf,\vartheta) \Big) \Pi_{\perp} +\sE_{0} \\ 
	& =  \Pi_{\perp} \Big(  \omega\,\pa_{\vf} + \pa_{\vartheta}\big( \tm_{1}(\alpha) + \fv_{\varepsilon}(\vf,\vartheta)  \big) + \pa_{\vartheta} \cR(0;0) + \cZ_{0,\alpha} \Big) \Pi_{\perp} +\sE_{0} \\
	& = : \Pi_{\perp} \big( \sL_{\rm Tr}+ \pa_{\vartheta} \cR(0;0) + \cZ_{0,\alpha} \big)  \Pi_{\perp} +\sE_{0} \,,
\end{align}
where the periodic transport operator is given by
\begin{equation}\label{transport.op}
	\sL_{\rm Tr} = \omega\,\pa_{\vf} + \pa_{\vartheta}\circ(\tm_{1}(\alpha) + \fv_{\varepsilon}(\vf,\vartheta)) \,.
\end{equation}

This operator encodes the main dynamical features of the system but still involves variable coefficients through the perturbation $\fv_{\varepsilon}$.
The first goal is to simplify its structure by means of a suitable change of variables. In particular, we aim to conjugate $\sL_{\rm Tr}$ to a constant-coefficient operator of the form $\omega\,\pa_{\vf} + \tm_{1,\varepsilon}\pa_{y}$. This reduction is the content of Lemma~\ref{almost.straight.lemma}.
Once this normalization is achieved, we exploit it to treat the full operator $\sL_{\perp}$. More precisely, we use the above transformation to conjugate $\sL_{\perp}$ to a diagonal operator acting on the normal directions, up to a smoothing remainder. This second step is carried out in Lemma~\ref{diffeo.conj.pert}.This scheme provides the fundamental ingredients needed to establish the almost invertibility of the operator $\sL_{\perp}$, which will be addressed in the next section.

\medskip

Our first main result concerns the conjugation of the transport operator $\sL_{\rm Tr}$ by a suitable symplectic transformation of the form
\begin{equation}
	\cS_{\varepsilon} := \bigl(1+\pa_{\vartheta}\beta_{\varepsilon}(\vf,\vartheta)\bigr)\circ \cB_{\varepsilon} \,,
\end{equation}
where the composition operator
\begin{equation}
	\cB_{\varepsilon}h(\vf,\vartheta) := h\bigl(\vf,\vartheta + \beta_{\varepsilon}(\vf,\vartheta)\bigr)
\end{equation}
is induced by a $\vf$-dependent diffeomorphism $y= \vartheta + \beta_{\varepsilon}(\vf,\vartheta)$ on the torus $\T_{\vartheta}$. Here, $$\beta_{\varepsilon}:\T_{\vf}\times\T_{\vartheta}\to\R
$$ is a small periodic function, assumed to be odd in $(\vf,\vartheta)$.
\\
We are now in a position to state the straightening result for the transport operator $\sL_{\rm Tr}$. This type of result is classical and has been implemented in various configurations with slight variations, see, for instance, \cite{FMT25,HR21}.

\begin{lem}\label{almost.straight.lemma}
	{\bf (Straightening of the transport operator).}
	There exists $\sigma=\sigma(\tau)\gg 1$  such that, for all $S>s_0+\sigma$, with $s_0$ as in \eqref{s0-embed}, there exists $\delta:=\delta(S,\tau) \in (0,1)$ such that, if \eqref{ansatz.approx} holds with $\mu_{0} \gg \sigma$ and if
	\begin{equation}
		\varepsilon \upsilon^{-1} \leq \delta\,,
	\end{equation}
	the following hold true:
	\\[1mm]
	\noindent $(i)$ There exist a map $(\omega,\alpha)\in \R\times [\alpha_{1},\alpha_{2}]\mapsto\tm_{1,\varepsilon}:= \tm_{1,\varepsilon}(\omega,\alpha):= \tm_{1,\varepsilon}(\omega,\alpha,\underline{\fu}_{\perp})\in\R$,  and a periodic odd function $(\vf,\vartheta)\in\mathbb{T}^2\mapsto\beta_{\varepsilon}(\vf,\vartheta)$, satisfying the estimates
	\begin{equation}
		|\tm_{1,\varepsilon}-\tm_{1}|^{5,\upsilon} \lesssim \varepsilon \,, \quad |\Delta_{12}\tm_{1,\varepsilon}| \lesssim \varepsilon \| \Delta_{12}\fu_{\perp} \|_{s_0+\sigma}\,, \label{tm.1e.est}
	\end{equation}
	and, for any $s_0 \leq s \leq S$ and any $s_1$ as in \eqref{s1},
	\begin{align}
		\| \beta_{\varepsilon}\|_{s}^{5,\upsilon}&  \lesssim_{S}\varepsilon\upsilon^{-1} \big( 1+ \| \underline{\fu}_{\perp} \|_{s+\sigma}^{5,\upsilon} \big) \,, \label{beta.eps.est} \\
		\| \Delta_{12} \beta_{\varepsilon}\|_{s_1} &  \lesssim_{s_1}\varepsilon\upsilon^{-1} \| \Delta_{12}\fu_{\perp} \|_{s_1+\sigma} \,; \label{beta.eps.est12}
	\end{align}
	\noindent $(ii)$ The diffeomorphism $\vartheta \mapsto \vartheta + \beta_{\varepsilon}(\vf,\vartheta)$ is invertible with inverse $y \mapsto y + \breve{\beta}_{\varepsilon}(\vf,y)$, where the function $\breve{\beta}_{\varepsilon}(\vf,y)$ satisfies bounds like in \eqref{beta.eps.est}, \eqref{beta.eps.est12};
	\\[1mm]
    \noindent $(iii)$ The operators $\cS_{\varepsilon}^{\pm 1}$ are reversibility preserving and satisfy,  for any $s_0\leq s\leq S-\sigma$ , any $s_1$ as in \eqref{s1} and for any $h\in H_{0}^s(\T^2)$,
	\begin{align}
		&	\| \cS_{\varepsilon}^{\pm 1}h \|_{s}^{5,\upsilon} \lesssim_{S} \| h \|_{s}^{5,\upsilon} +  \| \underline{\fu}_{\perp} \|_{s+\sigma}^{5,\upsilon}  \| h\|_{s_0}^{5,\upsilon} \,, \\
		&	\| (\cS_{\varepsilon}^{\pm 1}-{\rm Id}) h \|_{s}^{5,\upsilon} \lesssim_{S} \varepsilon\upsilon^{-1} \big( \|h \|_{s+1}^{5,\upsilon} + \| \underline{\fu}_{\perp} \|_{s+\sigma}^{5,\upsilon} \|h \|_{s_0+1}^{5,\upsilon}  \big) \,, \\
		& \| \Delta_{12} \cS_{\varepsilon}^{\pm 1} h \|_{s_1} \lesssim_{s_1} \varepsilon \upsilon^{-1}\| h \|_{s_1+1} \| \Delta_{12}\fu_{\perp} \|_{s_1+\sigma} \,;
	\end{align}
	\noindent $(iv)$ 
    For any $(\omega,\alpha)$ belonging to the set $ \tT\tC(\upsilon,\underline{\fu}_{\perp}) $ given by
	\begin{align}
		\tT\tC(\upsilon,\underline{\fu}_{\perp}) &:=  \bigcap_{(\ell,j)\in\Z^2 \setminus\{0\}} \big\{  (\omega,\alpha) \in \R \times [\alpha_{1},\alpha_{2}]  \,:\, |\omega\,\ell + \tm_{1,\varepsilon}\,j |\geq 4\upsilon|j|^{-\tau}
		\big\}\,,\label{nonres.transport}
	\end{align}
	the operator $\sL_{\rm Tr}$ in \eqref{transport.op} is conjugated to
	\begin{equation}
		\cS_{\varepsilon}^{-1} \sL_{\rm Tr} \cS_{\varepsilon} = \omega\,\pa_{\vf} + \tm_{1,\varepsilon} \pa_{y}\,.
	\end{equation}
	
\end{lem}
\begin{proof}
	The proof follows the same arguments presented in \cite[Lemma 9.8, Theorem C.2]{FMT25}, with an easy reduction from dimension 2 to the dimension 1 that we are considering here.  It is therefore omitted.
\end{proof}
The next aim is to conjugate the operator $\sL_{\perp}$ defined in \eqref{sL.perp}, which is localized on the normal modes. However, the transformation $\cS_{\varepsilon}$ acts on the full ambient space and is therefore not directly adapted to this setting. To make it suitable for our purposes, we restrict it to the normal directions by means of the projection $\Pi_{\perp}$.
More precisely, we introduce the localized transformation
\begin{align}\label{cS.eps.perp}
\cS_{\varepsilon,\perp}:= \Pi_{\perp} \cS_{\varepsilon}\, \Pi_{\perp}.
\end{align}
This operator captures the action of $\cS_{\varepsilon}$ on the normal component while preserving the underlying structure of the problem. 
In what follows, we collect several properties of the operator $\cS_{\varepsilon,\perp}$ that will be instrumental in the analysis.

	\begin{lem}\label{Lema-decompMon.eps}
	Let $M\in\N$ and let $S>s_0+\sigma$ be fixed as in Lemma $\ref{almost.straight.lemma},$  $s_0 \leq s \leq S-\sigma$, any $s_1$ as in \eqref{s1}, The following assertions hold:
    \\[1mm]
    \noindent $(i)$ The operator
	$\displaystyle{
		\cS_{\varepsilon,\perp} : H_{\perp}^{s}(\T^2) \to H_{\perp}^{s}(\T^2)}$
	is  invertible, reversibility preserving and  satisfies the following estimates:  for any $h \in H_{\perp}^{s}(\T^2)$,
	\begin{align*}
		&	\| \cS_{\varepsilon,\perp}^{\pm 1}h \|_{s}^{5,\upsilon} \lesssim_{S} \| h \|_{s}^{5,\upsilon} +  \| \underline{\fu}_{\perp} \|_{s+\sigma}^{5,\upsilon}  \| h\|_{s_0}^{5,\upsilon} \,,  \\
		\nonumber&	\| (\cS_{\varepsilon,\perp}^{\pm 1}-\Pi_{\perp}) h \|_{s}^{5,\upsilon} \lesssim_{S} \varepsilon\upsilon^{-1} \big(  \| h \|_{s+1}^{5,\upsilon} + \| \underline{\fu}_{\perp} \|_{s+\sigma}^{5,\upsilon}  \| h \|_{s_0+1}^{5,\upsilon} \big) \,,   \\
		\nonumber& \| \Delta_{12} \cS_{\varepsilon}^{\pm 1} h \|_{s_1} \lesssim_{s_1} \varepsilon \upsilon^{-1}  \| h \|_{s_1+1}\| \Delta_{12}\fu_{\perp} \|_{s_1+\sigma} \,; 
	\end{align*}
    \noindent $(ii)$
		We have the decompositions
        \begin{align}
            \cS_{\varepsilon,\perp}^{\pm1} = \cS_{\varepsilon}^{\pm1} + \tQ_{\pm 1}\,, 
        \end{align}
        where $\tQ_{\pm 1}$  are finite-rank operators satisfying, for any $m=0,1,...,M$ and  for any $h \in H_{\perp}^{s}(\T^2)$,
	\begin{align}
		&	\||{\rm D}|^{m} \tQ_{\pm1} h \|_{s}^{5,\upsilon} \lesssim_{S,m} \varepsilon\upsilon^{-1} \big(  \| h \|_{s}^{5,\upsilon} + \| \underline{\fu}_{\perp} \|_{s+\sigma+m}^{5,\upsilon}  \| h \|_{s_0}^{5,\upsilon}\big) \,,  \label{est.tQ} \\
	& \||{\rm D}|^{m} \Delta_{12}\tQ_{\pm1} h \|_{s_1} \lesssim_{s_1,m} \varepsilon \upsilon^{-1}  \| h \|_{s_1}\| \Delta_{12}\fu_{\perp} \|_{s_1+\sigma+m} \,. \label{est.tQ.12}
	\end{align}
\end{lem}

The proof of this lemma follows the same argument developed in  \cite[Lemma 6.3]{HR21}. We postpone its proof to Appendix \ref{app.diff.U}.

\medskip

We are now in a position to conjugate the full operator $\sL_{\perp}$ defined in \eqref{sL.perp}, see also \eqref{sL.perp.intro} and \eqref{transport.op}. To this end, we introduce the  operator
\begin{align}\label{L-perp-red}
\sL_{\perp,{\rm red}} := \cS_{\varepsilon,\perp}^{-1} \, \sL_{\perp} \, \cS_{\varepsilon,\perp}.
\end{align}
Here, the transformation $\cS_{\varepsilon,\perp}$ is defined in \eqref{cS.eps.perp}. 
By construction, the operator $\sL_{\perp,{\rm red}}$ is well-defined for all $(\omega,\alpha)\in\R\times[\alpha_{1},\alpha_{2}]$. This conjugation allows us to reduce $\sL_{\perp}$ to a diagonal Fourier multiplier, up to a smoothing remainder, as stated below.
\begin{lem}\label{diffeo.conj.pert}
	{\bf (Reduction of the transport term - perturbed)}.
	Under the assumptions of Lemma $\ref{Lema-decompMon.eps},$ 
    the operator $\displaystyle{
			\sL_{\perp,{\rm red}}  : H_{\perp,{\rm even}}^{s}(\T^2) \to H_{\perp,{\rm odd}}^{s-1}(\T^2)}
			\,,
		$ is  continuous and reversible. Moreover, for any $(\omega,\alpha)\in \tT\tC_{\infty}(\upsilon,\tau)$, defined by  \eqref{nonres.transport}, we have 
    \begin{align*}
    \sL_{\perp,{\rm red}}&  = \omega\,\pa_{\vf} + \sD_{\perp,{\rm red}}  + \sE_{\perp,{\rm red}}		\,,
    \end{align*}
    where the diagonal operator $\sD_{\perp,{\rm red}}$ is given by
	\begin{equation}\label{diag.D.perp}
		\begin{aligned}
			&\sD_{\perp,{\rm red}}  = \Pi_{\perp}\big( \tm_{1,\varepsilon} \pa_{y} +\pa_{y}\cR(0;0)+ \cZ_{0,\alpha} \big)\Pi_{\perp}= {\rm diag}_{j\in\Z_{\perp}} \big\{  \mu_{j,\varepsilon}(\omega,\alpha) \big\}  \,,  \\
			&	\begin{aligned}
				\mu_{j,\varepsilon}(\omega,\alpha)& := \im \big(  \tm_{1,\varepsilon}\, j - \tfrac{1}{2} {\rm sgn}(j) \big) + \tz_{j}(\alpha) \,, \\
				\mu_{-j,\varepsilon}(\omega,\alpha)&= -\overline{\mu_{j,\varepsilon}(\omega,\alpha)} \,,
			\end{aligned} \quad \ j\in\Z_{\perp} \,,
		\end{aligned}
	\end{equation}
	and the constant $\tm_{1,\varepsilon}$ is constructed in Lemma $\ref{almost.straight.lemma}$.
	In addition,  the smoothing and reversible operator $\sE_{\perp,{\rm red}}$ 
	satisfies the estimates, , for any $h\in H_{\perp.{\rm even}}^s(\T^2)$ and $m=0,1,...,M$,
	\begin{align}
		\| |{\rm D}|^{m}\sE_{\perp,{\rm red}}[h]  \|_{s}^{5,\upsilon} & \lesssim_{S,m} \varepsilon\upsilon^{-1} \big( \| h\|_{s}^{5,\upsilon}  + \| \underline{\fu}_{\perp}\|_{s+\sigma+m}^{5,\upsilon} \|  h \|_{s_0}^{5,\upsilon}  \big)  \,, \label{sE.perp.est} \\
		\| |{\rm D}|^{m}\Delta_{12}\sE_{\perp,{\rm red}}[h] \|_{s_1}& \lesssim_{s_1,m} \varepsilon \upsilon^{-1}  \| h \|_{s_1}  \| \Delta_{12}\fu_{\perp} \|_{s_1+\sigma}  \,.
	\end{align}
\end{lem}
\begin{proof}
The continuity and the reversibility  of the operator $\displaystyle{
			\sL_{\perp,{\rm red}}  : H_{\perp,{\rm even}}^{s}(\T^2) \to H_{\perp,{\rm odd}}^{s-1}(\T^2)}
			\,,
		$
follow directly from the corresponding properties of $\sL_{\perp}$ in Lemma \ref{expand.sL.pert}, together with the mapping properties of the transformation $\cS_{\varepsilon,\perp}^{\pm 1}$ established in Lemma \ref{Lema-decompMon.eps}-$(i)$.
    We now proceed to derive an explicit expansion of the operator $\sL_{\perp,{\rm red}}$. Using  \eqref{sL.perp.intro} and Lemma \ref{Lema-decompMon.eps}-$(i)$, we obtain
    \begin{align}
        \sL_{\perp,{\rm red}} & = \cS_{\varepsilon,\perp}^{-1} \sL_{\perp} \cS_{\varepsilon,\perp}  = \cS_{\varepsilon,\perp}^{-1} \LL_{\perp} \cS_{\varepsilon,\perp} + \cS_{\varepsilon,\perp}^{-1} \sE_{0} \cS_{\varepsilon,\perp}\,, \label{conj.eps1}
    \end{align}
    where we defined
    \begin{equation}
        \LL_{\perp}:= \sL_{\rm Tr} + \pa_{\vartheta} \cR(0;0) + \cZ_{0,\alpha} \,.
    \end{equation}
    By Lemma \ref{Lema-decompMon.eps}-$(ii)$, we deduce  that
    \begin{align}
        \cS_{\varepsilon,\perp}^{-1} \LL_{\perp} \cS_{\varepsilon,\perp} & = \big(\cS_{\varepsilon}^{-1} + \tQ_{-1} \big) \LL_{\varepsilon} \big( \cS_{\varepsilon}+ \tQ_{1} \big) \\
        & = \cS_{\varepsilon}^{-1} \LL_{\varepsilon} \cS_{\varepsilon} + \cS_{\varepsilon}^{-1} \LL_{\varepsilon} \tQ_{1} + \tQ_{-1} \LL_{\varepsilon} \cS_{\varepsilon,\perp} \,. \label{conj.eps2}
    \end{align}
    Consequently, an application of   Lemma \ref{almost.straight.lemma} yields
    \begin{align}
        \cS_{\varepsilon}^{-1} \LL_{\varepsilon} \cS_{\varepsilon}  & = \omega\,\pa_{\vf} + \pa_{y}(\tm_{1,\varepsilon,\alpha}+\cR(0;0)) + \cZ_{0,\alpha} \\
        & + \cS_{\varepsilon}^{-1}\big( \pa_{\vartheta} \cR(0;0) + \cZ_{0,\alpha} \big) \cS_{\varepsilon} - \pa_{y} \cR(0;0) - \cZ_{0,\alpha} \,.\label{conj.eps3}
    \end{align}
    Plugging \eqref{conj.eps2}, \eqref{conj.eps3} into \eqref{conj.eps1}, we obtain that $\sL_{\perp,{\rm red}}$ has the suitable form as indicated in the lemma, with $\sD_{\perp,{\rm red}}$ as in \eqref{diag.D.perp} and $\sE_{\perp,{\rm red}}$ explicitly given by
    \begin{align}
        \sE_{\perp,{\rm red}}  := &\Pi_{\perp}\big( \cS_{\varepsilon}^{-1}\big( \pa_{\vartheta} \cR(0;0) + \cZ_{0,\alpha} \big) \cS_{\varepsilon} - \pa_{y} \cR(0;0) - \cZ_{0,\alpha} \big) \Pi_{\perp} \\
        & + \Pi_{\perp} \big( \cS_{\varepsilon}^{-1} \LL_{\varepsilon} \tQ_{1} + \tQ_{-1} \LL_{\varepsilon} \cS_{\varepsilon,\perp} \big) \Pi_{\perp}  + \cS_{\varepsilon,\perp}^{-1} \sE_{0} \cS_{\varepsilon,\perp}\,. \label{sE.op.def}
    \end{align}
    It remains to prove the estimates \eqref{sE.perp.est}.
    First, we claim that the operators
    \begin{equation}\label{sE1sE2}
        \sE_{1}:=\Pi_{\perp}\big( \cS_{\varepsilon}^{-1} \pa_{\vartheta} \cR(0;0) \, \cS_{\varepsilon} - \pa_{y} \cR(0;0) \big) \Pi_{\perp} \,, \quad \sE_{2}:=\Pi_{\perp}\big( \cS_{\varepsilon}^{-1}\cZ_{0,\alpha} \,\cS_{\varepsilon} - \cZ_{0,\alpha} \big) \Pi_{\perp} \,,
    \end{equation}
    are smoothing, with estimates as in \eqref{sE.perp.est}. By Lemma \ref{lemma.cR.00}, we find  
    $$\pa_{\vartheta}\cR(0;0)[e^{\im j \vartheta}]= - \tfrac{\im}{2}{\rm sgn}(j)\, e^{\im j \vartheta}\,.
    $$ 
    By adapting the argument in \cite[Lemma 2.36]{BM20}, one has that 
    \begin{equation}
        \big( \cS_{\varepsilon}^{-1} \pa_{\vartheta} \cR(0;0) \, \cS_{\varepsilon} - \pa_{y} \cR(0;0)  \big)h(\vf,y) = \int_{\T} K_{\varepsilon}(\vf,y,\eta;\omega,\alpha) h(\vf,\eta) \wrt\eta \,,\label{K.Hilbert}
    \end{equation}
    where the kernel $K_{\varepsilon}$ is smooth and satisfies  the estimate
    \begin{equation}
        \| K_{\varepsilon} \|_{s}^{5,\upsilon} \lesssim_{s} \| \beta_{\varepsilon} \|_{s+1}^{5,\upsilon}\,.\label{K.Hilbert.est}
    \end{equation}
    Applying Lemma \ref{lem-productlaw}-$(ii)$, together with \eqref{K.Hilbert.est}, \eqref{beta.eps.est} in Lemma \eqref{almost.straight.lemma}-$(i)$ and Sobolev embedding, we obtain  
    \begin{align}
        \||{\rm D}|^{m}\sE_{1} h \|_{s} & \lesssim_{s} \| |{\rm D}|^{m}K_{\varepsilon} \|_{s_0}^{5,\upsilon} \| h\|_{s}^{5,\upsilon} + \| |{\rm D}|^{m}K_{\varepsilon} \|_{s}^{5,\upsilon} \| h\|_{s_0}^{5,\upsilon} \\
        & \lesssim_{s,m} \varepsilon\upsilon^{-1} \big( \| h\|_{s}^{5,\upsilon}  + \| \underline{\fu}_{\perp}\|_{s+\sigma+m}^{5,\upsilon} \| h\|_{s_0}^{5,\upsilon} \big) \,, \label{sE1.est}
    \end{align}
    for some $\sigma=\sigma(\tau)>0$. A similar argument applies to $\sE_{2}$ in \eqref{sE1sE2}. Indeed, by Proposition~\ref{prop.diag.red.equi}, the operator $\cZ_{0,\alpha}$ is a smoothing diagonal operator, and can therefore be represented as a convolution kernel operator, namely,
    \begin{equation}
        \cZ_{0,\alpha}h(\vf,\vartheta) = \int_{\T} h(\vf,\eta) K_{\cZ_{0,\alpha}}(\vartheta-\eta) \wrt \eta\,, \quad K_{\cZ_{0,\alpha}}(\eta) = \frac{\tz_{1}(\alpha)}{\pi} \cos(\theta)+ \frac{1}{2\pi}\sum_{|j|\geq 2} \tz_{j}(\alpha)e^{\im j \eta} \,,
    \end{equation}
    with $K_{\cZ_{0,\alpha}}$ being a smooth kernel by Lemma \ref{standard prop decay norm}-$(v)$ and \eqref{Y.gen.est} in Proposition \ref{prop.diag.red.equi}. Therefore, we obtain the desired estimate as in \eqref{sE1.est}, which concludes the proof of the claim.\\
    We now deduce the first estimate in \eqref{sE.perp.est} from \eqref{sE.op.def}, \eqref{sE1sE2}, \eqref{sE1.est}, Lemma~\ref{Lema-decompMon.eps}, and Lemma~\ref{expand.sL.pert}. The last estimate in \eqref{sE.perp.est} follows by similar arguments and  therefore its proof is omitted.
\end{proof}

\subsection{Invertibility of the linearized operator}\label{sect.almost.inv}
In this section, we address the invertibility properties of the linearized operator arising in the range equation. Building upon the reduction scheme developed in the previous sections, we exploit the almost diagonal structure obtained after conjugation to construct a full inverse with suitable tame estimates. 
The main difficulty lies in handling the small divisor effects generated by the transport dynamics. By combining the diagonalization procedure with smoothing properties of the remainders, we show that the operator is invertible on the normal directions. This result constitutes a key step in the implementation of the Nash–Moser scheme.

\smallskip

We first recall from Lemma \ref{diffeo.conj.pert}, \eqref{sL.perp} and \eqref{L-perp-red} that
	\begin{align}
		\sL_{\perp} & = \cS_{\varepsilon,\perp} \sL_{\perp,{\rm red}} \cS_{\varepsilon,\perp}^{-1} \\
		& = \cS_{\varepsilon,\perp}\big( \omega\,\pa_{\vf} + \sD_{\perp,{\rm red}} + \sE_{\perp,{\rm red}} \big)\cS_{\varepsilon,\perp}^{-1} \,.\label{sL.to.invert.full}
	\end{align}
We are now in a position to invert the operator $\sL_{\perp}$ and derive quantitative estimates for its inverse. When later applied in Section \ref{section.NASH}, its action will be truncated up to frequencies $|j|\leqslant N_{\tn}$, where $N_{\tn}\in\N$ denotes the cut-off at the $\tn$-th step of the Nash--Moser iteration:
this truncation will allow us to control the loss of derivatives at high frequencies and to construct approximated solutions in the nonlinear iteration. Having fully reduced the transport operator in Section \ref{sect.redu.eps}, the small divisors will depend on the $\tn$-th step of the iteration only through the approximate solution at which we are linearizing. This fact is key when imposing non-resonance conditions in Section \ref{section-measures}.

Recalling Lemma \ref{diffeo.conj.pert}, we write
\begin{equation}\label{split.op.red.full}
	\begin{aligned}
		\sL_{\perp,{\rm red}} & = \tL  + \sE_{\perp,{\rm red}}\,,\quad\textnormal{where} \quad   \tL := \omega\cdot\pa_{\vf} + \sD_{\perp,{\rm red}} \,. 
	\end{aligned}
\end{equation}
We define the set
\begin{align}\label{Cantor-n-2.full}
\t\Lambda(\upsilon,\underline{\fu}_{\perp}) &:=\bigcap_{(\ell,j)\in\Z \times \Z_{\perp} }\Big\{(\omega,\alpha)\in \R\times[\alpha_{1},\alpha_{2}]  \,: \,  | \omega\,\ell - \im \, \mu_{j,\varepsilon}(\omega,\alpha) | \geq \tfrac{2\upsilon}{ |j|^{\tau}}  \Big\}	 \,.
\end{align}
We first establish the invertibility of the operator $\sL_{\perp,{\rm red}}$ in \eqref{split.op.red.full},  and then recover the corresponding result for the original operator in \eqref{sL.to.invert.full} by conjugation.
The first result reads as follows.
\begin{pro}\label{prop-perp.full}
  There exists $M=M(\tau)\in \N$ such that the following assertions hold:
\\[1mm]
\noindent $(i)$
  There exists a real and reversible operator $\mathtt{T}$,  with  estimate
\begin{equation}\label{Est-T.full}
    \forall \,m =0,1,...,M\,, \ \forall \, s\in\R\,,\quad \||{\rm D}|^m\mathtt{T}h\|_{s}^{5,\upsilon}\lesssim_{s,m} \upsilon^{-1}\||{\rm D}|^{m+6\tau}h\|_{s}^{5,\upsilon}\,,
\end{equation}
such that,  for any $(\omega,\alpha)\in \t\Lambda(\upsilon,\underline{\fu}_{\perp})$ as in \eqref{Cantor-n-2.full},
one has
\begin{equation}\label{LT}
    \mathtt{L}\mathtt{T}={\rm Id}\,;
\end{equation}
\noindent $(ii)$ Assume \eqref{ansatz.approx}. For $\delta>0$ small enough and under the assumption
\begin{align}\label{small-second.full}
\varepsilon\upsilon^{-2}\leqslant \delta\,,
\end{align}
there exists a real and reversible operator $\mathcal{T}_{\perp,\textnormal{red}} $ such that, for  any
$(\omega,\alpha)\in\t\Lambda(\upsilon,\tau),$ one gets 
\begin{align} \label{inv.red.full}
		\sL_{\perp,{\rm red}} \mathcal{T}_{\perp,{\rm red}} & ={\rm Id}\,,
        \end{align}  
        satisfying the estimates
        \begin{align}
      \forall\, s \in [s_0,S-\sigma] \,, \quad  \| \mathcal{T}_{\perp,{\rm red}}[h]  \|_{s}^{5,\upsilon} 
      &\lesssim_{s}  \upsilon^{-1}\Big(\|h\|_{s+6\tau}^{5,\upsilon}+ \| \underline{\fu}_{\perp}\|_{s+\sigma}^{5,\upsilon} \|  h \|_{s_0+6\tau}^{5,\upsilon}\Big)\,. \label{ARI.red.est.full} 
	\end{align}
\end{pro}
\begin{proof}
We start with the proof of item $(i)$.
By the definition of $\tL$ in \eqref{split.op.red.full} and Lemma \ref{diffeo.conj.pert}, we have
\begin{equation}
 \forall\, (\ell,j)\in \Z \times \Z_{\perp}\,, \quad  \mathtt{L}{\bf e}_{\ell,j}=
 \im\big(\omega\,\ell-\ii\mu_{j,\varepsilon}\big) {\bf e}_{\ell,j} \,.
\end{equation}
We define the diagonal  operator  $\mathtt{T}$ by 
\begin{equation}\label{T.in.Fourier}
h=\sum_{(\ell,j)\in \Z \times \Z_{\perp}}h_{\ell,j}{\bf e}_{\ell,j},\quad\mathtt{T} h(\varphi,\theta) :=
-\im \sum_{ (\ell,j)\in \Z \times \Z_{\perp}}\frac{\chi\big((\omega\,\ell-\ii\mu_{j,\varepsilon})\upsilon^{-1}| j|^{\tau}\big)}{\omega\,\ell-\ii\mu_{j,\varepsilon}}h_{\ell,j}\,{\bf e}_{\ell,j}(\varphi,\theta)\,,
\end{equation}
where $\chi\in\mathscr{C}^\infty(\mathbb{R},[0,1])$ is an even positive cut-off function  such that 
\begin{equation}\label{properties cut-off function first reduction} 
   \chi(\xi) = \begin{cases}
        0 & \textnormal{if } \ |\xi|\leq \frac13 \,, \\
        1 & \textnormal{if } \ |\xi|\geq \frac12 \,.
    \end{cases}
\end{equation}
The operator $\tT$ in \eqref{T.in.Fourier} is real and reversible by \eqref{diag.D.perp} in Lemma \ref{diffeo.conj.pert}, Definition \ref{defi.rev} and \eqref{real.h}-\eqref{odd.h}.
Thus, for $(\omega,\alpha)\in\t\Lambda(\upsilon,\underline{\fu}_{\perp})$, one deduce \eqref{LT}.
Similarly to the proof of \cite[Theorem 9.1-(i)]{HHM21}, one   checks the estimate \eqref{Est-T.full}  directly from \eqref{T.in.Fourier}, recalling also the definition of the weighted topology in  \eqref{norm.weighted.def}.
This concludes the proof of item $(i)$.
\\[1mm]
\noindent We now prove item $(ii)$.
We define the operator
\begin{align}\label{T-red1}
\mathcal{T}_{\perp,{\rm red}}:=\big({\rm Id}+\mathtt{T}  \sE_{\perp,{\rm red}}\big)^{-1}\mathtt{T}\,.
\end{align}
To check the invertibility of ${\rm Id}+\mathtt{T}  \sE_{\perp,{\rm red}} $, we combine \eqref{Est-T.full} with $m=0$ and \eqref{sE.perp.est} from Lemma \ref{diffeo.conj.pert} with $m=6\tau$ to get
\begin{align*}
		\| \mathtt{T}   \sE_{\perp,{\rm red}} [h]  \|_{s}^{5,\upsilon} & \lesssim  \upsilon^{-1}\||{\rm D}|^{6\tau}  \sE_{\perp,{\rm red}}  [h]\|_{s}^{5,\upsilon}\\
        &\lesssim \varepsilon\upsilon^{-2} \big( \| h\|_{s}^{5,\upsilon}  + \| \underline{\fu}_{\perp}\|_{s+\sigma+6\tau}^{5,\upsilon} \|  h \|_{s_0}^{5,\upsilon}  \big)  \,.  
	\end{align*}
    Therefore, recalling the ansatz \eqref{ansatz.approx},
    \begin{align*}
		\| \mathtt{T}  \sE_{\perp,{\rm red}} [h]  \|_{s_0}^{5,\upsilon} 
        &\lesssim \varepsilon\upsilon^{-2}  \| h\|_{s_0}^{5,\upsilon}\,.
	\end{align*}
    Consequently, in view also of  the smallness condition \eqref{small-second.full}, we deduce, by a Neumann series argument, that the operator
\begin{equation}
    {\rm Id}+\mathtt{T}  \sE_{\perp,{\rm red}}  :W_\upsilon^5(\Lambda;H^s_{\perp,{\rm even}}(\T^2))\to {W_\upsilon^5(\Lambda;H^s_{\perp,{\rm even}}(\T^2))}
\end{equation}
is invertible, with
    \begin{align*}
		\| ({\rm Id}+\mathtt{T}  \sE_{\perp,{\rm red}})^{-1}[h]  \|_{s_0}^{5,\upsilon} 
        &\lesssim   \| h\|_{s_0}^{5,\upsilon}\,.
	\end{align*}
By a direct analysis based on a Neumann series argument together with \eqref{ansatz.approx}, \eqref{small-second.full}, we also obtain, for any $s\in [s_0,S-\sigma]$,
\begin{align}\label{neumann.tame.s.full}
		\| ({\rm Id}+\mathtt{T}  \sE_{\perp,{\rm red}})^{-1}[h]  \|_{s}^{5,\upsilon} 
        &\lesssim_{s}   \| h\|_{s}^{5,\upsilon} + \| \underline{\fu}_{\perp}\|_{s+\sigma+6\tau}^{5,\upsilon} \|  h \|_{s_0}^{5,\upsilon}\,.
	\end{align}
    Therefore, using \eqref{T-red1}, \eqref{Est-T.full} and \eqref{neumann.tame.s.full}, we conclude that 
    \begin{align}\label{T-EqP}
		\| \mathcal{T}_{\perp,{\rm red}}[h]  \|_{s}^{5,\upsilon} 
        &\lesssim_{s}   \| \mathtt{T}[h]\|_{s}^{5,\upsilon} + \| \underline{\fu}_{\perp}\|_{s+\sigma+6\tau}^{5,\upsilon} \|  \mathtt{T}[h] \|_{s_0}^{5,\upsilon}\\
      &\lesssim_{s} \upsilon^{-1}\Big(\|h\|_{s+6\tau}^{5,\upsilon}+ \| \underline{\fu}_{\perp}\|_{s+\sigma+6\tau}^{5,\upsilon} \|  h \|_{s_0+6\tau}^{5,\upsilon}\Big)\,,
	\end{align}
    which proves \eqref{ARI.red.est.full}.
    On the other hand, by \eqref{split.op.red.full}, \eqref{T-red1}, and \eqref{LT}, we infer that, for any $(\omega,\alpha)\in \t\Lambda(\upsilon,\underline{\fu}_{\perp})$,
    \begin{align}
        \sL_{\perp,{\rm red}}\mathcal{T}_{\perp,{\rm red}}&=\big(\mathtt{L}+  \sE_{\perp,{\rm red}}  \big)\big({\rm Id}+\mathtt{T}\sE_{\perp,{\rm red}}  \big)^{-1}\mathtt{T}\\
        &=\big(\mathtt{L}+\mathtt{L}\mathtt{T}\sE_{\perp,{\rm red}}  \big)\big({\rm Id}+\mathtt{T} \sE_{\perp,{\rm red}}  \big)^{-1}\mathtt{T} \\
        &=\mathtt{L}\mathtt{T} = {\rm Id} \,,
    \end{align}
which proves \eqref{inv.red.full}
       Finally, the operator $\cT_{\perp,{\rm red}}$ in \eqref{T-red1} is real and reversible by the fact that both $\tT$ and $\sE_{\perp,{\rm red}}$ are real and reversible (respectively by item $(i)$ and Lemma \ref{diffeo.conj.pert}). This concludes the proof of item $(ii)$ and the proof of the proposition.
\end{proof}
The goal of this final part is to establish the existence of the right inverse  for the operator $\sL_{\perp}$, building upon Proposition~\ref{prop-perp.full}. 
To this end, we first introduce a suitable Cantor-like set of parameters on which the required non-resonance conditions are satisfied:
\begin{align}\label{nonres.set.inv.full}
          \b\Lambda(\upsilon,\underline{\fu}_{\perp}) 
         := \tT\tC(\upsilon,\underline{\fu}_{\perp}) \cap \t\Lambda(\upsilon,\underline{\fu}_{\perp})\,,
    \end{align}
    where $\tT\tC(\upsilon,\underline{\fu}_{\perp})$ is introduced in \eqref{nonres.transport}, and $\t\Lambda (\upsilon,\underline{\fu}_{\perp})$ in \eqref{Cantor-n-2.full}. 
\begin{pro}\label{prop.full.inv}
    Assume \eqref{ansatz.approx} and \eqref{small-second.full}. 
There exist $\bar\sigma:=\bar\sigma(\tau)>0$ and a real and reversible operator $\mathcal{T}$ such that, for any
$(\omega,\alpha)\in\b\Lambda(\upsilon,\underline{\fu}_{\perp}),$ one has 
\begin{align}\label{full.inv.id}
		\sL_{\perp} \mathcal{T} & ={\rm Id}\,,
        \end{align} 
       where the following estimate hold,
        \begin{align}	
      \forall\, s\in [s_0,S-\bar\sigma] \,, \quad   \| \mathcal{T}[h]  \|_{s}^{5,\upsilon} 
      &\lesssim_{s}  \upsilon^{-1}\Big(\|h\|_{s+\bar\sigma}^{5,\upsilon}+ \| \underline{\fu}_{\perp}\|_{s+\bar\sigma}^{5,\upsilon} \|  h \|_{s_0+\sigma}^{5,\upsilon}\Big)\,.\label{cT.est.full}
	\end{align}
\end{pro}
\begin{proof}
  The relation between $\sL_{\perp}$ and $\sL_{\perp,{\rm red}}$ is given by \eqref{sL.to.invert.full}, namely
   \begin{align}
\cS_{\varepsilon,\perp}^{-1}\sL_{\perp} \cS_{\varepsilon,\perp}=\sL_{\perp,{\rm red}}\,.
        \end{align}
        In view of Proposition~\ref{prop-perp.full}-$(ii)$, it follows that 
        \begin{align}
		\cS_{\varepsilon,\perp}^{-1}\sL_{\perp} \cS_{\varepsilon,\perp}\mathcal{T}_{\tn,{\rm red}}=\sL_{\perp,{\rm red}}\mathcal{T}_{\perp,{\rm red}} =  {\rm Id}\,.
        \end{align}
        Conjugating back, we obtain the identity
        \begin{align}
		\sL_{\perp} \cS_{\varepsilon,\perp}\mathcal{T}_{\perp,{\rm red}}\cS_{\varepsilon,\perp}^{-1}
        & = \cS_{\varepsilon,\perp} \cS_{\varepsilon,\perp}^{-1} = {\rm Id}\,,
        \end{align}
        which proves \eqref{full.inv.id}  with $\cT$ defined as
        \begin{equation}\label{ARI.cT.n}
            \mathcal{T}:= \cS_{\varepsilon,\perp}\mathcal{T}_{\perp,{\rm red}}\cS_{\varepsilon,\perp}^{-1}\,.
        \end{equation}
        The corresponding estimate \eqref{cT.est.full} for $\mathcal{T}$  and the fact that $\cT$ is real and reversible follow from the properties  of  $\mathcal{T}_{\perp,{\rm red}}$ given in Proposition~\ref{prop-perp.full}-$(ii)$, together with those of $\cS_{\varepsilon,\perp}^{\pm1}$ established in Lemma~\ref{Lema-decompMon.eps}.  This ends the proof.
\end{proof}

\subsection{Nash-Moser iteration}\label{section.NASH}
In this section, we finally implement the Nash--Moser iteration scheme in order to solve the nonlinear range equation \eqref{range.eq} obtained after the Lyapunov--Schmidt reduction. 
The strategy developed here follows a rather classical Nash--Moser approach in the spirit of previous works on quasilinear Hamiltonian PDEs and KAM theory, see for instance  \cite{Berti-Bolle}. Starting from a sufficiently accurate approximate solution, we construct inductively a sequence of corrections obtained by solving suitable linearized equations. The crucial ingredient is the right inverse established in Section \ref{sect.almost.inv} through the reducibility and normal form analysis of Section \ref{sect.redu.eps}.
 The iterative estimates are carefully arranged so as to guarantee the rapid convergence of the approximate solutions toward an exact solution of the range equation. Moreover, the reversibility and symmetry properties of the problem are preserved throughout the scheme.
The outcome of the iteration is the existence of solutions for parameters satisfying suitable Melnikov-type non-resonance conditions. Combined with the measure estimates established in Section~\ref{section-measures}, this yields the resolution of the range equation on a Cantor-like set of asymptotically full measure.

\medskip

To formulate our result, we first introduce some preliminary material. In particular, we consider the following finite-dimensional subspaces
\begin{equation}
    \cH_{\tn} := \big\{ h \in L_{\perp}^2(\T^2) \, : \, h = \Pi_{\tn} h \big\}\,, \quad \tn \in \N_{0} \,,
\end{equation}
where $\Pi_{\tn}:= \Pi_{N_{\tn}} $ is the projection
\begin{equation}
    h(\vf,\vartheta) = \sum_{j \in \Z_{\perp}} h_j(\vf) e^{\im j\vartheta} \mapsto (\Pi_{\tn}h)(\vf,\vartheta) := \sum_{j\in \Z_{\perp} \atop |j|\leq N_{\tn}} h_{j}(\vf) e^{\im j\vartheta} \,.
\end{equation}
We also define $\Pi_{\tn}^\perp:= \Pi_{\perp} - \Pi_{\tn}$. The projections $\Pi_{\tn}$ and $\Pi_{\tn}^\perp$ satisfy the usual smoothing estimates
\begin{equation}\label{std.smoothing}
  \forall\,s,b\geq 0 \,, \quad   \| \Pi_{\tn}h \|_{s+b}^{5,\upsilon} \leq N_{\tn}^{b} \| h \|_{s}^{5,\upsilon} \,, \quad   \| \Pi_{\tn}^\perp h \|_{s}^{5,\upsilon} \leq N_{\tn}^{-b} \| h \|_{s+b}^{5,\upsilon} \,.
\end{equation}
For given parameters $\tau, N_{0}>1$, we introduce  the following quantities 
\begin{equation} \label{constants.NASH}
\begin{aligned}
& N_{\tn} := N_{0}^{\chi^\tn}, \quad \tn \geq 0, \quad N_{- 1} := 1\,, \quad \chi = 3/2, \quad  \kappa := 3\overline{\sigma} +12(\tau+1) \,, \\
&  \mathtt a_1 := {\rm max}\{ 6 \bar \sigma + 13\,,\, \chi^2(\tau + \tau^2 + 1) + \chi(2 \bar \sigma + 1) + 1 \} \,, \\
& \bar \tau := 2 \bar \sigma + 4 + \mathtt a_1 + \chi(\tau + \tau^2 + 1), \quad {\mathtt b}_1 := 2 \bar \sigma + 4  +\chi^{-1} ( \mathtt a_1 + \kappa)\,, \\
 & S:= s_0 + \bar\sigma + \tb_{1} \,.
\end{aligned}
\end{equation}
where $\bar \sigma > 1$ is given in Proposition \ref{prop.full.inv}. 

 \begin{rem}\label{remark.on.param.NM}
Let us now describe the parameters introduced in \eqref{constants.NASH} and explain their role in Proposition~\ref{NM.iteration}. 
The parameter $\bar\sigma>0$ measures the loss of derivatives in the approximate inverse of the linearized operator $\sL_{\perp}$ established in Proposition~\ref{prop.full.inv}. The parameter $\ta_{1}>0$ appears in the negative exponent of the estimate \eqref{stima.F.to.0} and quantifies the rate at which the low regularity norm at level $s_0$ of the nonlinear functional decreases along the Nash--Moser iteration. It also contributes in  the estimates \eqref{hn} governing the rapidly decaying corrections of the approximate solutions $\wt\fu_{\perp}^{(\tn)}$ in the regularity level $s_0+\bar\sigma$. 
As for the parameter $\tb_{1}>0$, it appears in the higher regularity index $s_0+\tb_{1}$ involved in the estimate \eqref{stima.alta.NM}. The corresponding high Sobolev norm is allowed to grow along the iteration, with growth rate controlled by the exponent $\kappa>0$.
 Finally, the parameter $\bar\tau>0$ enters the smallness condition \eqref{nash moser smallness condition}, which is crucial in order to close the inductive estimates and ensure the convergence of the Nash--Moser scheme.
\end{rem}
The tangential component $\fu_{\intercal}(\vf,\vartheta,0)$ introduced in \eqref{fu.tang.sol.NM} and constructed in Theorem~\ref{theo.FP.tangential} provides an approximate solution to the range equation $\wt\bG_{\perp}(\fu_{\perp})=0$ in \eqref{range.eq}. It therefore serves as a suitable initial approximation for the nonlinear Nash--Moser iteration leading later to the construction of a full solution. More precisely, we have the following result.
\begin{lem}\label{lemma.start.NM}
	{\bf (Initialization of the Nash--Moser iteration).}
	For any $s \geq 0$, there exists a constant $C(s)>0$ such that
	\begin{align}
		\|\widetilde\bG_{\perp}(0) \|_s^{5,\upsilon} \leq \varepsilon C(s) \big(\| \fu_{\intercal}(0)\|_{s}^{5,\upsilon}\big)^2\,.\label{approx.solution.normal}
		\end{align}
\end{lem}
\begin{proof}
    In view of \eqref{range.eq} and \eqref{LyaSch}, we deduce that
    \begin{align}
        \widetilde{\bG}_{\perp} (\fu_{\perp})  & = \bG_{\perp}\big( \fu_{\intercal}(\fu_{\perp}),\fu_{\perp} \big)  = \Pi_{\perp}\bG\big( \fu_{\intercal}(\fu_{\perp})+\fu_{\perp}  \big)  \,. 
    \end{align}
    Moreover, by \eqref{bG.equizero}, Proposition \ref{lemma.map.bG}-$(ii)$, Proposition \ref{prop.Malpha} and the fact that $\fu_{\intercal}=\Pi_{\intercal} \fu_{\intercal}$, 
    we have that
    \begin{equation}
        \bG(0) = 0 \,, \quad \Pi_{\perp} \di_{\fu} \bG(0) [\fu_{\perp}] = 0 \,.
    \end{equation}
    We deduce that
    \begin{align*}
      {\wt\bG_{\perp}(0)} & = {\Pi_{\perp}} {\bG}\big( \fu_{\intercal}(0) \big) \\
        & = \Pi_{\perp}\Big( \bG(0) + \di_{\fu}\bG(0)[\fu_{\intercal}(0)] +  \varepsilon \bQ\big(\fu_{\intercal}(0)\big) \Big) \\
         & = \varepsilon \Pi_{\perp} \bQ\big(\fu_{\intercal}(0)\big)\,,
    \end{align*}
    where
    \begin{equation}
        {\bQ}(\fu) := \varepsilon^{-1}\big( \bG(\fu) - \bG(0) - \di_{\fu} \bG(0)[\fu] \big)
    \end{equation}
    satisfies the quadratic estimate $\|\bQ(\fu) \|_{s}^{5,\upsilon}\lesssim_{s} \big( \| \fu\|_{s}^{5,\upsilon}\big)^2 $ by \eqref{G(u).def} and Lemma \ref{lem-productlaw}-$(v)$. 
\end{proof}
We now prove the following proposition, which constitutes the cornerstone in the construction of solutions to the range equation \eqref{range.eq}.
\begin{pro}\label{NM.iteration} 
{\bf (Nash-Moser)} 
With the notation in \eqref{constants.NASH},
there exist $ \delta \in (0, 1)$, $C_* > 0$ such that if
\begin{equation}  \label{nash moser smallness condition} 
	\begin{aligned}
		& N_{0}^{\overline \tau} \varepsilon  \leq \delta \,, \quad N_{0} := \upsilon^{- 1}\,, \quad \upsilon := \varepsilon^{\mathtt c}\,,  \quad \text{for some} \ \ \  0 < \mathtt c < \bar\tau\,^{-1} =: \tc_{0}  \,,
	\end{aligned}
\end{equation}
then the following properties hold for all $\tn \geq 0$: 
\\[1mm]
\noindent $({\mathcal P}1)_{\tn}$
There exists a $5$-times differentiable function 
$\wt\fu_{\tn}:=\wt\fu_{\perp}^{(\tn)}: \fT_{0} \to \cH_{\tn - 1} $, $\tn \geq 1$, 
$\wt\fu_{0}:=0$, $\cH_{- 1} := \{ 0 \}$, (with $\fT_{0}:=\t\Omega \times [\alpha_{1},\alpha_{2}]$ and $\t\Omega\subset \R $ as in \eqref{Omega-neighb})
satisfying
\begin{equation} \label{stima.bassa.NM}
\| \wt\fu_{\tn} \|_{s_0 + \bar\sigma}^{5,\upsilon,\fT_{0}} \leq C_* \,.
\end{equation}
If $\tn \geq 1$, 
the difference $ \wt\fh_{\tn} := \wt\fu_{\tn} - \wt\fu_{\tn-1}$ satisfies  
\begin{equation} \label{hn}
\| \wt\fh_{1} \|_{s_0 + \bar\sigma}^{5,\upsilon,\fT_{0}} 
\lesssim \varepsilon \upsilon^{-1}\,, \quad
\| \wt\fh_{\tn} \|_{s_0 + \bar\sigma}^{5,\upsilon,\fT_{0}} 
\leq C_*  N_{\tn-1}^{2 \overline \sigma + 1} N_{\tn-2}^{-\mathtt a_1} \,;
\end{equation}
\noindent $({\mathcal P}2)_{\tn}$
We define the sets (recalling \eqref{Omega-neighb} and $\b\Lambda(\upsilon_{\tn},\wt\fu_{\tn})$ defined in \eqref{nonres.set.inv.full})
\begin{align}\label{nonres.step.NM}
    \fT_{0}:= \t\Omega \times [\alpha_{1},\alpha_{2}] \,, \quad \fT_{\tn+1} &:= \fT_{\tn} \cap \b\Lambda(\upsilon_{\tn},\wt\fu_{\tn}) \\
    & = \fT_{\tn} \cap \tT\tC(\upsilon_{\tn},\wt\fu_{\tn}) \cap \t\Lambda(\upsilon_{\tn},\wt\fu_{\tn})\,, \quad \upsilon_{\tn}:=\upsilon(1+2^{-\tn})\,,  \quad \tn \geq 0 \,.
\end{align}
Then we have the estimate
\begin{equation} \label{stima.F.to.0}
\| \wt\bG_{\perp} (\wt\fu_{\tn})  \|_{s_{0}}^{5,\upsilon,\fT_{\tn}}
\leq C_* \varepsilon N_{\tn - 1 }^{- \mathtt a_1}\, ;
\end{equation}
\noindent $({\mathcal P}3)_{\tn}$
We have the estimate
\begin{equation}\label{stima.alta.NM}
	\| \wt\fu_{\tn} \|_{s_0 + \mathtt b_1}^{5,\upsilon,\fT_{0}}
	\leq C_*  N_{\tn-1}^{\kappa} \,.
\end{equation}
\end{pro}

\begin{proof}
To simplify notations, in this proof we write $\| \cdot \|_s$ instead of $\| \cdot \|_s^{5,\upsilon}$, without even specifying in the norms the set of parameter where the functions are defined, as it will be clear from the context.

\smallskip

\noindent
{\sc Proof of $({\mathcal P}1, 2, 3)_0$}.
By Lemma \ref{lemma.start.NM} with $s=s_0$,
$\| \wt\bG_{\perp}(0) \|_s  \leq \varepsilon C(s)$. Then \eqref{stima.bassa.NM}, \eqref{stima.F.to.0} and \eqref{stima.alta.NM} hold taking $ \tfrac12 C_0, \tfrac12C_*(s) \geq C(s)$ sufficiently large. In particular, we have
\begin{equation}\label{stima.w0}
	\| \fu_{0} \|_{s_0+\bar\sigma} \leq \tfrac12 C_0 \leq C_0\,.
\end{equation}

\smallskip

\noindent
{\sc Assume that $({\mathcal P}1,2,3)_{\tn}$ hold for some $\tn \geq 0$, 
and prove $({\mathcal P}1,2,3)_{\tn+1}$.}
By $({\mathcal P}1)_{\tn}$, one has $\| \wt\fu_{\tn} \|_{s_0 + \bar \sigma} \leq C_0$, for some $C_0 > 0$, independent of $\tn\in\N_0$.  
The assumption 
\eqref{nash moser smallness condition} implies 
the smallness condition 
$\varepsilon \upsilon^{- 1} \leq \delta$ 
of Lemma \ref{almost.straight.lemma}
by taking $\overline \tau (\tau, \nu)$ large enough and $S = s_0 + \mathtt b_1$. 
Then we apply Proposition \ref{prop.full.inv}  
to the linearized operator 
\begin{equation}\label{def.sL.n}
\sL_{\tn} := \sL_{\perp}((\omega,\alpha)) : = \di_{\fu_{\perp}} \wt\bG_{\perp} (\wt\fu_{\tn}(\omega,\alpha)) \,.
\end{equation}
This implies that, for any $\lambda:=(\omega,\alpha)\in \fT_{n}$, the operator $\sL_{n}(\lambda)$ admits an almost right inverse $\bT_{\tn}:=\bT_{\tn}(\lambda,\wt\fu_{\tn})$ satisfying the tame estimates, for any $s_0 \leq s \leq s_0 + \tb_{1}$,
\begin{align}
\|  \bT_{\tn} g \|_s 
&\lesssim_{s}  \upsilon^{-1} \big( \| g \|_{s+\bar\sigma}
+ \| \wt\fu_{n} \|_{s + \bar\sigma} \,  
\| g \|_{s_0 +\bar\sigma} \big)\,,\label{est.sTn.s} \\
\|  \bT_{\tn} g \|_{s_0} 
& \lesssim  \upsilon^{-1} \| g \|_{s_0 +\bar\sigma } \,. \label{est.sTm.s0} 
\end{align}
For all $\lambda \in \fT_{\tn+1} := \fT_{\tn} \cap \b\Lambda_{\tn+1}(\upsilon_{\tn},\tau,\wt\fu_{\tn})$ (see \eqref{nonres.step.NM}), we define the successive approximation 
\begin{equation}\label{approx.solutions}
\fu_{\tn + 1} := \wt\fu_{\tn} + \fh_{\tn + 1} \,, \quad 
\fh_{\tn + 1} :=  - \Pi_{\tn} \bT_{\tn}  \Pi_{\tn} \wt\bG_{\perp}(\wt\fu_{\tn}) 
\in {\mathcal H}_{n} \,.
\end{equation}
We now show that the iterative scheme in \eqref{approx.solutions} is rapidly converging. We write
\begin{equation}
    \wt\bG_{\tn}(\wt\fu_{\tn+1}) =  \wt\bG_{\tn}(\wt\fu_{\tn})  + \sL_{\tn} \fh_{\tn+1}+ \bQ_{\tn}\,,
\end{equation}
where
\begin{equation}\label{bQ.n}
    \bQ_{\tn} := \bQ(\wt\fu_{\tn},\fh_{\tn+1})\,, \quad  \bQ(\wt\fu_{\tn},\fh) := \wt\bG_{\perp}(\wt\fu_{\tn} + \fh) - \wt\bG_{\perp}(\wt\fu_{\tn}) - \sL_{\tn} \fh \,, \quad \fh \in \cH_{\tn} \,.
\end{equation}
Then, by the definition of $\fh_{\tn+1}$ in \eqref{approx.solutions}, we have
\begin{align}
      \wt\bG_{\tn}(\wt\fu_{\tn+1}) &  =  \wt\bG_{\tn}(\wt\fu_{\tn})  - \sL_{\tn} \Pi_{\tn} \bT_{\tn}  \Pi_{\tn} \wt\bG_{\perp}(\wt\fu_{\tn}) + \bQ_{\tn} \\
       &  =  \wt\bG_{\tn}(\wt\fu_{\tn})  - \sL_{\tn}  \bT_{\tn}  \Pi_{\tn} \wt\bG_{\perp}(\wt\fu_{\tn}) +  \sL_{\tn} \Pi_{\tn}^\perp \bT_{\tn}  \Pi_{\tn} \wt\bG_{\perp}(\wt\fu_{\tn}) + \bQ_{\tn}  \\
        &  \stackrel{\sL_{\tn}\bT_{\tn}={\rm Id}}{=}  \big({\rm Id} -\Pi_{\tn} \big) \wt\bG_{\tn}(\wt\fu_{\tn}) +  \sL_{\tn} \Pi_{\tn}^\perp \bT_{\tn}  \Pi_{\tn} \wt\bG_{\perp}(\wt\fu_{\tn}) + \bQ_{\tn} \\
     & = \Pi_{\tn}^\perp \wt\bG_{\tn}(\wt\fu_{\tn})  + \bR_{\tn} + \bQ_{\tn} \,, \label{approx.solutions.eval}
\end{align}
where
\begin{align}
  \bR_{\tn} := \sL_{\tn} \Pi_{\tn}^\perp \bT_{\tn}  \Pi_{\tn} \wt\bG_{\perp}(\wt\fu_{\tn})\,. \label{bRn.bPn}
\end{align}
First, we note that, for any $(\omega,\alpha)\in \t\Omega \times [\alpha_{1},\alpha_{2}]$, $s\geq s_0$, by the triangular inequality, \eqref{range.eq}, Proposition \ref{lemma.map.bG.solo}, Lemma \ref{expand.sL.pert}, \eqref{constants.NASH} and \eqref{nash moser smallness condition}, we have
\begin{align}
    \| \wt\bG_{\perp}(\wt\fu_{\tn}) \|_{s}  & \leq \| \wt\bG_{\perp}(0)\|_{s} + \| \wt\bG_{\perp}(\wt\fu_{\tn}) -\wt\bG_{\perp}(0) \|_{s} \lesssim \varepsilon \big( 1 + \| \wt\fu_{\tn} \|_{s+\bar\sigma}\big)\,, \label{est.bG.at.n.s}
\end{align}
and, by \eqref{nash moser smallness condition}, \eqref{stima.bassa.NM},
\begin{equation}\label{est.bG.at.n.s0}
     \| \wt\bG_{\perp}(\wt\fu_{\tn}) \|_{s_0} \leq 1 \,.
\end{equation}
We now want to prove the following claim: for all $\lambda \in \fI_{\tn+1}$, the following inductive Nash-Moser type estimates hold
\begin{align}
    \| \wt\bG_{\perp}(\wt\fu_{\tn+1}) \|_{s_0} & \lesssim N_{\tn}^{\kappa-\tb_{1}} \varepsilon\big( 1+ \|\wt\fu_{\tn} \|_{s_0+\tb_{1}} \big) + N_{\tn}^{\kappa} \| \wt\bG_{\perp}(\wt\fu_{\tn}) \|_{s_0}^2 \,, \label{NM.est.1} \\
    \| \fu_{\tn+1}\|_{s_0+\tb_{1}} & \lesssim N_{\tn}^{\kappa} \varepsilon\big( 1+ \|\wt\fu_{\tn} \|_{s_0+\tb_{1}} \big)\,.\label{NM.est.2}
\end{align}
where $\kappa$ and $\tb_{1}$ are as in \eqref{constants.NASH}. \\
We split the proof of \eqref{NM.est.1}, \eqref{NM.est.2} into several steps, where we estimate each contribution in \eqref{approx.solutions} and \eqref{approx.solutions.eval}.
\\[1mm]
\noindent $\blacktriangleright$ {\bf Estimates for $\fh_{\tn+1}$.} By \eqref{approx.solutions}, \eqref{std.smoothing}, \eqref{est.sTn.s}, \eqref{est.sTm.s0}, \eqref{est.bG.at.n.s}, \eqref{est.bG.at.n.s0} and \eqref{constants.NASH}, \eqref{nash moser smallness condition}, we have
\begin{align}
    \| \fh_{\tn+1} \|_{s_0+\tb_{1}} & \lesssim \upsilon^{-1} N_{\tn}^{\bar\sigma} \big( \|  \wt\bG_{\perp}(\wt\fu_{\tn}) \|_{s_0+\tb_{1}} + N_{\tn}^{\bar\sigma} \| \wt\fu_{\tn} \|_{s_0+\tb_{1}} \| \wt\bG_{\perp}(\wt\fu_{\tn}) \|_{s_0}  \big)  \label{esti.fh.n+1.1} \\
    & \lesssim \varepsilon N_{\tn}^{2\bar\sigma+1} \big( 1 + \| \wt\fu_{\tn} \|_{s_0+\tb_{1}}  \big) \,, \\
     \| \fh_{\tn+1} \|_{s_0} & \lesssim \upsilon^{-1} N_{\tn}^{\bar\sigma } \| \wt\bG_{\perp}(\wt\fu_{\tn}) \|_{s_0} \,. \label{esti.fh.n+1.2}
\end{align}
Moreover, by  \eqref{std.smoothing}, \eqref{approx.solutions}, \eqref{esti.fh.n+1.2} and Lemma \ref{lemma.start.NM}, we get
\begin{align}
    \| \fh_{1} \|_{s_0+\bar\sigma} & \lesssim \upsilon^{-1} \| \wt\bG_{\perp}(0)\|_{s_0+2\bar\sigma} \lesssim \varepsilon \upsilon^{-1}\,, \label{fh1} \\
     \| \fh_{\tn+1} \|_{s_0+\bar\sigma} &  \lesssim N_{\tn}^{\bar\sigma} \| \fh_{\tn+1} \|_{s_0} \lesssim \upsilon^{-1} N_{\tn}^{2\bar\sigma} \| \wt\bG_{\perp}(\wt\fu_{\tn}) \|_{s_0} \,.  \label{fhn+1}
\end{align}
\noindent $\blacktriangleright$ {\bf Estimate of $\bQ_{\tn}$.} By \eqref{bQ.n}, \eqref{approx.solutions}, \eqref{std.smoothing}, Lemma \ref{lem-productlaw}-$(v)$,
\eqref{esti.fh.n+1.2}, \eqref{constants.NASH}, \eqref{nash moser smallness condition}, we get, using that $\varepsilon\upsilon^{-2}\leq 1$,
\begin{align}
    \| \bQ_{\tn} \|_{s_0} \lesssim \varepsilon \| \fh_{\tn+1} \|_{s_0+1}^{2} \lesssim N_{\tn}^{2\bar\sigma+2} \|  \wt\bG_{\perp}(\wt\fu_{\tn})\|_{s_0}^2 \,. \label{bQn.est}
\end{align}
\noindent $\blacktriangleright$ {\bf Estimate of $\bR_{\tn}$.} We are left to estimate of $\bR_{\tn}$ in  \eqref{bRn.bPn}. By Lemma \ref{expand.sL.pert}, Proposition \ref{prop.Malpha}, \eqref{std.smoothing}, \eqref{est.sTn.s}, \eqref{est.sTm.s0}, \eqref{est.bG.at.n.s}, \eqref{est.bG.at.n.s0}, \eqref{constants.NASH}, \eqref{nash moser smallness condition}, we get
\begin{align}
    \| \bR_{\tn} \|_{s_0} & \lesssim \|  \Pi_{\tn}^\perp \bT_{\tn}  \Pi_{\tn} \wt\bG_{\perp}(\wt\fu_{\tn}) \|_{s_0+1} \lesssim N_{\tn}^{1-\tb_{1}} \|   \bT_{\tn}  \Pi_{\tn} \wt\bG_{\perp}(\wt\fu_{\tn}) \|_{s_0+\tb_{1}} \\
    & \lesssim N_{\tn}^{\bar\sigma+1-\tb_{1}}\upsilon^{-1} \big( \| \wt\bG_{\perp}(\wt\fu_{\tn}) \|_{s_0+\tb_{1}} + \| \wt\fu_{\tn} \|_{s_0+\tb_{1}} \| \wt\bG_{\perp}(\wt\fu_{\perp}) \|_{s_0}  \big) \\
    & \lesssim N_{\tn}^{2\bar\sigma+2-\tb_{1}}\varepsilon (1+ \| \wt\fu_{\tn} \|_{s_0+\tb_{1}}) \,. \label{bRn.est}
\end{align}

\noindent $\blacktriangleright$ {\bf Proof of the claimed estimates \eqref{NM.est.1}, \eqref{NM.est.2}.} By \eqref{approx.solutions.eval}, \eqref{bQ.n}, \eqref{bRn.bPn}, \eqref{std.smoothing}, \eqref{bQn.est}, 
\eqref{bRn.est}, we deduce \eqref{NM.est.1}, with $\kappa$ as in \eqref{constants.NASH}. Estimate \eqref{NM.est.2} for $\fu_{\tn+1}:=\wt\fu_{\tn}+\fh_{\tn+1}$ follows by \eqref{esti.fh.n+1.1}.

\medskip

We are now in the position to prove the induction estimates. First, note that, by \eqref{nonres.step.NM}, if $\lambda \in \fI_{\tn+1}$, then $\lambda\in \fI_{\tn}$, and so \eqref{stima.F.to.0}, \eqref{stima.alta.NM} hold. 

By Whitney extension Theorem (see for instance \cite[Th.m B.2]{FMT25}), we define a 5 times differentiable extension  $\wt\fh_{\tn+1}$ of $\fh_{\tn+1}|_{\fI_{\tn+1}}$ to the whole $\R\times [\alpha_{1},\alpha_{2}]$. We also define the function $\wt\fu_{\tn+1} := \wt\fu_{\tn}+ \wt\fh_{\tn+1}$, which is also a  5 times differentiable extension of $\fu_{\tn+1}|_{\fI_{\tn+1}}$ in \eqref{approx.solutions} to the whole $\R\times [\alpha_{1},\alpha_{2}]$. For convenience we consider the restrictions of $\wt\fh_{\tn+1}$ and $\wt\fu_{\tn+1}$ to $\fI_{0}:= \t\Omega\times [\alpha_{1},\alpha_{2}]$. The extensions $\wt\fh_{\tn+1}$ and $\wt\fu_{\tn+1}$ satisfy the same bounds for $\fh_{\tn+1}$ and $\fu_{\tn+1}$, respectively.

Then, \eqref{stima.F.to.0} at the step $\tn+1$ follows by \eqref{NM.est.1}, $(\cP 2)_{\tn}$, $(\cP 3)_{\tn}$, $\varepsilon \ll 1$ small enough and \eqref{constants.NASH}. For $\tn=0$, we also use \eqref{nash moser smallness condition}. The estimate \eqref{stima.alta.NM} at the step $\tn+1$ follows directly from \eqref{NM.est.2}, \eqref{stima.alta.NM}, \eqref{constants.NASH} and \eqref{nash moser smallness condition}, with $N_{0}\gg 1$ large enough. The estimates \eqref{hn} at the step $\tn+1$ follows from \eqref{fh1}, \eqref{fhn+1} and \eqref{stima.bassa.NM}. Finally, by \eqref{approx.solutions} and \eqref{hn}, we have
\begin{equation}
    \| \wt\fu_{\tn+1} \|_{s_0+\overline{\sigma}} \leq \sum_{k=1}^{\tn+1} \| \wt\fh_{k} \|_{s_0+\overline{\sigma}} \leq C_{*} \,,
\end{equation}
where the convergence of the series independently of $\tn$ follows by \eqref{constants.NASH}, \eqref{nash moser smallness condition}. This proves \eqref{stima.bassa.NM} at the step $\tn+1$ and concludes the proof the theorem.
\end{proof}

\subsection{Conclusion}\label{Concusion}
By \eqref{hn}, the sequence of functions $ (\wt\fu_{\perp}^{(\tn)})_{\tn\in\N_{0}}$ is a Cauchy sequence in $\| \,\cdot\, \|_{s_0+\bar\sigma}^{5,\upsilon}$ and it converges to
\begin{equation}
    \fu_{\perp}^{(\infty)} := \lim_{\tn\to\infty} \wt\fu_{\perp}^{(\tn)} \in H_\perp^{s_0+\bar\sigma}(\T^2) \,.
\end{equation}
We also deduce that
\begin{equation}
    \| \fu_{\perp}^{(\infty)} \|_{s_0 + \bar\sigma}^{5,\upsilon} 
\leq C_* \,, \quad \| \fu_{\perp}^{(\infty)} - \wt\fu_{\perp}^{(\tn)}\|_{s_0 + \bar\sigma}^{5,\upsilon} 
\leq C_*  N_{\tn}^{2 \overline \sigma + 1} N_{\tn-1}^{-\mathtt a_1} \,.
\end{equation}
Moreover, by \eqref{stima.F.to.0} in Proposition \ref{NM.iteration}-$(\cP2)_{\tn}$, we obtain that
\begin{equation}
    \wt\bG_{\perp}\big(\fu_{\perp}^{(\infty)}(\omega,\alpha);\omega,\alpha\big) = 0 \,,
\end{equation}
for all $(\omega,\alpha)$ belonging to the set of good parameters (recalling \eqref{nonres.step.NM}, \eqref{nonres.set.inv.full})
\begin{align}\label{final.set.nonres}
    \fT_{\infty} & := \bigcap_{\tn \geq 0} \fT_{\tn} = \fT_{0} \cap \bigcap_{\tn\geq 0} \b\Lambda(\upsilon_{\tn},\wt\fu_{\perp}^{(\tn)}) \\
    & = \fT_{0} \cap \bigg( \bigcap_{\tn\geq 0} \tT\tC(\upsilon_{\tn},\wt\fu_{\perp}^{(\tn)}) \bigg) \cap \bigg( \bigcap_{\tn\geq 0} \t\Lambda(\upsilon_{\tn},\wt\fu_{\perp}^{(\tn)}) \bigg) \,.
\end{align}
The fact that the set $\fT_{\infty}$ has large measure is left to be proven in the next section.

\section{Measure estimates}\label{section-measures}


In this section, we complete the proof of the main result by establishing suitable measure estimates for the set of admissible parameters. More precisely, our goal is to show that the non-resonance conditions required for the invertibility of the linearized operator along Nash-Moser scheme are satisfied on a large subset of the parameter space.
The main task is therefore to control the size of the excluded set and to prove that the remaining parameter set has asymptotically full measure as the perturbation parameter tends to zero. This requires a careful analysis of the dependence of the eigenvalues on the parameters, as well as sharp estimates on the measure of the resonant regions.
\subsection{Final Cantor set estimate}
We consider the set
\begin{equation}\label{Omega-neighb}
	\t\Omega := \big\{\omega \in \R \,:\, \inf_{\alpha\in[\alpha_{1},\alpha_{2}]} |\omega + \im\,\lambda_{\tJ}^{(\infty)}(\alpha)|<\delta \big\} \,, \quad 0<\delta\ll \alpha_{0} \,,
\end{equation}
be a $\delta$-neighbourhood of the curve ${-\im}\lambda_{\tJ}^{(\infty)}([\alpha_{1},\alpha_{2}])\subset \R$, where, recalling \eqref{alpha1alpha2}, 
\begin{equation}
	0<\tfrac14 \alpha_{0}=: \alpha_{1} < \alpha_{2} := \tfrac34 \alpha_{0}<\alpha_{0} \,.
\end{equation}
In Section \ref{section.normal}, and in particular in Subsection \ref{Concusion}, we have constructed the solution $\fu_{\perp}^{(\infty)}= \fu_{\perp}^{(\infty)}(\vf,\theta;\omega,\alpha,\varepsilon)$  via a Nash-Moser iteration. Now, we will  fix the frequency oscillation. For this aim, we 
return to Theorem \ref{theo.FP.tangential}  and produce
\begin{equation}
	\tw_{\infty}:= \tw_{\infty}(\omega,\alpha,\varepsilon) := -\varepsilon\braket{\tP_{\Theta}\big(\bY   \big(\,\cdot\,;\fu_{\perp}^{(\infty)}(\,\cdot\,;\omega,\alpha,\varepsilon),\omega,\alpha,\varepsilon\big),\fu_{\perp}^{(\infty)}(\,\cdot\,;\omega,\alpha,\varepsilon);\alpha,\varepsilon\big)}_{\vf}\,.
\end{equation}
Then, for any $(\omega,\alpha)\in \fT_{\infty}$, with $\fT_{\infty}$ as in \eqref{final.set.nonres},
we obtain a periodic solution of \eqref{mod.aa.sys} by Theorem \ref{theo.FP.tangential} with $\tw=\tw_\infty$. From the compatibility condition \eqref{tw.back}, the final frequency of oscillation $\bomega_{\tJ,\varepsilon}(\alpha)$ is implicitly given by the nonlinear equation
\begin{equation}\label{implicit.final.freq}
	\bomega_{\tJ,\varepsilon}(\alpha)=\omega= -\im\,\lambda_{\tJ}^{(\infty)}(\alpha)+ \tw_{\infty}(\omega,\alpha)= -\im\lambda_{\tJ}^{\infty}(\alpha)+ \tw_{\infty}(\bomega_{\tJ,\varepsilon}(\alpha),\alpha)\,.
\end{equation}
This equation is easily  solved by the Picard fixed point theorem, as the map $\omega \mapsto \tw_{\infty}(\omega,\alpha)$ is smooth and satisfies
\begin{equation}\label{hostel1}
	|\tw_{\infty}|^{5,\upsilon} \leq C \varepsilon \,. 
\end{equation}
Actually,  in view of Corollary \ref{cor.eigen.asympt}, one gets
		\begin{align}
			\bomega_{\tJ,\varepsilon}(\alpha) & = -\im\,\lambda_{\tJ}^{(\infty)}(\alpha) + \wt\tw_{\infty}(\alpha) \\
			&=   \big(\tfrac{1}{2} -\alpha^4 +  \tc_{1}(\alpha) \big)\tJ - \tfrac12 +\tfrac32 \alpha^4 \delta_{2,\tJ}   + \tr_{\tJ}(\alpha) + \wt\tw_{\tJ,\infty}(\alpha)\,,\label{ME.1} 
		\end{align}
        where $\tc_{1}(a)$, $\tr_{\tJ}(\alpha)$ satisfy the estimates \eqref{small.const.alpha} and where we denote
	\begin{equation}
		\wt\tw_{\tJ,\infty}(\alpha):= \tw_{\infty}(\bomega_{\tJ,\varepsilon}(\alpha),\alpha) \,, 
	\end{equation}
	which satisfies the estimate, using \eqref{hostel1},
	\begin{align}
		&\sup_{\alpha\in[\alpha_{1},\alpha_{2}]}|\pa_{\alpha}^{k}\wt\tw_{\tJ,\infty}(\alpha)| \leq C \varepsilon \upsilon^{-k} \,, \quad \forall \, |k|\leq 5 \,.\label{ME.2} 
	\end{align}
We now establish the main result of this section, showing that the final Cantor set has asymptotically full Lebesgue measure.
\begin{theo}{\bf (Measure estimates).}\label{meas.est.thm}
	Recalling \eqref{nash moser smallness condition}, let
	\begin{equation}\label{param.meas.est}
		\upsilon= \varepsilon^{\tc} \,, \quad 0<\tc<\tc_0 <1\,, \quad \tau > 7 \,.
	\end{equation}
	There exists $\varepsilon_0>0$ and $C>0$ such that for any  $\varepsilon\in (0,\varepsilon_{0})$, the Lebesgue measure of the set
    \begin{align}\label{set.to.measure}
        \Omega_{\varepsilon} &:= \big\{  \alpha\in[\alpha_{1},\alpha_{2}] \, : \, \big( \bomega_{\tJ,\varepsilon}(\alpha),\alpha \big) \in \fT_{\infty} \big\} 
    \end{align}
	satisfies $$|\Omega_{\varepsilon}|\geqslant  \alpha_{2}-\alpha_{1}- C\varepsilon^{\tc/4}\,.
    $$
\end{theo}
The proof will be carried out in several steps in the following sections via intermediate results.
\subsection{Cantor set decomposition}
In what follows, we provide a suitable decomposition of the final Cantor set described in \eqref{set.to.measure}. To this end, we recall some preliminaries. For any \( j \in \Z_{\perp} \), we introduce the notation, recalling Corollary \ref{cor.eigen.asympt} and Lemma \ref{diffeo.conj.pert},
    \begin{align}
			-\im \, \wt\mu_{j}^{(\tn)}(\alpha)&:= -\im\,\wt\mu_{j}^{(\tn)}\big(\bomega_{\tJ,\varepsilon}(\alpha),\alpha; \wt\fu_{\perp}^{(\tn)}\big)\\
			& = {- \im\,  \lambda_{j}^{(\infty)}(\alpha)  + \big(  \wt\tm_{1,\varepsilon}^{(\tn)}(\alpha)}-\big( \tfrac12 - \alpha^4 + \tc_{1}(\alpha) \big) \big)j \\
			& =   \wt\tm_{1,\varepsilon}^{(\tn)}(\alpha)j +{\rm sgn}(j)\big(- \tfrac12  + \tfrac32 \alpha^4 \delta_{2,|j|} \big)+ \tr_{j}(\alpha) \,, \label{ME.1.1}
		\end{align}
    where $\delta_{2, 2}:= 1$, $\delta_{2,|j|}:=0$ for $|j|\geq 3$. The functions   $\tc_{1}(a)$ and $\tr_{j}(\alpha)$ satisfy the estimates \eqref{small.const.alpha}. In addition
\begin{equation}
    \wt\tm_{1,\varepsilon}^{(\tn)}(\alpha) := \tm_{1,\varepsilon}\big(\bomega_{\tJ,\varepsilon}(\alpha),\alpha; \wt\fu_{\perp}^{(\tn)}\big) \,,
\end{equation}
 satisfies the estimate, by Lemma \ref{almost.straight.lemma},
\begin{equation}\label{diff-Hm}
 \sup_{\alpha\in[\alpha_{1},\alpha_{2}]}\big|\pa_{\alpha}^{k} \big(\wt\tm_{1,\varepsilon}^{(\tn)}(\alpha) -\big(\tfrac12 - \alpha^4 + \tc_{1}(\alpha) \big) \big)\big| \leq C \varepsilon \upsilon^{-k} \,,  \quad \forall\, |k|\leq 5 \,.   
\end{equation}
For any $\tn \geq 0$, recalling  \eqref{nonres.step.NM}, we denote by
\begin{equation}\label{set.n.to.measure}
     \Omega_{\tn} := \big\{  \alpha\in[\alpha_{1},\alpha_{2}] \, : \, \big( \bomega_{\tJ,\varepsilon}(\alpha),\alpha \big) \in \fT_{\tn} \big\} 
\end{equation}
and, recalling also \eqref{nonres.transport}, we define
\begin{equation}
    \fE_{\tn}^{(0)}  := \big\{  \alpha\in \Omega_{\tn} \, : \, \big( \bomega_{\tJ,\varepsilon}(\alpha),\alpha \big) \in \tT\tC (\upsilon_{\tn},\wt\fu_{\perp}^{(\tn)})  \big\} \,.
\end{equation}
 Note that $\Omega_{\tn}\setminus \fE_{\tn}^{(0)}$ is given by
\begin{align}
   \Omega_{\tn}\setminus \fE_{\tn}^{(0)} & = \bigcup_{(\ell,j)\in \Z^2\setminus\{0\}} R_{\ell,j}^{({\rm Tr},\tn)} \,, 
   \\
   R_{\ell,j}^{({\rm Tr},\tn)} & := \big\{  \alpha\in \Omega_{\tn} \, : \, \big|\bomega_{\tJ,\varepsilon}(\alpha)\,\ell + \wt\tm_{1,\varepsilon}^{(\tn)}(\alpha)j \big| <  2\upsilon_{\tn} \, |j|^{-\tau}  \big\}\,. \label{tappo1}
\end{align}
To estimate the measure of the final Cantor set \( \Omega_{\varepsilon} \) defined in \eqref{set.to.measure}, it suffices to analyze its complement, given by
\begin{equation}\label{compl.set.to.measure}
    [\alpha_{1},\alpha_{2}] \setminus \Omega_{\varepsilon} = \bigcup_{\tn \geq 0} \Omega_{\tn}\setminus \Omega_{\tn+1} \,. 
\end{equation}
By \eqref{set.n.to.measure}, \eqref{final.set.nonres}, \eqref{nonres.step.NM}, \eqref{nonres.set.inv.full}, \eqref{nonres.transport}, \eqref{Cantor-n-2.full}, we have, for any $\tn \geq 0$,
\begin{equation}
    \Omega_{\tn}\setminus \Omega_{\tn+1} \subseteq \big( \Omega_{\tn}\setminus \fE_{\tn}^{(0)} \big) \cup \fE_{\tn}^{(1)} \,, \quad  \fE_{\tn}^{(1)} := \bigcup_{\ell \in \Z\setminus\{0\} \atop
    j\in \Z_{\perp} } R_{\ell,j}^{(1,\tn)} \,, \label{tappo4}
\end{equation}
where the set $R_{\ell,j}^{(1,\tn)}$ is defined as
\begin{equation}
    R_{\ell,j}^{(1,\tn)} := \big\{  \alpha \in \fE_{\tn}^{(0)} \, : \, \big|\bomega_{\tJ,\varepsilon}(\alpha)\,\ell -\im\,\wt\mu_{j}^{(\tn)}(\alpha)\big| < 2 \upsilon_{\tn} \, |j|^{-\tau}   \big\}\,. \label{tappo5}
\end{equation}

\subsection{Trivial cases}
We aim to estimate the Lebesgue measure of each block \( \Omega_{\tn} \setminus \Omega_{\tn+1} \) in \eqref{compl.set.to.measure}. In the following two results below, we show that some of these sets are empty.

\begin{lem}\label{lemma.empty.n.TR}
{\bf (Trivial $R_{\ell,j}^{({\rm Tr},\tn)}$ sets).}
    There exist $\alpha_{0}, \varepsilon_{0} > 0$ sufficiently small such that, for all $\tn \in \mathbb{N}_0$, the following holds:
    \\[1mm]
    \noindent $(i)$ For any $\ell,j\in\Z\setminus\{0\}$, we have $R_{\ell,0}^{({\rm Tr},\tn)}=\emptyset$ and $R_{0,j}^{({\rm Tr},\tn)}=\emptyset$;
    \\[1mm]
    \noindent $(ii)$ If $\ell=\pm 1$ and $j\in\Z\setminus\{0\}$, with $j\neq \mp(\tJ-1)$, then  $R_{\ell,j}^{({\rm Tr},\tn)}=\emptyset$;
    \\[1mm]
    \noindent $(iii)$ If $\tJ\geq 3$, $|\ell|\geq 2$ and $j\in\Z$, with   $|\tJ \ell + j| \geq \tfrac52 |\ell|$, then $R_{\ell,j}^{({\rm Tr},\tn)}=\emptyset$;
    \\[1mm]
    \noindent $(iv)$ If $\tJ=2$, and either $|\ell|\geqslant\tfrac32|j|$ or $|j|\geqslant\tfrac32|\ell|$, with $(\ell,j)\neq(0,0),$ then $R_{\ell,j}^{({\rm Tr},\tn)}=\emptyset$.
\end{lem}
\begin{proof}
    By \eqref{ME.1}, \eqref{ME.1.1}, we have that
    \begin{align}
        \bomega_{\tJ,\varepsilon}(\alpha)\,\ell + \wt\tm_{1,\varepsilon}(\alpha) j
        & = -\im \lambda_{\tJ}^{(\infty)}(\alpha) \ell + \big( \tfrac12 -\alpha^4 + \tc_{1}(\alpha) \big)  \big)j \\
         & \quad + \big( \wt\tm_{1,\varepsilon}^{(\tn)}(\alpha) - \big( \tfrac12 -\alpha^4 + \tc_{1}(\alpha) \big)  \big)j + \wt\tw_{\tJ,\infty}(\alpha) \ell  \,,\label{smalldiv.TR}
    \end{align}
    where
    \begin{align}
    	 -\im \lambda_{\tJ}^{(\infty)}(\alpha) \ell & + \big( \tfrac12 -\alpha^4 + \tc_{1}(\alpha) \big)  \big)j  \\
         & = \big( \tfrac12 - \alpha^4 +\tc_{1}(\alpha)  \big)(\tJ \ell +j) - \tfrac12\ell + \tfrac32 \alpha^4 \delta_{2,\tJ} \ell + \tr_{\tJ}(\alpha)\ell\,. \label{smalldiv.unpert.TR}
    \end{align}
      {\sc Proof of $(i)$}. The claim easily follows by noting that, using Corollary \ref{cor.eigen.asympt} the estimates \eqref{ME.2}, \eqref{diff-Hm}, and $\tJ\geq 3$ by \eqref{tangential.mode}, 
      \begin{align}
          | \bomega_{\tJ,\varepsilon}(\alpha)\,\ell | & \geq |-\im \lambda_{\tJ}^{(\infty)}(\alpha)||\ell| - |\wt\tw_{\tJ,\varepsilon}(\alpha)||\ell| \\
          &  \geq \big( \tfrac14 \tJ - C\alpha^4 \tJ - C\varepsilon\big) |\ell| \geq \tfrac18 |\ell| \,, \\
          |\wt\tm_{1,\varepsilon}(\alpha)j| &  \geq \big| \tfrac12 - \alpha^4 + \tc_{1}(\alpha) \big||j| -\big| \wt\tm_{1,\varepsilon}^{(\tn)}(\alpha) - \big( \tfrac12 -\alpha^4 + \tc_{1}(\alpha) \big)  \big| |j| \\
          & \geq \big( \tfrac12 - C\alpha_0^4 -C\varepsilon \big) |j| \geq \tfrac14 |j| \,,
      \end{align}
      for $\alpha_0$ and $\varepsilon_{0}$ small enough.
      \\[1mm]
     {\sc Proof of $(ii)$}. For $\ell= \pm 1$, we begin by proving the following claim.
   \\[1mm]
   {\bf Claim 1.} {\it There exists $\alpha_{0}>0$ small enough such that, for any $\sigma\in \{-,+\}$}
    \begin{equation}\label{claim-1.TR}
    	\big| -\sigma \im \lambda_{\tJ}^{(\infty)}(\alpha)  + \big( \tfrac12 -\alpha^4 + \tc_{1}(\alpha) \big)  \big)j  \big| \geq \tfrac14\big| \sigma(\tJ-1)+j \big|  \quad \forall \, j \in\Z\setminus\{0\}\,, \ \ j\neq\sigma(1-\tJ) \,.
    \end{equation}   
To prove Claim 1, we introduce, for $\sigma=\{-,+\}$,
\begin{align}
  		f_j^{\sigma}(\alpha)& := -\sigma \im \lambda_{\tJ}^{(\infty)}(\alpha)  + \big( \tfrac12 -\alpha^4 + \tc_{1}(\alpha) \big)  \big)j \\
        & =  \big( \tfrac12 - \alpha^4 +\tc_{1}(\alpha)  \big)(\sigma \tJ  +j) -\sigma \tfrac12 +\sigma \tr_{\tJ}(\alpha)\,.\label{ell.pm1.TR}
  	\end{align}
Evaluating at $\alpha=0$, we obtain $f_{j}^{\sigma}(0) = \tfrac12 \big( \sigma (\tJ-1) + j \big)$.
Since we are assuming $j\in\Z\setminus\{0\}$ with $j\neq \sigma(1-\tJ)$, we have 
\begin{equation}
    |f_{j}^+(0)|  = \tfrac12 \big| \sigma(\tJ-1)+j \big|\geq \tfrac12 \,.
\end{equation}
Using that $|\sigma\tJ+j|\leq |\sigma(\tJ-1)+j|+1 \leq 2 |\sigma(\tJ-1)+j|$ for $j\neq \sigma(1-\tJ)$ and Corollary \ref{cor.eigen.asympt}, we deduce that
\begin{align}
    |f_{j}^\sigma(\alpha)| & \geq |f_{j}^\sigma(0)| - |-\alpha^4+\tc_{1}(\alpha)| |\sigma\tJ+j| - |\tr_{\tJ}(\alpha)| \\
     & \geq \tfrac12 \big| \sigma(\tJ-1)+j \big| - C\alpha_0^4|\sigma\tJ +j| -C\alpha_0^5 \\
    & \geq \big( \tfrac12 - 3 C \alpha_0^4  \big)\big| \sigma(\tJ-1)+j \big| \geq \tfrac14 \big| \sigma(\tJ-1)+j \big| \,,
\end{align}
for $\alpha_{0}$ small enough, where $C>0$ is a constant independent of $j$ and $\alpha$. This proves \eqref{claim-1.TR}.
  \\[1mm]
  Using \eqref{claim-1},  \eqref{ME.2} and \eqref{diff-Hm}, we estimate \eqref{smalldiv} with $\ell = \pm 1$  and $j\in \Z\setminus\{0\}$, $j\neq \mp(\tJ-1)$, as
  \begin{align}
  	 \big| \pm \bomega_{\tJ,\varepsilon}(\alpha) + \wt\tm_{1,\varepsilon}(\alpha)j \big| &  \geq \big|f_j^{\pm}(\alpha) \big| - \big| \wt\tm_{1,\varepsilon}^{(\tn)}(\alpha) - \big( \tfrac12 -\alpha^4 + \tc_{1}(\alpha) \big)  \big| |j| - |\wt\tw_{\tJ,\infty}(\alpha)| \\
  	 & \geq \tfrac14 \big| j \pm (\tJ-1) \big|- C\varepsilon |j| -C\varepsilon\,.
  \end{align}
  By the triangle inequality $|j| \leqslant \big|j \pm (\tJ-1)\big| + |\tJ-1|$, together with the fact that $\big|j \pm (\tJ-1)\big| \geqslant 1$ for $j\neq \mp(\tJ-1)$, we obtain
  \begin{align}
  	 \big| \pm \bomega_{\tJ,\varepsilon}(\alpha)  + \wt\tm_{1,\varepsilon}(\alpha)j\big| 
  	 & \geqslant \big(\tfrac14-C \varepsilon\big)\big| j \pm (\tJ-1) \big| -C\varepsilon |\tJ-1|-C\varepsilon \\
     &\geqslant  \big(\tfrac14-C \varepsilon-C\varepsilon |\tJ-1|-C\varepsilon \big)\geqslant  \tfrac{1}{8}\,,
  \end{align}
  where the last inequality holds for $\varepsilon>0$ small enough. This concludes the proof of item $(ii)$.
    	\\[1mm]
   	{\sc Proof of $(iii)$}. We have the following claim.
   	\\[1mm]
   	{\bf Claim 2.} {\it There exists $\alpha_{0}>0$ small enough such that}
   	\begin{equation}\label{claim-2.TR}
   			\big|  	-\im \lambda_{\tJ}^{(\infty)}(\alpha) \ell  + \big( \tfrac12 -\alpha^4 +\tc_{1}(\alpha)  \big) j  \big| \geq  \tfrac{3}{20} |\tJ\ell+j|\,,
   	\end{equation}
   	{\it occurs for $|\ell |\geq 2$ and $j\in \Z$ with  $|\tJ \ell + j|\geq 	\tfrac52 |\ell|$.}
    \\[1mm]
   	Let us prove this claim.
   		By \eqref{smalldiv.unpert.TR}, we have
   		\begin{align}
   				-\im \lambda_{\tJ}^{(\infty)}(\alpha) \ell  + \big( \tfrac12 -\alpha^4 +\tc_{1}(\alpha) \big) j   = \big( \tfrac12 - \alpha^4 +\tc_{1}(\alpha)  \big)(\tJ \ell +j) - \tfrac12\ell + \tr_{\tJ}(\alpha)\ell \,.
   		\end{align}
   		Since we are assuming that $|\tJ \ell + j| \geq \tfrac{5}{2}|\ell|$, it follows, using also \eqref{small.const.alpha} in Corollary \ref{cor.eigen.asympt}, that
        \begin{equation}
   			\begin{aligned}
   				\big|  -\im \lambda_{\tJ}^{(\infty)}(\alpha) \ell  + \big( \tfrac12 -\alpha^4 +\tc_{1}(\alpha) \big) j 	  \big| & \geq  \big|\tfrac12 -\alpha^4 + \tc_{1}(\alpha) \big| |\tJ\ell + j| -  \tfrac12 |\ell|- |\tr_{\tJ}(\alpha)| |\ell|  \\
   				& \geq  \big( \tfrac12 - C\alpha_0^4 - \tfrac15  - \tfrac25 C \alpha_0^5 \big) |\tJ\ell + j|  \\
   	            & \geq \tfrac{3}{20}|\tJ \ell + j| \,,
   			\end{aligned}
   		\end{equation}
   		for some  $C>0$ independent of $\ell,j,\alpha$ and for $\alpha_{0}$ sufficiently small. This ends the proof of \eqref{claim-3}.
        \\[1mm]
   	Using Claim 2, \eqref{ME.2}, and \eqref{diff-Hm}, we estimate \eqref{smalldiv.TR}, for $|\ell|\geq 2$ and $j\in\Z$ with  $|\tJ \ell + j|\geq \tfrac52 |\ell|$, by
   	\begin{align}
   		\big| \bomega_{\tJ,\varepsilon}(\alpha)  \ell + \wt\tm_{1,\varepsilon}(\alpha)j \big| &  \geq \big|  -\im \lambda_{\tJ}^{(\infty)}(\alpha) \ell  + \big( \tfrac12 -\alpha^4 +\tc_{1}(\alpha) \big) j 	  \big| \\
        &-  \big| \wt\tm_{1,\varepsilon}^{(\tn)}(\alpha) - \big( \tfrac12 -\alpha^4 + \tc_{1}(\alpha) \big)  \big| |j| - |\wt\tw_{\tJ,\infty}(\alpha)| |\ell|  \\
   		& \geq \big( \tfrac{3}{20} - C\varepsilon\big(\tfrac{2\tJ}{5}+1\big)  - \tfrac25 C\varepsilon \big)|\tJ\ell+j| \geq \tfrac{3}{100}|\ell|  \,,
   	\end{align}
   	where the last inequality holds for $\varepsilon>0$ small enough, and where we have used in the last line the fact that
    $$
    |j|\leqslant \tJ|\ell|+|\tJ\ell+j|\leqslant \big(\tfrac{2\tJ}{5}+1\big)|\tJ\ell+j|\,.
    $$
    This concludes the proof of item $(iii)$.
    \\[1mm]
    \noindent {\sc Proof of $(iv)$.} When $\tJ=2$, by \eqref{smalldiv.TR}, \eqref{smalldiv.unpert.TR} we have
    \begin{align}
        \bomega_{2,\varepsilon}(\alpha)\,\ell + \wt\tm_{1,\varepsilon}(\alpha) j & = \tfrac12 (\ell + j) + \big( -\alpha^4 + \tc_1(\alpha) \big) (2\ell + j) + \big( \tfrac32 \alpha^4 + \tr_{2}(\alpha) +\wt\tw_{\tJ,\infty}(\alpha) \big)\ell  \\
        &  + \big( \wt\tm_{1,\varepsilon}^{(\tn)}(\alpha) - \big( \tfrac12 -\alpha^4 + \tc_{1}(\alpha) \big)  \big)j \,. 
    \end{align}
    If $|\ell|\geq \frac32 |j|> |j|$, we estimate, by \eqref{small.const.alpha} in Corollary \ref{cor.eigen.asympt}, \eqref{ME.2} and \eqref{diff-Hm},
    \begin{align}
        | \bomega_{2,\varepsilon}(\alpha)\,\ell + \wt\tm_{1,\varepsilon}(\alpha) j | & \geq \tfrac12 ( |\ell| - |j| ) - C \alpha_{0}^4 ( 2|\ell| +|j| ) - \big( C\alpha_0^4 + C\varepsilon \big)|\ell| - C\varepsilon |j| \\
        & \geq \tfrac{7}{16} |\ell| - \tfrac{9}{16}  |j| \geq \big( \tfrac{21}{32} - \tfrac{9}{16} \big) |j| = \tfrac{3}{32} |j|\,,
    \end{align}
    assuming $\alpha_0$ and $\varepsilon_0$ small enough. The analogous estimate holds also in the case when $|j|\geq  \tfrac32 |\ell|>|\ell|$. This concludes the proof of item $(iv)$ and of the lemma.
\end{proof}

\begin{lem}\label{lemma.empty.n}
{\bf (Trivial $R_{\ell,j}^{(1,\tn)}$ sets).}
    There exist $\alpha_{0}, \varepsilon_{0} > 0$ sufficiently small such that, for all $\tn \in \mathbb{N}_0$, the following holds:
    \\[1mm]
    \noindent $(i)$ If $\ell=0$ or $|\ell|=1$, and $j\in\Z_{\perp}$, then $R_{\ell,j}^{(1,\tn)}=\emptyset$;
    \\[1mm]
    \noindent $(ii)$  If $\tJ\geq 3$, $|\ell|\geq 2$ and $|j|=2$, then $R_{\ell,j}^{(1,\tn)}=\emptyset$;
    \\[1mm]
    \noindent $(iii)$ If $\tJ\geq 3$, $|\ell|\geq 2$ and $j\in\Z_{\perp}$, with $|j|\geq 3$,  satisfy $|\tJ \ell + j| \geq \tfrac52 |\ell|$, then $R_{\ell,j}^{(1,\tn)}=\emptyset$;
    \\[1mm]
    \noindent $(iv)$  If $\tJ=2$, $|\ell|\geqslant 2$  and $j\in \Z_{\perp} $, either satisfying $|\ell|\geqslant \tfrac32|j-{\rm sgn}(j)|$ or  $ |j-{\rm sgn}(j)|\geqslant \tfrac32|\ell|$, then $R_{\ell,j}^{(1,\tn)}=\emptyset$.
\end{lem}
\begin{proof}
    By \eqref{ME.1}, \eqref{ME.1.1}, we have that
    \begin{align}
        \bomega_{\tJ,\varepsilon}(\alpha)\,\ell  -\im\,\wt\mu_{j}^{(\tn)}(\alpha)  
        & = -\im \big( \lambda_{\tJ}^{(\infty)}(\alpha) \ell + \lambda_{j}^{(\infty)}(\alpha) \big)  \\
         & \quad + \big( \wt\tm_{1,\varepsilon}^{(\tn)}(\alpha) - \big( \tfrac12 -\alpha^4 + \tc_{1}(\alpha) \big)  \big)j + \wt\tw_{\tJ,\infty}(\alpha) \ell  \,,\label{smalldiv}
    \end{align}
    where
    \begin{align}
    	-\im \big( \lambda_{\tJ}^{(\infty)}(\alpha) \ell + \lambda_{j}^{(\infty)}(\alpha) \big)   & = \big( \tfrac12 - \alpha^4 +\tc_{1}(\alpha)  \big)(\tJ \ell +j) - \tfrac12(\ell +{\rm sgn}(j)) \\
    	&\quad+ \tr_{\tJ}(\alpha)\ell + {\rm sgn}(j)\tfrac32\alpha^4 \delta_{2,|j|} + \tr_{j}(\alpha)\,. \label{smalldiv.unpert}
    \end{align}
     {\sc Proof of $(i)$}. The case $\ell=0$ easily follows by noting that, using \eqref{smalldiv.unpert}, Corollary \ref{cor.eigen.asympt} and the estimate \eqref{diff-Hm}, for any $j\in \Z_{\perp}$,
      \begin{align}
         |-\im \wt\mu_{j}^{(\tn)}(\alpha)| & \geq |-\im \lambda_{j}^{(\infty)}(\alpha)|  -\big| \wt\tm_{1,\varepsilon}^{(\tn)}(\alpha) - \big( \tfrac12 -\alpha^4 + \tc_{1}(\alpha) \big)  \big| |j| \\
          & \geq \tfrac12 \big| j- {\rm sgn}(j) \big| - \big| \big(-\alpha^4 +\tc_{1}(\alpha) \big) j + {\rm sgn}(j) \tfrac32 \alpha^4 \delta_{2,|j|} \big| - C\varepsilon |j| \\
          & \geq \big( \tfrac14 - C \alpha_0^4 - C \varepsilon \big) |j| \geq \tfrac18 |j|\,,
      \end{align}
      for $\alpha_0$ and $\varepsilon_{0}$ small enough. We now consider the cases $\ell= \pm 1$. We begin by proving the following claim.
   \\[1mm]
   {\bf Claim 1.} {\it There exists $\alpha_{0}>0$ small enough such that}
    \begin{equation}\label{claim-1}
    	\big| -\im \big( \pm \lambda_{\tJ}^{(\infty)}(\alpha) + \lambda_{j}^{(\infty)}(\alpha) \big)  \big| \geq \tfrac18|j\pm \tJ| \quad \forall \, j \in\Z_{\perp} \,.
    \end{equation} 
To prove Claim 1, we introduce
\begin{align}
  		f_j^{\pm}(\alpha)&:=-\im \big( \pm \lambda_{\tJ}^{(\infty)}(\alpha)+ \lambda_{j}^{(\infty)}(\alpha) \big)  \\
        & = \big( \tfrac12 - \alpha^4 +\tc_{1}(\alpha)  \big)(\pm \tJ  +j) - \tfrac12(\pm 1 +{\rm sgn}(j)) \\
  		& \quad \pm \tr_{\tJ}(\alpha) + \tfrac32\alpha^4 \delta_{2,j} + \tr_{j}(\alpha)\,.\label{ell.pm1}
  	\end{align}
Evaluating at $\alpha=0$, we obtain
$$
2 \,f_j^{\pm}(0)=  \pm\tJ  + j\mp1-{\rm sgn}(j)\,.
$$
We distinguish two cases:  $j\geqslant2$ and $j\leqslant-2.$ In the first case, we get
$$
2 \,f_j^{+}(0)=  \tJ  + j-2\quad\hbox{and}\quad 2 \,f_j^{-}(0)= - \tJ  + j\,.
$$
Since $\tJ\geqslant2$ and $j\geqslant2$, $j\neq \tJ$, we have $\frac12(\tJ+j)\geq 2$, which implies that
$$
2|f_j^{\pm}(0)|\geqslant \tfrac12|\pm \tJ+j|\,.
$$
For  the case $j\leqslant-2,$ we have
$$
2 \,f_j^{+}(0)=  \tJ  + j\quad\hbox{and}\quad 2 \,f_j^{-}(0)= - \tJ  + j+2\,.
$$
Since $\tJ\geq 2$ and $j\leq -2$, $j\neq -\tJ$ we have $\frac12|-\tJ+j|=\frac12(\tJ+|j|)\geq 2$, which implies that
$$
2|f_j^{\pm}(0)|\geqslant \tfrac12|\pm \tJ+j|\,.
$$
Summing up, we deduce that
\begin{align}\label{fj-low}
\forall \,|j|\geqslant 2,\quad |f_j^{\pm}(0)|\geqslant \tfrac14|j\pm \tJ|\,.
\end{align}
Combining \eqref{fj-low} with \eqref{small.const.alpha} from Corollary \ref{cor.eigen.asympt}, we obtain for all $j \in \mathbb{Z}_\perp$, having $|j \pm \tJ| \geqslant 1,$
\begin{align}
    |f_j^\pm(\alpha)|& \geqslant |f_j^\pm(0)|-C\alpha_0^4|j\pm \tJ|-C\alpha_0^5\\
   & \geqslant |j\pm \tJ|\big(\tfrac14-2C\alpha_0^4\big) \,,
\end{align}
where  $C$ is a constant independent of $j$ and $\alpha.$ This proves \eqref{claim-1}  for $\alpha_0$ small enough.
  \\[1mm]
  Using \eqref{claim-1},  \eqref{ME.2} and \eqref{diff-Hm}, we estimate \eqref{smalldiv} with $\ell = \pm 1$  and $j\in \Z_{\perp}$ as
  \begin{align}
  	 \big| \pm \bomega_{\tJ,\varepsilon}(\alpha)  -\im\,\wt\mu_{j}^{(\tn)}(\alpha) \big| &  \geq \big|f_j^{\pm}(\alpha) \big| - \big| \wt\tm_{1,\varepsilon}^{(\tn)}(\alpha) - \big( \tfrac12 -\alpha^4 + \tc_{1}(\alpha) \big)  \big| |j| - |\wt\tw_{\tJ,\infty}(\alpha)| \\
  	 & \geq \tfrac18|j\pm \tJ| - C\varepsilon |j| -C\varepsilon \,.
  \end{align}
  Applying the triangle inequality $|j| \leqslant |j \pm \tJ| + \tJ$, together with the fact that $|j \pm \tJ| \geqslant 1$ since $|j|\neq \tJ$, we obtain
  \begin{align}
  	 \big| \pm \bomega_{\tJ,\varepsilon}(\alpha)  -\im\,\wt\mu_{j}^{(\tn)}(\alpha) \big| 
  	 & \geqslant \big(\tfrac18-C \varepsilon\big)|j\pm \tJ| -C\varepsilon \tJ-C\varepsilon\\
     &\geqslant  \big(\tfrac18-C \varepsilon-C\varepsilon \tJ-C\varepsilon\big)\geqslant  \tfrac{1}{16}\,,
  \end{align}
  where the last inequality holds for $\varepsilon>0$ small enough. This concludes the proof of item $(i)$.
  \\[1mm]
   {\sc Proof of $(ii)$}. For $|\ell |\geq 2$ and $j=\pm 2$, we have the following claim.
  \\[1mm]
  {\bf Claim 2.} {\it There exists $\alpha_{0}>0$ small enough such that}
  \begin{equation}\label{claim-2}
  	\big|  	-\im\big(  \lambda_{\tJ}^{(\infty)}(\alpha) \ell  \pm \lambda_{2}^{(\infty)} (\alpha) \big)   \big| \geq \tfrac14 |\ell|\,, \quad \forall \,|\ell| \geq 2 \,.
  \end{equation}
Let us prove this claim. 
  By \eqref{smalldiv.unpert}, we have
  	\begin{equation}
  		-\im\big(  \lambda_{\tJ}^{(\infty)}(\alpha) \ell  \pm \lambda_{2}^{(\infty)} (\alpha) \big)  
  		 = \tfrac12 \tJ \ell - \tfrac12 (\ell\pm 1) - \alpha^4 \big( \tJ \ell \pm\tfrac12 \big) + \tc_{1}(\alpha) \big( \tJ \ell \pm 2 \big) + \tr_{\tJ}(\alpha)\ell \pm \tr_{2}(\alpha) \,.
  	\end{equation}
  	We note that
  	\begin{align}
  		\big|  \tfrac{\tJ \ell}{2} - \tfrac12 (\ell\pm 1)\big| \geq \tfrac{|\ell|}{2}  \big( |\tJ|- \tfrac{|\ell \pm 1|}{|\ell|} \big) \geq \tfrac{|\ell|}{2}   \big( |\tJ| - 2 \big) \geq \tfrac{|\ell|}{2}  \,,
  	\end{align}
  	since we are assuming $\tJ \geq 3$, see \eqref{tangential.mode}. Moreover, we have that
  	\begin{equation}
  	|\tJ||\ell|\geqslant6\quad\hbox{and}\quad	\big|  \tJ \ell \pm \tfrac12  \big| \leqslant 2|\tJ| |\ell| \,, \quad \big|  \tJ \ell \pm 2 \big| \leqslant 2|\tJ| |\ell| \,.
  	\end{equation}
  	Consequently, we obtain
  	\begin{align}
  		\big|  	-\im\big(  \lambda_{\tJ}^{(\infty)}(\alpha) \ell  \pm \lambda_{2}^{(\infty)} (\alpha) \big)   \big| &\geq  \big|   \tfrac12 \tJ \ell - \tfrac12 (\ell\pm 1)  \big|  - \alpha^4 	\big|  \tJ \ell \pm \tfrac12  \big| \\& - |\tc_{1}(\alpha)| \big| \tJ\ell \pm 2 \big| - |\tr_{\tJ}(\alpha)| |\ell| - |\tr_{2}(\alpha)| \\
  		& \geq  |\ell| \big( \tfrac12 - 2 \alpha_{0}^4  |\tJ|  -C \alpha_{0}^5   \big) \\
  		& \geq \tfrac14 |\ell| \,,
  	\end{align}
  	for some $C>0$ independent of $\alpha$ and for $\alpha_{0}>0$ small enough. This proves \eqref{claim-2}.
    \\[1mm]
   Now, using Claim 2 and  \eqref{ME.2}, we estimate \eqref{smalldiv} for $|\ell| \geq 2$ and $j= \pm 2$ by
   	 \begin{align}
   		\big| \bomega_{\tJ,\varepsilon}(\alpha)  \ell \mp \im\,\wt\mu_{2}^{(\tn)}(\alpha) \big| &  \geq \big| -\im \big( \lambda_{\tJ}^{(\infty)}(\alpha)  \ell \pm \lambda_{2}^{(\infty)}(\alpha) \big) \big| \\
   		& \quad- 2 \big| \wt\tm_{1,\varepsilon}^{(\tn)}(\alpha) - \big( \tfrac12 -\alpha^4 + \tc_{1}(\alpha) \big)  \big|  - |\wt\tw_{\tJ,\infty}(\alpha)| |\ell|  \\
   		& \geq \big( \tfrac14 - 3C\varepsilon   \big)|\ell| \geq \tfrac18 |\ell| \,,
   	\end{align}
   	where the last inequality holds for $\varepsilon \ll 1$ small enough. This concludes the proof of item $(ii)$.
   	\\[1mm]
   	{\sc Proof of $(iii)$}. We have the following claim.
   	\\[1mm]
   	{\bf Claim 3.} {\it There exists $\alpha_{0}>0$ small enough such that}
   	\begin{equation}\label{claim-3}
   			\big|  	-\im\big(  \lambda_{\tJ}^{(\infty)}(\alpha) \ell  + \lambda_{j}^{(\infty)} (\alpha) \big)   \big| \geq  \tfrac{1}{17} |\tJ\ell+j|\,,
   	\end{equation}
   	{\it occurs for $|\ell |\geq 2$ and $j\in \Z_{\perp}$ with $|j|\geq3$ and $|\tJ \ell + j|\geq 	\tfrac52 |\ell|$.}
    \\[1mm]
   	Let us prove this claim.
   		By \eqref{smalldiv.unpert}, we have
   		\begin{align}
   			-\im\big(  \lambda_{\tJ}^{(\infty)}(\alpha) \ell & + \lambda_{j}^{(\infty)} (\alpha) \big)   = \big( \tfrac12 - \alpha^4 +\tc_{1}(\alpha)  \big)(\tJ \ell +j) - \tfrac12(\ell +{\rm sgn}(j)) + \tr_{\tJ}(\alpha)\ell + \tr_{j}(\alpha)\,.
   		\end{align}
   		Since we are assuming that $|\tJ \ell + j| \geq \tfrac{5}{2}|\ell|$, it follows, using also \eqref{small.const.alpha} in Corollary \ref{cor.eigen.asympt} and the fact that $|\ell\pm1|\leq \tfrac32 |\ell|$ for $|\ell|\geq 2$, that
        \begin{equation}
   			\begin{aligned}
   				\big|  	-\im\big(  \lambda_{\tJ}^{(\infty)}(\alpha) \ell & + \lambda_{j}^{(\infty)} (\alpha) \big)   \big| \geq  \big|\tfrac12 -\alpha^4 + \tc_{1}(\alpha) \big| |\tJ\ell + j| -  \tfrac12 |\ell+{\rm sgn}(j)|- |\tr_{\tJ}(\alpha)| |\ell| - |\tr_{j}(\alpha)| \\
   				& \geq  \tfrac{3}{8} |\tJ\ell + j| - \tfrac12 |\ell \pm 1| - C \alpha_{0}^5 |\ell| - C \alpha_{0}^5 \\
   				& \geq \tfrac{3}{8} |\tJ\ell + j| - \tfrac34 |\ell | - C \alpha_{0}^5 |\ell| - C \alpha_{0}^5 \\
   				& \geq |\tJ\ell + j|\big(\tfrac38-\tfrac{3}{10}-C\alpha_0^5  \big)\geq\tfrac{1}{17}|\tJ\ell + j|\,.
   			\end{aligned}
   		\end{equation}
   		for some  $C>0$ independent of $\ell,j,\alpha$ and for $\alpha_{0}$ sufficiently small. This ends the proof of \eqref{claim-3}.
        \\[1mm]
   	Using Claim 3, \eqref{ME.2}, and \eqref{diff-Hm}, we estimate \eqref{smalldiv}, for $|\ell|\geq 2$ and $j\in\Z_{\perp}$ with $|j|\geq 3$ and $|\tJ \ell + j|\geq \tfrac52 |\ell|$, by
   	\begin{align}
   		\big| \bomega_{\tJ,\varepsilon}(\alpha)  \ell + \im\,\wt\mu_{j}^{(\tn)}(\alpha) \big| &  \geq \big| -\im \big( \lambda_{\tJ}^{(\infty)}(\alpha)  \ell \pm \lambda_{j}^{(\infty)}(\alpha) \big) \big| \\
        &-  \big| \wt\tm_{1,\varepsilon}^{(\tn)}(\alpha) - \big( \tfrac12 -\alpha^4 + \tc_{1}(\alpha) \big)  \big| |j| - |\wt\tw_{\tJ,\infty}(\alpha)| |\ell|  \\
   		& \geq \big( \tfrac{1}{17} - C\varepsilon\big(\tfrac{2\tJ}{5}+1\big)  - \tfrac25 C\varepsilon  \big)|\tJ\ell+j| \geq \tfrac{1}{8}|\ell|  \,,
   	\end{align}
   	where the last inequality holds for $\varepsilon>0$ small enough, and where we have used in the last line the fact that
    $$
    |j|\leqslant \tJ|\ell|+|\tJ\ell+j|\leqslant \big(\tfrac{2\tJ}{5}+1\big)|\tJ\ell+j|\,.
    $$
    This concludes the proof of item $(iii)$.
    	\\[1mm]
   	{\sc Proof of $(iv)$}. When $\tJ=2$, by \eqref{smalldiv}, \eqref{smalldiv.unpert}, we have
    \begin{align}
        \bomega_{2,\varepsilon}(\alpha)\,\ell - \im \wt\mu_{j}^{(\tn)}(\alpha) & = \tfrac12 (\ell + j- {\rm sgn}(j)) + \big( -\alpha^4 + \tc_1(\alpha) \big) (2\ell + j) \\
        & + \big( \tfrac32 \alpha^4 + \tr_{2}(\alpha) +\wt\tw_{\tJ,\infty}(\alpha) \big)\ell   + \tr_{j}(\alpha) + \big( \wt\tm_{1,\varepsilon}^{(\tn)}(\alpha) - \big( \tfrac12 -\alpha^4 + \tc_{1}(\alpha) \big)  \big)j \,. 
    \end{align}
    If $|\ell|\geq \frac32 |j-{\rm sgn}(j)|> |j-{\rm sgn}(j)|$, noting also that $\tfrac12 |j| \leq |j-{\rm sgn}(j)|\leq 2|j|$,
    we estimate, by \eqref{small.const.alpha} in Corollary \ref{cor.eigen.asympt}, \eqref{ME.2} and \eqref{diff-Hm},
    \begin{align}
        |   \bomega_{2,\varepsilon}(\alpha)\,\ell - \im \wt\mu_{j}^{(\tn)}(\alpha) | & \geq \tfrac12 ( |\ell| - |j-{\rm sgn}(j)| ) - C \alpha_{0}^4 ( 2|\ell| +|j| ) - \big( C\alpha_0^4 + C\varepsilon \big)|\ell| - C \alpha_0^5 - C\varepsilon |j| \\
        & \geq \tfrac{7}{16} |\ell| - \tfrac{9}{16} |j-{\rm sgn}(j)| \geq \big( \tfrac{21}{32} - \tfrac{9}{16} \big) |j-{\rm sgn}(j)| \geq  \tfrac{3}{64} |j|\,,
    \end{align}
    assuming $\alpha_0$ and $\varepsilon_0$ small enough. The analogous estimate holds also in the case when $|j-{\rm sgn}(j)|\geq  \tfrac32 |\ell|>|\ell|$. This concludes the proof of item $(iv)$ and of the lemma.
\end{proof}

The next step is to provide a more precise description of the set in \eqref{tappo4}, based on Lemma \ref{lemma.empty.n.TR} and Lemma \ref{lemma.empty.n}.
First, we observe that, in the case $\tJ\geq 3$, if $|\tJ \ell + j| < \tfrac52 |\ell|$, then, by the triangle inequality, we infer that
\begin{equation}
	\tJ |\ell| -|j|  \leq |\tJ \ell + j | < \tfrac52 |\ell| \quad \Rightarrow \quad |\ell| < C_0 |j|\quad\hbox{with}\quad C_0:=\big( \tJ -\tfrac52 \big)^{-1}\,,
\end{equation}
which is justified as  $\tJ\geqslant 3$ by \eqref{tangential.mode}. Therefore,  Lemma \ref{lemma.empty.n.TR} and Lemma \ref{lemma.empty.n} imply that, when $\tJ\geq 3$, the set $\Omega_{\tn}\setminus\Omega_{\tn+1}$ in \eqref{tappo4} is actually given by
\begin{align}\label{compl.set.less}
	\Omega_{\tn}\setminus\Omega_{\tn+1}& =  \bigg(\bigcup_{ \ j\in \Z\backslash\{0\} \atop 1 \leqslant  |\ell| \leqslant C_0 |j| } R_{\ell,j}^{({\rm Tr},\tn)}\bigg)  \bigcup \bigg(\bigcup_{ \ j\in \Z_\perp   \atop 2\leqslant  |\ell| \leqslant C_0 |j| } R_{\ell,j}^{(1,\tn)}\bigg) \,,
\end{align}
and, when $\tJ= 2$, by
\begin{align}\label{compl.set.less.2}
	\Omega_{\tn}\setminus\Omega_{\tn+1}& \subseteq \bigg(\bigcup_{ j \in \Z\setminus\{0\} \atop \frac23|j| < |\ell| < \frac32 |j| } R_{\ell,j}^{({\rm Tr},\tn)}\bigg)  \bigcup \bigg(\bigcup_{  j\in \Z_\perp,  \atop \frac23|j-\rm{sgn}(j)|<  |\ell| < 3 |j| } R_{\ell,j}^{(1,\tn)}\bigg) \,.
\end{align}

The next goal is to provide a precise description of the sets introduced in \eqref{tappo1} and \eqref{tappo4}.
\begin{lem}\label{lemma.empty.shell.n}
    For $\tJ\geqslant 2$ and for any $\tn \geq 0$, the following inclusions hold.
    \begin{equation}
       (i): \   \Omega_{\tn}\setminus {\fE}_{\tn}^{{{0}}} \subseteq \bigcup_{ \ \ell,j\neq 0 \atop   |j| \geqslant N_{\tn-1} } R_{\ell,j}^{({\rm Tr},\tn)}\,; \quad \quad (ii): \ \fE_{\tn}^{(1)} \subseteq \bigcup_{\ell \neq0 , \
    j\in \Z_{\perp}  \atop  |j| \geq N_{\tn-1} } R_{\ell,j}^{(1,\tn)} \,. \label{tappo3}
    \end{equation}
    As a consequence, for $\tJ\geq 3$, we have
    \begin{align}\label{compl.set.less1}
	\Omega_{\tn}\setminus\Omega_{\tn+1}& \subseteq \bigg(\bigcup_{ \ |j|\geqslant N_{\tn-1} \atop 1 \leqslant  |\ell| \leqslant C_0 |j| } R_{\ell,j}^{({\rm Tr},\tn)}\bigg)  \bigcup \bigg(\bigcup_{ \ j\in \Z_\perp,\, |j|\geqslant N_{\tn-1}  \atop 2\leqslant  |\ell| \leqslant C_0 |j| } R_{\ell,j}^{(1,\tn)}\bigg) \,,
\end{align}
and, for $\tJ=2$,
\begin{align}\label{compl.set.less-J2}
	\Omega_{\tn}\setminus\Omega_{\tn+1}& \subseteq \bigg(\bigcup_{ \ |j|\geqslant N_{\tn-1} \atop \frac23|j| < |\ell| < \frac32 |j| } R_{\ell,j}^{({\rm Tr},\tn)}\bigg)  \bigcup \bigg(\bigcup_{ \ j\in \Z_\perp,\, |j|\geqslant N_{\tn-1}  \atop \frac23|j-\rm{sgn}(j)|<  |\ell| < 3 |j| } R_{\ell,j}^{(1,\tn)}\bigg) \,.
\end{align}
\end{lem}
\begin{proof}
    We begin with the proof of item $(i)$. For $\tn = 0$, we have, in view of \eqref{constants.NASH}, that $N_{-1} = 1$: the result then follows directly from \eqref{set.n.to.measure}. Now let $\tn \geq 1$ and $\alpha\in \Omega_{\tn}$. By \eqref{set.n.to.measure}, \eqref{nonres.step.NM} and \eqref{nonres.set.inv.full}, we have that $(\b\omega_{\tJ,\varepsilon}(\alpha),\alpha) \in \tT\tC(\upsilon_{\tn-1},\tau,\wt\fu_{\perp}^{(\tn-1)})
    $,  which means that it satisfies (recall \eqref{nonres.transport})
    \begin{equation}\label{tappo2}
        \big| \b\omega_{\tJ,\varepsilon}(\alpha)\ell +\wt\tm_{1,\varepsilon}^{(\tn-1)}(\alpha)j \big| \geqslant \tfrac{4\upsilon_{\tn-1}}{|j|^{\tau}} \,, \quad \forall\, \ell,j \in\Z\setminus\{0\} \,.
    \end{equation}
    Therefore, for any $\ell\in \Z\setminus\{0\}$ and for any $j\in\Z\setminus\{0\}$ such that $|j|\leq N_{\tn-1}$, we compute, by \eqref{tappo2}, \eqref{beta.eps.est12} in Lemma \ref{almost.straight.lemma}-$(i)$ and \eqref{hn} in Proposition \ref{NM.iteration},
    \begin{align}
        \big| \b\omega_{\tJ,\varepsilon}(\alpha)\ell +\wt\tm_{1,\varepsilon}^{(\tn)}(\alpha)j \big| &\geq \big| \b\omega_{\tJ,\varepsilon}(\alpha)\ell +\wt\tm_{1,\varepsilon}^{(\tn-1)}(\alpha)j \big| - \big| \wt\tm_{1,\varepsilon}^{(\tn)}(\alpha) - \wt\tm_{1,\varepsilon}^{(\tn-1)}(\alpha) \big| |j| \\
        & \geq \tfrac{4\upsilon_{\tn-1}}{|j|^{\tau}} - C \varepsilon \| \wt\fu_{\perp}^{(\tn)}-\wt\fu_{\perp}^{(\tn-1)} \|_{s_0+\bar\sigma}^{5,\upsilon} |j| \\
        & \geq \tfrac{4\upsilon_{\tn-1}}{|j|^{\tau}} - C C_{*} \varepsilon N_{\tn-1}^{2\bar\sigma+1} N_{\tn-2}^{-\ta_{1}} |j| \geq \tfrac{4\upsilon_{\tn}}{|j|^{\tau}} \,,
    \end{align}
   provided  that
    \begin{equation}
         C C_{*} \varepsilon N_{\tn-1}^{2\bar\sigma+1} N_{\tn-2}^{-\ta_{1}} |j| 
         \leq \tfrac{4\,2^{-\tn}}{|j|^{\tau}}\,, 
    \end{equation}
    which holds for any $|j|\leq N_{\tn-1}$ by \eqref{constants.NASH}, assuming that  $N_{0}$ taken  sufficiently large, independently of $\tn$. Therefore, we deduce that  $R_{\ell,j}^{({\rm Tr},\tn)}=\emptyset$ for any $\ell\in \Z \setminus \{0\}$ and $1\leqslant |j|\leq N_{\tn-1}$. \\
    The proof of the inclusion in $(ii)$ follows by similar arguments, recalling \eqref{ME.1.1}, \eqref{tappo4}, \eqref{tappo5}, and therefore is omitted. \\The  claims \eqref{compl.set.less1} and \eqref{compl.set.less-J2} follow from the foregoing results together with \eqref{compl.set.less} and \eqref{compl.set.less.2}.
    This ends the proof of the lemma.
\end{proof}
\subsection{Transversality and proof of Theorem \ref{meas.est.thm}}
To estimate the measure of each set in \eqref{compl.set.less}, and achieve the proof of Theorem \ref{meas.est.thm}, we need to control the associated small divisors by quantifying their transversality with respect to the parameter $\alpha$. 
Relying on a result of R\"ussmann \cite{Russ}, we show that these quantities satisfy a uniform non-degeneracy condition, in the sense that a finite number of $\alpha$-derivatives cannot vanish simultaneously, with a lower bound proportional to $|j|$. This property is crucial for the forthcoming measure estimates.
\\
The precise statement is given below.
\begin{lem}{\bf (Transversality).}\label{transversality.lemma}
	There exists $\varepsilon_0>0$ small enough such that, for any  $\varepsilon\in[0,\varepsilon_{0})$  and for all $\tn\geq 0$, $\alpha\in [\alpha_{1},\alpha_{2}]$,  the following hold:
    \\[1mm]
    \noindent $(i)$ When $\tJ\geq 3$, we have
		\begin{align}
             \big| \pa_{\alpha}^4 \big( \bomega_{\tJ,\varepsilon}(\alpha)\,\ell + \wt\tm_{1,\varepsilon}^{(\tn)}(\alpha) j  \big| \geq \tfrac{1}{16\tJ}  |j|
			\,, \quad  \forall\, 1\leqslant |\ell|\leqslant C_0 |j|\,, \quad j \in \Z\setminus\{0\}\,, 
            \end{align}
            and
            \begin{align}
            \big| \pa_{\alpha}^4 \big( \bomega_{\tJ,\varepsilon}(\alpha)\,\ell -\im\,\wt\mu_{j}^{(\tn)}(\alpha) \big)  \big| \geq \tfrac{1}{32\tJ}   |j|
			\,, \quad\forall\,2\leqslant |\ell|\leqslant C_0 |j|\,,\quad  j\in \Z_{\perp} \,, 
		\end{align}
        with $C_0= \big( \tJ -\tfrac52 \big)^{-1}$;
        \\[1mm]
        \noindent $(ii)$ When $\tJ=2$, we have
        \begin{align}
			   \big| \pa_{\alpha}^4 \big( \bomega_{2,\varepsilon}(\alpha)\,\ell + \wt\tm_{1,\varepsilon}^{(\tn)}(\alpha) j  \big| \geq \tfrac32  |j|
			\,, \quad   \forall\, \tfrac23|j|< |\ell|< \tfrac32 |j|\,, \quad j \in \Z \setminus \{0\} \,,  
            \end{align}
            and
            \begin{align}
			  \big| \pa_{\alpha}^4 \big( \bomega_{2,\varepsilon}(\alpha)\,\ell -\im\,\wt\mu_{j}^{(\tn)}(\alpha) \big)  \big| \geq  \tfrac38  |j|
			\,, \quad  \forall \,2 \leq \tfrac23|j-{\rm{sgn}}(j)|< |\ell|< \tfrac32 |j-{\rm sgn}(j)|\,,\quad |j|\geqslant3 \,. 
		\end{align}
\end{lem}
\begin{proof}
We prove only the second estimates of both items $(i)$ and $(ii)$, as the first ones follow by similar arguments and are therefore omitted. Since we are interested in estimating the sets $ R_{\ell,n}^{(1,\tn)} $ in \eqref{compl.set.less}, Lemma \ref{lemma.empty.shell.n} shows that only the regimes on $j$ and $\ell$ stated therein are nontrivial.
\\[1mm]
\noindent $\blacktriangleright$ {\sc Proof of $(i)$.}
	We split the proof into two steps. In the first part, we show the claim in the unperturbed case, namely $\varepsilon=0$, then we will extend it to the perturbed one.
     We introduce the function
\begin{equation}\label{fj-alpha}
     \tf_{\varepsilon,\ell,j}(\alpha) := \bomega_{\tJ,\varepsilon}(\alpha) \ell -\im \, \wt\mu_{j}^{(\tn)}(\alpha) \,.
	\end{equation}
\noindent {\bf Step 1.}	
We first show that for all $\alpha\in [\alpha_{1},\alpha_{2}]$, we have
    \begin{equation}\label{transv.unpert}
        \big| \pa_{\alpha}^4 \tf_{0,\ell,j} \big| \geqslant \tfrac{1}{16\tJ}  |j|
		\,, \quad \forall\, 2\leqslant |\ell|\leqslant C_0 |j|\,, \quad j \in \Z_{\perp} \,. 
	\end{equation}
    For $\varepsilon=0,$ we have, recalling Corollary \ref{cor.eigen.asympt},
	\begin{align}
			\tf_{0,\ell,j}(\alpha) & = -\im\big( \lambda_{\tJ}^{(\infty)}(\alpha)\,\ell + \lambda_{j}^{(\infty)}(\alpha)\big) \\
			& =	\big(\tfrac12 - \alpha^4+\tc_{1}(\alpha) \big)(\tJ \ell + j) - \tfrac12 (\ell +{\rm sgn}(j) ) +  \tr_{\tJ}(\alpha) \ell + \tr_{j}(\alpha)		\,.\label{f-o0}
	\end{align}
	 We note that
     \begin{align}\label{f-000}
			\tf_{0,\ell,j}(0) 
			& =	\tfrac12 (\tJ \ell + j) - \tfrac12 (\ell +{\rm sgn}(j))
	\end{align}
    and 
     \begin{equation}
		\pa_{\alpha}^4\tf_{0,\ell,j}(\alpha) =\big(- 24 + \pa_{\alpha}^4\tc_{1}(\alpha) \big)(\tJ \ell +j) +\pa_{\alpha}^4 \tr_{\tJ}(\alpha) \ell + \pa_{\alpha}^4 \tr_{j}(\alpha)\,.
	\end{equation}
     Let $\delta>0$ that will be fixed later. We discuss the following three cases:
     \begin{equation}\label{abc.cases}
          ({\rm a})\quad  \tJ\ell+j=0\,,\quad ({\rm b}) \quad \, 1\leqslant |\tJ\ell+j|<\delta |j|\,, \quad ({\rm c})\quad |\tJ\ell+j|\geqslant\delta |j|\,.
     \end{equation}
     In the first case $({\rm a})$, having $|\ell\pm 1|\geq \tfrac12$ for $|\ell|\geqslant2$, we deduce
     \begin{align}
			|\tf_{0,\ell,j}(0)| 
			& =	  \tfrac12 |\ell +{\rm sgn}(j)|  \geq \tfrac14 |\ell| = \tfrac{1}{4 \tJ}|j|\,. \label{Ineq-tran1}
	\end{align}
    In the second case $({\rm b})$, by \eqref{f-000}, we have
    \begin{align*}
			|\tf_{0,\ell,j}(0)| 
			& \geqslant \tfrac12\big(|\ell|-1-\delta |j|\big)\,.
	\end{align*}
    From the triangle inequality, we infer
    \begin{equation}
        (1+\delta)|j| \geq |j| + |\tJ \ell + j|  \geqslant \tJ|\ell|\geqslant |j|-|\tJ\ell+j|\geqslant (1-\delta)|j| \,,
    \end{equation}
    implying that
    \begin{align}\label{ell-to-j}
    \tfrac{1+\delta}{\tJ}|j|\geqslant |\ell|\geqslant \tfrac{1-\delta}{\tJ} |j|\,.\end{align}
    Therefore, we have
\begin{align*}
			|\tf_{0,\ell,j}(0)| 
			& \geqslant \tfrac12\Big(\tfrac{1-\delta(1+\tJ)}{\tJ}|j|-1\Big)\,.	
	\end{align*}   
    Now, if we choose $\delta>0$ such that $ \delta\leqslant \tfrac{1}{2(1+\tJ)}$,
    then
    \begin{align*}
			|\tf_{0,\ell,j}(0)| 
			& \geqslant \tfrac12\big(\tfrac{1}{2\tJ}|j|-1\big)\,.	
	\end{align*}  
    When $|j|\geqslant 4\tJ$,  we obtain
    \begin{align}
			|\tf_{0,\ell,j}(0)| 
			& \geqslant \tfrac12\big(\tfrac{1}{2\tJ}-\tfrac{1}{4\tJ}\big) \geqslant \tfrac12\big(2|j|-\tfrac{|j|}{4\tJ}\big) \geqslant\tfrac{1}{8\tJ}|j| \,.
	\end{align}
    On the other hand, for $|j|\leqslant 4\tJ$, one should get from the constraint $({\rm b})$
   \begin{align}\label{ine-X1}
   1\leqslant \delta |j|\leqslant 4\delta \tJ\,.
   \end{align}
  By  fixing  $\delta$ as below
   \begin{align}\label{delta-fix}
   \delta:=\tfrac{1}{8\tJ} \leq \tfrac{1}{2(1+\tJ)} \,,
   \end{align}
   we obtain that \eqref{ine-X1} is impossible. Thus in this regime, we necessarily get
   \begin{align}
			|\tf_{0,\ell,j}(0)| 
    &\geqslant \tfrac{1}{8\tJ}|j|\,. \label{Ineq-tran2}
	\end{align} 
    We note that, in both cases $({\rm a})$ and $({\rm b})$ we obtain from \eqref{Ineq-tran1} and \eqref{Ineq-tran2}
    \begin{align}
			|\tf_{0,\ell,j}(0)| 
    &\geqslant \tfrac{1}{8\tJ}|j|\,. \label{Ineq-tran4}
	\end{align}
    Coming back to \eqref{f-o0},
    by the triangle inequality, Corollary \ref{cor.eigen.asympt} and \eqref{ell-to-j}, we find that, in both cases $({\rm a})$ and $({\rm b})$ in \eqref{abc.cases},
    \begin{align}
			|\tf_{0,\ell,j}(\alpha)| &  \geqslant |\tf_{0,\ell,j}(0)|- |-\alpha^4 +\tc_{1}(\alpha)||\tJ\ell+j| +|\tr_{\tJ}(\alpha)||\ell| + |\tr_{j}(\alpha)| \\ 
            & \geq \tfrac{1}{8\tJ}|j| - \delta C\alpha_{0}^4 |j| - \tfrac{1+\delta}{\tJ} C\alpha_{0}^5 |j| - C \alpha_{0}^5 \\
            &\geqslant\big(\tfrac{1}{8\tJ}- \tfrac{1}{\tJ}C\alpha_{0}^4\big)|j| \geqslant\tfrac{1}{16\tJ}|j| \,. \label{Ineq-tran4}
	\end{align}
    We now consider the last case $({\rm c})$ in \eqref{abc.cases}. By the triangle inequality and Corollary \ref{cor.eigen.asympt},
  we  find  that, since  $|\tJ\ell+j|\geq \delta |j|$ and $2\leqslant |\ell|\leqslant C_0|j|$,
    \begin{align*}
		|\pa_{\alpha}^4\tf_{0,\ell,j}(\alpha)| &\geqslant ( 24 \delta-C\alpha_{0}) |j|-C\alpha_{0} |\ell|-C\alpha_{0} \\
         &\geqslant \big( 24 \delta-3C C_0\alpha_{0}\big) |j|\,.
\end{align*}
Thus, for $0<\alpha_{0}\ll 1$ small enough, we deduce, using \eqref{delta-fix}, that
\begin{align}\label{Ineq-tran3}
		|\pa_{\alpha}^4\tf_{0,\ell,j}(\alpha)| 
         &\geqslant 12 \delta |j|=\tfrac{3}{2\tJ}|j|\,.
\end{align}
 Putting together \eqref{Ineq-tran4} and \eqref{Ineq-tran3}  yields
 \begin{align}
        |\pa_{\alpha}^4 \tf_{0,\ell,j}(\alpha)| 
         \geqslant \tfrac{1}{16\tJ}|j|\,.
\end{align}
This concludes the proof of the claimed estimate in \eqref{transv.unpert}.
\\[1mm]
     \noindent {\bf Step 2.}	
     By \eqref{fj-alpha}, \eqref{ME.1.1}, \eqref{transv.unpert}, \eqref{ME.2} and \eqref{diff-Hm}, we get, in view of $2\leqslant|\ell|\leqslant C_0|j|,$ 
	\begin{align}
            |\pa_{\alpha}^4\tf_{\varepsilon,\ell,j}(\alpha)|& \geq 
            |\pa_{\alpha}^4 \tf_{0,\ell,j}(\alpha)| - 
            |\pa_{\alpha}^4 \wt\tw_{\tJ,\infty}(\alpha)| |\ell| \\&-
            \big| \pa_{\alpha}^4 \big( \wt\tm_{1,\varepsilon}(\alpha) - \tfrac12 + \alpha^4 - \tc_{1}(\alpha) \big) \big| |j| \\
			& \geq \big(\tfrac{1}{16\tJ}   - C C_0 \varepsilon \upsilon^{-4}  - C \varepsilon \upsilon^{-4}  \big) |j|\,.
		\end{align}
	Therefore, for  $0< \varepsilon_{0} \ll 1$ small enough, we obtain 
    $$
    |\pa_{\alpha}^4\tf_{\varepsilon,\ell,j}(\alpha)|\geq\tfrac{1}{32 \tJ}|j| \,.
    $$
    This concludes the proof of the result of item $(i)$ in the case $\tJ\geqslant3.$
    \\[1mm]
    \noindent $\blacktriangleright$ {\sc Proof of $(ii)$.}
We now consider the case $\tJ=2$. As before, we split the proof in two steps. Let
\begin{equation}\label{fj-alpha-2}
     \tf_{\varepsilon,\ell,j}(\alpha) := \bomega_{2,\varepsilon}(\alpha) \ell -\im \, \wt\mu_{j}^{(\tn)}(\alpha) \,.
	\end{equation}
\noindent {\bf Step 1.}	
We first show that for all $\alpha\in [\alpha_{1},\alpha_{2}]$, we have
    \begin{equation}\label{transv.unpert-2}
        \big| \pa_{\alpha}^4 \tf_{0,\ell,j} \big| \geqslant \tfrac34  |j|
		\,, \quad  \forall \,2 \leq \tfrac23|j-{\rm{sgn}}(j)|< |\ell|< \tfrac32 |j-{\rm sgn}(j)|\,,\quad |j|\geqslant3 \,. 
	\end{equation}
For $\varepsilon=0$, we have, recalling Corollary \ref{cor.eigen.asympt},
	\begin{align}
			\tf_{0,\ell,j}(\alpha) & = -\im\big( \lambda_2^{(\infty)}(\alpha)\,\ell + \lambda_{j}^{(\infty)}(\alpha)\big) \\
			& =	\big(\tfrac12 - \alpha^4+\tc_{1}(\alpha) \big)(\tJ \ell + j) - \tfrac12 (\ell +{\rm sgn}(j) ) +\tfrac32 \alpha^4 \ell +  \tr_{\tJ}(\alpha) \ell + \tr_{j}(\alpha)		\,, \label{f-ooo}
	\end{align}
and that
\begin{equation}\label{der4.ooo}
	\pa_{\alpha}^4\tf_{0,\ell,j}(\alpha)
	=
	\big( -24 + \pa_{\alpha}^4\tc_{1}(\alpha) \big)(2\ell+j)
	+ \tfrac32 \cdot 24 \ell+\pa_{\alpha}^4 \tr_{\tJ}(\alpha)\,\ell
	+\pa_{\alpha}^4 \tr_{j}(\alpha)\,.
\end{equation}
We split the discussion into two cases: $|j|\geqslant 4$ and $|j|=3$. In the first case, the triangle inequality and Corollary \ref{cor.eigen.asympt} give
\begin{align}
    |\pa_{\alpha}^4\tf_{0,\ell,j}(\alpha)| \geq \big| - 12 \ell - 24 j \big| + C\alpha_{0} |\ell| + C \alpha_0 \,.
\end{align}
Therefore, having for any $|\ell|\geq 2$ and $|j|\geq 4$ that
\begin{equation}\label{bound.ooo}
   \tfrac23(|j|-1)<\tfrac23|j-\operatorname{sgn}(j)|< |\ell|< \tfrac32 |j-{\rm sgn}(j)| < \tfrac32(|j|+1) \,,
\end{equation}
we infer that, using also that $|j|-3 \geq \tfrac14 |j|$ for $|j|\geq 4$,
\begin{align*}
	|\pa_{\alpha}^4\tf_{0,\ell,j}(\alpha)|
	&\geqslant 12\big( 2|j|- |\ell|\big) - C\alpha_0|\ell| \\
	&\geqslant 12 \big( 2|j| - \tfrac32 |j| - \tfrac32 \big) - 3 C\alpha_0 |j| \\
    & \geq 6 (|j|-3) - 3C \alpha_0 |j| \\
    & \geq \tfrac32 |j| - 3C \alpha_0 |j| \geq \tfrac34 |j| \,,
\end{align*}
for $\alpha_0$ small enough. On the other hand, when $|j|=3$, we note that
\begin{equation}
    \ell + 2 j = 0 \,, \quad  j= \pm 3 \quad \Rightarrow \quad \ell = \mp 6 \,.
 \end{equation}
 But
 \begin{equation}
    \tfrac43 =\tfrac23 |3 - 1| < |\ell| < \tfrac32 (3-1)= 3  \quad \Rightarrow \quad |\ell| = 2 \,.
 \end{equation}
 We deduce that $\ell + 2j\neq 0$ in this regime. Consequently, we simply get from \eqref{der4.ooo} that, for $|j|=3$ and $|\ell|=2$,
 \begin{align}
     |\pa_{\alpha}^4\tf_{0,\ell,j}(\alpha)| & \geq 12 |\ell + 2 j| - C\alpha_0 \geq 6 = 2|j| \,.
 \end{align}
 This concludes the proof of the claimed estimate in  \eqref{transv.unpert-2}.
 \\[1mm]
     \noindent {\bf Step 2.}	
     By \eqref{fj-alpha}, \eqref{ME.1.1}, \eqref{transv.unpert-2}, \eqref{ME.2} and \eqref{diff-Hm}, we get, in view of \eqref{bound.ooo},
	\begin{align}
            |\pa_{\alpha}^4\tf_{\varepsilon,\ell,j}(\alpha)|& \geq 
            |\pa_{\alpha}^k \tf_{0,\ell,j}(\alpha)| - 
            |\pa_{\alpha}^4 \wt\tw_{2,\infty}(\alpha)| |\ell| \\&-
            \big| \pa_{\alpha}^4 \big( \wt\tm_{1,\varepsilon}(\alpha) - \tfrac12 + \alpha^4 - \tc_{1}(\alpha) \big) \big| |j| \\
			& \geq \big( \tfrac34   - C C_0 \varepsilon \upsilon^{-4}  - C \varepsilon \upsilon^{-4}  \big) |j|\,.
		\end{align}
	Therefore, for  $0< \varepsilon_{0} \ll 1$ small enough, we obtain 
    $$
    |\pa_{\alpha}^4\tf_{\varepsilon,\ell,j}(\alpha)|\geq \tfrac38 |j| \,.
    $$
    This concludes the proof of the result of item $(ii)$ in the case $\tJ= 2$ and the proof of the lemma.
\end{proof}

\begin{lem}\label{est.Relln}{\bf (Estimates of the resonant sets).} Under the assumptions on $\ell,j$ detailed in Lemma \ref{transversality.lemma}, we get
	\begin{equation}
	|R_{\ell,j}^{({\rm Tr},\tn)}| \leqslant C (\upsilon |j| ^{-(\tau+1)})^{\frac14}\quad\hbox{and}\quad   |R_{\ell,j}^{(1,\tn)}| \leqslant C (\upsilon |j| ^{-(\tau+1)})^{\frac14} \,,
	\end{equation}
	where the sets   $R_{\ell,j}^{({\rm Tr},\tn)}$ and  $R_{\ell,j}^{(1,\tn)}$ are defined in \eqref{tappo1} and \eqref{tappo5}, respectively.
\end{lem}
\begin{proof}
	We shall give the proof for $\tJ\geqslant3$ and the remaining case $\tJ=2$ is similar. We introduce the function
	\begin{equation}
		\wt\tf_{\varepsilon,\ell,j}(\alpha) := \tfrac{1}{|j|} \tf_{\varepsilon,\ell,j}(\alpha) := \tfrac{1}{|j|} \big( \bomega_{\tJ,\varepsilon}(\alpha) \ell -\im \, \wt\mu_{j}^{(\tn)}(\alpha) \big) \,.
	\end{equation}
	Then, we rewrite the set $R_{\ell,j}^{(1,\tn)}$ in \eqref{tappo5} as
	\begin{equation}
		R_{\ell,j}^{(1,\tn)} = \big\{ \alpha\in[\alpha_{1},\alpha_{2}]  \, : \, |\wt\tf_{\varepsilon,\ell,j}(\alpha)| < 2\upsilon |j|^{-\tau-1} \big\} \,.
	\end{equation}
	By Lemma \ref{transversality.lemma}, we know that
	\begin{equation}
		\max_{0\leq k \leq 4} |\pa_{\alpha}^k \wt\tf_{\varepsilon,\ell,j}( \alpha) | \geq  |\pa_{\alpha}^4 \wt\tf_{\varepsilon,\ell,j}( \alpha) |  \geq \tfrac{1}{32 \tJ} \,.
	\end{equation}
	Moreover, as $2\leqslant |\ell|\leqslant |j|$ and $\varepsilon \upsilon^{-5}$ is small enough, we get from straightforward computations that the following uniform estimate 
	$$ 	\max_{0\leq k \leq 5} |\pa_{\alpha}^k \wt\tf_{\varepsilon,\ell,j}( \alpha) |  \leq C\,, \quad \forall \, \alpha\in [\alpha_{1},\alpha_{2}]\,.
    $$
	Therefore  \cite[Theorem 17.1]{Russ}  applies and the claim follows. The proof for $|R_{\ell,j}^{({\rm Tr},\tn)}|$ uses similar arguments and therefore it is omitted.
\end{proof}

\begin{proof}[Proof of Theorem \ref{meas.est.thm} completed.]
 We shall carry out the proof for $\tJ\geqslant 3$: the case $\tJ=2$ can be treated in a similar way and it is therefore omitted.  We now estimate the measure of the set $\Omega_{\varepsilon}$ in \eqref{set.to.measure}. By \eqref{compl.set.to.measure}, we have to estimate the measure of the set $\Omega_{\tn}\setminus\Omega_{\tn+1}$ for any $\tn\geq 0$. By \eqref{tappo4}, we have
    \begin{equation}\label{tappo5}
        |\Omega_{\tn}\setminus\Omega_{\tn+1}| \leq |\Omega_{\tn}\setminus \fE_{\tn}^{(0)}| + |\fE_{\tn}^{(1)}| \,.
    \end{equation}
    Then, applying Lemma \ref{lemma.empty.n}, \eqref{compl.set.less}, Lemma \ref{lemma.empty.shell.n} and Lemma \ref{est.Relln}, we have that
    \begin{align*}
    |\Omega_{\tn}\setminus \fE_{\tn}^{(0)}|+|\fE_{\tn}^{(1)}| & \leq  \sum_{1\leq |\ell|\leq C_{0}|j| \atop 
        |j|\geq N_{\tn-1}} |R_{\ell,j}^{({\rm Tr},\tn)}|+ \sum_{2\leq |\ell|\leq C_{0}|j| \atop 
       j\in\Z_\perp , \ |j|\geq N_{\tn-1}} |R_{\ell,j}^{(1,\tn)}| \\
        &\lesssim \upsilon^{\frac14} \sum_{|j|\geq N_{\tn-1}}  |j|^{-\frac{\tau+1}{4}+1}\,.
        \end{align*}        According to  \eqref{param.meas.est}, we have ${\tau>7}$   and therefore  the series   $ \sum_{j\in \Z\setminus\{0\}}  |j|^{-\frac{\tau-3}{4}}$ is  convergent.\\ By \eqref{tappo5}, it implies that, for any $\tn\geq 0$,
    \begin{equation}
     |\Omega_{\tn}\setminus\Omega_{\tn+1}| \lesssim \upsilon^{\frac14} N_{\tn-1}^{\frac{7-\tau}{4}} \,.
    \end{equation}
    Therefore, recalling \eqref{compl.set.to.measure}, 
    \begin{equation}
        |[\alpha_{1},\alpha_{2}]\setminus\Omega_{\varepsilon}| \leq \sum
        _{\tn\geq 0} |\Omega_{\tn}\setminus\Omega_{\tn+1}|  \lesssim \upsilon^{\frac14} \sum
        _{\tn\geq 0} N_{\tn-1}^{\frac{7-\tau}{4}} \lesssim \upsilon^{\frac14}\,,
    \end{equation}
    since the geometric series $\sum_{\tn\geq 0} N_{\tn-1}^{\frac{7-\tau}{4}}$ is convergent. The claimed estimate then follows using that $\upsilon=\varepsilon^{\tc}$ as in \eqref{param.meas.est}. The proof of the theorem is concluded.
 \end{proof}

\section*{ Acknowledgements.}  R. Montalto  
is supported by the ERC STARTING GRANT 2021 ``Hamiltonian Dynamics, 
Normal Forms and Water Waves'' (HamDyWWa), Project Number: 101039762. 
Views and opinions expressed are however those of the authors only and do not necessarily reflect those of the European Union or the European Research Council. Neither the European Union nor the granting authority can be held responsible for them. \\
L. Franzoi is supported by ``GNAMPA - INdAM'', CUP E53C25002010001 and he has been supported by the ERC STARTING GRANT 2021 ``Hamiltonian Dynamics, 
Normal Forms and Water Waves'' (HamDyWWa), Project Number: 101039762 until January 2026. \\ T. Hmidi is supported by Tamkeen under the NYU Abu Dhabi Research Institute grant.  

 \section*{Statements and declarations}

\medskip

\noindent
 All the authors state that there is no conflict of interest and certify that they have no affiliations or involvement with any organization or entity with any
financial interest, or non-financial interest in the subject matter or materials discussed in this manuscript.

Moreover, data sharing is not applicable to this article as no datasets were generated or analyzed during the
current study.

\appendix

\section{Useful trigonometric identities}\label{app.trigo}

For the convenience of the reader, we collect in this appendix elementary trigonometric identities that we used in several computations of the paper.
\begin{enumerate}
	\item $\cos(k(\theta+\eta))\sin(\eta) = \tfrac12 \big( \sin(k\theta+(k+1)\eta) - \sin(k\theta+(k-1)\eta) \big)$;
	\item $\cos^2(\theta) = \tfrac12\big( \cos(2\theta)+1 \big)$;
	\item We have that
	\begin{equation}
		\big( \cos(2\theta)-\cos(2\eta) \big) {\rm cotg} \big( \tfrac{\theta-\eta}{2} \big) = - 2 \sin(\theta+\eta) \big( \cos(\theta-\eta)+1 \big) \,.
	\end{equation}
	It follows by computing
	\begin{equation}
		\begin{aligned}
				\big( \cos(2\theta)&-\cos(2\eta) \big) {\rm cotg} \big( \tfrac{\theta-\eta}{2} \big)  = -2 \sin(\theta+\eta) \sin(\theta-\eta) {\rm cotg} \big( \tfrac{\theta-\eta}{2} \big) \\
				& = - 4\sin(\theta+\eta) \cos \big( \tfrac{\theta-\eta}{2} \big)  \sin \big( \tfrac{\theta-\eta}{2} \big)  {\rm cotg} \big( \tfrac{\theta-\eta}{2} \big)\\
				& =- 4\sin(\theta+\eta)  \cos^2 \big( \tfrac{\theta-\eta}{2} \big) \\
				& = - 2 \sin(\theta+\eta) \big( \cos(\theta-\eta)+1 \big) ;
		\end{aligned}
	\end{equation}
	\item $\sin(k\eta)- \sin(k\theta) = 2 \cos\big(\tfrac{k}{2}(\eta+\theta)\big) \sin\big(\tfrac{k}{2}(\eta-\theta)\big) $ for any $k\in\N$;
	\item
	For $k \in \N$, we have
	\begin{equation}\label{sink}
		\begin{aligned}
			\sin(k \alpha) 
			&=\begin{cases}
				\sin(\alpha) & k=1 \,, \\
				2 \sin(\alpha) \sum_{m=0}^{ \frac{k}{2} -1 }\cos((k-1-2m)\alpha) & k \ {\rm even} \,, \\
				2 \sin(\alpha) \sum_{m=0}^{ \frac{k-3}{2}  }\big( \cos((k-1-2m)\alpha) +\tfrac12 \big)& k \ {\rm odd} \,.
			\end{cases}
		\end{aligned}
	\end{equation}
	It follows by computing
	\begin{align}
	\sin(k\alpha) & = \frac{1}{2\im}e^{\im k \alpha} \big( 1- e^{-2\im k\alpha} \big)  = \frac{1}{2\im} e^{\im k\alpha} (1-e^{2\im\alpha}) \sum_{m=0}^{k-1} e^{-2\im m\alpha} \\
	& = \sin(\alpha)  \sum_{m=0}^{k-1} e^{\im(k-1-2m)\alpha} \\
	& = \begin{cases}
		\sin(\alpha) & k=1 \,, \\
		2 \sin(\alpha) \sum_{m=0}^{\lfloor \frac{k}{2} \rfloor -1 }\cos((k-1-2m)\alpha) & k \ {\rm even} \,, \\
		2 \sin(\alpha) \sum_{m=0}^{\lfloor \frac{k}{2} \rfloor -1}\big( \cos((k-1-2m)\alpha) +\tfrac12 \big)& k \ {\rm odd} \,;
	\end{cases}
	\end{align}
	\item 
Combining the previous two results, we have that
	\begin{equation}
		\begin{aligned}
			&	(\sin(k\eta)-\sin(k\vartheta)) {\rm cotg} \big( \tfrac{\eta-\vartheta}{2} \big) = 2 \cos\big(\tfrac{k}{2}(\eta+ \vartheta)\big) \sin\big(\tfrac{k}{2}(\eta-\vartheta)\big) {\rm cotg} \big(\tfrac12 (\eta-\vartheta)\big) \\
			&\quad  = \begin{cases}
				2 \cos\big(\tfrac{1}{2}(\eta+ \vartheta)\big) \cos\big(\tfrac{1}{2}(\eta- \vartheta)\big)  & k=1 \,, \\
				4 \cos\big(\tfrac{k}{2}(\eta+ \vartheta)\big)\cos\big(\tfrac{1}{2}(\eta- \vartheta)\big)\sum_{m=0}^{\frac{k}{2}-1} \cos\big(\tfrac{k-1-2m}{2}(\eta- \vartheta)\big) & k \ {\rm even} \,,\\
				4 \cos\big(\tfrac{k}{2}(\eta+ \vartheta)\big) \cos\big(\tfrac{1}{2}(\eta- \vartheta)\big)\Big( \sum_{m=0}^{\frac{k-3}{2}} \cos\big(\tfrac{k-1-2m}{2}(\eta- \vartheta)\big) +\tfrac12 \Big) & k \ {\rm odd} \,;
			\end{cases}
		\end{aligned}
	\end{equation}
	\item For any $k\in\N$, we have that
	\begin{equation}
		\begin{aligned}
			\cos\big(\tfrac{k}{2}(\eta+ \vartheta)\big)& \cos\big(\tfrac{1}{2}(\eta- \vartheta)\big) = \tfrac12 \Big( \cos\big( \tfrac{k+1}{2}\eta + \tfrac{k-1}{2}\vartheta \big)  + \cos\big(\tfrac{k-1}{2}\eta + \tfrac{k+1}{2}\vartheta\big)\Big) \,, \\
			\cos\big(\tfrac{k}{2}(\eta+ \vartheta)\big)& \cos\big(\tfrac{1}{2}(\eta- \vartheta)\big)\cos\big(\tfrac{k-1-2m}{2}(\eta- \vartheta)\big) \\
			& = \tfrac14 \Big( \cos\big( (k-1-m)\eta + (m+1)\vartheta \big)  + \cos(m\eta+ (k-m)\vartheta) \\
			& \quad \ \,+ \cos((k-m)\eta+m\vartheta) + \cos((m+1)\eta+(k-1-m)\vartheta) \Big)\,.
		\end{aligned}
	\end{equation}
\end{enumerate}

\section{Proofs of Lemma \ref{Lema-decompMon} and Lemma \ref{Lema-decompMon.eps}}\label{app.diff.U}
In this appendix we prove both Lemma \ref{Lema-decompMon} and Lemma \ref{Lema-decompMon.eps}. Both proofs follow by adapting the idea developed in the proof  in \cite[Lemma 6.3]{HR21}.

\begin{proof}[Proof of Lemma \ref{Lema-decompMon}]
The proof of item $(i)$ is standard and we omit the details, see for instance \cite[Lemma 6.3]{HR21}. We now prove item $(ii)$. 
\\[1mm]
Let $h\in H_{\circ}^s(\T^2)$ and we look for a decomposition in the form (recalling \eqref{basis-1})
\begin{equation}\label{IAC1}
    \cS h = g  + \rho \,, \quad g = \Pi_{\rm ph} g\,, \quad \rho = \rho_{-1} \bs_{1} \,.
\end{equation}
By projecting \eqref{IAC1} on the mode $\bs_{1}$, we obtain
\begin{equation}\label{IAC2}
    \rho_{-1}:= \cS_{0}[h] := \braket{\cS h , \bs_{1}}_{L^2(\T)} =\braket{ h ,  \cB \bs_{1}}_{L^2(\T)} = \braket{ h ,\big( \cB -{\rm Id}  \big)\bs_{1}}_{L^2(\T)} \,.
\end{equation}
Inserting \eqref{IAC2} into \eqref{IAC1}, we deduce that
\begin{equation}
    g = \cS_{\rm ph} h  = \cS h - \cS_{0}[h] \bs_{1} \,,
\end{equation}
with $\cS_{0}[h]$ as in \eqref{IAC2}.
\\[1mm]
Now, given $g\in H_{\circ}^s(\T^2)$, we want to find $h\in H_{\circ}^s(\T^2)$ such that
\begin{equation}
    \cS_{\rm ph} f = \Pi_{\rm ph} \cS \Pi_{\rm ph} h = g \,.
\end{equation}
We extend it to an equation in $H_{0}^{s}(\T^2)$ by solving
	\begin{equation}\label{IAC3}
		 \cS f = g + h \,,
	\end{equation}
	with $f,g,h \in H_{0}^{s}(\T^2)$ such that, recalling \eqref{proj.1s}
	\begin{equation}\label{IAC4}
		\Pi_{1,\bs}f = 0 \,, \quad \Pi_{1,\bs} g = 0 \,, \quad  \Pi_{1,\bs} h = h \,.
	\end{equation}
	Since $\cS$ is invertible by Lemma \ref{lemma.reparam.0}-$(ii)$, we solve \eqref{IAC3} and  we get
	\begin{equation}\label{IAC5}
		f= \cS^{-1}\big( g+ h  \big)\,.
	\end{equation}
	We now impose, by \eqref{IAC4}, that $\Pi_{1,\bs}f = 0$, that is 
    \begin{equation}
        0 = \braket{\cS^{-1}\big( g + h \big), \bs_{1}}_{L^2(\T)}  = \braket{ g+h, \cB^{-1}\bs_{1}}_{L^2(\T)}  \,.
    \end{equation}
    Looking for $h=\Pi_{1,\bs} h$, we deduce that
    \begin{equation}\label{IAC6}
        h(\vf,\theta) = - \cS_{1}[g(\vf,\,\cdot\,)] \bs_{1}(\theta)\,, \quad \cS_{1}[g(\vf,\,\cdot\,)]:=\frac{\braket{g(\vf,\,\cdot\,),\cB^{-1}\bs_{1}}_{L^2(\T)}}{\braket{\bs_{1},\big(\cB^{-1}-{\rm Id}\big)\bs_{1}}_{L^2(\T)}} \,.
    \end{equation}
    Inserting \eqref{IAC6} into \eqref{IAC5}, we obtain
    \begin{equation}
        f = \Pi_{\rm ph} g = \cS_{\perp}^{-1} g = \cS^{-1}g - \cS_{-1}[g] \cS^{-1}\bs_{1} \,.
    \end{equation}
    with $\cS_{1}[g]$ as in \eqref{IAC6}. This concludes the proof of the representations in \eqref{representations}. The estimates \eqref{est.repre} then follow by \eqref{IAC2}, \eqref{IAC6}, Lemma \ref{lemma.kernelHS}-$(ii)$ and Lemma \ref{lemma.reparam.0}-$(ii)$.
\end{proof}

\begin{proof}[Proof of Lemma \ref{Lema-decompMon.eps}] The proof of item $(i)$ follows as in \cite[Lemma 6.3]{HR21} and we omit the details. Also, we focus here on proving the decomposition
\begin{equation}\label{Q-1.dec}
    \cS_{\varepsilon,\perp}^{-1} = \cS_{\varepsilon}^{-1} + \tQ_{-1}\,, 
\end{equation}
as the decomposition $\cS_{\varepsilon,\perp} = \cS_{\varepsilon} + \tQ_{1}$ follows by similar and easier arguments, see for instance \cite[Lemma 6.3-$(i)$]{HR21}. Therefore, we now prove \eqref{Q-1.dec}.
\\[1mm]
	Given $g\in H_{\perp}^s(\T^2)$, we want to find $f \in H_{\perp}^s(\T^2)$ such that
	\begin{equation}
		\cS_{\varepsilon,\perp} f = \Pi_{\perp}  \cS_{\varepsilon} \Pi_{\perp} f = g \,.
	\end{equation}
	We extend it to an equation in $H_{0}^{s}(\T^2)$ by solving
	\begin{equation}\label{sun1}
		 \cS_{\varepsilon}f = g + h \,,
	\end{equation}
	with $f,g,h \in H_{0}^{s}(\T^2)$ such that
	\begin{equation}\label{moon1}
		({\rm Id}_{L_0^2} - \Pi_{\perp})f = 0 \,, \quad ({\rm Id}_{L_0^2} - \Pi_{\perp}) g = 0 \,, \quad ({\rm Id}_{L_0^2} - \Pi_{\perp}) h = h \,.
	\end{equation}
	Since $\cS_{\varepsilon}$ is invertible by Lemma \ref{almost.straight.lemma}, we solve \eqref{sun1} and  we get
	\begin{equation}\label{sun3}
		f= \cS_{\varepsilon}^{-1}\big( g+ h  \big)\,.
	\end{equation}
	We now require, by \eqref{moon1}, that $({\rm Id}_{L_0^2} - \Pi_{\perp})f = 0$, that is (recalling \eqref{basis-1})
    \begin{align}
        0 & = \braket{\cS_{\varepsilon}^{-1}\big( g + h \big), \bs_{1}}_{L^2(\T)}  = \braket{ g+h, \cB_{\varepsilon}^{-1}\bs_{1}}_{L^2(\T)}  \,, \\
        0 & = \braket{\cS_{\varepsilon}^{-1}\big( g+h\big), \bc_{1}}_{L^2(\T)}  = \braket{g+h , \cB_{\varepsilon}^{-1}\bc_{1}}_{L^2(\T)}   \,,  \\
        0 & = \braket{\cS_{\varepsilon}^{-1}\big( g+h\big), \be_{\pm \tJ}}_{L^2(\T)}  = \braket{g+h, \cB_{\varepsilon}^{-1}\be_{\pm\tJ}}_{L^2(\T)}  \,.
    \end{align}
	Looking for $h = ({\rm Id}_{L_0^2} - \Pi_{\perp})h$, we deduce that
	\begin{equation}
	     h(\vf,\vartheta) = \tA_{\varepsilon}\big[g(\vf,\,\cdot\,) \big](\vartheta) \,,
	\end{equation}
    where $\tA_{\varepsilon}$ is the finite rank operator given by
	\begin{equation}\label{tA-eps.FR}
		\begin{footnotesize}
			\tA_{\varepsilon}[g](\vartheta)  := \begin{pmatrix}
			    \tA_{\varepsilon}^{(-1)}[g] \\ \tA_{\varepsilon}^{(1)}[g]
			\end{pmatrix} \cdot \begin{pmatrix}
			   \bs_{1}(\theta) \\  \bc_{1}(\theta)
			\end{pmatrix}     +  \begin{pmatrix}
			    \tA_{\varepsilon}^{(\tJ)}[g] \\ \tA_{\varepsilon}^{(-\tJ)}[g]
			\end{pmatrix}\cdot \begin{pmatrix}
			    \be_{\tJ}(\theta) \\ \be_{-\tJ}(\theta)
			\end{pmatrix}  \,,
		\end{footnotesize}
	\end{equation}
    with
	\begin{equation}\label{tA.def}
		\begin{footnotesize}
			\begin{pmatrix}
				\tA_{\varepsilon}^{(-1)}[g]	 \\
                \tA_{\varepsilon}^{(1)}[g]	 \\
				\tA_{\varepsilon}^{(\tJ)}[g]\\
				\tA_{\varepsilon}^{(-\tJ)}[g]
			\end{pmatrix} = -\tM_{\varepsilon}^{-1} \begin{pmatrix}
            \braket{g,\cB_{\varepsilon}^{-1}\bs_{1}} \\ 
				\braket{g,\cB_{\varepsilon}^{-1}\bc_{1}} \\ 	\braket{g,\cB_{\varepsilon}^{-1}\be_{\tJ}}\\	\braket{g,\cB_{\varepsilon}^{-1}\be_{-\tJ}}
			\end{pmatrix}\,,
		\end{footnotesize}
	\end{equation}
    and
    \begin{equation}\label{matrix.tM}
        \begin{footnotesize}
            \tM_{\varepsilon} := \begin{pmatrix}
            \braket{\bs_{1},\cB_{\varepsilon}^{-1}\bs_{1}} & 	\braket{\bc_{1},\cB_{\varepsilon}^{-1}\bs_{1}} & 	\braket{\be_{\tJ},\cB_{\varepsilon}^{-1}\bs_{1}}  & 	\braket{\be_{-\tJ},\cB_{\varepsilon}^{-1}\bs_{1}} \\
			 \braket{\bs_{1},\cB_{\varepsilon}^{-1}\bc_{1}} & 	\braket{\bc_{1},\cB_{\varepsilon}^{-1}\bc_{1}} & 	\braket{\be_{\tJ},\cB_{\varepsilon}^{-1}\bc_{1}}  & 	\braket{\be_{-\tJ},\cB_{\varepsilon}^{-1}\bc_{1}} \\
			 \braket{\bs_{1},\cB_{\varepsilon}^{-1}\be_{\tJ}} & 	\braket{\bc_{1},\cB_{\varepsilon}^{-1}\be_{\tJ}} & 	\braket{\be_{\tJ},\cB_{\varepsilon}^{-1}\be_{\tJ}}  & 	\braket{\be_{-\tJ},\cB_{\varepsilon}^{-1}\be_{\tJ}} \\
			\braket{\bs_{1},\cB_{\varepsilon}^{-1}\be_{-\tJ}} & 	\braket{\bc_{1},\cB_{\varepsilon}^{-1}\be_{-\tJ}} & 	\braket{\be_{\tJ},\cB_{\varepsilon}^{-1}\be_{-\tJ}}  & 	\braket{\be_{-\tJ},\cB_{\varepsilon}^{-1}\be_{-\tJ}} 
			\end{pmatrix} \,.
        \end{footnotesize}
    \end{equation}
    Inserting \eqref{tA-eps.FR} into \eqref{sun3}, we obtain
    \begin{equation}
        f  = \Pi_{\perp}f := \cS_{\varepsilon,\perp}^{-1} g = \cS_{\varepsilon}^{-1}g + \tQ_{-1}g\,, 
    \end{equation}
    where we deduce that $\tQ_{-1}$ is the finite rank operator
    \begin{equation}\label{tQ-1.FR}
		\begin{footnotesize}
			\tQ_{-1}[g](\vf,\vartheta) := \begin{pmatrix}
			    \tA_{\varepsilon}^{(-1)}[g] \\ \tA_{\varepsilon}^{(1)}[g]
			\end{pmatrix} \cdot \begin{pmatrix}
			   \cS_{\varepsilon}^{-1}[\bs_{1}](\vf,\vartheta) \\  \cS_{\varepsilon}^{-1}[\bc_{1}](\vf,\vartheta)
			\end{pmatrix}     +  \begin{pmatrix}
			    \tA_{\varepsilon}^{(\tJ)}[g] \\ \tA_{\varepsilon}^{(-\tJ)}[g]
			\end{pmatrix}\cdot \begin{pmatrix}
			    \cS_{\varepsilon}^{-1}[\be_{\tJ}](\vf,\vartheta) \\  \cS_{\varepsilon}^{-1}[\be_{-\tJ}](\vf,\vartheta)
			\end{pmatrix}  \,,
		\end{footnotesize}
	\end{equation}
     with coefficients as in \eqref{matrix.tM}. The estimates \eqref{est.tQ}, \eqref{est.tQ.12} then follow by \eqref{tQ-1.FR}, \eqref{matrix.tM} and  Lemma \ref{almost.straight.lemma}-$(iii)$. This concludes the proof.
\end{proof}

 \end{document}